%% file: main.tex
\documentclass{amsart}
\usepackage[a4paper,margin=3.2cm]{geometry}
\usepackage{graphicx} 
\usepackage{xcolor,amsmath,amssymb,amsthm,amsfonts,hyperref,faktor,mathabx,tikz-cd,stmaryrd,pifont,fdsymbol}

\hypersetup{
	colorlinks=true,
	linkcolor=blue,
	hypertexnames=false,
	filecolor=black,      
	urlcolor=blue,
	citecolor=blue
}
\usepackage[nameinlink]{cleveref}

\theoremstyle{plain}
\newtheorem{theorem}{Theorem}[section]
\newtheorem*{theorem*}{Theorem}
\newtheorem{proposition}[theorem]{Proposition}
\newtheorem{lemma}[theorem]{Lemma}
\newtheorem{corollary}[theorem]{Corollary}

\newtheorem*{conjecture}{Conjecture}

\theoremstyle{definition}
\newtheorem{definition}[theorem]{Definition}
\newtheorem*{convention}{Convention}
\theoremstyle{remark}
\newtheorem{remark}[theorem]{Remark}
\newtheorem{example}[theorem]{Example}
\newtheorem*{example*}{Example}

\definecolor{Blue}{RGB}{0, 122, 255}
\definecolor{Green}{RGB}{12, 163, 2}
\definecolor{Orange}{RGB}{245, 154, 35}
\definecolor{DarkBlue}{RGB}{5, 70, 143}
\definecolor{Pink}{RGB}{255, 0, 162}
\definecolor{Red}{RGB}{209, 0, 0}
\definecolor{Purple}{RGB}{130, 33, 139}
\definecolor{Teal}{RGB}{51, 186, 150}
\definecolor{Yellow}{RGB}{255, 201, 0}

\usepackage{tikz,subcaption}
\usetikzlibrary{positioning, calc}
\usetikzlibrary{arrows.meta}
\usetikzlibrary{decorations.markings}
\usetikzlibrary{shapes.geometric}
\tikzstyle{mutable}=[circle,draw=black,fill=black,inner sep=0pt,minimum size=6]
\tikzstyle{frozen}=[draw=black,fill=black,inner sep=0pt,minimum size=6]
\tikzstyle{mutableBig}=[circle,draw=black,fill=black,inner sep=0pt,minimum size=10]
\tikzstyle{frozenBig}=[draw=black,fill=black,inner sep=0pt,minimum size=10]
\tikzstyle{arrow}=[->,very thick]
\tikzstyle{barrow}=[->,very thick,>=Triangle]
\tikzstyle{xarrow}=[->,very thick,color=Blue,>={Triangle[open]}]
\tikzstyle{halfarrow}=[->,very thick,dashed]
\tikzstyle{angleindicator}=[very thick,color=Yellow]
\tikzstyle{carrow}=[->>,very thick,color=Green]
\tikzstyle{box}=[rectangle,draw=black]
\tikzstyle{line}=[thick]
\tikzstyle{dline}=[thick,dashed]
\tikzstyle{qarrow}=[very thick,decoration={
    markings,
    mark=at position 0.6 with {\arrow{>}}}]
\tikzstyle{sarrow}=[thick,decoration={
    markings,
    mark=at position 0.7 with {\arrow{>}}}]
\tikzstyle{straightangle}=[black,thick]
\tikzstyle{crossangle}=[Blue, thick,dashdotted]
\tikzstyle{vertex}=[circle,inner sep=0pt,minimum size=0]
\tikzstyle{halfcircle} =[circle, draw=black, 
        path picture={\fill[black] (path picture bounding box.north west) rectangle (path picture bounding box.south);}]
\tikzstyle{flag}=[circle,draw=black,fill=black,inner sep=0pt,minimum size=4]
\tikzset{>=stealth, node distance = 40}
\newcommand{\distance}{1.2}
\newcommand{\distL}{1.4}
\usepackage{junkcommands}
\tikzset{
 NCnode/.pic={
    \node[circle,draw,line width=0.8mm,minimum size=6mm] (-circle) {};
    \draw[thick] (-0.6,0) -- (-0.3,0);
    \draw[thick] (0.6,0) -- (0.3,0);
    \filldraw[-] (0,0) -- (90:0.35) arc (90:270:0.35) -- cycle;
    \draw[thick] (0,-0.7) -- (0,0.7);
  }
}
\tikzset{
 NCnodeHalfR/.pic={
    \node[circle,draw,line width=0.8mm,minimum size=6mm] (-circle) {};
    \draw[thick] (0.6,0) -- (0.3,0);
    \draw[thick] (0,-0.7) -- (0,0.7);
    \filldraw[-] (0,0) -- (90:0.3) arc (90:270:0.3) -- cycle;
  }
}
\tikzset{
 NCnodeHalfL/.pic={
    \node[circle,draw,line width=0.8mm,minimum size=6mm] (-circle) {};
    \draw[thick] (-0.6,0) -- (-0.3,0);
    \draw[thick] (0,-0.7) -- (0,0.7);
    \filldraw[-] (0,0) -- (90:0.3) arc (90:270:0.3) -- cycle;
  }
}

\newcommand{\CL}{\mathrm{Cl}}
\newcommand{\R}{\mathbb{R}}
\newcommand{\Z}{\mathbb{Z}}
\newcommand{\Q}{\mathbb{Q}}
\newcommand{\C}{\mathbb{C}}
\newcommand{\In}{\mathrm{In}}
\newcommand{\Out}{\mathrm{Out}}
\newcommand{\SP}{\mathrm{Sp}}
\newcommand{\SL}{\mathrm{SL}}
\newcommand{\GL}{\mathrm{GL}}
\newcommand{\SO}{\mathrm{SO}}
\newcommand{\Spin}{\mathrm{Spin}}
\newcommand{\sau}{\sigma\hspace{-.8mm}\tau}

\newcommand{\SPA}{\SP_2(A,\sigma)}
\newcommand{\weaviso}[1]{
\begin{tikzpicture}[#1]
	\node[vertex](s) at (0,-0.6ex)	 [label = above: $\sim$]{};
	\draw (-0.8ex,-0.3ex) -- (0.8ex,0.3ex);
	\draw (-0.8ex,0.3ex) -- (0.8ex,-0.3ex);
\end{tikzpicture}}
\newcommand{\pruned}{\mathrm{pr}}

\newcommand{\quiverAlgebra}{\mathcal{Q}}
\newcommand{\seedAlgebra}{\mathcal{S}}
\newcommand{\clusterAlgebra}{\mathcal{A}}

\newcommand{\AngleSkeleton}{\mathrm{AS}}
\newcommand{\angleParity}{\ell}
\newcommand{\canonAngle}[1]{\Delta_{#1}}
\newcommand{\canonLabel}{L}

\newcommand{\CP}{\C{P}}
\newcommand{\Bl}{\mathrm{Bl}}

\renewcommand{\bar}[1]{\overline{#1}}
\renewcommand{\epsilon}{\varepsilon}
\renewcommand{\tilde}[1]{\widetilde{#1}}
\newcommand{\keyword}[1]{\textbf{#1}}

\title[noncommutative Polygonal Cluster Algebras]{noncommutative Polygonal Cluster Algebras}

\author{Zachary Greenberg}
\address{Max Planck Institute for Mathematics in the Sciences\\
Inselstr. 22\\04103 Leipzig, 
Germany }
\email{greenberg@mis.mpg.de}

\author{Dani Kaufman}
\address{University of Copenhagen \\
         Department of Mathematical Sciences \\
             2100 Copenhagen \o, Denmark }
\email{dk@math.ku.dk}

\author{Merik Niemeyer}
\address{Max Planck Institute for Mathematics in the Sciences\\
Inselstr. 22\\04103 Leipzig, 
Germany }
\email{niemeyer@mis.mpg.de}

\author{Anna Wienhard}
\address{Max Planck Institute for Mathematics in the Sciences\\
Inselstr. 22\\04103 Leipzig, 
Germany }
\email{wienhard@mis.mpg.de}

\thanks{
D.K. was supported by the Danish National Research Foundation (CPH-GEOTOP-DNRF151) and the Alexander von Humboldt Foundation. Z.G., M.N. and A.W. were supported by the European Research Council under ERC-Advanced Grant 101018839. A.W. thanks the Hector Fellow Academy for support. \\
}

\date{}

\begin{document}

\begin{abstract}
    We define a new family of noncommutative generalizations of cluster algebras called polygonal cluster algebras. These algebras generalize the noncommutative surfaces of Berenstein-Retakh, and are inspired by the emerging theory of $\Theta$-positivity for the groups $\Spin(p,q)$. They are generated by mutations of quivers which we call ST-compatible, and which encode the order of the products that appear in the exchange relations. We show that these ST-compatible quivers can be represented by tilings of surfaces by polygons, a generalization of the description of surface type cluster algebras. As examples, we construct tilings which produce ST-compatible versions of the Del Pezzo quivers and the quivers first described by Le for Fock-Goncharov coordinates for Lie groups of type $B$. We show that polygonal cluster algebras have natural evaluations in Clifford algebras, which we use to produce noncommutative generalizations of the Somos sequences and to parameterize the $\Theta$-positive semigroup of $\Spin(2,n)$. We indicate how this will be done for the semigroup in $\Spin(p,q)$ and how one will give coordinates for general $\Theta$-positive representations into $\Spin(p,q)$.
\end{abstract}

\maketitle

\begin{figure}[!h]
    \centering
    \begin{subfigure}[b]{0.48\textwidth}
	\centering
	\begin{tikzpicture}
     	\node[regular polygon,draw,regular polygon sides=5, minimum size=2.8*2.3cm, very thick] (p) at (0,0) {};
     	\draw[angleindicator] (p.side 1) -- (p.corner 4);
      	\draw[angleindicator] (p.side 2) -- (p.corner 5);
    	\draw[angleindicator] (p.side 3) -- (p.corner 1);
    	\draw[angleindicator] (p.side 4) -- (p.corner 2);
    	\draw[angleindicator] (p.side 5) -- (p.corner 3);
        \orientedarc (p.corner 1:p.corner 2)
        \orientedarc (p.corner 2:p.corner 3)
        \orientedarc (p.corner 3:p.corner 4)
        \orientedarc (p.corner 4:p.corner 5)
        \orientedarc (p.corner 5:p.corner 1)
        \centerarc[straightangle](p.corner 1)(216:324:0.4)
        \centerarc[straightangle](p.corner 2)(288:396:0.4)
        \centerarc[straightangle](p.corner 3)(0:108:0.4)
        \centerarc[straightangle](p.corner 4)(72:180:0.4)
        \centerarc[straightangle](p.corner 5)(144:252:0.4)
        \node[regular polygon,draw,regular polygon sides=5, minimum size=2.8*2.3cm, thick] (p) at (0,0) {};

        \node[circle,line width=0.8mm,minimum size=6mm] (-circle) at (p.side 3) {};
	\end{tikzpicture}
\end{subfigure}
    \begin{subfigure}[b]{0.48\textwidth}
 \scalebox{1}{
        \begin{tikzpicture}[]
        \node[regular polygon,draw,regular polygon sides=5, minimum size=2.8*2.3cm, very thick,dotted] (p) at (0,0) {};        
        \pic[name=A,rotate=270-2*72] at (p.side 1) {NCnodeHalfL};  
        \pic[name=B,rotate=270-72] at (p.side 2) {NCnodeHalfL};        
        \pic[name=C,rotate=270] at (p.side 3) {NCnodeHalfL};
        \pic[name=D,rotate=270+72] at (p.side 4) {NCnodeHalfL};
        \pic[name=E,rotate=270+2*72] at (p.side 5) {NCnodeHalfL};

        \draw[barrow] (C-circle.120) to[out=120,in=360-12] (B-circle.360-12);
        \draw[barrow] (B-circle.48) to[out=48,in=276] (A-circle.276);
        \draw[barrow] (A-circle.336) to[out=336,in=204] (E-circle.204);
        \draw[barrow] (E-circle.264) to[out=264,in=132] (D-circle.132);
        \draw[barrow] (D-circle.192) to[out=192,in=60] (C-circle.60);

        \draw[xarrow] (C-circle.90+60) to[out=90+60,in=306-60] (A-circle.306-60);
        \draw[xarrow] (B-circle.18+60) to[out=18+60,in=234-60] (E-circle.234-60);
        \draw[xarrow] (A-circle.6) to[out=306+60,in=162-60] (D-circle.162-60);
        \draw[xarrow] (E-circle.234+60) to[out=234+60,in=90-60] (C-circle.90-60);
        \draw[xarrow] (D-circle.162+60) to[out=162+60,in=18-60] (B-circle.18-60);
    \end{tikzpicture}
            }
\end{subfigure}

    \caption{A decorated polygon and associated ST-compatible quiver.}
    \label{fig:coolTitlePicture}
\end{figure}
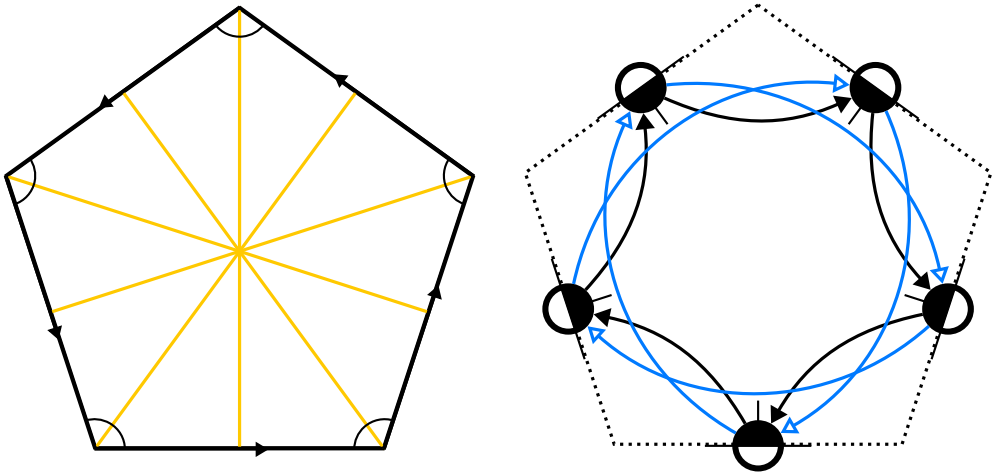

\setcounter{tocdepth}{2}
\makeatletter
\def\l@subsection{\@tocline{2}{0pt}{2.5pc}{5pc}{}}
\makeatother

\pagebreak
\tableofcontents

\section{Introduction}

\subsection{Motivation from Teichm\"uller Theory}\label{sec:introTM}

Cluster algebras were introduced by Fomin and Zelevinsky in the early 2000's \cite{fomin2002cluster} and are a class of constructively defined commutative rings. Roughly speaking, they can be thought of as growing from a seed consisting of a quiver and a set of (commuting) variables via a process known as mutation.

The simplest family of cluster algebras are those of ``surface type'' originally defined in \cite{fomin2008clusterTri}. In this case the seeds are in correspondence with ideal triangulations of a surface with punctures and mutations are captured by flipping arcs of the triangulation. The quiver for each seed is easily constructed from the triangulation by inscribing a clockwise cycle in each triangle, see \Cref{fig:PentagonTriangulationQuiver}.

\input{figurePentagonTriangulationQuiver}

These cluster algebras have an important geometric realization: their positive real valued points can be identified with Penner's decorated Teichmüller space \cite{penner1987decorated}. Points in this space parameterize discrete and faithful representations $\pi_1(\Sigma) \to \SL_2(\R)$ along with some extra data associated to the punctures of $\Sigma$.

Fock, Goncharov and Shen \cite{fock2006moduli, goncharov2019quantum} greatly generalized this theory in the context of higher Teichm\"uller spaces for split real Lie groups, see also \cite{le2019cluster, gilles2021fock} for further work. 

These spaces can be seen as positive real valued points of a cluster algebra built out of an ideal triangulation of $\Sigma$ using a more complicated quiver inscribed in each ideal triangle.
The inscribed quiver depends on the Dynkin type of $G$, with the $A_1$ case $(G=\SL_2(\R))$ reproducing the quiver used for the original surface type cluster algebra.

Berenstein and Retakh introduced a noncommutative version of a surface type cluster algebra \cite{berenstein2018noncommutative}. 
A geometric realization similar to Penner's decorated Teichm\"uller space of these noncommutative algebras was observed in the construction of cluster coordinates for maximal representations $\pi_1(\Sigma) \to \SP_{2n}(\R)$ \cite{alessandrini2019noncommutative}, which form another family of higher Teichm\"uller spaces associated to Hermitian Lie groups of tube type.  

The theory of $\Theta$-positivity introduced in \cite{guichard2022generalizing} now provides a common framework for both these families of higher Teichm\"uller spaces, as well as two new families. 
Whenever a semisimple Lie groups admits a $\Theta$-positive structure, the set of $\Theta$-positive representations form connected components of the representation variety consisting of discrete and faithful representations 
\cite{guichard2021, beyrer2021, beyrer2024}, thus they give rise to higher Teichm\"uller spaces. 

Simple Lie groups admitting a $\Theta$-positive structure are classified \cite{guichard2022generalizing}, and in each case (other than the split real case) the $\Theta$-positive structure encodes a new root system which we call the type of the structure.  
\begin{enumerate}
    \item (Split real): $G$ is split real.
    \item (Type $A_1$): $G$ is a real form with type $C_n$ restricted root system.
    \item (Type $B_n$): $G$ is a real form with type $B_{n+1}$ restricted root system.
    \item (Type $G_2$): $G$ is a real form with type $F_4$ restricted root system.
\end{enumerate}

The symplectic group $\SP_{2n}(\R)$ appears twice on this list, once as a split real group and as a group with a type $A_1$ $\Theta$-positive structure. This corresponds to the Hitchin component and the components of maximal representations. Thus, the results of \cite{alessandrini2019noncommutative,rogozinnikov2020symplectic} state that Berenstein and Retakh's algebra is the algebra whose geometric realization parameterizes  type $A_1$ $\Theta$-positive representations. This suggests there should be analogous cluster structures whose geometric realizations parameterize types $B_n$ and $G_2$. Our construction will provide a framework to realize type $B_n$ $\Theta$-positive representations.

\subsection{Main Results}

Our overarching goal for this paper is to describe a family of noncommutative algebras with a cluster-like structure generalizing Berenstein-Retakh's algebra. Our algebras are ``grown'' from an initial seed associated to decorated quivers which we require to be \emph{ST-compatible}. This is an algebraic condition that leads to strong restrictions on the structure of the quiver (\Cref{thm:GSTRestrictions}), allowing us to associate a decorated polygonal tiling of a surface to each ST-compatible quiver $Q$, where the tiles correspond to locally transitive tournaments in $Q$, see \Cref{fig:coolTitlePicture} for an example. 

Moreover we associate a non-commutative algebra called the seed algebra, $\seedAlgebra(Q)$, to an ST-compatible quiver $Q$. This seed algebra satisfies the triangle relations of the Berenstein-Retakh algebra with some additional centrality conditions which essentially imply that the square of each cluster variable is central. The seed algebra replaces the ambient field of fractions used in the usual definition of a cluster algebra.

We extend quiver mutation to ST-compatible quivers and call mutations admissible if they return new ST-compatible quivers and hence new decorated tilings. 
\begin{theorem*}[\ref{thm:seedalgebra_iso}]
    Let $Q$ be an ST-compatible quiver.  An admissible mutation of a quiver $Q$ at $k$ induces an isomorphism of seed algebras $\seedAlgebra_Q \cong \seedAlgebra_{\mu_k(Q)}$. 
\end{theorem*}

The subalgebra $\clusterAlgebra_Q$ generated by images of cluster variables is a noncommutative analogue of a cluster algebra. In fact if all the cluster variables are assumed to commute, $\clusterAlgebra_Q$ recovers the classical commutative cluster algebra. This contrasts with the noncommutative surface cluster algebra which only recovers the cluster algebra localized at every cluster variable. Moreover $\clusterAlgebra_Q$ satisfies a Laurent phenomenon:
\begin{theorem*}[\ref{thm:laurent}]
    In the cluster algebra $\clusterAlgebra_Q$ generated from a ST-compatible quiver $Q$, every cluster variable is a noncommutative Laurent polynomial in the initial cluster variables.
\end{theorem*}

The family of ST-compatible quivers is much larger than the family of quivers associated to triangulations of a surface. In particular it contains quivers whose undecorated counterparts are those quivers associated to toric Del Pezzo surfaces and the initial quivers for representations of type $B_n$ given in \cite{le2019cluster}. 

We call an ST-compatible quiver fruitful if every mutation is admissible. We show that not every ST-compatible quiver is fruitful. Nevertheless there are arbitrarily large ST-compatible quivers (\Cref{thm:large_fuitefull_forks}). We also prove the initial quiver for $B_2$ representations is fruitful.

Just as the usual cluster algebra can be evaluated in positive real numbers, $\clusterAlgebra_Q$ can be evaluated in the Clifford algebra $\CL(\R^{1,n})$ (\Cref{thm:clifford_evaluations_general}).
\begin{example*}
    The Somos-4 and 5 sequences can be generated by an infinite mutation sequence on the quiver given in \Cref{fig:Del Pezzo Quivers}. This quiver for Somos-4 is associated to a decorated tiling of a punctured disk with digons and triangles. Somos-5 is given by a tiling of a punctured torus with one boundary component by triangles and squares.
    The Somos mutation sequences are admissible and provide a noncommutative generalization of the Somos-4 and 5 sequences in $\CL(\R^{1,n})$.
\end{example*}

Another family of fruitful quivers are those associated to the usual triangulations of marked surfaces $\Sigma$. We call the associated cluster algebras $\mathcal{B}_1(\Sigma)$ as we expect positive points in $\CL(\R^{1,n-1})$ to parameterize $\Theta$-positive representations of type $B_1$ into $\Spin(2,n)$. While we do not prove this in general here, we prove the following initial result:
 
 \begin{theorem*}[\ref{prop:so2n_warmpu}]
     Let $\Sigma$ be a disk with four marked points. A subset of $\CL(\R^{1,n-1})$-points of $\mathcal{B}_1(\Sigma)$ parameterize an open subset of the group $\Spin(2,n)$. Moreover the positive such points parameterize the $\Theta$-positive semigroup $\Spin(2,n)_\Theta^{>0}$.
 \end{theorem*}
 
As the initial quivers of type $B_{p-1}$ given in \cite{le2019cluster,goncharov2019quantum} have ST-compatible counterparts, we conjecture that their positive $\CL(\R^{1,q-p+1})$-points parameterize $\Theta$-positive representations into $\Spin(p,q)$.

\begin{remark}
    There are several other noncommutative generalizations of cluster algebras. A noncommutative cluster-like algebra of type $A_n$ was given in \cite{goncharov2021spectral}. In this case mutation is only allowed at vertices that are ``surface type'' while our algebra will allow more general mutations.
    Super cluster algebras \cite{os-ClusterAlgebrasWithGrassmanVariables,shemyakova-superClusterAlgebras} add additional anticommuting variables. These extended quivers provide a different noncommutative generalization of Somos sequences. 
\end{remark}

\subsection*{Structure of the paper}

First we describe the algebra $\mathcal{B}_1(\Sigma)$  in \Cref{sec:warmup}. This provides a self contained version of the theory which parallels that of Berenstein-Retakh and does not use ST-compatible quivers. 
In \Cref{sec:2colQuivers}, we will introduce 2-colored quivers in the commutative setting. We will describe them in terms of 2-colored exchange matrices and explain how to understand their combinatorics through an appropriate unfolding. We endow the 2-colored quivers with additional decoration in \Cref{sec:ncQuivers} and use this to define ST-compatible quivers. In \Cref{sec:Polygons} we define decorated tilings, which provide a combinatorial model for ST-compatible quivers. \Cref{sec:Algebra} contains our main construction of noncommutative polygonal cluster algebras, a proof of their Laurent phenomenon, and a description of their evaluation in Clifford algebras.
After that, in \Cref{sec:examples}, we will provide a wide variety of examples, before briefly considering a connection of our construction with higher Teichm\"uller theory in \Cref{sec:TM}.

\section{Warm-Up: The Algebra \texorpdfstring{$\mathcal{B}_1(P_n)$}{B1 of Pn}}\label{sec:warmup}
In this section we review the noncommutative cluster algebras introduced by Berenstein-Retakh \cite{berenstein2018noncommutative}, which can be considered as a noncommutative $A_1$ theory. We then discuss a generalization which can be considered as a noncommutative $B_1$ theory. 

\subsection{Review of the cluster-like structure on noncommutative surfaces}
To construct Berenstein-Retakh's algebra \cite{berenstein2018noncommutative} one considers a regular $n$-gon, $P_n$, with vertices labeled $1,...,n$ and the set of all its oriented chords. Associate to the chord $i\rightarrow j$ a variable $x_{ij}$ and consider the $\mathbb{Q}$-algebra $\mathcal{A}_1(P_n)$ generated by the set $\{x_{ij},x_{ij}^{-1}|i\neq j\in [1,n]\}$, subject to the
\begin{itemize}
	\item{\keyword{triangle relations}: for any $i,j,k\in [1,n]$
		\begin{equation}\label{eq:triangleRel}
			T_j^{ki}:=x_{kj}^{-1}x_{ki}x_{ji}^{-1}=x_{ij}^{-1}x_{ik}x_{jk}^{-1}=: T_j^{ik}\,,
		\end{equation}}
	\item{and the \keyword{exchange relations}: for any cyclic $(i,j,k,\ell)$ in $[1,n]$
		\begin{equation}\label{eq:exchRel}
			x_{\ell j}=x_{\ell k}x_{ik}^{-1}x_{ij}+x_{\ell i}x_{ki}^{-1}x_{kj}\,.
		\end{equation}}
\end{itemize}

We call the expression $T_j^{ki}$ defined in the triangle relation the \keyword{noncommutative angle} of the triangle $\{ijk\}$ based at $j$. As the variables are associated to oriented chords of the $n$-gon, we see that the triangle relation indeed only contains variables associated to (orientations) of the triangle with corners $i,j,k$, while the exchange relation is related to a flip of the diagonal in the 4-gon with corners $i,j,k,\ell$, as depicted in \Cref{fig:FlipSurfaceCase}.

\input{figureFlipSurfaceCase.tex}

The exchange relation is equivalent to the assertion that the angles add in a natural way: $$ T_k^{j\ell}= x_{jk}^{-1}x_{j\ell}x_{k\ell}^{-1} = x_{jk}^{-1}x_{ji}x_{ki}^{-1}+x_{ik}^{-1}x_{i\ell}x_{k\ell}^{-1} = T_k^{ji}+T_k^{i\ell} $$

The algebra carries a cluster(-like) structure with clusters formed by sets of variables which correspond to triangulations of the $n$-gon and mutations given by the flips in these triangulations where the variables mutate according to the exchange relations. In fact, we can generate the entire algebra starting with the variables associated with a single triangulation subject to the triangle conditions. After a flip of triangulation the new variable is expressed via the exchange relation, and the new triangle conditions are implied by the old. 

This is essentially a noncommutative version of the surface type cluster algebra for a disk with $n$ marked points on the boundary, the only difference is that we have imposed that all of the cluster variables are invertible. Indeed, this commutative algebra can be reobtained by setting $x_{ij}=x_{ji}$, and assuming that all variables commute, which causes the triangle relations to become trivial.

\begin{remark}\label{rem:A1angleMonodromy}
   The angles based at different vertices of the same triangle can be related to one another in a simple way: 
\begin{equation} \label{eqn:BRAngleTransition}
    T_j^{ki} = x_{kj}^{-1}x_{ki}x_{ji}^{-1} = x_{ij}^{-1} x_{ij} x_{kj}^{-1}x_{ki}x_{ji}^{-1}   = ( x_{ji} T_i^{kj} x_{ij})^{-1}. 
\end{equation} 
Thus it suffices to consider one angle from each triangle The other angles are recovered as images under compositions of maps for each edge $p_{ji}$ defined  by $p_{ji}(T) = (x_{ji}Tx_{ij})^{-1}$. We will see this structure in the generalizations defined in \Cref{sec:ncQuivers,sec:Polygons}.
\end{remark}

\subsection{Geometric realization of noncommutative surfaces}
To obtain a geometric interpretation of $\mathcal{A}_1(P_n)$ we look at its points in some noncommutative ring $A$ endowed with an anti-involution $\sigma$. We note that $\mathcal{A}_1(P_n)$ has a natural anti-involution, which we also denote by $\sigma$, that is defined by $ \sigma(x_{ij})= x_{ji}$. 
\begin{definition}
    An $(A,\sigma)$-point of the algebra $\mathcal{A}_1(P_n)$ is an algebra homomorphism
    \begin{equation*}
        \phi \colon \mathcal{A}_1(P_n) \rightarrow A \text{ such that } \phi\circ\sigma = \sigma\circ\phi\,.
    \end{equation*}
\end{definition}
Note that for $\phi$ to be a homomorphism, $\phi$ must satisfy that 
\begin{enumerate}
    \item $\phi(x_{ij})$ is invertible in $A$
    \item $\phi(T_j^{ik})$ is fixed by $\sigma$
\end{enumerate}

Given a noncommutative ring with involution $(A,\sigma)$, one can generalize the symplectic group $\SP_{2n}(\R)$ to the group $\SPA$ \cite{alessandrini2022-SP2Asigma}. In particular $\SP_{2n}(\R)$ can be realized as $\SPA$ with $A$ being the ring $M_n(\R) $ of $n\times n$ matrices over $\R$ and $\sigma$ given by matrix transpose. The group $\SPA$ should be thought of as a group of $2\times2$ matrices over the noncommutative ring $A$. We can use the algebra $\mathcal{A}_1(P_4)$ to parameterize this group as proved in \cite{GreenbergEtAl2024MathrmSL_2}:

\begin{proposition}\label{prop:A1pointsParametermizeSPA}
    The points of $\mathcal{A}_1(P_4)$ which send the variables $x_{12},x_{34}$ to $1\in A$ parameterize elements of the group $\SPA$ with invertible entries.
\end{proposition}

In general it can be hard to define $(A,\sigma)$-points because of invertibility. However when $A$ has an additional \emph{positive structure}, a point is defined by considering a single triangulation of $P_n$. Informally, a positive structure is a subset $A^+\subset A$ of \emph{positive} elements, that is closed under inverse, addition, and symmetric conjugation $a \to xa\sigma(x) $ by elements of $x \in A^+$. If the image of $x_{ij}$ for each arc of the triangulation is invertible and the image of each angle $T_j^{ik}$ is positive then the collection obtained by mutation will be an $(A,\sigma)$-point. This is analogous to the situation in commutative cluster algebras where specifying a positive real number in an initial seed guarantees all other cluster variables are positive. 

Interest in $\SPA$ comes from the fact that groups carrying a $\Theta$-positive structure of type $A_1$ are of this form for some ring $A$ with the positive structure required above. The following proposition was shown for $\SP(2n,\R)$ in \cite{alessandrini2019noncommutative} and extended to more general $(A,\sigma)$ in \cite{rogozinnikov2020symplectic,kineider2022-FramedLocalSystems}:
\begin{proposition}
Let $\Sigma$ be an oriented surface with boundary and punctures. The $\Theta$-positive decorated representations from $\pi_1(\Sigma)$ to $\SPA$ are parameterized by the positive $(A,\sigma)$-points of $\mathcal{A}_1(\Sigma)$.
\end{proposition}

However, it has been noted in \cite[Section 10]{rogozinnikov2020symplectic} that the groups $\SO(2,n)$, which carry a $\Theta-$positive structure of type $B_1$, are not realized directly by the group $\SPA$ for any choice of $A$. Instead they are found as a subgroup of $\SPA$ when $A$ is a Clifford algebra. Our next goal is to construct an algebra $\mathcal{B}_1(P_n)$ such that
\begin{enumerate}
    \item there is a natural map $\mathcal{A}_1(P_n) \to \mathcal{B}_1(P_n)$, and
    \item homomorphisms $\mathcal{B}_1(P_4) \to A$ such that $x_{12},x_{34}$ are mapped to $1\in A$ parameterize the group $\Spin(2,n)$ in analogy with \Cref{prop:A1pointsParametermizeSPA}.
\end{enumerate}

\subsection{The Algebra \texorpdfstring{$\mathcal{B}_1(P_n)$}{B1 of Pn}}\label{sec:AlgebraB1Pn}

Again consider $P_n$ a regular $n$-gon with vertices labeled $1,...,n$.
We now associate to each chord two variables $x_{ij},y_{ij}$ and consider the $\mathbb{Q}$-algebra $\mathcal{D}_2(P_n)$ generated by the set $\left\{x_{ij}^{\pm1},y_{ij}^{\pm1}|i,j\in [1,n]\right\}$, subject to the triangle relations and exchange relations in both of the $x$ and $y$ variables independently. We may think of this ring as a noncommutative tensor product of two copies of the ring $\mathcal{A}_1(P_n)$. The notation $\mathcal{D}_2$ is suggestive of the fact that the $D_2$ Dynkin diagram is $A_1\times A_1$ and folds to the Dynkin diagram $B_1$. 

Now we define the ring $\mathcal{B}_1(P_n)$ as the ring obtained from $\mathcal{D}_2(P_n)$ by imposing the  \keyword{symmetry conditions} and \keyword{centrality conditions} for each arc:
\begin{align}
    \quad & x_{ij}y_{ji} = y_{ij}x_{ji} \tag{Symmetry}\\
    \quad &x_{ij}y_{ji} \quad \text{is central} \tag{Centrality}
\end{align}
and localizing by the center. 

The algebra $\mathcal{B}_1(P_n)$ comes with two commuting anti-involutions we will denote by $\sigma$ and $\tau$. $\sigma$ is defined as before by $$\sigma(x_{ij})=x_{ji} \qquad \sigma(y_{ij})=y_{ji}$$ and $\tau$ is defined by $$\tau(x_{ij}) = y_{ji}$$ 

\begin{definition}\label{def:normMap}
    We define the \keyword{norm} of an element $x\in \mathcal{B}_1(P_n)$ by $N(x) = x\tau(x)$ and we say that $x$ has \keyword{central norm} if $N(x)$ is central and fixed by $\sigma$. 
    Linearizing the norm, we define a symmetric pairing $b: \mathcal{B}_1(P_n) \times \mathcal{B}_1(P_n) \to \mathcal{B}_1(P_n)$ by $$b(\alpha,\beta) = N(\alpha+\beta) - N(\alpha)-N(\beta) =  \alpha\tau(\beta) + \beta\tau(\alpha)$$
\end{definition}

Since $x$ is invertible, $N(x)=N(\tau(x))$, and we see that the symmetry and centrality conditions are equivalent to all four variables associated to a single arc having the same central norm.

Using the centrality of the norm, we can rewrite the triangle and exchange relations in $\mathcal{B}_1(P_n)$ in a way which does not involve any inverses of elements

\begin{definition}
    We call $\Delta_j^{ki}:= y_{jk}x_{ki}y_{ij}$ the \keyword{canonical angle} of the triangle $ijk$ based at $j$. In $\mathcal{B}_1(P_n)$ we have the
\begin{itemize}
	\item{\keyword{canonical triangle relations}: for any $i,j,k\in [1,n]$
		\begin{equation}\label{eq:triangleRel_noinverse}
			\Delta_j^{ki}=y_{jk}x_{ki}y_{ij}=y_{ji}x_{ik}y_{kj}= \Delta_j^{ik}\,,
		\end{equation}}
	\item{and the \keyword{canonical exchange relations}: for any cyclic $(i,j,k,\ell)$ in $[1,n]$
		\begin{equation}\label{eq:exchRel_noinverse}
			N(x_{ik})x_{\ell j}=x_{\ell k}y_{ki}x_{ij}+x_{\ell i}y_{ik}x_{kj}\,.
		\end{equation}}
\end{itemize}
\end{definition}

These canonical triangle and exchange relations are clearly equivalent to the usual triangle and exchange relations.

\begin{remark}
    The canonical exchange relations are equivalent to the following additivity property of the canonical angles  
    \begin{equation}\label{eq:canonical_adding_warmup}
        N(x_{ki})\Delta_k^{j\ell}=N(x_{k\ell})\Delta_k^{ji} + N(x_{kj})\Delta_k^{i\ell}\,
    \end{equation}
\end{remark}

\begin{definition}
    The cluster algebra $\mathcal{C}_n \subset \mathcal{B}_1(P_n)$ is the subalgebra generated by $\{x_{ij},y_{ij}\}$ for each arc in $P_n$
\end{definition}

Note that $\mathcal{A}_1(P_n)$ does not naturally contain an analogue of the usual cluster algebra in it, due to the fact that the exchange relations cannot be written without any inverses. On the other hand, if we pass to a commutative version of the algebra $\mathcal{C}_n\subset\mathcal{B}_1(P_n)$ by setting $x_{ij}=y_{ij}$ and making all variables commute, we recover the classical cluster algebra.

\begin{remark}\label{rem:non-central}
    The centrality conditions of only the variables associated to arcs in a single triangulation \emph{do not} imply the centrality conditions for the rest of the arcs.
\end{remark}

If $\alpha$ and $\beta$ both have central norm, then clearly their product does as well. However, their sum does if and only if $b(\alpha,\beta)$ is central. Note that the exchange relations explicitly produce a new element with central norm out of the sum of two elements of central norm. Thus, it is natural to consider exactly when the sum of elements with central norm has central norm. To help understand the structure of these centrality conditions, we will first consider some maps from $\mathcal{B}_1(P_n)$ to Clifford algebras.

\subsection{Geometric realization of \texorpdfstring{$\mathcal{B}_1(P_n)$}{B1 of Pn}}

Our goal is to study the space of algebra homomorphisms $\mathcal{B}_1(P_n) \to (B,\sigma,\tau)$ where $B$ is a noncommutative ring with two commuting anti-involutions. The key difficulty in constructing such a map are the centrality conditions, otherwise such a map would consist of two maps from $\mathcal{A}_1(P_n)$ related by $\tau$. 

Let $\mathbb{F}$ be a field, and let $V$ be an $\mathbb{F}$-vector space with a quadratic form $q$. We write $\CL(V,q)$ for the Clifford algebra over $(V,q)$ (we refer to \Cref{app:Clifford} for details and notations about Clifford algebras). For the rest of the following we fix a vector $e\in V$ with $q(e)\neq 0$. We define two involutions $\sigma_e,\tau$ on $\CL(V,q)$ as follows 
\begin{equation}
    \tau(x)= x^T \quad \quad\quad\sigma_e(x) = \frac{1}{q(e)}ex^Te
\end{equation}
where $x\in\CL(V,q)$ and $x^T$ is the transpose map.

\begin{definition}
    We write $\CL(V,q)^{\sigma_e}$ for the vector space of elements which are fixed by $\sigma_e$. We also write $\CL(V,q)^\tau$ for the subgroup of elements $x \in \CL(V,q)^\times$ which satisfy that $x\tau(x) $ is central.
\end{definition}
 
The vector $e$ gives a canonical embedding of $V$ into $\CL(V,q)^{\sigma_e}$ by $v \to ev$ for $v\in V$, we write $eV$ for this image. We also note that the Clifford group $\Gamma(V,q)$ is a subgroup of $\CL(V,q)^\tau$.

We can now define a natural family of homomorphisms $\mathcal{B}_1(P_n) \to \CL(V,q)$.
\begin{definition}
    A \emph{$\CL(V,q,e)-$point} of the algebra $B_n$ is an algebra homomorphism $\phi:\mathcal{B}_1(P_n) \to \CL(V,q)$ such that 
    \begin{enumerate}
        \item $\phi$ intertwines the two anti-involutions defined on  $\mathcal{B}_1(P_n)$ and $\CL(V,q,e)$.
        \item $\phi(x_{ij}) \in \Gamma(V,q)$ for all arcs.
        \item $\phi(\Delta_{j}^{ik}) \in eV$ for all angles.
    \end{enumerate}
\end{definition}

Fix $\Delta$ a triangulation of the polygon $P_n$.
\begin{theorem}\label{thm:Clifford_warmup}
    A $\CL(V,q,e)-$point of $\mathcal{B}_1(P_n)$ $\phi$ may be constructed by picking a collection of elements $\{x_{ij}\} \subset \CL(V,q)$ for $\{ij\} \in \Delta$ such that 
    \begin{enumerate}
        \item $x_{ji} =\sigma_e(x_{ij})$
        \item $x_{ij} \in \Gamma(V,q)$.
        \item $x_{ji}^{T}x_{jk}x_{ik}^{T}\in eV$ for all triangles $\{ijk\} \in \Delta$ 
    \end{enumerate}
    and defining the map by abuse of notation $\phi(x_{ij})= x_{ij}$, $\phi(y_{ji})=(x_{ij})^T$, provided that each exchange relation produces a new invertible element of $\CL(V,q)$. 
\end{theorem} 
\begin{proof}
    Clearly this choice determines the images of the $x_{ij}$ and  $y_{ij}$ in $\mathcal{B}_1(P_n)$ for all arcs. We need to check the triangle and centrality conditions. Consider a single flip of the arc $13$ in a square labeled ${1,2,3,4}$. The exchange relation for a flip of a triangulation implies the triangle relations in the new triangulation because of the additivity of angles. Since the angles in our initial triangulation all live in $eV$, their sum does as well. Thus, the new angles after a flip are elements of $eV$ and therefore are elements of $\Gamma(V,q)$ as well, provided that they are invertible. One easily sees that invertibility of the new coordinate $x_{24}$ is equivalent to invertibility of the angle $x_{23}^{T}x_{24}x_{34}^{T}$. Thus by assumption this angle is invertible and is hence in the Clifford group. This further implies that $x_{24}$ is in the Clifford group as well since both $x_{23}$ and $x_{34}$ are assumed to be. This implies that $x_{24}(x_{24})^T$ is central.

    This argument applied to each flip implies that the variables associated to each triangulation satisfy the centrality conditions.

\end{proof}

 Let $(V,q)\cong \R^{1,n-1}$ be an $n$-dimensional $\R$-vector space and $q$ a quadratic form of signature $(1,n-1)$ on $V$. We pick $e \in V$ to be a vector of norm $+1$. We denote by $V^+$ the proper cone of vectors with positive norm which contains $e$. This cone forms the analogue of a positive structure for the Clifford algebra. 

\begin{definition}\label{def:positive_points}
    A $\CL\left(\R^{1,n-1},e\right)$-point of $\mathcal{B}_1(P_n)$ is \keyword{positive} if all of the angles are mapped into $eV^+$ and all of the arcs are mapped to elements of $\Gamma\left(\R^{1,n-1}\right)$ with positive real norm.
\end{definition}

\begin{proposition}
    A positive $\CL\left(\R^{1,n-1},e\right)$-point of $\mathcal{B}_1(P_n)$ can be given simply by defining the images of the arcs in an initial triangulation according to \Cref{thm:Clifford_warmup}, the invertibility of the images of the new arcs is implied by positivity. 
\end{proposition}
\begin{proof}
    This follows simply because $ eV^+$ is a proper convex cone consisting of invertible elements which is fixed by symmetric conjugation by elements in $\Gamma(\R^{1,n-1})$ with positive real norm. 
\end{proof}

\begin{proposition}\label{prop:so2n_warmpu}
    The $\CL(\R^{1,n-1},e)-$points of $\mathcal{B}_1(P_4)$ which send $x_{12},x_{34}$ to $1$ parameterize an open subset of the group $\Spin(2,n)$ while the \emph{positive} such points parameterize the $\Theta-$positive semigroup $\Spin(2,n)^{>0}$.
\end{proposition}
\begin{proof}
    In \cite{rogozinnikov2020symplectic} a generic element of $\SP(\CL(V,q,e))$ is factored as
    \[\begin{pmatrix}1 & \\ y & 1\end{pmatrix}\begin{pmatrix}x & \\  & \sigma(x)^{-1}\end{pmatrix}\begin{pmatrix}1 & z\\  & 1\end{pmatrix}\]
    for $x,y,z\in \Gamma(V,q,e)$, $y,z$ fixed by $\sigma$. In \cite[Theorem 2.10.20]{rogozinnikov2020symplectic} they use this factorization to provide an isomorphism between this group and $\Spin(2,n)$.\\
    By \cite[Theorem 7.24]{GreenbergEtAl2024MathrmSL_2} such a matrix is parameterized by $\mathcal{A}_1(P_4)$ with \[x_{14} =  x z \hspace{2pc} x_{13}= x \hspace{2pc} x_{23}=y x \hspace{2pc} x_{12}=x_{34} = 1\]
    However the assumptions on $\mathcal{A}_1(P_4)$ only show that the elements of the matrix are units of the Clifford algebra, not necessarily elements of the Clifford group. The assumptions on $\mathcal{B}_1(P_4)$ force each element to be in the Clifford group and for each angle to be in $eV$.
    
    The conditions on the $\Theta-$positive semi-group defined in \cite{guichard2022generalizing} translate directly to the positive conditions of \Cref{def:positive_points}.
\end{proof}

\subsection{Structure of centrality conditions}

The centrality conditions of $\mathcal{B}_1(P_n)$ are given for every chord of $P_n$. We now give an alternate definition of the algebra based on a single triangulation that we will generalize in \Cref{sec:Algebra}.

First consider $\Delta$  the ``fan triangulation'' of polygon $P_n$ based at the vertex 1, i.e the triangulation consisting of the arcs $13,14,\dots, 1(n-1)$, cf. \Cref{fig:PentagonTriangulationQuiver}.
\begin{proposition}
    The algebra defined as a quotient of $\mathcal{D}_2(P_n)$ by the relations that 
    \begin{itemize}
        \item $\left\{x_{1j},x_{i,i+1}\,|\, 2\leq i\leq n-1, 2\leq j\leq n\right\}$ have central norm
        \item $b\left(\Delta_1^{i,i+1},\Delta_1^{j,j+1}\right)$ is central 
    \end{itemize}
    and localized by the center, is isomorphic to $\mathcal{B}_1(P_n)$.
\end{proposition}
\begin{proof}
    Since $b\left(\Delta_1^{i,i+1},\Delta_1^{i+1,j+2}\right)$ is central, the sum $\Delta_1^{i,i+2}$ has central norm and inductively $\Delta_1^{i,j}$ has central norm for any $i,j$. Thus any arc $x_{ij}$ has central norm as well. Conversely, if all arcs $x_{ij}$ have central norm, any sum of consecutive angles has central norm, and thus the symmetric pairing between any two angles has to be central.
\end{proof}

Note that the centrality condition for $x_{ij}$ immediately implies the condition for $x_{ji}$, so that there are at most $ {\binom{n}{2}}$ independent centrality conditions. In fact these are all independent.
\begin{corollary}
        All of the  $ {\binom{n}{2}}$ centrality conditions are independent. 
\end{corollary}
\begin{proof}
It is easier to use the centrality conditions on the fan triangulation. We see  $2n-3$ centrality conditions on the arcs of the fan triangulation are independent by first considering a single triangle with vertices $1,2,3$. The only dependence that could arise among the three centrality conditions for the three arcs are due to the triangle relations and symmetry conditions. In \cite[Theorem 2.7]{berenstein2018noncommutative} they show the algebra generated by $x_{ij}$ with the triangle relation is rank 5. Thus including the $y_{ij}$ variables we obtain a free algebra of rank 7 generated by $x_{12},x_{23},x_{31},x_{21},y_{21},y_{13},y_{32}$ after imposing the 3 symmetry conditions. It is then clear that $x_{12}y_{21},x_{31}y_{13},x_{23}y_{23}$ can replace $y_{21},y_{13},y_{32}$ as independent generators. Thus the centrality condition on these three elements are independent.

 Once we add these elements to the center of $\mathcal{B}_1(P_n)$, we can consider the vector space over the fraction field of the center generated by the $\Delta_1^{i,i+1}$ along with the (now) inner product $b$. This inner product is determined by a symmetric $(n-2)\times(n-2)$-matrix. The diagonal entries are simply products of norms of variables in $\Delta$, while the off diagonal entries are pairs. 

Thus in total we have $2n-3 + {\binom{n-2}{2}} = {\binom{n}{2}} $ centrality relations. 

\end{proof}

We would like to extend this formulation to an arbitrary triangulation. In the fan triangulation we could take every angle to be based at the same vertex, so it made sense to add them. In an arbitrary triangulation we need a procedure to move all the angles to a single vertex where we impose the centrality conditions.

As in \Cref{rem:A1angleMonodromy} the canonical angles can be related by maps along the arcs of the triangulation. For a directed arc $ij$ we can define a map $p_{ij}: \mathcal{B}_1(P_n)^\times \to \mathcal{B}_1(P_n)^\times$ by \[p_{ij}(\alpha) = \frac{1}{N(x_{ij})} \sigma\tau(x_{ji}\alpha x_{ij}) = \frac{1}{N(x_{ij})} y_{ji} (\sigma\tau)(\alpha) y_{ij},\]
As before we have that $p_{ij}(\Delta_i^{jk}) = \Delta_j^{ki}$. Note that since we replaced inversion with $\tau$, this map is linear in contrast to the situation in $\mathcal{A}_1(P_n)$.

Fix a general triangulation $T$ of $P_n$ with triangles $T_i$, and fix a vertex $1$. In the 1-skeleton of $T$ chose a path $\gamma_i$ from from a vertex $v_i$ of $T_i$ to $1$. Let $\Delta_i$ be the angle for $T_i$ based at $v_i$. For each $\gamma_i$ we define a map $m_i$ given by composing $p_{e}$ for each edge $e$ of $\gamma_i$. If the choice of paths is generic enough we have the following:

\begin{proposition}
    The algebra $\mathcal{B}_1(P_n)$ can be equivalently defined as the quotient of $\mathcal{D}_2(P_n)$ by the relations
    \begin{itemize}
        \item $x_e$ has central norm for every edge $e$ of $T$
        \item $\Delta_i = \sigma(\Delta_i)$ for each triangle $T_i$
        \item $b(m_i\Delta_i,m_j\Delta_j)$ is central for every pair of triangles $T_i,T_j$.
    \end{itemize}
    The vector space over the center of $\mathcal{B}_1(P_n)$ generated by the images of all of the angles of triangles in $\Delta$ mapped to the vertex 1 along paths in $T$ is dimension $n-2$ and is independent of the choice of $\Delta$ and $T$.  
\end{proposition}
This will be a consequence of \Cref{thm:finiteCentrConds}. The key idea is that  the monodromy around a triangle $ijk$ of the maps $p_{ij}$ acts via a reflection formula
\begin{equation*}
    v \to -v + 2\frac{b\left(v,\Delta_i^{jk}\right)}{b\left(\Delta_i^{jk},\Delta_i^{jk}\right)}\Delta_i^{jk} \,.
\end{equation*} 
This allows the centrality assumptions at 1 to propagate to any sum of angles as needed.

\subsection{Towards generalizations}
The next step is to construct an algebra $\mathcal{B}_{p}(\Sigma)$ whose points parameterize representations into $\Spin(p+1,q)$. From \cite{guichard2022generalizing} these representations should be of type $B_{p}$ which is a Dynkin diagram with roots of two different sizes when $p>1$. \Cref{fig:FGQuivers} shows the quivers first described in \cite{le2019cluster} of types $B_2$ and $B_3$. The \emph{short} root in $B_p$ is represented by the \emph{big} nodes in the quiver. This is described by a skew-symmetrizable exchange matrix where big nodes are assigned 2 as their skew-symmetrizer. We remark the subquiver of big nodes is of surface type, but this no longer holds after mutations at the small nodes.
\input{figureFGQuivers}

In $B_1$ all the nodes are big nodes. This corresponds to the two ``layers'' of isomorphic subalgebras in $\mathcal{B}_1(P_n)$, one generated by $x_{ij}$ and the other generated by $y_{ji}$. To capture this structure with quivers and mutations we consider two disconnected copies of the same quiver associated to $\mathcal{A}_1(P_n)$. This makes sense as $B_1$ should be thought of as ``folding'' the two tails of a $D_2$ Dynkin diagram (by definition the same as an $A_1\times A_1$ diagram). At the level of algebras the two disconnected layers of quivers are the $D_2$ type cluster algebra and the folding is the quotient to establish the centrality and symmetry conditions. 

To parameterize $\Theta$-positive representations of type $B_p$, the big nodes take values in a Clifford algebra and the small nodes take values in the real numbers. Thus we will assign commutative variables to small nodes and noncommutative variables to big nodes. Furthermore in the exchange relation for small nodes, the variables for big nodes always appear squared. We replace the squared term with the norm (\Cref{def:normMap}) to generalize this exchange relation. 

We also notice that mutation at a small node adds pairs of parallel and crossing arrows between the big nodes. Thus we consider pairs of parallel and crossing arrows to be a commutative phenomenon. Moreover these mutations connect the two layers corresponding the big nodes (\Cref{fig:MixedMutation}). To capture these two kinds of arrows we will need the notion of 2-colored quivers.

\begin{figure}[!hb]
    \centering
    \begin{tikzpicture}
        \begin{scope}[shift={(0,0)}]
            \node[mutable]		(1a) at (0,0.3*\distL)		[label = left: $1$]	{};
            \node[mutable]		(1b) at (0,-0.3*\distL)		[label = left: $1'$]	{};
            \node[mutable]		(2a) at (\distL,0.3*\distL)	[label = right: $2$]	{};
            \node[mutable]		(2b) at (\distL,-0.3*\distL)[label = right: $2'$]	{};
            \node[mutable]      (3)  at (\distL/2,-0.9*\distL) [label = right:$3$] {};
            \draw[arrow]		(1a) to (2a);
            \draw[arrow]		(1b) to (2b);
            \draw[arrow]        (3) to (1a);
            \draw[arrow]        (3) to (1b);
            \draw[arrow]        (2a) to (3);
            \draw[arrow]        (2b) to (3);
        \end{scope}
        \begin{scope}[shift={(1.8*\distL,0)}]
            \node[vertex] (v1) at (0,-0.5) []{};
            \node[vertex] (v2) at (0.5*\distL,-0.5) []{};
            \draw[<->] (v1) -- (v2)
            node[above,midway,text centered, text width=2cm]{$\mu_3$};
        \end{scope}
        \begin{scope}[shift={(3*\distL,0)}]
            \node[mutable]		(1a) at (0,0.3*\distL)		[label = left: $1$]	{};
            \node[mutable]		(1b) at (0,-0.3*\distL)		[label = left: $1'$]	{};
            \node[mutable]		(2a) at (\distL,0.3*\distL)	[label = right: $2$]	{};
            \node[mutable]		(2b) at (\distL,-0.3*\distL)[label = right: $2'$]	{};
            \node[mutable]      (3)  at (\distL/2,-0.9*\distL) [label = right:$3$] {};
            \draw[arrow]		(2a) to (1b);
            \draw[arrow]		(2b) to (1a);
            \draw[arrow]        (3) to (2a);
            \draw[arrow]        (3) to (2b);
            \draw[arrow]        (1a) to (3);
            \draw[arrow]        (1b) to (3);
        \end{scope}
    \end{tikzpicture}
    \caption{Mutation at 3 ``mixes'' layers of quiver.}
    \label{fig:MixedMutation}
\end{figure}

\section{Two-colored quivers}\label{sec:2colQuivers}
In this section, we will introduce the notion of 2-colored quivers, for now still in a commutative setting. The colored edges will describe a canonical unfolding related to the classical folding procedure for quivers \cite{felikson2012cluster}.

\subsection{Reminder: (One-colored) quivers}

Recall the following definitions, e.g. from \cite{fomin2017introduction}.
\begin{definition}
    A \keyword{quiver} is an oriented graph, which has neither 2-cycles nor self-loops. The vertices are called \keyword{nodes}, the edges \keyword{arrows}.
\end{definition}

For a quiver $Q$, we will denote the set of nodes by $Q_0$, the set of arrows by $Q_1$. Let $Q$ be a quiver, $Q_0=[1,n]$. We can produce a new quiver by \keyword{mutation} at $k\in Q_0$:
\begin{definition}\label{def:CommutativeQuiverMutation}
    The quiver $\mu_k(Q)$ after mutation at $k$ is produced by these steps:
    \begin{enumerate}
        \item For every sequence $i\rightarrow k\rightarrow j$ introduce an arrow $i\rightarrow j$.
        \item Reverse any arrow connected to $k$.
        \item Cancel any resulting 2-cycles.
    \end{enumerate}
\end{definition}

\begin{proposition}\label{prop:quiverMutInvolution}
    Let $Q$ be a quiver and $k\in Q_0$. Mutation at $k$ is an involution, i.e.
    \begin{equation*}
        \mu_k^2(Q)=Q\,.
    \end{equation*}
\end{proposition}

This can also be encoded in terms of the exchange matrix, which allows for a more general version of quivers, which we will need subsequently.

\begin{definition}
	Let $Q$ be a quiver with $Q_0=[1,n]$. Its \keyword{exchange matrix} $B(Q)=(b_{ij})_{i,j\in [1,n]}$ is the integer matrix with entries
	\begin{equation*}
		b_{ij}=|\{\alpha\in Q_1|\alpha:i\to j\}|-|\{\alpha\in Q_1|\alpha:j\to i\}|\,.
	\end{equation*}
\end{definition}
Note that at least one of the terms vanishes because $Q$ is 2-acyclic and that $B(Q)$ is a skew-symmetric matrix. On these exchange matrices we can perform \keyword{matrix mutation}.
\begin{definition}
	Let $B=(b_{ij})\in\mathbb{Z}^{n\times n}$ be skew-symmetric and $1\leq k\leq n$. The matrix $\mu_k(B)=(b'_{ij})$ after mutation at $k$ has entries
	\begin{equation}
		b'_{ij}=\begin{cases}
			-b_{ij} & k\in\{i,j\}\\
			b_{ij}+b_{ik}^+b_{kj}^+-b_{ik}^-b_{kj}^- & k\notin\{i,j\}
		\end{cases}
	\end{equation}
	where $z^\pm=\pm\mathrm{max}\{0,\pm z\}$ for $z\in\mathbb{Z}$ denotes the usual positive, respectively negative, part.
\end{definition}

This provides two equivalent descriptions of mutation in the following sense.
\begin{proposition}\label{prop:classicMutEquiv}
	For a quiver $Q$ with a node $k$, we have
	\begin{equation*}
		B(\mu_k(Q))=\mu_k(B(Q))\,.
	\end{equation*}
\end{proposition}

Matrix mutation generalizes to the larger class of skew-symmetrizable matrices:
\begin{definition}\label{def:skewSym}
	A matrix $B\in\mathbb{Z}^{n\times n}$ is called \keyword{skew-symmetrizable} if there exists a diagonal matrix $D$ with diagonal entries in $\mathbb{Z}_{>0}$ such that $DB$ is skew-symmetric.
\end{definition}

\begin{lemma}\label{lemma:Dmut}
	For $B$ and $D$ as in the previous definition, $D\mu_k(B)$ is skew-symmetric for any $k\in [1,n]$.
\end{lemma}

These skew-symmetrizable exchange matrices can be encoded either by \emph{valued quivers}\cite{marsh2013lectureOnClusterAlgebras} with weighted edges or by \emph{quivers with weighted nodes} \cite{zickert2020RankTwoLieGroups}. 
\begin{example}\label{ex:SkewSymmetrizableQuivers}
    Consider the following skew symmetrizable matrix represented as a valued quiver and a quiver with weighted nodes.
    \begin{equation*}
        \hspace*{2pc}
        \begin{pmatrix}
            1 & & \\
            & 1 & \\
            & & 2
        \end{pmatrix}
        \begin{pmatrix}
           0 & 1 &0\\
          -1 & 0 & -2 \\
           0 & 1 & 0
        \end{pmatrix}
        \hspace{2pc}
        \vcenter{\begin{tikzpicture}
            \begin{scope}[shift={(0,0)}]
                \node[mutable]	(X1) at (0,0) [label = above: $1$]	{};
                \node[mutable]	(X2) at (\distL,0) [label = above: $2$] {};
                \node[mutable]	(X3) at (2*\distL,0) [label = above: $3$]	{};
                \draw[arrow] (X1) to (X2);
                \draw[arrow] (X3) to (X2);
                \node[] (e) at (1.5*\distL,0) [label= below: {$(2,1)$}] {};
            \end{scope}
            \begin{scope}[shift={(3*\distL,0)}]
                \node[mutable]	(X1) at (0,0) [label = above: $1$]	{};
                \node[mutable]	(X2) at (\distL,0) [label = above: $2$] {};
                \node[mutableBig]	(X3) at (2*\distL,0) [label = above: $3$]	{};
                \draw[arrow] (X1) to (X2);
                \draw[arrow] (X3) to (X2);
            \end{scope}
        \end{tikzpicture}}
    \end{equation*}
    The corresponding valued quiver has an arrow $i \xrightarrow{|b_{ij}|~|b_{ji}|} j$ if $b_{ij} > 0$ with two weights  corresponding to $b_{ij}$ and $b_{ji}$. By convention, an unlabeled arrow corresponds to $b_{ij} = -b_{ji} = 1$. \\
    The quiver with weighted nodes instead encodes the matrix $D$ as the weights of each node. Here the first two small nodes have weight 1 while the last large node has weight 2. The entry $b_{ij}$ is recovered by $d_{j}/\mathrm{gcd}(d_i,d_j)$ multiplied by the number of arrows from $i$ to $j$ . 
\end{example}

In many cases, including every case in this paper, such a matrix can also be encoded by \emph{folding} a larger quiver representing a skew symmetric matrix.

\begin{definition}
    A \keyword{folding} of a quiver is a choice of disjoint sets of nodes, called \keyword{folding groups}, which satisfies the following \keyword{folding conditions}:
    \begin{itemize}
        \item nodes in a such a folding group are not connected by any arrows, and
        \item this is satisfied after any number of group mutations.
    \end{itemize}
    A \keyword{group mutation} is carried out by performing quiver mutation on every node in the chosen group. We call a folding \keyword{valid} if it satisfies the folding conditions.
\end{definition}
\begin{example}
    The following quiver has a folding with groups $\{1\},\{2\}, \{3,3'\}$. 
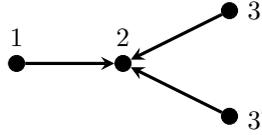
    \begin{equation*}
        \begin{tikzpicture}
            \node[mutable]	(X1) at (0,0) [label = above: $1$]	{};
            \node[mutable]	(X2) at (\distL,0) [label = above: $2$] {};
            \node[mutable]	(X3a) at (2*\distL,0.5*\distL) [label = right: $3$]	{};
            \node[mutable]	(X3b) at (2*\distL,-0.5*\distL) [label = right: $3'$]	{};
            \draw[arrow] (X1) to (X2);
            \draw[arrow] (X3a) to (X2);
            \draw[arrow] (X3b) to (X2);
        \end{tikzpicture}
    \end{equation*}
    This folding represents the same matrix in \Cref{ex:SkewSymmetrizableQuivers} as it has the same combinatorics, i.e. group mutation corresponds exactly to the skew symmetrizable matrix mutation. 
    Note that the size of the group corresponds exactly to the weight of each node.
\end{example}

Generally, it can be difficult to tell if a folding is valid. However, there is a special case in which this can be guaranteed: 
\begin{lemma}[\cite{kaufman2023special}]
    Let $\epsilon$ be an order two automorphism of $Q$. Then folding of $Q$ into orbits of $\epsilon$ is a valid folding. Moreover each quiver group mutation equivalent to $Q$ also has $\epsilon$ as an automorphism.
\end{lemma}

If a skew symmetrizable matrix $B$ can be encoded by a folding, the larger quiver is called an  \emph{unfolding} of $B$.  In general such an unfolding is not unique. In the next section we will explain how to add data to the quiver which encodes the unfolding in the case where the folding is by an order 2 automorphism.

\subsection{Two-colored quivers and mutation}
We will now introduce 2-colored quivers and their mutation. We will first describe this intrinsically and then explain how it can be understood in terms of an unfolding, where the coloring captures information that is usually lost in the folding.

\begin{definition}\label{def:2coloredQuivers}
    A \keyword{2-colored quiver} is a quiver $Q$ with a partition
    \begin{equation*}
    	Q_1 = Q_1^\shortparallel\sqcup Q_1^\times
    \end{equation*}
    of the arrows into two disjoint sets. The arrows in these are called \keyword{parallel} and \keyword{crossing}, respectively.
\end{definition}

\begin{remark}
	A 2-colored quiver is the data $\left\{Q,Q_1^\shortparallel\right\}$ but we will often suppress the partition of arrows in the notation and simply write $Q$. It is usually clear whether we are referring to the 2-colored quiver or its underlying quiver.
\end{remark}

In figures we will use black arrows (with a solid tip) as parallel arrows and blue arrows (with an empty tip) as crossing arrows. Pairs of parallel and crossing arrows will play an important role in our theory. Thus, we make the following definition for later convenience.

\begin{definition}
	Let $Q$ be a 2-colored quiver. The \keyword{pruned quiver} $Q^\pruned$ is obtained from $Q$ by removing all pairs of parallel and crossing arrows between any two nodes.
\end{definition}

We extend the notion of mutation to 2-colored quivers:

\begin{definition}
	Let $Q$ be a 2-colored quiver with a node $k\in Q_0$. The \keyword{mutation at $k$} produces a new 2-colored quiver, $\mu_k(Q)$, through the following steps
    \begin{enumerate}
        \item{ For every sequence $i\rightarrow k\rightarrow j$, introduce an arrow $i\rightarrow j$, whose color is determined by the \keyword{coloring rules}
	        \begin{align*}
	    		\textcolor{black}{\rightarrow} + \textcolor{black}{\rightarrow} = \textcolor{black}{\rightarrow}\\
	    		\textcolor{Blue}{\rightarrowtriangle} + \textcolor{Blue}{\rightarrowtriangle} = \textcolor{black}{\rightarrow}\\
	    		\textcolor{Blue}{\rightarrowtriangle} + \textcolor{black}{\rightarrow} = \textcolor{Blue}{\rightarrowtriangle}
			\end{align*}}
        \item reverse any arrow connected to $k$,
        \item cancel any resulting 2-cycles of the same color.
    \end{enumerate}
\end{definition}

\begin{remark}
    For a slight generalization of 2-colored quivers, one can allow 2-cycles in a 2-colored quiver $Q$ under the condition that the arrows have different colors. In step 3 of the mutation procedure we only cancel 2-cycles of the same color, meaning that the resulting quiver after mutation may contain such a 2-cycle with different colors. This is related to the notion of a special folding \cite{kaufman2023special}, and we will describe a case where this is useful in \Cref{sec:puncDigon}. For now we will simply assume that we only perform mutations which do not lead to 2-cycles.
\end{remark}

\begin{example}
	An example of a mutation on a simple 2-colored quiver is shown in \Cref{fig:colMutEx}.
	\input{figureColoredMutation.tex}
\end{example}

\begin{proposition}
	Mutation of 2-colored quivers is an involution.
\end{proposition}

\begin{proof}
    Simply combine \Cref{prop:quiverMutInvolution} with the coloring rules. 
\end{proof}

\subsection{Two-colored exchange matrices}
We will now give a description of two-colored quivers and their mutation in terms of an exchange matrix $C(Q)$. To do so, we need to pass to an appropriate ring.

\begin{definition}
	The \keyword{split-complex numbers} are the ring $\mathbb{R}[\varepsilon]/(\varepsilon^2-1)$. We denote by $\mathbb{Y}=\mathbb{Z}[\varepsilon]/(\varepsilon^2-1)$ the integer points, and by $\mathbb{Y}_{>0}=\{a+\varepsilon b\neq 0|a,b\in\mathbb{Z}_{\geq 0}\}$ the positive points. The \keyword{real} and \keyword{imaginary part} functions on $\mathbb{R}[\varepsilon]/(\varepsilon^2-1)$ are defined by
	\begin{align*}
		\mathfrak{R}(a+b\varepsilon)&=a\,,\\
		\mathfrak{I}(a+b\varepsilon)&=b\,.
	\end{align*}
\end{definition}

\begin{definition}
	Let $Q$ be a 2-colored quiver with $Q_0=[1,n]$. Its \keyword{2-colored exchange matrix} $C(Q)\in\mathbb{Y}^{n\times n}$ has entries
	\begin{align*}
		c_{ij}=|\{\alpha &\in Q_1^\shortparallel|\alpha:i\to j\}|-|\{\alpha\in Q_1^\shortparallel|\alpha:j\to i\}|\\
				&+\varepsilon\cdot\big(|\{\alpha\in Q_1^\times|\alpha:i\to j\}|-|\{\alpha\in Q_1^\times|\alpha:j\to i\}|\big)\,.
	\end{align*}
\end{definition}
In other words, we think of parallel arrows as having weight 1 while crossing arrows have weight $\varepsilon$.

\begin{proposition}
	Let $(Q,Q_1^\shortparallel)$ be a 2-colored quiver. Then
	\begin{equation*}
		C(Q)|_{\varepsilon=1}=B(Q)\,,
	\end{equation*}
	i.e. the coloring of the quiver can be ``forgotten''.
\end{proposition}

The positive and negative part functions on $\mathbb{Z}$ naturally extend to $\mathbb{Y}$ as
\begin{align*}
    (a+b\varepsilon)^+=a^++b^+\varepsilon\,,\\
    (a+b\varepsilon)^-=a^-+b^-\varepsilon\,
\end{align*}
which we can use to extend matrix mutation to 2-colored matrices.

\begin{definition}
	Let $C\in\mathbb{Y}^{n\times n}$ be a skew-symmetric matrix and $k\in [1,n]$. The \keyword{matrix mutation} at $k$ yields the matrix $\mu_k(C)$ with entries
	\begin{equation}
		c'_{ij}=\begin{cases}
			-c_{ij} & k\in\{i,j\}\\
			c_{ij}+c_{ik}^+c_{kj}^+-c_{ik}^-c_{kj}^- & k\notin\{i,j\}
		\end{cases}
	\end{equation}
\end{definition}

\begin{proposition}
	For a 2-colored quiver $Q$ with $k\in Q_0$, we have
	\begin{equation*}
		C(\mu_k(Q))=\mu_k(C(Q))\,.
	\end{equation*}
\end{proposition}

\begin{proof}
	We only need to check that the mutation of the exchange matrix produces arrows according to the coloring rule of quiver mutation. But this is true, simply because $1\cdot\varepsilon=\varepsilon\cdot 1=\varepsilon$ and $1\cdot 1=\varepsilon\cdot\varepsilon=1$.
\end{proof}

\begin{remark}
    One can perform a more general construction of exchange matrices over a group ring, $\Z[G]$, for any group $G$. Under this interpretation, the construction above is for $\Z/2\Z$, the cyclic group of order 2.
\end{remark}

\begin{definition}
	A matrix $C\in\mathbb{Y}^{n\times n}$ is \keyword{skew-symmetrizable} (over $\mathbb{Y}$) if there exists a diagonal matrix $D$ with diagonal entries in $\mathbb{Y}_{>0}$ such that $DC$ is skew-symmetric.
\end{definition}

As in the one-colored case, we can extend matrix mutation to skew-symmetrizable matrices, and we have an analogue of \Cref{lemma:Dmut} in this case:

\begin{lemma}\label{lemma:skewSymmMut}
	For $C$ and $D$ as in the previous definition, $D\mu_k(C)$ is skew-symmetric for any $k\in [1,n]$.
\end{lemma}

\subsection{Unfolding of two-colored quivers}
Recall that sometimes an unfolding can be used to describe mutation of skew-symmetrizable matrices in terms of quivers in the one-colored case. If we pass to 2-colored quivers, which by definition have skew-symmetric exchange matrix over $\mathbb{Y}$, a similar unfolding procedure can be used to interpret the different types of arrows. This is motivated by quivers with skew-symmetrizable exchange matrix that arise in higher Teichm\"uller theory, see \Cref{sec:mixed2colQuivers} and \Cref{sec:TM} for more details.

\begin{definition}\label{def:Unfolding2ColQuivers}
	Let $Q$ be a 2-colored quiver with $Q_0=[1,n]$. The \keyword{unfolded quiver} is the quiver $\overline{Q}$ with vertex set $\overline{Q}_0=\{i,i'\,|\,i\in [1,n]\}$ and a pair of parallel, resp. crossing, arrows for every arrow in $Q_1^\shortparallel$, resp. $Q_1^\times$, as shown in \Cref{fig:arrowTypes}.
\end{definition}

\begin{figure}[ht]
	\begin{center}
		\begin{tikzpicture}[]
			\begin{scope}
				\node[mutableBig]	(1) at (0,0)			[label = above: $1$]		{};
				\node[mutableBig]	(2) at (\distL,0)	[label = above: $2$]		{};
				\draw[barrow]		(1) to (2);
			\end{scope}
			\begin{scope}[shift={(1.5*\distL,0)}]
				\node[vertex]	(1) at (0,-0.2)			[label = $\mathrm{=}$]		{};
			\end{scope}
			\begin{scope}[shift={(2*\distL,0)}]
				\node[mutable]		(1a) at (0,0.3*\distL)		[label = above: $1$]	{};
				\node[mutable]		(1b) at (0,-0.3*\distL)		[label = below: $1'$]	{};
				\node[mutable]		(2a) at (\distL,0.3*\distL)	[label = above: $2$]	{};
				\node[mutable]		(2b) at (\distL,-0.3*\distL)[label = below: $2'$]	{};
				\draw[arrow]		(1a) to (2a);
				\draw[arrow]		(1b) to (2b);
			\end{scope}
		\end{tikzpicture}\hspace{6pc}%
        \begin{tikzpicture}[]
			\begin{scope}
				\node[mutableBig]	(1) at (0,0)			[label = above: $1$]		{};
				\node[mutableBig]	(2) at (\distL,0)	[label = above: $2$]		{};
				\draw[xarrow]		(1) to (2);
			\end{scope}
			\begin{scope}[shift={(1.5*\distL,0)}]
				\node[vertex]	(1) at (0,-0.2)			[label = $\mathrm{=}$]		{};
			\end{scope}
			\begin{scope}[shift={(2*\distL,0)}]
				\node[mutable]		(1a) at (0,0.3*\distL)		[label = above: $1$]	{};
				\node[mutable]		(1b) at (0,-0.3*\distL)		[label = below: $1'$]	{};
				\node[mutable]		(2a) at (\distL,0.3*\distL)	[label = above: $2$]	{};
				\node[mutable]		(2b) at (\distL,-0.3*\distL)[label = below: $2'$]	{};
				\draw[arrow]		(1a) to (2b);
				\draw[arrow]		(1b) to (2a);
			\end{scope}
		\end{tikzpicture}
	\end{center}
	\caption{Parallel and crossing arrows and their unfoldings.}
    \label{fig:arrowTypes}
\end{figure}

Clearly, the 2-colored quiver $Q$ is determined by its unfolding $\bar{Q}$, and in fact the unfolding can be used to describe the mutation of 2-colored quivers.

\begin{proposition}
	Let $Q$ be a two-colored quiver with a node $k$. The unfolding $\overline{\mu_k(Q)}$ can be obtained from $\overline{Q}$ by group mutation at $\{k,k'\}$.
\end{proposition}

\begin{proof}
	It suffices to check that the composition of arrows under group mutation of the unfolded quiver matches the coloring rules for arrows introduced by 2-colored mutation. This follows since composition of two pairs of crossing arrows introduces a pair of parallel arrows etc.
\end{proof}

\subsection{Isomorphisms of two-colored quivers}
For 2-colored quivers there are two natural notions of isomorphism:
\begin{definition}
	A map $\xi:Q\to Q'$ between 2-colored quivers is an \keyword{isomorphism} if it is an isomorphism of the underlying quivers which preserves the partition of arrows, i.e. $\xi\left(Q_1^\shortparallel\right)=(Q'_1)^\shortparallel$ and $\xi\left(Q_1^\times\right)=(Q'_1)^\times$. If there exists an isomorphism between $Q$ and $Q'$ we write $Q\cong Q'$.
\end{definition}
A weaker version is provided by considering the unfoldings.
\begin{definition}
	A map $\xi:Q\to Q'$ between 2-colored quivers is a \keyword{weaving isomorphism} if it is an isomorphism of the underlying quivers and there exists a quiver isomorphism $\overline{\xi}:\overline{Q}\to\overline{Q'}$ such that the diagram commutes
    \begin{center}
		\begin{tikzcd}
            \overline{Q} \arrow[r,"\overline{\xi}"] \arrow[d] & \overline{Q'} \arrow[d]\\
            Q \arrow[r,"\xi"] & Q'
		\end{tikzcd}
	\end{center}
 where the vertical arrows correspond to the folding. If $Q$ and $Q'$ are weaving isomorphic we write $Q\weaviso{}Q'$.
\end{definition}
The name is motivated by the fact that a weaving isomorphism can move nodes between the two levels of the unfolding. Thus it preserves the folding groups, and consequently, while it does not preserve the coloring, it is not completely unaware of it, as it can only change colors in compatible ways.
\begin{example}
	Let $Q$ be a 2-colored quiver. The easiest example of a weaving isomorphism is provided by the \keyword{weaving at $k$} for a node $k$: Define the 2-colored quiver $w_k(Q)$ as having the same underlying quiver as $Q$ but switch the colors of all arrows adjacent to $k$. This is weaving-isomorphic to $Q$ because
	\begin{align*}
		w_k:\overline{Q}&\to\overline{w_k(Q)}\\
		i&\mapsto\begin{cases}
			i & i\notin \{k,k'\}\\
			k' & i=k\\
			k & i=k'
		\end{cases}\\
		(i\to j)&\mapsto\begin{cases}
			(i\to j) & i,j\notin\{k,k'\}\\
			(k'\to j) & i=k\\
			(k\to j) & i=k'\\
			(i\to k') & j=k\\
			(i\to k) & j=k'
		\end{cases}
	\end{align*}
	defines an isomorphism of the unfolded quivers which covers the identity. \Cref{fig:weavingIso} shows 2-colored quivers which are weaving-isomorphic via $w_3$ but not isomorphic as 2-colored quivers.
	
	\input{figureWeavingIso.tex}
\end{example}

\begin{proposition}\label{prop:weaving}
	Let $Q$ be a 2-colored quiver. Any weaving-isomorphic 2-colored quiver is obtained by applying an isomorphism composed with $w_{k_m}\circ\dotsb\circ w_{k_1}$ for some $k_1,\dots,k_m\in Q_0$.
\end{proposition}
\begin{proof}
    We need to consider quivers $Q'$ for which the unfolded quiver $\overline{Q'}$ is isomorphic to $\overline{Q}$ via an isomorphism which preserves the groups of the unfolding. Since an isomorphism can be seen as a simple relabeling of nodes, it is clear that we can at most relabel entire groups (isomorphism of 2-colored quiver) and exchange the two nodes which constitute a group ($w_k$ for some $k$).
\end{proof}

\subsection{Mixed 2-colored quivers}\label{sec:mixed2colQuivers}
As described before, much of this work is motivated by higher Teichm\"uller theory which we will return to in \Cref{sec:TM}. The quivers that arise there have a skew-symmetrizable exchange matrix $B$, which can be diagonalized by a matrix $D$ s.t. $d_{ii}\in\{1,2\}$ for all $i$. Observe that this restriction is stable under mutation due to \Cref{lemma:Dmut}. The nodes of these quivers are of two different types: Small nodes labeled by $\{i\,|\,d_{ii}=1\}$, and big nodes labeled by $\{i\,|\,d_{ii}=2\}$. We will now consider quivers which arise by `coloring' the arrows between big nodes.
\begin{definition}
    A \keyword{mixed 2-colored quiver} is a quiver with two types of nodes, \keyword{big} and \keyword{small}, and the structure of a 2-colored quiver on the full subquiver of big nodes. We define the \keyword{pruned quiver} $Q^\pruned$ by first removing all small nodes, then pruning the remaining quiver.
\end{definition}

\begin{definition}
    The mutation of these quivers is given by mutation of weighted quivers where big nodes have weight $1+\epsilon$ and small nodes have weight 1 and edges between small nodes have weight $\frac{1+\epsilon}{2}$. 
\end{definition}
Equivalently matrix mutation naturally extends to an exchange matrix $C(Q)$, defined as follows: For the full subquiver of big nodes the entries are the same as before, while an arrow from a big node to a small node as in \Cref{fig:mixedUnfold} corresponds to entries $c_{12}=1$ and $c_{21}=-(1+\varepsilon)$ with opposite signs for an arrow going from a small node to a big node. Finally, an arrow between small nodes, $i\rightarrow j$ yields entries $c_{ij}=\frac{1+\epsilon}{2}=-c_{ji}$. Observe that these entries are not integer in the split complex numbers in the sense we described before but become integers if we set $\varepsilon=1$.
\begin{remark}
    We observe that $\frac{1+\epsilon}{2}^2 = \frac{1+\epsilon}{2}$. Thus after any number of mutations the weight of arrows between small nodes is a multiple of $\frac{1+\epsilon}{2}$ and this class of quivers is preserved under mutation.
\end{remark}

 A mixed 2-colored quiver unfolds in such a way that the small nodes sit between the two layers and connect to both of them:
\begin{definition}
    Let $Q$ be a mixed 2-colored quiver. The \keyword{unfolded quiver} $\overline{Q}$ is obtained by unfolding the subquiver of big nodes as in \Cref{def:Unfolding2ColQuivers}, and adding a copy of the subquiver of small nodes which is connected as shown in \Cref{fig:mixedUnfold}.

    \begin{figure}[ht]
    	\begin{center}
	       	\begin{tikzpicture}[]
			     \begin{scope}
				    \node[mutableBig]	(1) at (0,0)			[label = above: $1$]		{};
				    \node[mutable]	    (2) at (\distL,0)	    [label = above: $2$]		{};
				    \draw[arrow]		(1) to (2);
			     \end{scope}
			     \begin{scope}[shift={(2*\distL,0)}]
				    \node[vertex]	(1) at (0,-0.2)			[label = $\mathrm{=}$]		{};
			     \end{scope}
			     \begin{scope}[shift={(3*\distL,0)}]
				    \node[mutable]		(1a) at (0,0.3*\distL)		[label = above: $1$]	{};
				    \node[mutable]		(1b) at (0,-0.3*\distL)		[label = below: $1'$]	{};
				    \node[mutable]		(2) at (\distL,0)	        [label = above: $2$]	{};
				    \draw[arrow]		(1a) to (2);
				    \draw[arrow]		(1b) to (2);
			     \end{scope}
		    \end{tikzpicture}
	   \end{center}
	   \caption{Unfolding arrows in mixed quivers.}
        \label{fig:mixedUnfold}
    \end{figure}
\end{definition}

\begin{proposition}
    The quiver $\overline{\mu(Q)}$ is isomorphic to the quiver obtained by group mutation on $\overline{Q}$ under an isomorphism preserving the layers. Explicitly within each group the plain and primed nodes are mapped onto each other.
\end{proposition}

\begin{proof}
    We have already seen that mutations of the subquiver defined by the large nodes agree with mutations of the unfolded quiver. We see in the unfolded quiver mutation at a small group adds a pair of crossing and parallel arrows between big groups. This corresponds exactly to adding $1+\epsilon$ arrows between big nodes, which is the defined mutation of mixed two colored quivers.
\end{proof}

\begin{proposition}
    Weaving isomorphism classes of mixed 2-colored quivers are equivalent to the data of a regular quiver folded by an order two automorphism. 
\end{proposition}

\begin{proof}
    The quiver $\overline{Q}$ has an automorphism $\epsilon$ given by swapping the elements of all groups $\{i,i'\}$ coming from big nodes. It is clear that a weaving of $Q$  at node $i$ produces the same unfolding as this corresponds to switching the nodes in group $i$ of $\overline{Q}$. Weaving at every node is the same as applying $\epsilon$ to $\overline{Q}$. 
\end{proof}

\begin{example}
	An example of a mutation on a small mixed 2-colored quiver is shown in \Cref{fig:mixedColMutEx}. Observe the coloring rules and the introduction of pairs of arrows due to a mutation at a small node.
	\input{figureMarkedColoredMutation.tex}
\end{example}

\section{Ordered and ST-compatible Quivers}\label{sec:ncQuivers}
Our ultimate goal is to construct a new class of noncommutative cluster algebras. In order to do so, we will capture the underlying combinatorial data in two different ways: Via quivers, which are the focus of this section, and via polygons, which we will see in the next section. We will show that these descriptions are equivalent.

\subsection{Reminder: Cluster algebras}
A cluster algebra is typically defined by starting with a \keyword{seed} $s=(Q,x)$, which consists of a quiver $Q$ and a set of variables $x=(x_1,\dots,x_n)\in\mathbb{Q}(x_1,\dots,x_n)$, one for each node of the quiver. Let $k\in Q_0$. The quiver can be mutated as described before to obtain $\mu_k(Q)$, and the variables get updated at the same time according to the rule
\begin{align}\label{eq:classicalVarMut}
    \mu_k(x_i)=
    \begin{cases}
        x_i & i\neq k\\
        \frac{1}{x_k}\left(\prod_{k\leftarrow j}x_j+\prod_{k\rightarrow j}x_j\right) & i=k
    \end{cases}\,.
\end{align}
Now, the \keyword{cluster variables} are the variables obtained by performing all possible mutations. The \keyword{cluster algebra} $\mathcal{A}_s$ is defined to be the subalgebra of $\mathbb{Q}(x_1,\dots,x_n)$ generated by the cluster variables.

\subsection{Ordered quivers}\label{subsec:quivNCCAs}
One problem we immediately encounter in the construction of noncommutative cluster algebras is to find the correct analogue of \Cref{eq:classicalVarMut} in this setting. If we simply assume that the variables in the initial seed lie in some skew field, we need to decide on an order for the products. Moreover, we need to have a way to choose a consistent order for the variables in every possible mutation. To keep track of the order, we add some additional decoration to 2-colored quivers.

\begin{definition}
    An \keyword{ordered (2-colored) quiver} consists of a 2-colored quiver $Q$ together with a splitting of the incoming and outgoing arrows in $Q^\pruned$ at every node $k\in Q_0$ into two ordered tuples, denoted $\In_f(k)$ and $\Out_f(k)$ for $f\in\{0,1\}$. The remaining pairs of parallel and crossing arrows are called \keyword{commutative} and live in $\In^c(k)$ and $\Out^c(k)$ respectively. We will sometimes call the arrows in $Q^\pruned$ \keyword{noncommutative arrows}.
\end{definition}

Note that the tuples $\In_f(k)$ and $\Out_f(k)$ are allowed to be empty. 

\begin{convention}
    Each node now has two sides gives by the index $f=0$ or $f = 1$. We visualize the two sides of the nodes by the fill as in \Cref{fig:decNode}. Thus we use the subindices  $\circ,\bullet$ to represent the two sides. In computation we consider $\circ=0$ and $\bullet=1$ as elements of $\Z/2\Z$. To visualize commutative arrows, we use a green arrow with a double tip for a pair of parallel and crossing arrows, and indicate its weight, i.e. the number of pairs, by a number along the arrow.
\end{convention}

\begin{figure}[ht]
	\begin{center}
		\begin{tikzpicture}
			\node[vertex]	(A) at (-1.5,1)		[label = left: $\In_\bullet(k)$]	{};
			\node[vertex]	(B) at (-1.5,-0.8)	[]{};
			\node[vertex]	(C) at (-1.1,-1.1)	[]{};
			\node[vertex]	(AB) at ($0.5*(B)+0.5*(C)+(0,-0.2)$) [label = left: $\Out_\bullet(k)$]	{};
			\node[vertex]	(D) at (1.5,1)		[label = right: $\Out_\circ(k)$]{};
			\node[vertex]	(F) at (1.5,-1)		[label = right: $\In_\circ(k)$]	{};

			\pic[name=E] at (0,0) {NCnode};
			\node[vertex] (X) at (0,0.8)		[label = above: $k$]		{};
            \node[vertex]	(OutC) at (3.5,0.0)		[label = right: $\Out^c(k)$]	{};

			\draw[barrow] (A.-20) to [out=-20,in=140] (E-circle.140);
			\draw[xarrow] (E-circle.220) to [out=220,in=10] (B.10);
			\draw[barrow] (E-circle.240) to [out=250,in=60] (C.60);
			\draw[xarrow] (E-circle.45) to [out=45,in=-180] (D.-180);
			\draw[barrow] (F.180) to [out=180,in=-45] (E-circle.-45);
            \draw[carrow] (0,0) to[out=15,in=195] (OutC);
            
		\end{tikzpicture}
	\end{center}
 	\caption{Decoration of big nodes in an ordered quiver with partition of the arrows.}
	\label{fig:decNode}
\end{figure}

Visually the ordering on $\In_\bullet, \Out_\bullet,\In_\circ,\Out_\circ$ is given by the counterclockwise order on the boundary of the node. Note that each node is required to have $\In_f$ appear before $\Out_f$ in this order.

\subsection{Isomorphisms of ordered quivers}
Every ordered quiver consists of the data of a 2-colored quiver, the In and Out tuples and the filling of the nodes. We will see that the most important aspect is the order of the tuples, and accordingly we have different types of isomorphisms, depending on the amount of the data which is preserved:

\begin{definition}
	Let $Q,Q'$ be ordered quivers, $\xi:Q\to Q'$ a map of the underlying 2-colored quivers. Then
    \begin{enumerate}
        \item $\xi$ is an \keyword{isomorphism} if $\xi$ is an isomorphism of 2-colored quivers and for all $k\in Q_0, f\in\{\circ,\bullet\}$ we have $\xi(\In_f(k))=\In_f(\xi(k))$, $\xi(\Out_f(k))=\Out_f(\xi(k))$ as ordered tuples.
        \item $\xi$ is a \keyword{switching isomorphism} if $\xi$ is an isomorphism of 2-colored quivers an for all $k\in Q_0, f\in\{\circ,\bullet\}$ there exists a $g\in\{\circ,\bullet\}$ s.t. $\xi(\In_f(k))=\In_g(\xi(k)),\xi(\Out_f(k))=\Out_g(\xi(k))$ as ordered tuples.
        \item  $\xi$ is a \keyword{weaving isomorphism} if $\xi$ is a weaving-isomorphism of the underlying 2-colored quiver and for all $k\in Q_0, f\in\{\circ,\bullet\}$ there exists a $g\in\{\circ,\bullet\}$ s.t. $\xi(\In_f(k))=\In_g(\xi(k))$, $\xi(\Out_f(k))=\Out_g(\xi(k))$ as ordered tuples.
    \end{enumerate}
\end{definition}

The amount of data that is preserved by the different types of isomorphisms defined above, is summarized in \Cref{tab:orderIsom}.
\begin{table}[h]
    \centering
    \begin{tabular}{c|ccc}
         & isomorphism & switching isom. & weaving isom.\\
         \hline
        order of arrows  & \ding{51} & \ding{51} & \ding{51} \\
        color of arrows & \ding{51} & \ding{51} & \ding{55} \\
        type of groups & \ding{51} & \ding{55} & \ding{55} \\
    \end{tabular}
    \caption{Data preserved by different types of isomorphism for ordered quivers.}
    \label{tab:orderIsom}
\end{table}

\begin{example}
	An elementary example showing that the notion of switching isomorphism for ordered quivers is slightly richer than just isomorphisms of the underlying quivers, is the following \keyword{switch at $k$}: Let $Q$ be an ordered quiver, $k$ a node. Denote by $\sigma_k(Q)$ the ordered quiver which has the same underlying 2-colored quiver but with the fill of the node $k$ reversed. Phrased differently, this means
	\begin{align*}
		\In_f(k,\sigma_k(Q))&=\In_{f+1}(k,Q)\\
		\Out_f(k,\sigma_k(Q))&=\Out_{f+1}(k,Q)
	\end{align*}
	for $f\in\Z/2\Z$.
\end{example}

\begin{example}
	The \keyword{weaving at $k$} defined for 2-colored quivers extends to the noncommutative case: Simply recolor all arrows incident to a node $k$ in $Q$ without changing any of the other decoration. We will still denote the resulting ordered quiver by $w_k(Q)$.
\end{example}

\subsection{Quiver algebras and angles}
With the quivers in place, we now turn our attention to the variables we use in the construction of our cluster algebras. Here our approach draws inspiration from Berenstein-Retakh's noncommutative surfaces again: The variables will not lie in any skew field, instead there will be some relations among them which should be thought of as a generalization of the triangle relations. We need to make some algebraic constructions to prepare this.

Let $Q$ be an ordered quiver.
\begin{definition}
	The \keyword{quiver algebra} $\quiverAlgebra_Q$ is the $\Q$-algebra generated by
	\begin{equation*}
		\left\{X_i^{(f)},Y_i^{(f)} \Big|\,i\in Q_0, f\in \mathbb{Z}/2\mathbb{Z}\right\}
	\end{equation*}
	subject to the relations
	\begin{align*}
        & X_i^{(f)}Y_i^{(f+1)} = Y_i^{(f)}X_i^{(f+1)} \\
		& X_i^{(f)}Y_i^{(f+1)} \quad \text{is central}
	\end{align*}
	for any $i,j\in Q_0, f \in\Z/2\Z$, localized by the center.
\end{definition}

\begin{remark}
	Observe that the quiver algebra does not actually depend on the quiver $Q$ but only on its set of nodes $Q_0$.
\end{remark}

One should think that we associate a 4-tuple of variables  to each node such that the (noncommutative) variables
pair in different ways to give central elements in $\quiverAlgebra_Q$. To lighten the notation, we will introduce a pair $\sigma,\tau$ of anti-automorphisms on $\quiverAlgebra_Q$. On the generators these are given by
\begin{align*}
	&\begin{matrix}
		\sigma: & X_j^{(f)} & \mapsto & X_j^{(f+1)}
	\end{matrix}\\
	&\begin{matrix}
		\tau: & X_j^{(f)} & \mapsto & Y_j^{(f+1)}
	\end{matrix}\,.
\end{align*}

\begin{convention}
    To reduce overhead the symbol $\sau$ is used to denote the composition $\sigma\circ\tau$ throughout the paper. 
\end{convention}

\begin{lemma}
	$\sigma$ and $\tau$ commute and are trivial on $\mathrm{Z}(\quiverAlgebra_Q)$.
\end{lemma}

\begin{proof}
    $\sigma$ and $\tau$ commute by definition. To see that they are trivial on the center of the quiver algebra, we first observe that it is generated by expressions of the form $X_i^{(f)}Y_i^{(f+1)}$ (and their inverses). These expressions are fixed by $\sigma$ due to the symmetry relation and by $\tau$ since
    \begin{equation*}
        X_i^{(f)}X_i^{(f)}Y_i^{(f+1)} = X_i^{(f)}Y_i^{(f+1)}X_i^{(f)}
    \end{equation*}
    which implies due to invertiblility of $X_i^{(f)}$
    \begin{equation*}
       X_i^{(f)}Y_i^{(f+1)} = Y_i^{(f+1)}X_i^{(f)}
    \end{equation*}
\end{proof}

\begin{definition}
	The \keyword{norm} of an element $X\in\quiverAlgebra_Q$ is
	\begin{equation*}
		N(X):=X\tau(X)\,.
	\end{equation*}
\end{definition}

\begin{definition}   
	We say $X \in \quiverAlgebra_Q$ has \keyword{central norm} if $N(X)\in \mathrm{Z}(\quiverAlgebra_Q)$ and is fixed by $\sigma$.
\end{definition}

Clearly all the variables $X_i$ have central norm.

\begin{lemma}\label{lemma:BQinversion}
	For $X\in\quiverAlgebra_Q$ with central norm, we have
	\begin{equation*}
		X^{-1}=\tau(X)N(X)^{-1}\in\quiverAlgebra_Q\,.
	\end{equation*}
\end{lemma}
Observe that this inverse is contained in $\quiverAlgebra_Q$ since we have localized by the center.

\medskip

We can now generalize the triangle invariants. Again, let $Q$ be an ordered quiver, $k$ some node. We represent a general arrow from to $i$ to $j$ as
\begin{equation*}
	i \xrightarrow{f~t~g} j
\end{equation*}
 where $f,g\in\Z/2\Z$ are the fill of the side of the node at $i,j$ respectively. If the arrow is parallel then $t = 0$, and for a crossing arrow $t=1$.
\begin{definition}\label{def:angles}
	The \keyword{angles} at $k$ are the following elements in $\quiverAlgebra_Q$:
	\begin{align}\label{eq:defAngles}
        \Delta_{g}(k) =\hspace{-3mm}\prod_{\substack{\alpha\colon\: i\xrightarrow{f~t~g} k\\ \alpha \in \In_{g}(k)}}\hspace{-3mm}\sigma^f\sau^{t+1}(X_i)\cdot \sigma^g(X_k)\cdot\hspace{-3mm}\prod_{\substack{\alpha\colon\: k\xrightarrow{g~t~f} j\\ \alpha\in\Out_{g}(k)}}\hspace{-3mm}\sigma^f\sau^{t+1}(X_j)
	\end{align}
 where the products are understood to be ordered.
\end{definition}
The angles are easy to read from the ordered quiver. There are two angles per node, associated to the filled and empty side of the node. For each side, read the neighboring nodes following the cyclic order starting from the incoming arrows. For each arrow there are four possible choices of associated variable, which is determined by the color of the arrow and fill of the corresponding node. The fill of the neighbor determines the power of $\sigma$ while the power of of $\sau$ is determined by the color.

\begin{remark}
    The convention of reading the incoming arrows first is an arbitrary choice. By making this choice we fix a cyclic order on the set of variables on each side of the node. We will later impose an ``angle relation'' which relates this choice with the correct angle given by reading the nodes in the reverse order starting with the outgoing nodes.
\end{remark}

\begin{example}
    An example of an ordered quiver is shown in \Cref{fig:angleEx}.
    \input{figureAngleEx}

    The angles at node $6$ are given by
    \begin{align*}
        \Delta_\circ(6)&=\sigma\tau(X_5)\cdot X_6\cdot \sigma\tau(X_4)\\
        \Delta_\bullet(6)&=\tau(X_1)\cdot \sigma(X_6)\cdot\sigma(X_2)\tau(X_3)\,.
    \end{align*}
    Observe that we can turn the underlying quiver into one associated to a surface by removing node $2$ and the commutative arrow from 4 to 1. In that case the angles become
    \begin{alignat*}{2}
        \Delta_\circ(6)&=\sigma\tau(X_5)\cdot X_6\cdot \sigma\tau(X_4)&&={N(X_4X_5)}\sigma(X_5)^{-1}\cdot X_6\cdot \sigma(X_4)^{-1}\\
        \Delta_\bullet(6)&=\tau(X_1)\cdot \sigma(X_6)\cdot\tau(X_3)&&={N(X_1X_3)}X_1^{-1}\cdot\sigma(X_6)\cdot X_3^{-1}\,.
    \end{alignat*}
    Up to central elements, i.e. norms, these are the same as the angles defined by Berenstein and Retakh \cite{berenstein2018noncommutative}, see also \Cref{sec:warmup}.
\end{example}

From the construction of the angles we clearly have:
\begin{proposition}
	The ordered quiver $Q$ is determined by its set of vertices and the angles.
\end{proposition}

The appropriate generalization of triangle relations is now given by the \keyword{angle relations}
\begin{equation*}
    \Delta_f(k)=\sigma(\Delta_f(k))
\end{equation*}
where $\{\Delta_f(k)\}\subset\quiverAlgebra_Q$ is the set of angles associated to $Q$. We would like to simply quotient $\quiverAlgebra_Q$ by the ideal generated by these relations and define mutation in the resulting algebra. Unfortunately, this is not possible since the resulting new cluster variables would not have central norm in this algebra, which is a crucial feature of our construction. Thus, we must first restrict to a class of quivers where it is possible to guarantee that the new variables have central norm.

Before we do that, let us note that weaving isomorphic quivers are equivalent for the above constructions.
\begin{proposition}\label{prop:weavingIsoQuiverAlgebras}
    Suppose $Q\weaviso{}Q'$. Then $\quiverAlgebra_Q\sim\quiverAlgebra_{Q'}$ via an isomorphism that sends angles to angles.
\end{proposition}

\begin{proof}
    A full isomorphism of the ordered quivers clearly preserves angles, so it suffices to consider switching and weaving at a node $k$. For the switch at $k$, we get an isomorphism $\sigma_k$ of the corresponding quiver algebras induced by
    \begin{align*}
        X_i\mapsto\begin{cases}
            X_i' & i\neq k\\
            \sigma(X_k') & i=k
        \end{cases}
    \end{align*}
    where we extend this to a $\sigma$- and $\tau$-equivariant map. \Cref{eq:defAngles} clearly shows that this sends angles to angles, i.e.
    \begin{align*}
        \sigma_k(\Delta_f(i))=\begin{cases}
            \Delta'_f(i) & i\neq k\\
            \Delta'_{f+1}(k) & i=k
        \end{cases}
    \end{align*}
    where $\Delta_f(i),\Delta'_f(i)$ denote the angles for $Q$ and $Q'$, respectively.

    For weaving at $k$ the argument is completely analogous, we define the isomorphism $w_k:\quiverAlgebra_Q\to\quiverAlgebra_{Q'}$ by $w_k(X_k)=\sigma\tau(X_k')$ and $w_k(X_i)=X_i'$ for all $i\neq k$.
\end{proof}

\subsection{ST-compatible quivers}\label{sec:ST}
Recall that in the case of noncommutative surfaces the three triangle relations associated to each triangle are equivalent, in the sense that imposing one of them in the algebra automatically implies the other ones. Our next definition enforces this behavior on an ordered quiver. We will see this innocuous condition imposes strong structural constraints on the ordered quiver.

\begin{definition}
	An ordered quiver $Q$ is \keyword{source-target-compatible} or \keyword{ST-compatible} if for every arrow $i \xrightarrow{f~t~g} j$ in the pruned quiver $Q^\pruned$ the angle relations at the source and target are equivalent, i.e.
    \begin{align*}
        \Delta_f(i)=\sigma(\Delta_f(i))\Leftrightarrow \Delta_g(j)=\sigma(\Delta_g(j))\,.
    \end{align*}
\end{definition}

\begin{example}
    The quiver in \Cref{fig:angleEx} is ST-compatible.
\end{example}

\begin{lemma}\label{lemma:elemEquiv}
	Two angle relations for an ordered quiver $Q$ are equivalent if and only if the corresponding angles can be related by a sequence of the following operations:
	\begin{enumerate}
		\item apply $\sigma\tau$
		\item $AB\leftrightarrow\sigma\tau(B)A$ for monomials $A,B\in\quiverAlgebra_Q$ (cyclic rotation).
	\end{enumerate}
\end{lemma}
\begin{proof}
    Suppose we have two angles $\Delta$ and $\Delta'$ such that their angle relations are equivalent. Observe that $\Delta,\Delta'$ are monomials in the generators of $\quiverAlgebra_Q$ which contain no inverses of the generators. Now $\Delta=\sigma(\Delta)$ implies $\Delta'=\sigma\left(\Delta'\right)$ if and only if
    \begin{equation*}
        \Delta'=c\cdot\sigma^\delta\tau^\varepsilon(A\Delta\sigma(A))
    \end{equation*}
    for a monomial $A$ with no inverses, a central element $c$ and $\delta,\varepsilon\in\Z/2\Z$. Since $\Delta'$ is an angle it cannot contain symmetric terms as these would correspond to 2-cycles in $Q$. Therefore we find $\Delta=\tau(A)B$ or $\Delta=B\sau(A)$ for an inverse-free monomial $B$. In the first case we consequently find
    \begin{equation}\label{eq:angleEquivIntermediate}
        \Delta'=c\cdot N(A)\sigma^\delta\tau^\varepsilon(B\sigma(A))\,.
    \end{equation}
    As this last part is an inverse-free monomial and the angle only depends on the pruned quiver all central terms should equal 1. Suppose that $\Delta$ is an angle at a node $i$, $\Delta'$ at $k$. Since variables appearing in both are associated to the same sets of nodes by \Cref{eq:angleEquivIntermediate}, $i$ appears in $\Delta'$ and vice versa. As the angles are such that In variables appear before Out variables, $i$ and $k$ need to appear in the same order in both $\Delta$ and $\Delta'$, forcing $\delta=\varepsilon$ in \Cref{eq:angleEquivIntermediate}. In conclusion, we have $\Delta'=\sau^\delta(B\sigma(A))$. The argument is completely analogous for the second case, i.e. $\Delta=B\sau(A)$.
\end{proof}

Being ST-compatible severely restricts the possible shape of the ordered quiver as we will show next. Let us first decompose ST-compatible quivers into building blocks:

\begin{definition}
    Let $Q$ be an ST-compatible ordered quiver. An \keyword{ST-tile} in $Q$ is a subquiver consisting of all arrows and nodes in $Q^\pruned$ which contribute to a set of equivalent angle relations.
\end{definition}

In other words, these subquivers are obtained by starting with a node and adding all arrows that lie on one side of it. We then add all nodes which are connected to our original node via these arrows. Next we add all the arrows which are on the same side of the new nodes as the original arrows. Now repeat this process until it terminates.

\begin{example}
    Consider again \Cref{fig:angleEx}. This ST-compatible quiver decomposes into two ST-tiles: The full subquiver with nodes $1,2,3$ and $6$, and the full subquiver with nodes $4,5,6$.
\end{example}

Note that an ST-tile is not necessarily the full subquiver formed by the set of nodes associated to an angle, since there might be arrows between these nodes, which are used in the second angle. 

\begin{lemma}\label{lem:st-tile_decomposition}
    Let $Q$ be an ST-compatible quiver. Then $Q^\pruned$ decomposes into finitely many overlapping ST-tiles, with each node contained in at most two tiles. 
\end{lemma}
\begin{proof}
    Every node generates at most two ST-tiles corresponding to the two sides of the node. Call the variable at this node $X$.  It is not possible for the angles at both sides of the node $X$ to generate the same ST-tile, as the expression for the angle at a given node can only contain $X$ once, but the angle written at other nodes in this tile would contain $X$ twice. Therefore, the only way a node can appear in only one ST-tile is for it to have no arrows on one side of it.
\end{proof}

\begin{definition}
	Let  $Q$ be a 2-colored quiver, $k$ a node. $Q$ is a \keyword{3-cycle} if it is weaving-isomorphic to the quiver shown in \Cref{fig:3cycles}(a). $Q$ is a \keyword{3-pseudocycle} based at $k$ if it is weaving-isomorphic to the quiver shown in \Cref{fig:3pseudocycles}(a), where the corresponding isomorphism of the unfolded quivers sends $k$ to $1$.
	\input{figurePseudoCycles.tex}
\end{definition}

Recall that a \keyword{tournament} is a directed complete graph, and that this is called \keyword{transitive} if it is totally ordered. It is called \keyword{locally transitive} if for every node the full subquivers of all nodes connected by an incoming, resp. outgoing, arrow are totally ordered.

With this preparation, we can now completely describe ST-tiles of an ST-compatible quiver:
\begin{theorem}\label{thm:GSTRestrictions}
    For an ST-compatible quiver
    \begin{enumerate}
    	\item Every ST-tile is a locally transitive tournament. 
    	\item Every ST-tile is colored in such a way that every set of three nodes either forms a 3-cycle or a 3-pseudocycle.
    	\item The equivalent angles have the same cyclic order for their variables.
    \end{enumerate}
\end{theorem}

\begin{remark}
	These restrictions affect every level of the ordered quiver: (1) is a statement about the underlying quiver, (2) about its coloring, and (3) about the order of the arrows. (3) relies on our convention to write the angles such that the incoming variables appear first. Applying $\sigma$ clearly reverses the cyclic order. 
\end{remark}

\begin{proof}
    Let $Q$ be an ordered quiver which is an ST-tile.

    \textit{Statement (3):} This is an immediate consequence of \Cref{lemma:elemEquiv}.
    
    \textit{Statement (1):} The first part is due to the fact that equivalence of the angle relations at all the nodes makes it necessary that the same variables appear with the same multiplicity. Thus the group of nodes needs to form a tournament.

    In order to prove that this is locally transitive, we pick a node, say $k$, and prove that all the nodes which are connected to this via incoming arrows on one side, $\In_f(k)$, form a transitive tournament (the same argument then applies for all the nodes connected by outgoing arrows). The angle relation associated to $k$ is of the form
    \begin{equation*}
        \Delta_f(k)=\sigma^{s_1}\tau^{t_1}(X_{i_1})\cdots \sigma^{s_a}\tau^{t_a}(X_{i_a})\,\sigma^f(X_k)\,A\,,
    \end{equation*}
    for some $s_1,t_1,\dots,s_a,t_a\in\Z/2\Z$ and $A\in\quiverAlgebra_Q$ a monomial. Thus, $\In_f(k)=(X_{i_1},\dots, X_{i_a})$. Due to statement (3), the cyclic order of the angle at $X_{i_1}$ is the same, so that $X_{i_1}$ is connected to $X_{i_l}$ for $2\leq l\leq a$ by an outgoing arrow since it is connected to $X_k$ by an outgoing arrow. This means that $X_{i_1}$ is a source of $\In_f(k)$. Proceeding inductively, we see that $\In_f(k)$ is totally ordered, i.e. a transitive tournament.\medskip

    \textit{Statement (2):} This is another consequence of how we obtain equivalent angle relations: Suppose $i,j,k$ are three nodes, which are part of the same group, and that $i$ and $k$ are connected to $j$ as
    \begin{equation*}
        i\xrightarrow{f~t~g} j\xrightarrow{g~t'~f'} k\,.
    \end{equation*}
    Then $j$ has angle
    \begin{equation*}
        \Delta_g(j)=\sigma^f(\sigma\tau)^{t+1}(X_i)\,\sigma^g(X_j)\,\sigma^{f'}(\sigma\tau)^{t'+1}(X_k)\,,
    \end{equation*}
    where we have set all other variables to $1$. Then the possible angles at $i$ are
    \begin{align*}
        &\Delta_f(i)=\sigma^f(X_i)\,\sigma^g(\sigma\tau)^{t+1}(X_j)\,\sigma^{f'}(\sigma\tau)^{t'+t}(X_k)\,\mathrm{, or}\\
        &\Delta_f(i)=\sigma^{f'}(\sigma\tau)^{t'+t+1}(X_k)\sigma^f(X_i)\,\sigma^g(\sigma\tau)^{t+1}(X_j)\,\,.
    \end{align*}
    Either choice determines the 2-colored noncommutative subquiver consisting of these three nodes. The first possibility corresponds to a 3-pseudocycle, the second corresponds to a 3-cycle.\medskip
\end{proof}

\begin{remark}
    We have seen that an ST-tile consists of a locally transitive tournament, an admissible coloring (i.e. every group of 3-nodes forms a cycle or pseudocycle), a compatible cyclic order for the nodes, and the filling of the nodes.
\end{remark}

Let us focus on the first two parts of the data and make some observations which will be useful later:
\begin{lemma}
    Any coloring of the arrows incident at any given node in a locally transitive tournament uniquely determines an admissible coloring. In particular, any locally transitive tournament admits an admissible coloring.
\end{lemma}

\begin{proof}
    Suppose Q is a locally transitive tournament, $i\in Q_0$. Fix any coloring of the arrows connected to $i$. Now for $j\to k\in Q_1$ with $j,k\neq i$, the color is uniquely determined by the condition that $i,j,k$ form a cycle or pseudocycle. We need to show that the unique coloring obtained in this way is admissible, i.e. that any group $\{j,k,l\}\notni i$ of nodes forms a cycle or pseudocycle. Since $Q$ is a locally transitive tournament, we only need to consider the situations shown in \Cref{fig:LTTcoloring}, which all have the desired property.

    \input{figureLTTcolor}
\end{proof}

\begin{corollary}\label{cor:LTTcoloringWeavingUnique}
    Any two admissible colorings of the same locally transitive tournament are weaving isomorphic.
\end{corollary}

\begin{proof}
    Since any admissible coloring is uniquely determined by the color of the arrows connected to a fixed node $i\in Q_0$, two admissible colorings can be related by weaving at the nodes in $Q_0\setminus\{i\}$, which allows us to independently change the colors of the arrows connected to $i$.
\end{proof}

\begin{proposition}
    Let $Q$ be an admissibly colored locally transitive tournament with an arbitrary choice for the fill of the nodes. Then $Q$ is an ST-tile with ordering induced by the tournament.
\end{proposition}
\begin{proof}
    By \cite{brouwer1980enumeration} every node in a locally transitive tournament induces the same cyclic order given by the total order on incoming arrows, followed by $k$, and then the total order on outgoing arrows (see also the proof of \Cref{thm:decPolyElemSTCorrespondence}). It suffices to show that the two angles associated to nodes $i \xrightarrow{f~t_1~g} j$ in the cyclic order are equivalent.
    We split the remaining elements of the tournament into 4 ordered tuples $IO,II,OI,OO$, where $IO=\In_f(i)\cap\Out_g(j)$ etc. Then the angles at $i$ and $j$ are of the form  
    \begin{align*}
        \sau^{t_1+1}\Delta_f(i) =& \Phi_1(IO)\Phi_2(II)\sau^{t_1+1}\sigma^f(X_i)\Phi_3(OI)\sigma^g(X_j)\Phi_4(OO) \\
        \Delta_g(j) =& \Psi_2(II)\sau^{t_1+1}\sigma^f(X_i)\Psi_3(OI)\sigma^g(X_j)\Psi_4(OO)\Psi_1(IO)
    \end{align*}
    where we use $\Phi_1(IO),\dots$ as a shorthand to denote an ordered product over the variables $X_k$ for $k\in IO$ with some combination of $\sigma$ and $\tau$ applied to every individual factor.
    
    Let $k\in Q_0\setminus\{i,j\}$. Let $t_2$ be the color of the arrow between $i$ and $k$ and $t_3$ be the color of the arrow between $j$ and $k$. Let $h$ be the fill at $k$. From the angle formula $\Phi_r = \sigma^h\sau^{t_1+t_2}$ and $\Psi_r = \sigma^h\sau^{t_3+1}$. Every node $k$ in $II,OI,OO$ is in a pseudocycle with $i,j$. Thus for the coloring to be admissible the colors satisfy $t_1+t_2+t_3 = 1$ and so $\Phi_r = \Psi_r$ for $r=2,3,4$.
    Every node $k$ in $IO$ is in a cycle with $i$ and $j$. Thus for the coloring to be admissible $t_1+t_2+t_3=0$ and $\Phi_1 = \sau \Psi_1$. Therefore the two angles are related by the operations in \Cref{lemma:elemEquiv} and are equivalent.

\end{proof}
\begin{corollary}\label{cor:TilingsImplySTcompatible}
    If $Q$ is a quiver such that $Q^\pruned$ can be divided into a collection of admissibly colored locally transitive tournaments then $Q$ is ST-compatible.
\end{corollary}

\subsection{Mutation of ST-compatible quivers}\label{sec:MutationOfSCQuivers}

Let $Q$ be an ST-compatible ordered quiver and $k$ a node of $Q$. By \Cref{thm:GSTRestrictions} $Q^\pruned$ can be divided into ST-tiles which are locally transitive tournaments and there are no arrows in the pruned quiver between ST-tiles. To define mutation we make the following technical \textbf{assumption}: the node $k$ belongs to two distinct ST-tiles.
\begin{definition}\label{def:quiverMut}
	In this setting, \keyword{(quiver) mutation at $k$} gives the ordered quiver $\mu_k(Q)$ obtained by the following procedure:
	\begin{enumerate}
		\item Perform mutation at $k$ on the underlying 2-colored quiver.
        \item Reestablish the ST-tiles. The only tiles which change are the two tiles containing $k$, which become $\{\In_\bullet(k,Q),k,\Out_\circ(k,Q)\}$ and $\{\In_\circ(k,Q),k,\Out_\bullet(k,Q)\}$.
		\item Establish the cyclic order in each new ST-tile. The order of the (reversed) arrows connected to $k$ is given by
			\begin{align*}
				\In_f(k,\mu_k(Q))&=\Out_{f+1}(k,Q)\,,\\
				\Out_f(k,\mu_k(Q))&=\In_f(k,Q)\,,
			\end{align*}
			for $f\in\Z/2\Z$. This should establish the cyclic order in the new ST-tiles.\\
           For any node not incident to $k$ in the pruned quiver the order is preserved by mutation. 
	\end{enumerate}
     A mutation is \keyword{admissible} if the new ST-tiles satisfy the conclusion of \Cref{thm:GSTRestrictions}, i.e. form a well colored locally transitive tournament and the order of the tournament matches the order defined at $k$.
     We call the ordered quiver $Q$ \keyword{fruitful} if all mutation sequences are admissible.
\end{definition}

As a first check that admissible mutations exist we observe the following:
\begin{lemma}
    The new ST-tiles after mutation at $k$ are disconnected in the pruned quiver.
\end{lemma}
\begin{proof}
    This follows from \Cref{thm:GSTRestrictions} restriction on coloring of subsets of three nodes. Each arrow between $\In^\pruned_{f}(k)$ and $\Out^\pruned_{f}(k)$ is either canceled or become a pair of parallel and crossing arrows by this assumption.
\end{proof}

\begin{remark}
    Admissible mutations result in new ST-compatible quivers by \Cref{cor:TilingsImplySTcompatible}.
\end{remark}

Recall that mutation of (2-colored) quivers is of order 2. For noncommutative quivers we have.
\begin{proposition}\label{prop:MutationOrder4}
	\begin{equation*}
		\mu_k^2(Q)=\sigma_k(Q)\,.
	\end{equation*}
\end{proposition}
\begin{proof}
	After two mutations, the underlying 2-colored quiver comes back to the original quiver but we have
	\begin{align*}
	    \In_f(k,\mu_k^2(Q)) =& \Out_{f+1}(k,\mu_k(Q)) = \In_{f+1}(k,Q)\\
        \Out_f(k,\mu_k^2(Q) =& \In_f(k,\mu_k(Q)) = \Out_{f+1}(k,Q)     
	\end{align*}
    This is exactly a switch at $k$ as claimed.
\end{proof}

\begin{example}\label{ex:NCmutation}
	An example of a mutation of an ST-compatible quiver with a single mutable node is shown in \Cref{fig:mutationEx}.
	
	\input{figureMutationExample.tex}
\end{example}

\begin{remark}
    When mutation at $k$ is not admissible, the mutation rule given does not naturally give an order on the nodes incident to $k$. One can define an explicit mutation which maintains an order at each node throughout the process. Under this rule, a mutation that is not admissible will result in a quiver with a two-cycle that cannot be canceled due to the order.
\end{remark}

\subsection{ST-compatible mixed quivers}
Recall that in \Cref{sec:mixed2colQuivers} we described a generalization of 2-colored quivers which includes small nodes.  We will now explain how to include an additional layer of small nodes. The general philosophy is that the big nodes contribute noncommutatively while the small nodes only contribute commutatively.

The additional structure for ordered quivers only lives on the pruned quiver:
\begin{definition}
    An \keyword{ordered mixed quiver} is a mixed 2-colored quiver $Q$ with the structure of an ordered quiver on $Q^\pruned$.
\end{definition}

The quiver algebra $\quiverAlgebra_Q$ for such a quiver is constructed by adding one central variable for every small node to $\quiverAlgebra_{Q^\pruned}$. We associate angles only to the big nodes and the notion of ST-compatibility extends naturally: A mixed ordered quiver is ST-compatible if the pruned quiver is.

If we consider mutation of ST-compatible mixed ordered quivers, there are two separate cases to consider: The mutation at big nodes works the same as before but it might add arrows between the small nodes. Mutation at small nodes on the other hand typically adds pairs of parallel and crossing arrows between big nodes, see \Cref{fig:mixedColMutEx}. This makes the interplay of the commutative and noncommutative part of the quiver interesting and nontrivial, in particular not every mutation at a small node is admissible and these mutations interact with each other as the example in \Cref{fig:mixedColMutEx2} shows: Starting with the quiver in the middle, mutation at both $4$ and $5$ is admissible, but once we perform either of these mutations, the other one is not admissible anymore. This behavior makes fruitfulness even harder to prove for ordered mixed quivers.

\input{figureSmallNodeMutNonAdmissible}

\section{Decorated Polygonal Tilings}\label{sec:Polygons}

In this section, we will construct a combinatorial model for ST-compatible quivers, which we will use to define mutations, seed algebras and ultimately our polygonal cluster algebras. This draws inspiration from triangulated surfaces which capture the combinatorics of Berenstein-Retakh's noncommutative surfaces, but instead of triangles we will glue bigger polygons. This leads to a new aspect: In a triangle, there is a canonical choice of opposite angle for every side, which is connected to the choice of a triangle invariant (and noncommutative angle) associated to every node in a surface type cluster algebra. Moving from triangles to polygons there is no canonical choice, which we will address by adding some decoration. Adding this allows for a second quiver to be encoded by a triangle (we have already seen these pseudocycles).

\subsection{Decorated polygons}

Let us first present the main definitions and result: Denote by $P_n$ the polygon in $\R^2\cong\C$ with vertex set $\left\{e^{ik2\pi/n}|1\leq k\leq n\right\}$. The group $\Z/n\Z$ acts naturally by rotation by $2\pi/n$ on $\C$, preserving $P_n$.

\begin{definition}\label{def:decPoly}
    A \keyword{decorated polygon} consists of the following data
    \begin{itemize}
        \item a polygon $P_n$,
        \item an orientation of its sides,
        \item an admissible partition of its vertices into two disjoint sets, called \keyword{parallel} and \keyword{crossing angles},
        \item a choice of vertex for every side such that the straight lines connecting the midpoint of each side to its corresponding vertex, which we call \keyword{angle indicators}, intersect pairwise where we allow intersection at the endpoints.
    \end{itemize}
    A decorated polygon is called \keyword{singular} if all angle indicators end at the same vertex, otherwise \keyword{regular}.
    
    A decoration of a regular odd-sided, respectively even-sided, polygon is called \keyword{admissible} if it has an odd, respectively even, number of parallel angles.

\end{definition}

Decorated polygons can be depicted graphically, as shown in \Cref{fig:decPolyEx}.

\input{figureDecPoly}

Recall that an ST-tile is an ordered quiver such that the underlying 2-colored quiver is a locally transitive tournament with arrows colored in such a way that each subquiver determined by three nodes is either a 3-cycle or a 3-pseudocycle, which is a constraint on the coloring. We will now establish a connection between these and decorated polygons:

\begin{theorem}\label{thm:decPolyElemSTCorrespondence}
    There is a 1:1 correspondence between decorated polygons up to rotation and isomorphism classes of ST-tiles.
\end{theorem}

The correspondence is explicit and the construction we will explain in the proof is shown for some examples in \Cref{fig:decPolyToQuiverEx}.

\input{figureDecPolyToQuiver}

\begin{proof}
    The proof has three steps as outlined in \Cref{tab:correspondenceLevels}.
    \begin{table}[h!]
        \centering
        \begin{tabular}{l|l}
            polygons + angle indicators & locally transitive tournaments\\
            + color of angles & + color of arrows\\
            + orientation of sides & + filling of nodes\\
            \hline
            decorated polygons & ST-tiles
        \end{tabular}
        \caption{Correspondence between levels of decorated polygons and ST-tiles.}
        \label{tab:correspondenceLevels}
    \end{table}

    \textit{Step 1:}
     We associate a quiver with cyclically ordered nodes to a polygon with angle indicators in the following way: For every side of the polygon, we get a node, the cyclic order is simply given by the cyclic order of the sides of the polygon. The angle indicator determines whether the other nodes are connected via an incoming arrow (it lies on the right side of the line) or via an outgoing arrow (left side), see \Cref{fig:decPolyToQuiver}.

    In this way we clearly obtain a tournament, i.e. an orientation of a complete graph. To show that this is locally transitive, pick any node and consider the subtournament given by all nodes connected via incoming (respectively outgoing) arrows. This is described by the above rule by collapsing the polygon, see \Cref{fig:decPolyCollapse}. The result is a decoration where all angle indicators end at the same vertex, which clearly corresponds to a totally ordered quiver, and thus a transitive tournament.

    \input{figureDecPolyTournProof}

    For the converse direction, we start with a locally transitive tournament $T$ and take a vertex $v$. Denoting by $T_{in}(v)$, respectively $T_{out}(v)$ the subtournament of vertices connected to $v$ by an incoming, respectively outgoing, arrow, we know that both of these are transitive. Therefore we obtain a total order on the vertex set, which we can denote
    \begin{equation*}
        T_{in}(v)\rightarrow v\rightarrow T_{out}(v)\,.
    \end{equation*}
    We get such a total order for any vertex and we aim to show that these are the same up to cyclic permutation. As mentioned before this was stated in \cite{brouwer1980enumeration} but we will give a full proof for completeness: Pick a vertex $v'\neq v$, without loss of generality $v'\in T_{in}(v)$. We get
    \begin{equation*}
        \underbrace{T_{in}(v)\cap T_{in}(v')}_{=:u}\rightarrow v' \rightarrow \underbrace{T_{in}(v)\cap T_{out}(v')}_{=:u'}\rightarrow v \rightarrow T_{out}(v)\,.
    \end{equation*}
    Furthermore, let us define
    \begin{align*}
        w&:=T_{out}(v)\cap T_{out}(v')\,,\\
        w'&:=T_{out}(v)\cap T_{in}(v')\,.
    \end{align*}
    Since $T$ is by assumption locally transitive we can deduce that the tournament looks schematically as shown in \Cref{fig:tournProof}, where the orange arrows are the ones whose orientation we deduced from local transitivity.

    \input{figureTournamentProof}

    Thus, the orders as defined by $v$ and $v'$ are
    \begin{align*}
        &u\rightarrow v'\rightarrow u'\rightarrow v\rightarrow w\rightarrow w'\,,\\
        &w'\rightarrow u\rightarrow v'\rightarrow u'\rightarrow v\rightarrow w\,,
    \end{align*}
    and indeed are just cyclically permuted. Therefore, we can associate a polygon with angle indicators by inverting the above construction. The angle indicators intersect pairwise, as a direct consequence of local transitivity, see \Cref{fig:decPolyCollapse}.\medskip

    \textit{Step 2:} Next we will show how the coloring of angles and arrows determine each other. Starting with the polygon with colored angles, we color the corresponding locally transitive tournament by coloring the \keyword{peripheral arrows} between nodes which are adjacent in the cyclic order, according to the color of the angle between the corresponding sides of the polygon. The only exception is the distinguished angle in a singular decorated polygon, i.e. the one where all the angle indicators end. In this case we color the corresponding arrow in the opposite color of the angle.
    
    We need to prove that this rule uniquely determines an admissible coloring of the entire quiver. To do so let $i,j$ be two sides of the polygon. Denote by $(i,j)$ all angles between $i$ and $j$ such that the angle indicators of $i,j$ do not end at angles in this interval, and set
    \begin{equation*}
        c(x):=\begin{cases}
            0 & x \text{ is parallel}\\
            1 & x \text{ is crossing}
        \end{cases}
    \end{equation*}
    for $x$ either an angle in the polygon or an arrow in the quiver. Fix the coloring of the quiver as follows: The color of the arrow $i\to j$, for two sides of the polygon, respectively nodes of the associated quiver is
    \begin{equation}\label{eq:arrowColorFromAngles}
        c(i,j)=\sum_{\alpha\in(i,j)}\big(1+c(\alpha)\big)+1\in\Z/2\Z\,.
    \end{equation}
    Observe that this indeed assigns the color of the angles to the corresponding peripheral arrows. Moreover, the quiver is consistently colored: There are two cases to consider which are schematically illustrated in \Cref{fig:CorrespondenceColorProof}.

    \input{figureCorrespondenceColorProof}
    
    If $(i,j,k)$ form an oriented cycle then
    \begin{align*}
        c(i,j)+c(j,k)+c(i,k)=\sum_{\alpha}\big(1+c(\alpha)\big)+1=0
    \end{align*}
    because the sum runs over all angles and we can thus use the parity condition for decorated polygons. This equality implies exactly that the cycle is colored consistently. On the other hand, when $(i,j,k)$ form a totally ordered subquiver with source $k$ and sink $j$, we have $(j,k)=(i,j)\sqcup (i,k)$ and therefore
    \begin{align*}
        c(i,j)+c(j,k)+c(i,k)=2\left(\sum_{\alpha\in(j,k)}\big(1+c(\alpha)\big)\right)+1=1\,,
    \end{align*}
    so that $(i,j,k)$ are the nodes of a pseudocycle in the 2-colored quiver.

    Uniqueness follows from \Cref{cor:LTTcoloringWeavingUnique} since the only weaving isomorphism which fixes the color of all peripheral arrows is the identity. For the converse direction simply color the angles of the decorated polygon according to the color of the peripheral arrows. To see that this is admissible we use induction on the number of nodes: Clearly then decorated polygon corresponding to an ST-tile with two nodes has one parallel and one crossing angle. If we add a node to an ST-tile it is not hard to see that the two new peripheral arrows have colors $c_1, c_2$ such that $c_1+c_2=c+1\in\Z/2\Z$ where $c$ is the color of the arrow that is no longer peripheral in the bigger ST-tile. But this is enough to inductively conclude that the angle colors satisfy the parity condition.\medskip

    \textit{Step 3:}
    To complete the translation of decorated polygons to ST-tiles, we need to decide the fill of the nodes. We choose these such that the empty ($\circ$) side is on the right-hand side of the corresponding side of the polygon when following the orientation.
\end{proof}

This will help us to obtain a better understanding of the combinatorics of ST-compatible quivers. In particular we make the following observation about decorated polygons:
\begin{proposition}
    Any decorated polygon can be obtained from a standard decorated polygon through a series of side collapses.
\end{proposition}

\begin{definition}
    A \keyword{standard decorated polygon} is an odd-sided polygon with the decoration chosen such that the angle indicators connect a side to the opposite angle, and such that all vertices are parallel angles.
\end{definition}

The first few standard decorated polygons are depicted in \Cref{fig:stdDecPolyEx}.
\input{figureStandardDecPoly}

\begin{definition}
    A \keyword{side collapse} at a side $S$ of a decorated $n$-gon yields a decorated $(n-1)$-gon by removing $S$ and identifying the two vertices that bound it. The angle indicators in the resulting polygon end at the image under the collapse of their original endpoint. The colors of the angles away from the edge which is collapsed stay the same, while the new angle is colored so that the new polygon is admissible. 

    A \keyword{vertex blowup} at vertex $v$ of a decorated polygon adds a new side between the two sides touching $v$ in such a way so that the collapse of the new side returns the old polygon. 
\end{definition}
There are multiple choices for a vertex blowup, since there are possibly many different polygons which can collapse to the original polygon. These choices are determined by choosing where the new angle indicator ends, and by choosing the colors of the two new angles. If $n$ angle indicators end at $v$ then there are $n+1$ valid choices for where to place the angle indicator of the new side; this determines the endpoints of all the other indicators. There are only two choices for admissible coloring of the new angles; a crossing angle can blow up into a pair of matching angles while a parallel angle blows up into distinct angles.

\begin{proof}
    First, we can embed any polygon in an odd sided polygon with standard angle indicators but possibly non-standard angle colors simply by blowing up each vertex which has $n$ angle indicators ending at it $n-1$ times, choosing the new angle indicators so that each new vertex only has one indicator ending at it. 

    Now assume we have a $2n+1$ sided polygon with standard angle indicators and some even number $k>0$ of crossing angles.
    We can fix the crossing angle indicators in the following way:
    \begin{enumerate}
        \item Label the vertices $1,2,\dots,2n+1$ with vertex 1 a crossing angle and at least one other crossing angle in the range $\{2,\dots,n\}$
        \item If vertex $n$ or $n+1$ is a crossing angle, then blow up both vertex 1 and the other crossing angle to make a $2n+3$ sided polygon with standard angle indicators and $k-2$ crossing angles. Go back to step 1.
        \item Otherwise blow up vertex $n$ to make an edge with new vertex $n'$ between $n$ and $n+1$ with vertex $n$ now a crossing angle and new angle indicator ending at vertex 1.
        \item Next blow up vertex 1 to create a new vertex labeled $1'$ between $1$ and $2n+1$ so that we have an odd sided polygon with standard angle indicators and angle $n$ crossing, but angles $1,1',n'$ straight. 
        \item Now, in the same fashion, blow up vertices $n$ and $2n+1$. If vertex $2n+1$ is crossing, then go back to step 1 with 2 less crossing angles, otherwise make vertex $(2n+1)'$ between $2n+1$ and $1'$ straight and vertex $2n+1$ crossing and vertex $n$ and $n''$ both straight. 
        \item Go back to step 1 with vertex $2n+1$ labeled as vertex 1. 
    \end{enumerate}
    This process terminates since we either reduce the number of crossing angles or we note that after renumbering in the last step we have a $2(n+2)+1$-gon with the vertex which was labeled $n-1$ now labeled $n+2$ and there is at least one crossing vertex between $\{5,\dots,n+2\}$, so that after at most $n$ iterations we will be able to cancel a crossing vertex directly at step 2. 
\end{proof}

\subsection{Decorated tilings}
\begin{definition}
    A \keyword{decorated tiling} of an oriented surface $\Sigma$ with marked points consists of the following data
    \begin{itemize}
        \item a set of non-intersecting oriented arcs between marked points on $\Sigma$, this means that the complementary regions are polygons,
        \item an admissible decoration for every such polygon s.t. the orientation of the edges agree with the orientation of the arcs on the surface,
        \item a choice of arrows connecting the arcs such that any arrow connecting two sides $A$ and $B$ of the same sub-polygon has the angle indicator of $A$ (equivalently $B$) on its right-hand side as given by the orientation of the arrow (again up to homotopy).
    \end{itemize}
\end{definition}

Our correspondence between ST-tiles and decorated polygons now extends to establish a connection with ST-compatible quivers:
\begin{theorem}\label{thm:polygonCorrespondence}
    There is a 1:1 correspondence between ST-compatible quivers and decorated tilings.
\end{theorem}

Using this, we will denote by $Q(T)$ the quiver corresponding to the tiling $T$, and by $T(Q)$ the tiling corresponding to the quiver $Q$. Sometimes we will also just use the quiver and the tiling interchangeably.

\begin{proof}
   We know from \Cref{thm:GSTRestrictions} that any ST-compatible quiver decomposes into ST-tiles and pairs of parallel and crossing arrows. Combined with \Cref{thm:decPolyElemSTCorrespondence} this gives the desired correspondence.
\end{proof}

\begin{remark}
    There are some technicalities when considering tilings of surfaces with nontrivial topology. Firstly, since the commutative arrows do not rely on the order or sides of the node they can ``move across the topology'' meaning that they can cancel with other arrows in surprising ways. Moreover, the decorated tiling may glue two sides of the same decorated polygon, which is a generalization of the self-folded triangles that can appear in triangulations of surfaces, cf. \Cref{sec:puncDigon}. We cannot currently deal with this phenomenon for bigger polygons, and will continue to require that mutation is only performed at nodes which are part of two distinct ST-tiles, equivalently arcs which border two different polygons.
\end{remark}

\subsection{Mutation of decorated tilings}\label{sec:MutationDecoratedTiling}

\begin{definition}\label{def:decoratedTilingFlip}
    Given a decorated tiling containing two decorated polygons $P_1$ and $P_2$ which are glued along a side $A$, a \keyword{flip} at $A$ is performed as follows:
    \begin{enumerate}
        \item Replace the angle indicators at $A$ with a chord of the polygon connecting their endpoints, the side $A'$. The orientation of $A'$ is chosen such that it is rotated in the clockwise direction with respect to $A$.
        \item Replace the original side $A$ by angle indicators at $A'$ with the same endpoints.
        \item The new angle indicators at all the other sides of $P_1$ and $P_2$ facing $A$ are such that they end at one of the endpoints of $A'$ and intersect the corresponding angle indicators at $A'$. This determines them uniquely.
        \item The color $c'$ of an angle at an endpoint of the angle indicators for $A'$ is given by $c'=c_1+c_2\in\Z/2\Z$ were $c_1,c_2$ are the colors of the two angles that were separated by $A$ before the flip. The color of the angle between $A'$ and an adjacent side $B$ is $c(A,B)$ from the original decorated polygon as given by \Cref{eq:arrowColorFromAngles}. All other angles remain unchanged.
        \item Introduce and cancel commutative arrows.
	\end{enumerate}
\end{definition}

Steps (1) and (2) just mean that we flip the colors between the side $A$ and its angle indicators. Step (5) requires some additional explanation: Suppose we perform a flip at a side $A$. For every side $B$ on the right of the angle indicator at $A$ add a commutative arrow to every side $C$ to the left of both the angle indicator of $A$ and of $B$, see \Cref{fig:commArrowRules}(a).

Moreover, if there is a commutative arrow connected to $A$, it gets reversed. If it has $A$ as its source, resp. target, we add a commutative arrow from its target, resp. source, to every side on the right, resp. left, of an angle indicator for $A$.

Furthermore commutative arrows follow the concatenation rules of quiver mutation when we mutate at their endpoints but since we think of them as arrows with weight $1+\varepsilon$, if there are $b$ commutative arrows from $B$ to $A$ and $c$ such arrows from $A$ to $C$, we introduce $2bc$ commutative arrows from $B$ to $C$ when mutating at $A$.

Finally, for cancellation, we need to look for commutative arrows which connect two sides, say $A$ and $B$, of the same decorated polygon in the following way: The angle indicator of $A$ meets an endpoint of $B$ on the left of the commutative arrow as given by its orientation, and vice versa, see \Cref{fig:commArrowRules}(c).

\input{figureCommArrowRules}

\begin{remark}
    We will give an easier description of the new angle color in terms of a $\Z/2\Z$-valued distance in \Cref{sec:Algebra}.
\end{remark}

Because of the extra decoration by commutative arrows, not every flip leads to another decorated tiling.
\begin{definition}
	A flip on a decorated tiling is \keyword{admissible} if the result is a decorated tiling, i.e. if the commutative arrows are admissible.
\end{definition}

\begin{remark}
    This is a generalization of the flips which relate different triangulations: Given a triangulation, choose the decoration where each angle indicator connects a side to the opposite angle. Then we have no pseudocycles and the mutation rule gives a flip in the triangulation.
\end{remark}

The mutation rule, including the introduction of commutative arrows, is illustrated in \Cref{fig:decPentaFlip}. If one follows the rules to mutate back, the rules for cancellation of commutative arrows need to be applied. Moreover, the decorated tilings in the polygons correspond exactly to the ST-compatible quivers in \Cref{fig:mutationEx}.

\input{figureDecPentaFlip}

This is of course a general observation:
\begin{proposition}
    The flip on decorated tilings is equivalent to mutation of the corresponding quiver. A flip at a side is admissible if and only if mutation at the corresponding node is.
\end{proposition}

\begin{proof}
    Simply translate the rules given above to rules on the quiver and observe that our flip respects the cyclic order since it is modeled on polygons. Concerning admissibility, our criterion on commutative arrows for decorated tilings is equivalent to the admissibility criteria on quivers: A non-admissible commutative arrow would correspond to a two-cycle in the quiver which, if canceled, would lead to a quiver with a tile that is not an ST-tile with the prescribed order.
\end{proof}

\begin{lemma}\label{lem:MutationPreservesAngleColors}
    In polygons that do not border $k$, mutation at $k$ can only change the angle indicators. In particular the colors of all angles are unchanged.
\end{lemma}
\begin{proof}
    The sides of a polygon $P$ not incident to $k$ can only be connected to $k$ commutatively. Thus mutation at $k$ can only change $P$ by introducing commutative arrows inside $P$. This can change the angle indicators. However to change the color of an angle we look at the quiver. This color would be changed if we introduce a pair of arrows which flips a peripheral arrow. But this is only possible if the mutation introduces enough commutative arrows to make the nodes connected by this arrow a source and sink in the ST-tile, respectively. In other words, the resulting decorated polygon would be singular. In this case our coloring convention is to use the opposite color of the arrow for the angle. 
\end{proof}

\section{Polygonal Cluster Algebras}\label{sec:Algebra}
We now use decorated tilings to define a \emph{seed algebra}, a quotient of the quiver algebra, which will serve as a replacement for the ambient field containing all cluster variables in the classical theory of cluster algebras. The key feature is that the seed algebra will contain a discrete family of vector spaces called \emph{angle spaces} over the polygonal tiling in which angles can be added while preserving symmetry and central norm. 

\subsection{Angle skeleton and transition maps}

\begin{definition}
    The \keyword{angle skeleton} of a polygonal tiling $T$ is a the graph, $\AngleSkeleton(T)$ given by associating two vertices to every arc of the tiling, one near each marked point. There is an oriented edge between two vertices on the same arc of a polygon with the same orientation as the side, and an edge between vertices near the same marked point on neighboring edges (See \Cref{fig:AngleSkeleton}).\\
    We classify the edges of the angle skeleton into three kinds, \keyword{sides}, \keyword{parallel angles}, or \keyword{crossing angles} depending on the feature of the tiling that they trace.
\end{definition}

We will associate two invariants to paths in the angle skeleton, a parity in $\Z/2\Z$ and a label in $\quiverAlgebra_{Q(T)}$. The label associated to a loop around a polygon will recover the angle associated with that polygon.
\begin{definition}
    Let $\gamma$ be a path in the angle skeleton. The \keyword{parity} of $\gamma$, $\angleParity(\gamma)$ is the additive function on paths given by $\angleParity(e) = 1$ if $e$ is a side or crossing angle and $\angleParity(e) = 0$ if $e$ is a parallel angle. 
\end{definition}
\begin{example}
    The parity of the clockwise path from $v_0$ to $v_2$ in \Cref{fig:AngleSkeleton} is 1 as $\gamma$ contains an odd number of sides and crossing angles.
\end{example}

\input{figureAngleSkeleton}

The parity condition for decorated polygons immediately translates to the following:
\begin{lemma}
    The parity of any loop around a polygon in an admissible tiling is odd. 
\end{lemma}
\begin{remark}
    The parity of a path $\gamma$ only depends on the homotopy class of $\gamma$ with fixed endpoints in the angle skeleton. Furthermore, we have $\angleParity(\overline{\gamma})=\angleParity(\gamma)$ for the reversed path $\overline{\gamma}$.
\end{remark}

Recall that the coloring of a quiver $Q(T)$ is determined by the corresponding decorated tiling $T$. We have seen in the proof of \Cref{thm:decPolyElemSTCorrespondence}:
\begin{lemma}\label{lemma:ParityDefinesArrowColor}
    Let $P$ be a polygon of a decorated tiling and $X,Y$ be two sides of $P$. The color of the arrow from $X$ to $Y$ is
    \begin{align*}
        c(X,Y)=\angleParity(\gamma)
    \end{align*}
    where $\gamma$ is the path in the angle skeleton around $P$ from the end of $X$ closest to $Y$ to the end of $Y$ closest to $X$ that does not cross the endpoints of the angle indicators for $X,Y$.
\end{lemma}
\begin{example}
    The color of the arrow from $A_1$ to $A_3$ in the quiver corresponding to the decorated tiling in \Cref{fig:AngleSkeleton} is $c(A_1,A_3)=0$, i.e. it is a parallel arrow.
\end{example}
\begin{proof}
    The parity of the path $\gamma$ coincides with the formula given in \Cref{eq:arrowColorFromAngles}.
\end{proof}

\begin{definition}
Let $Q$ be a ST-compatible quiver and $\gamma$ a path in the angle skeleton. The \keyword{canonical label} associated to $\gamma$ is the element in $\quiverAlgebra_Q$ defined inductively by
\begin{itemize}
    \item The canonical label associated to the empty path, straight angle, or crossing angle is the identity of $\quiverAlgebra_Q$.
    \item The canonical label associated to a side labeled $A$ oriented along $\gamma$ is $\sigma\tau(A)$. If the orientation of $A$ is opposite to $\gamma$ assign $\tau(A)$.
    \item The canonical label associated to the concatenation of paths $\gamma *\delta$ is given by
        \begin{equation*}
            \canonLabel(\gamma * \delta) = \canonLabel(\gamma) \sau^{\angleParity(\gamma)}(\canonLabel(\delta))
        \end{equation*}
\end{itemize}
When $P$ is a polygon of $Q$ and $\gamma$ is the clockwise oriented loop around $P$ based at $v_0$ we call $\canonLabel(\gamma)$ the \keyword{canonical angle} and write $\canonAngle{v_0,P}$ or $\canonAngle{v_0}$ when the polygon is clear.
\end{definition}

\begin{corollary}
    Let $P_1$ and $P_2$ be two polygons with the same angle skeleton. Then the canonical angles for $P_1$ are equal to the canonical angles for $P_2$, i.e. they do not depend on the angle indicators.
\end{corollary}

The canonical angles are closely related to the angles in the ST-compatible quiver.
\begin{proposition}
    Let $k$ be a side of a polygon $P$. The angle $\Delta_f(k)$ at $k$ in the corresponding ST-compatible quiver is
    \begin{equation*}
        \Delta_f(k)=\sau^{\angleParity(\gamma)}(\canonAngle{v})\,,
    \end{equation*}
    where $v$ is a vertex that is adjacent to the endpoint of the angle indicator of $k$ in $P$, and $\gamma$ is the clockwise path around $P$ starting at $v$ which passes through $k$ and ends at the vertex of $\AngleSkeleton(P)$ at the end of the side $k$, see \Cref{fig:angleVsCanonicalAngle}.
\end{proposition}

\input{figureAngleVsCanonicalAngle}

\begin{proof}
    First, observe that the order of the variables in $\alpha:=\sau^{\angleParity(\gamma)}(\canonAngle{v})$ and $\beta:=\Delta_f(k)$ are the same. Furthermore, we may assume without loss of generality that all sides of $P$ are oriented clockwise, since changing the orientation of a side $i$ amounts to replacing $X_i$ with $\sigma(X_i)$ in both $\alpha$ and $\beta$, and $\sigma$ is applied to the variables in both expressions in the same way as the fill $g$ of a node $X_i$ corresponds to the orientation of the corresponding side, where $g=0$ if $k$ is oriented clockwise around the polygon, and $g=1$ otherwise. So we only need to prove that the power $d_i$ of $\sigma\tau$ that is applied to a variable $X_i$ is the same in $\alpha$ and $\beta$.

    For $i\neq k$, in $\beta$ we have $d_i=t_i+1$ by \Cref{eq:defAngles} where $t_i$ is the color of the arrow connecting $i$ and $k$ in the quiver. Using \Cref{lemma:ParityDefinesArrowColor}, we find that $t_i=\angleParity(\gamma)+\angleParity(\gamma_i)$, where $\gamma_i$ is the clockwise path from $v$ to the vertex at the beginning of the side $i$. We have
    \begin{equation*}
        \alpha=(\sigma\tau)^{\angleParity(\gamma)}(\canonAngle{v})=(\sigma\tau)^{\angleParity(\gamma)}(\canonLabel(\gamma_i)(\sigma\tau)^{\angleParity(\gamma_i)}(\sigma\tau(X_i)\dotsc))
    \end{equation*}
    proving our claim. Finally, $X_k$ appears without $\sigma\tau$ in $\beta$. Since we have $\angleParity(\gamma_k)=\angleParity(\gamma)+1$, the above formula tells us that the same is true in $\alpha$, which concludes the proof.
\end{proof}

\begin{lemma}\label{lemma:CanonLabelReverse}
    The canonical labels associated to a path $\gamma$ and its reverse $\overline{\gamma}$ are related by
    \begin{equation*}
        \canonLabel(\overline{\gamma}) = \sau^{\angleParity(\gamma)+1}(\sigma(\canonLabel(\gamma))).
    \end{equation*}
\end{lemma}
\begin{proof}
    Suppose that $\gamma$ is a path from $v$ to $w$ in the angle skeleton. We may assume that the sides along that path are oriented according to $\gamma$ and that they are labeled $A_1,\dots,A_m$. The power of $\sigma\tau$ applied to $A_i$ in $\canonLabel(\gamma)$ is $d_i=\angleParity(\gamma|_{[v,v_i]})$, the parity of the path from $v$ to the vertex $v_i$ of $A_i$ farthest from $v$. Similarly, the power appearing in $\canonLabel(\overline{\gamma})$ is $e_i=\angleParity(\overline{\gamma}|_{[w,w_i]})$ with $w_i$ the vertex at the other end of $A_i$. Thus $\angleParity(\gamma)+1 = d_i + e_i$, and therefore $e_i=d_i+\angleParity(\gamma)+1$. Finally, as the edges are oriented in opposite directions we pick up a factor of $\sigma$ at each edge.
\end{proof}
\begin{corollary}\label{cor:canonLabelTrivPath}
    The canonical label of $\gamma * \overline{\gamma}$ is $N(\canonLabel(\gamma))$
\end{corollary}
\begin{proof}
    This is a simple computation combining \Cref{lemma:CanonLabelReverse} and the formula for concatenation of paths. 
\end{proof}
\begin{corollary}\label{cor:normCanonLabel}
    The canonical label for any path $\gamma$ in $\AngleSkeleton(T)$ is a monomial in $\quiverAlgebra_{Q(T)}$ and as such has central norm $N(\canonLabel(\gamma))$. For the reversed path $\overline{\gamma}$ we have
    \begin{align*}
        N(\canonLabel(\overline{\gamma}))=N(\canonLabel(\gamma))\,.
    \end{align*}
    For concatenation of paths we similarly find
    \begin{equation*}
        N(\canonLabel(\gamma*\delta))=N(\canonLabel(\gamma))\cdot N(\canonLabel(\delta))\,.
    \end{equation*}
\end{corollary}

\begin{example}
    The canonical angle at $v_0$ in \Cref{fig:AngleSkeleton} is
    \begin{align*}
        \canonAngle{v_0}=\sigma\tau(A_1)\sigma\tau(A_2)A_3\sigma\tau(A_4)\,,
    \end{align*}
    while the angle at $v_1$ is
    \begin{align*}
        \canonAngle{v_1}=A_2\sigma\tau(A_3)A_4\sigma\tau(A_1)\,.
    \end{align*}
    This means that the canonical angles are related via
    \begin{align*}
        \canonAngle{v_1}=\frac{1}{N(A_1)}\tau(A_1)\canonAngle{v_0}\sigma\tau(A_1)\,,
    \end{align*}
    which is an analogue of \Cref{eqn:BRAngleTransition}.
\end{example}

The canonical label can be used to relate canonical angles based at different vertices of a polygon in the angle skeleton, like we just saw in the example.
\begin{definition}\label{def:AngleTransitionMaps}
    The \keyword{transition map} associated to a path $\gamma$ in the angle skeleton is the map $m_\gamma \colon \quiverAlgebra_Q \rightarrow \quiverAlgebra_Q$ given by 
    \begin{equation*}
        m_\gamma(x):=\frac{1}{N(\canonLabel(\gamma))}\sau^{\angleParity(\gamma)}\Big(\tau(\canonLabel(\gamma))\,x\, \sigma\tau(\canonLabel(\gamma))\Big)
    \end{equation*}
\end{definition}
Note that if $\gamma$ is a parallel angle then $m_\gamma$ is the identity, and if $\gamma$ is a single crossing angle then $m_\gamma = \sigma\tau$.

\begin{lemma}\label{lemma:pathConcatenationTransition}
    For any two composable paths $\gamma,\delta$ in $\AngleSkeleton(T)$ we have
    \begin{align*}
        m_{\gamma*\delta}=m_\delta\circ m_\gamma\,.
    \end{align*}
    In particular
    \begin{align*}
        m_{\overline{\gamma}} \circ m_\gamma=\mathrm{id}_{\quiverAlgebra_{Q(T)}}.
    \end{align*}
\end{lemma}

\begin{proof}
We use \Cref{cor:normCanonLabel} and definition of the canonical label to compute:
\begin{align*}
    m_{\gamma*\delta} &= \frac{1}{N(\canonLabel(\gamma*\delta))}\sau^{\angleParity(\gamma*\tau)}\Big(\tau(\canonLabel(\gamma*\delta))\,x\,\sau(\canonLabel(\gamma*\delta))\Big)\\
    &= \frac{1}{N(\canonLabel(\gamma))N(\canonLabel(\delta))}\sau^{\angleParity(\gamma)+\angleParity(\delta)}{\bigg (}\tau{\Big(}\canonLabel(\gamma)\sau^{\angleParity(\gamma)}(\canonLabel(\delta)){\Big)}\,x\,\sau{\Big(}\canonLabel(\gamma)\sau^{\angleParity(\gamma)}(\canonLabel(\delta)){\Big)}{\bigg)}\\
    &= \frac{1}{N(\canonLabel(\delta))}\sau^{\angleParity(\delta)}{\bigg(}\tau(\canonLabel(\delta))\frac{1}{N(\canonLabel(\gamma))}\sau^{\angleParity(\gamma)}{\Big(}\tau(\canonLabel(\gamma))\,x\,\sau(\canonLabel(\gamma)){\Big)}\sau(\canonLabel(\gamma)){\bigg)}\\
    &= m_\delta(m_\gamma(x))\,.
\end{align*}
To prove the second claim, we simply use the concatenation rule we just proved.
\begin{align*}
    m_{\overline{\gamma}}(m_\gamma(x))=m_{\gamma*\overline{\gamma}}(x)=\frac{1}{N(\canonLabel(\gamma))}\sau^{2\angleParity(\gamma)}\Big(\tau(\canonLabel(\gamma*\overline{\gamma}))\,x\,\sau(\canonLabel(\gamma*\overline{\gamma}))\Big)=x
\end{align*}
because \Cref{cor:canonLabelTrivPath} tells us that $\canonLabel(\gamma*\overline{\gamma})=N(\canonLabel(\gamma))$\,.
\end{proof}

\begin{corollary}
    The transition map $m_\gamma$ only depends on $\gamma$ up to homotopy with fixed endpoints in the angle skeleton.
\end{corollary}

The transition maps relate the different canonical angles for the same polygon.
\begin{lemma}\label{lemma:anglePolygonInvariant}
    Let $\gamma$ be a path in the angle skeleton along the boundary of a polygon $P$ from $v_0$ to $v_1$. Then 
    \begin{equation}
        \canonAngle{v_1} = m_{\gamma}(\canonAngle{v_0})\,.
    \end{equation}
\end{lemma}
\begin{proof}
    Let $M_0$ be the path around $P$ starting at $v_0$ and let $M_1$ be the analogous path at $v_1$. Then $\overline{\gamma}*M_0*\gamma = \overline{\gamma}*\gamma *M_1$ or $M_1 * \overline{\gamma}*\gamma$ depending on whether $\gamma$ goes clockwise or counterclockwise around $P$. Using the properties of the canonical label we have
    \begin{align*}
        \canonLabel(\overline{\gamma}*M_0*\gamma) =~& \canonLabel(\overline{\gamma})\sau^{\angleParity(\gamma)}(\canonLabel(M_0)) \sau^{\angleParity(\gamma)+\angleParity(M_0)}(\canonLabel(\gamma))\\
        =~& \sau^{\angleParity(\gamma)+1}(\sigma(\canonLabel(\gamma)))\sau^{\angleParity(\gamma)}(\canonAngle{v_0}) \sau^{\angleParity(\gamma)+1}(\canonLabel(\gamma))\,.
    \end{align*}
    On the other hand
    \begin{equation*}
        \canonLabel(M_1*\overline{\gamma}*\gamma ) = \canonLabel(\overline{\gamma}*\gamma *M_1) = N(\canonLabel(\gamma)) \canonLabel(M_1) = N(\canonLabel(\gamma)) \canonAngle{v_1}\,.
    \end{equation*}
    So $\canonAngle{v_1} = \frac{1}{N(\canonLabel(\gamma))}\,\sau^{\angleParity(\gamma)}\left(\tau(\canonLabel(\gamma)) \canonAngle{v_0} \sau(\canonLabel(\gamma)) \right)$ as needed.\\
\end{proof}

The transition maps from \Cref{def:AngleTransitionMaps} define a discrete local system on the tiled surface with a copy of $\quiverAlgebra_Q$ over each vertex of the angle skeleton. 

\begin{corollary}\label{cor:AngleMonodromy}
    The monodromy of the local system, $m_{P,v}$ around a polygon $P$ starting at a vertex $v$ fixes the canonical angle at $v$.
 
    If moreover $\sigma(\canonAngle{v})=\canonAngle{v}$, then
    \begin{align}
        m_{P,v}^2 &= \mathrm{Id}_{\quiverAlgebra_Q}\\
        m_{P,v}(x) &= -x+2\frac{b(x,\canonAngle{v})}{b(\canonAngle{v},\canonAngle{v})}\canonAngle{v} \label{eqn:AngleMonodromyReflection}
    \end{align}
    where $b(u,v) = u\tau(v) + v\tau(u)$.
\end{corollary}
\begin{proof}
This is a simple computation using \Cref{def:AngleTransitionMaps} and the fact that the path around $P$ is assumed to have parity $1$. 
\end{proof}
\begin{remark}
The function $b \colon \quiverAlgebra_Q \times \quiverAlgebra_Q \rightarrow \quiverAlgebra_Q$ from \Cref{cor:AngleMonodromy} is symmetric, bilinear and invariant under $\tau$. Furthermore if the image was a field it would be an inner product and the monodromy around a polygon would be a reflection. 
\end{remark}
Additionally, $b(u,v)$ is the obstruction for $u+v$ to have central norm when both $u$ and $v$ have central norm.
\begin{lemma}
If $u,v \in \quiverAlgebra_Q$ have central norm then  $b(u,v)$ is central if and only if $u+v$ has central norm.
\end{lemma}

\subsection{Seed algebra}
With our discrete local system in place, we are now in a position to construct the \emph{seed algebra} in which mutation of variables will make sense. This algebra should be constructed in such a way that all angles are $\sigma$-invariant and can be added whenever we use transition maps to `move them to the same vertex' of the angle skeleton.

For this suppose that $Q$ is an ST-compatible quiver with a corresponding decorated tiling $T$ of a surface $\Sigma$ by polygons $P_i$ for $i\in I$.
\begin{definition}\label{def:seedAlgebra}
    Let $\{v_i|i\in I\}$ be a set of vertices of $\AngleSkeleton(T)$ such that $v_i$ lies in $P_i$. Denote $\canonAngle{i}:=\canonAngle{v_i,P_i}$. The \keyword{seed algebra} $\seedAlgebra_Q$ of $Q$ is obtained in two steps:
    \begin{enumerate}
        \item quotient $\quiverAlgebra_Q$ by the relations $\canonAngle{i}=\sigma(\canonAngle{i})$ for all $i\in I$ and $b(m_{\eta_i}\canonAngle{i},m_{\eta_j}\canonAngle{j})$ is central and fixed by $\sigma$ for all vertices $p$ of $\AngleSkeleton(T)$, and all pairs of paths $\eta_{i}: v_i \to p$ and $\eta_j:v_{j}\to p$ 
        \item localize the intermediate algebra by its center\,.
    \end{enumerate}
\end{definition}

The seed algebra is well-defined and in fact the construction can be simplified slightly.
\begin{proposition}\label{prop:seedAlgebraChoiceIndep}
    The seed algebra is independent of the choice of vertices $\{v_i\}$. Moreover, we can fix a base point $p$ in the construction and still obtain the same algebra.
\end{proposition}

The proof makes use of the following lemma:
\begin{lemma}[Sliding Lemma]\label{lemma:sliding}
    Suppose $v,w$ lie in (some quotient of) $\quiverAlgebra_Q$ and are such that $b(v,w)$ is central and fixed by $\sigma$. Then
    \begin{align*}
        b(m_\delta(v),m_\delta(w))=b(v,w)
    \end{align*}
    for any path $\delta$ in $\AngleSkeleton(T)$.
\end{lemma}
\begin{proof}
    Assuming $b(v,w)$ is central and fixed by $\sigma$, and setting $l:=\angleParity(\delta)$ and $x:=\canonLabel(\delta)$, we have
    \begin{align*}
        b(m_\delta(v),m_\delta(w)) &= m_\delta(v)\tau(m_\delta(w))+m_\delta(w)\tau(m_\delta(v))\\
        &= \frac{1}{N(x)^2}\sau^l \left(\tau(x)v\sau(x)\sigma(x)\tau(w)x+\tau(x)w\sau(x)\sigma(x)\tau(v)x\right)\\
        &= \frac{1}{N(x)}\sau^l \left(\tau(x)b(v,w)x\right) =b(v,w)\,.
    \end{align*}
    Now \Cref{lemma:pathConcatenationTransition} lets us reverse $\delta$ to prove the reverse implication.
\end{proof}

\begin{proof}[Proof of \Cref{prop:seedAlgebraChoiceIndep}]
    First, suppose that we replace a vertex $v_i$ in the construction of the seed algebra $\seedAlgebra_Q$ with some other vertex $v_i'$ of $\AngleSkeleton(T)$ which lies in the same polygon $P_i$. Then the two can be joined by a path $\gamma:v_i'\to v_i$ along the boundary of $P_i$. By \Cref{lemma:anglePolygonInvariant} we have $m_\gamma(\canonAngle{v_i'})=\canonAngle{v_i}$, and thus using \Cref{lemma:pathConcatenationTransition} we find
    \begin{align*}
        b(m_{\eta_i}\canonAngle{v_i},m_{\eta_j}\canonAngle{v_j})=b(m_{\gamma*\eta_i}\canonAngle{v_i'},m_{\eta_j}\canonAngle{v_j})\,,
    \end{align*}
    for $\eta_i,\eta_j$ as above. Using the same calculation with the roles of $v_i$ and $v_i'$ exchanged, we find that the intermediate algebras agree, and therefore we obtain the same seed algebra for both choices.\medskip

    For the second statement, we fix some base point $p$. Any other point $p'$ can be connected to it by some path $\gamma:p'\to p$. Consider two paths $\eta_{i} : v_i \to p'$ and $\eta_{j} : v_{j} \rightarrow p'$. Composing with $\gamma$ gives paths to our fixed base point $p$. By assumption the Sliding Lemma (\Cref{lemma:sliding}) applies and
    \begin{align*}
        b(m_{\eta_i}\canonAngle{i},m_{\eta_j}\canonAngle{j})=b(m_\gamma m_{\eta_i}\canonAngle{i},m_\gamma m_{\eta_j}\canonAngle{j})\,.
    \end{align*}
    Therefore the compatibility of paths ending at an arbitrary point $p'$ are implied by the compatibility of paths ending at $p$ as needed.
\end{proof}

In certain cases the construction of the seed algebra can be simplified considerably as we will be able to write down a \emph{finite} set of relations.

\begin{theorem}\label{thm:finiteCentrConds}
    Let $\{v_i\}$, $\canonAngle{i}$ be as in \Cref{def:seedAlgebra}, and let $p$ be a fixed base point in $\AngleSkeleton(T)$. Denote by $\delta_i$ the clockwise loop around the boundary of the polygon $P_i$ based at $v_i$. Suppose that there exists a set of paths $\{\gamma_i:v_i\to p\,|\,i\in I\}$ such that $\{\overline{\gamma_i}*\delta_i*\gamma_i\,|\,i\in I\}$ is a generating set for $\pi_1(\AngleSkeleton(T),p)$. Then $\seedAlgebra_Q$ is obtained by
    \begin{enumerate}
        \item taking the quotient $\quiverAlgebra_Q$ by the relations that $\canonAngle{i}=\sigma(\canonAngle{i})$ for all $i\in I$ and $b(m_{\gamma_i}\canonAngle{i},m_{\gamma_j}\canonAngle{j})$ is central and fixed by $\sigma$ for all $i,j\in I$
        \item localize the intermediate algebra by its center.
    \end{enumerate}
\end{theorem}

\begin{proof}
    Set $\delta'_i:=\overline{\gamma_i}*\delta_i*\gamma_i$ for all $i\in I$. These paths are loops based at $p$ which by assumption represent a set of generators for $\pi_1(\AngleSkeleton(T),p)$. Suppose now that $i,j\in I$ and that $\eta_{i,j}:v_{i,j}\to p$ are two paths. We aim to show that
    \begin{align*}
        b(m_{\eta_i}\canonAngle{v_i},m_{\eta_j}\canonAngle{v_j})=b(m_{\eta_i}m_{\overline{\gamma_i}}m_{\gamma_i}\canonAngle{v_i},m_{\eta_j}m_{\overline{\gamma_j}}m_{\gamma_j}\canonAngle{v_j})=:x
    \end{align*}
    is central and fixed by $\sigma$ in the algebra $\seedAlgebra'_Q$ defined in the theorem, which would allow us to conclude that $\seedAlgebra'_Q=\seedAlgebra_Q$ (since the reverse direction is clearly true).

    To do so, consider
    \begin{align*}
        b(m_{\gamma_i}\canonAngle{i},m_{\gamma_i}m_{\overline{\eta_i}}m_{\eta_j}m_{\overline{\gamma_j}}\canonAngle{j})=:y\,.
    \end{align*}
    We use \Cref{lemma:pathConcatenationTransition} and the assumption to see
    \begin{align*}
        m_{\gamma_i}m_{\overline{\eta_i}}m_{\eta_j}m_{\overline{\gamma_j}}=m_{\overline{\gamma_j}*\eta_j*\overline{\eta_i}*\gamma_i}=\prod_k m_{\delta'_k}\,.
    \end{align*}
    A quick calculation shows that $\canonLabel(\delta'_k)=m_{\gamma_k}\canonAngle{k}$ which is fixed by $\sigma$ in $\seedAlgebra'_Q$. This means that \Cref{cor:AngleMonodromy} applies:
    \begin{align*}
        m_{\delta'_k}m_{\gamma_j}\canonAngle{j}=-m_{\gamma_j}\canonAngle{j}+2\frac{b(m_{\gamma_k}\canonAngle{k},m_{\gamma_j}\canonAngle{j})}{b(m_{\gamma_k}\canonAngle{k},m_{\gamma_k}\canonAngle{k})}m_{\gamma_k}\canonAngle{k}\,.
    \end{align*}
    By induction on $k$ we find that $y$ is central and fixed by $\sigma$, which allows us to apply the Sliding Lemma (\Cref{lemma:sliding}) to conclude $x=y$.
\end{proof}

\begin{corollary}
    If the quiver $Q$ corresponds to a decorated tiling $T$ of the disk with marked points in the boundary, the seed algebra $\seedAlgebra_Q$ can be constructed as in \Cref{thm:finiteCentrConds}.
\end{corollary}

\begin{proof}
    Fix any base point $p$ in $\AngleSkeleton(T)$. We have $\pi_1(\AngleSkeleton(T),p)=\ast^n\Z$, i.e. the free group on $n$ generators, where $n=|I|$ is the number of polygons in the tiling. Let us now consider the dual graph of the tiling, with vertices corresponding to polygons which are connected by an edge whenever they are adjacent. Since we assume that the underlying surface is the disk with marked points in the boundary, this dual graph is a tree. For every edge $i$---$j$ in the tree we choose an undirected path $\gamma_{ij}:v_i\to v_j$ in $\AngleSkeleton(T)$, which can then be concatenated and oriented to give paths $\gamma'_i:v_i\to v_1$. Now any path $\eta:v_1\to p$ allows us to take $y_i:=\gamma'_i*\eta$, which are as in the assumption of \Cref{thm:finiteCentrConds}.
\end{proof}

The above constructions also allow us to associate a vector space to every vertex of the angle skeleton:
\begin{definition}
    The \keyword{angle space} over a vertex $v$ of the angle skeleton is the $\mathrm{Z}(\seedAlgebra_Q)$-vector space generated by $\{m_{\eta_i}\canonAngle{i}\,|\,\eta_i:v_i\to p, i\in I\}$.
\end{definition}

\begin{remark}
    Typically, the angle space is infinite-dimensional. Whenever \Cref{thm:finiteCentrConds} applies, the space is finite-dimensional.
\end{remark}

\begin{proposition}
    The bilinear map $b$ restricts to an inner product on any angle space. Moreover, for any path $\gamma:v\to v'$ the map $m_\gamma$ restricts to an isomorphism of the angle spaces at $v$ and $v'$ which preserves $b$.
\end{proposition}

\begin{proof}
    On any angle space the bilinear form takes values in the field of fractions of $\mathrm{Z}(\seedAlgebra_Q)$. The corresponding quadratic form is clearly nondegenerate. The second part is proved analogously to \Cref{prop:seedAlgebraChoiceIndep}.
\end{proof}

\begin{proposition}\label{prop:weavIsoSeedAlgebras}
    If $Q,Q'$ are weaving isomorphic ST-compatible quivers, then $\seedAlgebra_Q\cong\seedAlgebra_{Q'}.$
\end{proposition}

\begin{proof}
    We only need to observe that the isomorphism of the quiver algebras $\quiverAlgebra_Q\cong\quiverAlgebra_{Q'}$ from \Cref{prop:weavingIsoQuiverAlgebras} descends to an isomorphism of the intermediate algebras. This is due to the fact that it sends canonical angles to canonical angles, and more generally maps labels for the same path to each other. This can be seen in the same way as in \Cref{prop:weavingIsoQuiverAlgebras}. Therefore the ideals in the construction of the intermediate algebra for both quivers get mapped to each other.
\end{proof}

\begin{remark}
    The construction of the seed algebra also extends naturally to mixed quivers. The angle skeleton of the pruned quiver defines all the relations to quotient the quiver algebra by. The resulting algebra is simply the seed algebra of the pruned quiver with additional central variables for the small nodes.
\end{remark}

\subsection{Mutation of variables}
In \Cref{sec:MutationOfSCQuivers} and \Cref{sec:MutationDecoratedTiling} we defined mutation of ST-compatible quivers and decorated tilings. We will now build a noncommutative cluster algebra from a quiver by choosing the variables $X_i$ from the seed algebra $\seedAlgebra_Q$ as initial variables for $Q$. We will see that $\seedAlgebra_Q$ was constructed to contain all possible cluster variables obtained via admissible mutation.

Let us define the formula to mutate a variable $X$ attached to an arc between two adjacent polygons. Consider \Cref{fig:GeneratlMutation} where we see a schematic local picture of the tiling and its angle skeleton before and after mutation (commutative arrows are suppressed). Here each outer edge represents a sequence of edges of the actual polygonal tiles.

\input{figureMutationVariables}

We label the corners of the polygon incident to the edge $X$ and its angle indicators clockwise $t,u,v,w$. So $X$ is the arc from $t$ to $v$ and will mutate to the arc $X'$ from $u$ to $w$ under the decorated tiling mutation described in \Cref{sec:MutationDecoratedTiling}. By slight abuse of notation let $t$ and $v$ be the vertices of the angle skeleton of $X$ near the corresponding corners before mutation, similar for $u,w$ after mutation.

Moreover, let $\gamma_A$ be the path in the angle skeleton connecting $u_1$ and $v_1$. Similarly $\gamma_B$ connects $u_2$ and $t_2$, $\gamma_C$ connects $t_1$ and $w_1$, and $\gamma_D$ connects $v_2$ and $w_2$ as shown in \Cref{fig:GeneratlMutation}. Note that $X$ has an arrow in from every edge of $\gamma_A$ and $\gamma_C$, and an arrow out to every edge of $\gamma_B$ and $\gamma_D$. The important part is that their parity will tell us whether an angle is parallel or crossing. Finally, we label the relevant angles before and after mutation as shown in \Cref{fig:GeneratlMutation}.

We now express the mutation formula using the canonical labels of these paths:
\begin{align}\label{eqn:FullMutationFormula_withnorm}
\begin{split}
    N(X)\mu(X) = \prod\limits_{Y\in \In^c(X)}\hspace{-2mm}N(Y)\cdot &\sau^{\angleParity(\gamma_A*\alpha_v)}\canonLabel(\gamma_A*\alpha_v*\overline{X}*\beta_t*\gamma_C)\\
    &+\prod_{Z\in \Out^c(X)}\hspace{-2mm}N(Z)\cdot \sau^{\angleParity(\gamma_B*\alpha_t)}\canonLabel(\gamma_B*\alpha_t*X*\beta_v*\gamma_D)\,.
\end{split}
\end{align}
To simplify this, we observe
\begin{align*}
    \sau^{\angleParity(\gamma_A*\alpha_v)}\canonLabel(\gamma_A*\alpha_v*\overline{X}*\beta_t*\gamma_C) &= \sau^{\angleParity(\gamma_A*\alpha_v)}\Big
    (\canonLabel(\gamma_A*\alpha_v)\Big)\canonLabel(\overline{X}*\beta_t*\gamma_C)\\
    &=\sau^{\angleParity(\gamma_A*\alpha_v)}\Big(\canonLabel(\gamma_A*\alpha_v)\Big)\,\tau(X)\,\sau\Big(\canonLabel(\beta_t*\gamma_C)\Big)
\end{align*}
and similarly for the second term. Now recall that $\tau(X)/N(X)=X^{-1}$ if $X$ has central norm to arrive at the following compact version of the mutation formula:
\begin{equation}\label{eqn:FullMutationFormula}
    \mu(X) = \prod_{Y\in \In^c(X)}N(Y)\cdot AX^{-1}C + \prod_{Z\in \Out^c(X)}N(Z)\cdot B\,\sigma(X)^{-1}D
\end{equation}
where
\begin{equation*}
    \begin{array}{lll}
        A= \sau^{\angleParity(\gamma_A*\alpha_v)}(\canonLabel(\gamma_A)) & \hspace{2pc} & B= \sau^{\angleParity(\gamma_B*\alpha_t)}(\canonLabel(\gamma_B))\\
        C= \sau(\canonLabel(\beta_t*\gamma_C)) & \hspace{2pc} & D=\sau(\canonLabel(\beta_v*\gamma_D))\,.
    \end{array}
\end{equation*}

\begin{remark}
    Note that when $\gamma_A,\gamma_B,\gamma_C,\gamma_D$ all contain one edge and $X$ has no commutative arrows this mutation rule reduces to the usual noncommutative mutation of triangulated surfaces as given \Cref{eq:exchRel}. Consequently, the paths used in the mutation formula \Cref{eqn:FullMutationFormula_withnorm} can be understood as the familiar zigzag paths which connect the two opposite vertices $u$ and $w$ and go through $X$.
\end{remark}

\begin{lemma}
    The mutation rule is equivalently read from the corresponding ordered quiver using colored arrows:
    \begin{equation*}
        X\xrightarrow{f~t~g}Z \mapsto \sigma^{f+g}\sau^t(Z) \hspace{2pc} Y\xrightarrow{f~t~g}X \mapsto \sigma^{f+g+1}\sau^t(Y)\,.
    \end{equation*}
    So that the mutation formula can be expressed as
    \begin{align*}
        \mu(X)=\prod\limits_{Y\in \In^c(X)}\hspace{-2mm}&N(Y)\cdot \hspace{-3mm}\prod_{\substack{\alpha\colon\: A\xrightarrow{f~t~\bullet} X\\ \alpha \in \In_\bullet(X)}}\hspace{-3mm}\left(\sigma^{f}\sau^t(A)\right)\cdot X^{-1}\cdot\hspace{-3mm}\prod_{\substack{\alpha\colon\: C\xrightarrow{f~t~\circ} X\\ \alpha \in \In_\circ(X)^*}}\hspace{-3mm}\left(\sigma^{f+1}\sau^t(C)\right)\\
        &+\prod\limits_{Z\in \Out^c(X)}\hspace{-2mm}N(Z)\cdot \hspace{-3mm}\prod_{\substack{\alpha\colon\: X\xrightarrow{\bullet~t~f} B\\ \alpha \in \Out_\bullet(X)^*}}\hspace{-3mm}\left(\sigma^{f+1}\sau^t(B)\right)\cdot \sigma(X)^{-1}\cdot\hspace{-3mm}\prod_{\substack{\alpha\colon\: X\xrightarrow{\circ~t~f} D\\ \alpha \in \Out_\circ(X)}}\hspace{-3mm}\left(\sigma^{f}\sau^t(D)\right)
    \end{align*}
    where $\In_f(X)^*,\Out_f(X)^*$ denote the respective tuples with the reversed order.
\end{lemma}
\begin{proof}
    Suppose that an edge of the angle skeleton in the path $\gamma_A$ is labeled $A$ and is oriented the same as $\gamma_A$, i.e. $\gamma_A=\gamma_A'*A*\gamma_A''$. This assumption is equivalent to assuming the fill at $X$ and $A$ are different. By additivity of the parity we find $\angleParity(\gamma_A)=\angleParity(\gamma_A')+\angleParity(\gamma_A'')+1$ and since we are working over $\Z/2\Z$, we have equivalently $\angleParity(\gamma_A'')=\angleParity(\gamma_A)+\angleParity(\gamma_A')+1$. Now the color $t$ of the arrow between the nodes corresponding to $A$ and $X$ according to \Cref{lemma:ParityDefinesArrowColor} is given by
    \begin{align*}
        t=\angleParity(\gamma_A'')+\angleParity(\alpha_v)=\angleParity(\gamma_A)+\angleParity(\gamma_A')+\angleParity(\alpha_v)+1\,.
    \end{align*}
    According to the formula from the Lemma, $A$ should appear as $\sau^t(A)$ in $\mu(X)$. On the other hand, the mutation formula \Cref{eqn:FullMutationFormula} includes the variable $A$ in
    \begin{align*}
        \sau^{\angleParity(\gamma_A*\alpha_v)}(\canonLabel(\gamma_A)) &= \sau^{\angleParity(\gamma_A*\alpha_v)}(\canonLabel(\gamma_A'))\sau^{\angleParity(\gamma_A*\alpha_v)}\sau^{\angleParity(\gamma_A')}(\sau(A))\cdots\\
        &= \cdots \sau^{\angleParity(\gamma_A)+\angleParity(\gamma_A')+\angleParity(\alpha_v)+1}(A)\cdots = \cdots\sau^t(A)\cdots\,
    \end{align*}
    as claimed. Note that switching the fill at one of these nodes is equivalent to flipping the orientation of the corresponding side of the polygon and thus to applying $\sigma$ to the variable. Therefore the formula for the power of $\sigma$ agrees for all sides appearing in $\gamma_A$. Analogous computations work for all other paths noting that for paths out from $X$, i.e. $\gamma_B,\gamma_D$, the fills aligning with the paths disagree with the fill at $X$.
\end{proof}

\begin{example}
    For the decorated tiling shown in \Cref{fig:AngleSkeleton}, the mutation at $A_2$ gives the new variable
    \begin{equation*}
        A_2'=\tau(A_6)A_2^{-1}\tau(A_1)+A_5\sigma(A_2)^{-1}A_3\sau(A_4)\,.
    \end{equation*}
\end{example}

We will now present the main theorem of this subsection.
\begin{theorem}\label{thm:seedalgebra_iso}
    Let $Q$ and $Q'$ be two ST-compatible quivers related by admissible mutation at the node $k$. Then $\seedAlgebra_{Q'}$ is isomorphic to $\seedAlgebra_Q$ via the map $\mu_k$ which is given on the generators of $\seedAlgebra_{Q'}$ by
    \begin{align*}
        X'_i\mapsto\begin{cases}
            X_i & i\neq k\\
            \mu(X_k) & i=k
        \end{cases}
    \end{align*}
    with $\mu(X_k)$ as given in \Cref{eqn:FullMutationFormula}. We will call this map a \keyword{mutation isomorphism}.
\end{theorem}

\begin{proof}
    The above map is a priori only defined between the quiver algebras. To show that it descends to the seed algebras, i.e. is well-defined we must show that
    \begin{itemize}
        \item the image of canonical angles are fixed by $\sigma$
        \item the image of canonical angles are in the angle spaces
        \item the image of $X_k'$ has central norm
    \end{itemize} 
    First we lighten notation by taking $X = X_k$ and $X' = X_k'$.  We can now simply choose base points for the canonical angles at all polygons that do not border $k$ away from polygons bordering $k$. Since admissible mutation does not change the canonical angles for the polygons, these angles of $Q'$ map exactly onto angles of $Q$ and thus are fixed by $\sigma$ and in the angle vector spaces.
    
    We will show that the remaining angles in the new polygon can be written as a linear combination with central coefficients of angles in the old polygon. Thus they are sums of $\sigma$-fixed elements and belong to the angle vector spaces as needed.
    
    Consider again \Cref{fig:GeneratlMutation}. For the polygons bordering $X_k$ in the original tiling, base the canonical angles at $t$. Let $\blacktriangle_t$ be the angle for the top polygon and $\triangledown_t$ be the angle for the bottom polygon. Then
    \begin{align*}
        \blacktriangle_t &= \canonLabel(\alpha_t*\overline{\gamma_B}*\theta_u*\gamma_A*\alpha_v*\overline{X})\\
        \triangledown_t &= \canonLabel(X*\beta_v*\gamma_D*\theta_w*\overline{\gamma_C}*\beta_t)\,.
    \end{align*}
    
    In the mutated tiling the vertex $t$ of the angle skeleton does not exist. So we instead compute the angle $\triangleleft_{t_2}$ at $t_2$ and use the transition map $m_{\alpha_t}$ to move it to $t$ once we compute the image in $\quiverAlgebra_Q$. We have
    \begin{align*}
        \triangleleft_{t_2} &= \canonLabel(\overline{\gamma_B}*\beta_u*X'*\beta_w *\overline{\gamma_C}*\theta_t)\\
        &= \canonLabel(\overline{\gamma_B})\sau^{\angleParity(\gamma_B)+\angleParity(\beta_u)+1}(X'\,\canonLabel(\beta_w*\overline{\gamma_C}))\\
        &= \canonLabel(\overline{\gamma_B})\sau^{\angleParity(\alpha_t)}(X'\,\canonLabel(\beta_w*\overline{\gamma_C}))\,,
    \end{align*}
    where we used the additivity of the parity, the fact that the canonical label of angles is 1, and the mutation rule for the angle colors in decorated polygons, which implies $\angleParity(\beta_u)=\angleParity(\gamma_B)+\angleParity(\alpha_t)+1$. To see this combine \Cref{def:decoratedTilingFlip} with \Cref{lemma:ParityDefinesArrowColor}.

    Therefore we get the following after applying the transition map
    \begin{align*}
        m_{\alpha_t}(\triangleleft_{t_2}) &= \sau^{\angleParity(\alpha_t)}(\triangleleft_{t_2})=\sau^{\angleParity(\alpha_t)}\big(\canonLabel(\overline{\gamma_B})\big)X'\canonLabel(\beta_w*\overline{\gamma_C}) = \canonLabel(\alpha_t*\overline{\gamma_B})X'\canonLabel(\beta_w*\overline{\gamma_C})\\
        \intertext{and plugging in the exchange relation, \Cref{eqn:FullMutationFormula}, yields}
        &= \canonLabel(\alpha_t*\overline{\gamma_B}){\bigg(}\frac{I_X}{N(X)}\cdot\sau^{\angleParity(\gamma_A*\alpha_v)}\Big(\canonLabel(\gamma_A*\alpha_v*\overline{X})\Big)\sau(\canonLabel(\beta_t*\gamma_C))\\
        & \hphantom{= \canonLabel(\alpha_t*\overline{\gamma_B})\bigg(+}+\frac{O_X}{N(X)}\cdot \sau^{\angleParity(\gamma_B*\alpha_t)}\Big(\canonLabel(\gamma_B*\alpha_t)\Big)\canonLabel(X*\beta_v*\gamma_D){\bigg)}\canonLabel(\beta_w*\overline{\gamma_C})
    \end{align*}
    with $I_X=\prod_{Y\in \In^c(X)}N(Y)$ and $O_X=\prod_{Z\in \Out^c(X)}N(Z)$. Next, we observe that the parity condition for decorated polygon implies that $\angleParity(\gamma_A)+\angleParity(\alpha_v)=\angleParity(\alpha_t)+\angleParity(\gamma_B)+\angleParity(\theta_u)$, which allows us to use the concatenation rule for canonical labels to further simplify the above expression to
    \begin{align*}
        N(X)\cdot m_{\alpha_t}(\triangleleft_{t_2}) &= I_X\cdot \canonLabel(\alpha_t*\overline{\gamma_B}*\theta_u*\gamma_A*\alpha_v*\overline{X})\sau^{\angleParity(\beta_t)+1}(\canonLabel(\gamma_C))\sau^{\angleParity(\beta_w)}(\canonLabel(\overline{\gamma_C}))\\
        & \hphantom{= \canonLabel(\alpha_t*\overline{\gamma_B})\bigg(} + O_X\cdot \canonLabel(\alpha_t*\overline{\gamma_B}*\gamma_B*\alpha_t)\canonLabel(X*\beta_v*\gamma_D)\sau^{\angleParity(\beta_w)}(\canonLabel(\overline{\gamma_C}))
    \end{align*}
    Now we already see $\blacktriangle_t$ appearing in the first summand. To simplify the second summand in a similar way, use again the mutation rule for angle colors to get $\angleParity(\beta_w) = \angleParity(\beta_v) + \angleParity(\gamma_D) + \angleParity(\theta_w) +1$, and plug this in to get
    \begin{align*}
        N(X)\cdot m_{\alpha_t}(\triangleleft_{t_2}) &= I_X\cdot\blacktriangle_t \sau^{\angleParity(\beta_t)+1}(\canonLabel(\gamma_C))\sau^{\angleParity(\beta_w)+1}(\canonLabel(\overline{\gamma_C}))+O_X\cdot N(\canonLabel(\gamma_B))\cdot\triangledown_t
    \end{align*}
    For the final simplification, we use that $\angleParity(\beta_w)=\angleParity(\beta_t)+\angleParity(\gamma_C)+1$ and \Cref{lemma:CanonLabelReverse}, to arrive at
    \begin{align}\label{eq:mutationAngleSum}
        N(X)\cdot m_{\alpha_t}(\triangleleft_{t_2}) = I_X\cdot N(\gamma_C)\cdot \blacktriangle_t + O_X \cdot N(\gamma_B)\cdot \triangledown_t
    \end{align}
    Since transition maps preserve symmetry and the angle vector spaces, this shows that the image of the left angle is valid. A similar computation shows the other angle is essentially the sum of $\blacktriangle_v$ and $\triangledown_v$ as needed.

    Finally, since the image of $X'$ can be written using $\triangleleft_{t_2}$ as a product of elements with central norm in $\seedAlgebra_Q$ it has central norm as well.\medskip
    
    To show that the map is an isomorphism recall that two mutations at the same node $k$ in an ST-compatible quiver $Q$ give us the original quiver after switching is applied at $k$, which we denoted by $\sigma_k(Q)=:Q''$. We want to prove that
    \begin{equation*}
        \mu_k^2(X_k)=\sigma(X_k)\,.
    \end{equation*}
    More precisely, this is saying that the map between the seed algebras induced by double mutation is exactly the isomorphism $\sigma_k:\seedAlgebra_{Q''}\to\seedAlgebra_Q$ from \Cref{prop:weavIsoSeedAlgebras}, which proves that the homomorphism induced by a single mutation at $k$ is actually an isomorphism.

    To prove the above equation we consider $X'':=\mu_k^2(X_k)$. We have
    \begin{align*}
        X''=I_{X'} ~\cdot~ & \sau^{\angleParity(\gamma_D*\alpha_w)}(\canonLabel(\gamma_D))\,X'^{-1}\,\sau(\canonLabel(\beta_u*\gamma_B))\\
        & +O_{X'}\cdot\sau^{\angleParity(\overline{\gamma_A}*\alpha_u)}(\canonLabel(\overline{\gamma_A}))\,\sigma(X')^{-1}\,\sau(\canonLabel(\beta_w*\overline{\gamma_C}))\,.
    \end{align*}
    We know that all the angles are fixed by $\sigma$ in $\seedAlgebra_{Q'}$, which implies that
    \begin{align*}
        \sau^{\angleParity(\overline{\gamma_A}*\alpha_u)}(\canonLabel(\overline{\gamma_A}))\,\sigma(X')^{-1}\,\sau(\canonLabel(\alpha_w*\overline{\gamma_D}))
    \end{align*}
    is fixed by $\sigma$. We use this to write
    \begin{equation}\label{eq:doubleMutationProofX''}
        X''=\sau^{\angleParity(\gamma_D*\alpha_w)}(\canonLabel(\gamma_D))\,X'^{-1}\,M
    \end{equation}
    with
    \begin{align*}
        M=I_{X'} ~\cdot~& \sau(\canonLabel(\beta_u*\gamma_B))\\
        &+O_{X'}\cdot \sau^{\angleParity(\alpha_u)+1}(\canonLabel(\gamma_A))\sau(\canonLabel(\alpha_w*\overline{\gamma_D})^{-1})\sau(\canonLabel(\beta_w*\overline{\gamma_C}))
    \end{align*}
    where we have already used the familiar rules for parities and canonical labels to simplify the expression. A very similar computation yields
    \begin{equation}\label{eq:doubleMutationProofX'}
        X'=W\,\sigma(X)^{-1}\,\sau(\canonLabel(\beta_v*\gamma_D))
    \end{equation}
    with
    \begin{align*}
        W=I_X\sau^{\angleParity(\gamma_A*\alpha_v)}(\canonLabel(\gamma_A))\sau^{\angleParity(\theta_w*\overline{\gamma_D}*\beta_v)}(\canonLabel( & \theta_w*\overline{\gamma_D}*\beta_v)^{-1})\tau(\canonLabel(\beta_t*\gamma_C))\\
        &+ O_X\sau^{\angleParity(\gamma_B*\alpha_t)}(\canonLabel(\gamma_B))\,.
    \end{align*}
    To finish the proof, note that $I_{X'}=O_X$ and $O_{X'}=I_X$ and using the same rules as before to see $M=W$. Plug \Cref{eq:doubleMutationProofX'} into \Cref{eq:doubleMutationProofX''} to conclude $X''=\sigma(X)$ as claimed.
\end{proof}

In the proof we have seen:
\begin{proposition}\label{prop:doubleMutationIsomorphism}
    For two mutations at the same node of an ST-compatible quiver $Q$ the concatenation of mutation maps is the isomorphism $\sigma_k:\seedAlgebra_{\sigma_k(Q)}\to\seedAlgebra_Q$.
\end{proposition}

\begin{remark}
    On the polygonal tiling our mutation rule can be seen as a clockwise flip. \Cref{prop:doubleMutationIsomorphism} implies that threefold mutation at the same mutable node $k$ can be understood as a counterclockwise flip. We will sometimes denote the corresponding mutation as $\overline{\mu_k}$. Observe that $\overline{\mu_k}(\sigma_k(Q))=\mu_k(Q)$.
\end{remark}

\subsection{Polygonal cluster algebras}
\begin{definition}\label{def:clusterAlgebra}
    Let $Q$ be an ST-compatible quiver. We call the variables associated to the nodes of any quiver obtained via a sequence of admissible mutations \keyword{cluster variables}. The \keyword{(polygonal) cluster algebra}, $\clusterAlgebra_Q$, is the subalgebra of $\seedAlgebra_Q$ generated by the image of all cluster variables under the mutation isomorphisms.
\end{definition}

\begin{theorem}\label{thm:mutationClusterAlgebraIsomorphism}
    Let $Q,Q'$ be ST-compatible quivers. If $Q'$ is obtained from $Q$ by a series of admissible mutations, then $\clusterAlgebra_Q$ is isomorphic to $\clusterAlgebra_{Q'}$.
\end{theorem}
\begin{proof}
    By definition we use the mutation isomorphisms to map all cluster variables to a fixed seed algebra. If we choose a different but mutation-equivalent quiver as the starting point of this process, the resulting cluster algebras are isomorphic to each other under the mutation isomorphism relating $\seedAlgebra_Q$ to $\seedAlgebra_{Q'}$.
\end{proof}

\begin{proposition}
If $Q\weaviso{}Q'$ are weaving isomorphic ST-compatible quivers, then $\clusterAlgebra_Q\cong\clusterAlgebra_{Q'}$. More precisely, the cluster algebras are isomorphic via the isomorphism of the seed algebras from \Cref{prop:weavIsoSeedAlgebras}. Thus the cluster algebra generated by any isomorphism class of ordered quivers is unique.
\end{proposition}
\begin{proof}
    In the case of an isomorphism, i.e. a permutation of the nodes, the statement is obvious.\medskip
    
    For the switch at $k$, we first observe that $\sigma_k$ commutes with any mutation at a node $j\neq k$, i.e. $\sigma_k\circ\mu_j=\mu_j\circ\sigma_k:\seedAlgebra_{\mu_j(Q)}\to\seedAlgebra_{\sigma_k(Q)}$. Consequently, the cluster variable $\mu_j(X_i)\in\seedAlgebra_Q$ gets mapped by $\sigma_k$ to $\sigma_k(\mu_j(X_i))=\mu_j(\sigma_k(X_i))\in\seedAlgebra_{\sigma_k(Q)}$ which is another cluster variable.
    
    For mutation at $k$ on the other hand, we find $\sigma_k\circ\mu_k=\overline{\mu_k}:\seedAlgebra_{\sigma_k(\mu_k(Q))}\to\seedAlgebra_{\sigma_k(Q)}$. This is clear if we think of $\mu_k$ as flipping the orientation of the arc labeled $k$ in corresponding decorated tiling. Similarly to above we find that $\sigma_k(\mu_k(X_i))$ is a cluster variable in $\seedAlgebra_{\sigma_k(Q)}$.

    Inductively, we find that the isomorphism $\sigma_k:\seedAlgebra_Q\to\seedAlgebra_{\sigma_k(Q)}$ induces a bijection of cluster variables.\medskip

    Finally, for weaving at $k$, we find that $w_k\circ\mu_j=\mu_j\circ w_k:\seedAlgebra_{\mu_j(Q)}\to\seedAlgebra_{w_k(Q)}$ for any $j$, including $k$. Just like above this means that the cluster variables in the respective seed algebras are in bijection via the map $w_k$.
\end{proof}

\begin{remark}
     Admissible mutations of mixed quivers also induce mutation isomorphisms between the corresponding seed algebras. Admissible mutations at small nodes do not change the angle skeleton by the same argument as given in the proof of \Cref{thm:seedalgebra_iso}. Then the mixed cluster algebra, $\clusterAlgebra_Q$, is defined in the same way as in \Cref{def:clusterAlgebra}.
\end{remark}

\subsection{Laurent Phenomenon}
Berenstein and Retakh generalized the Laurent phenomenon of commutative cluster algebras to noncommutative surface cluster algebras \cite[Theorem 3.30]{berenstein2018noncommutative}. We establish an analogous Laurent phenomenon for polygonal cluster algebras. 
\begin{theorem}\label{thm:laurent}
    Let $Q$ be a ST-compatible quiver and $\clusterAlgebra_{Q}$ the polygonal cluster algebra. Each cluster variable in  $\clusterAlgebra_{Q}$ is a noncommutative Laurent polynomial in the initial cluster variables, i.e. a sum of monomials in the initial seed algebra $\seedAlgebra_Q$. 
\end{theorem}
\begin{proof}
    We have seen that mutation is equivalent to the following addition formula for angles (\Cref{eq:mutationAngleSum}):
    \begin{equation*}
    N(X)m_{\alpha_t}(\triangleleft_t) =
        I_X\cdot N(\gamma_C)\cdot \blacktriangle_t + O_X\cdot N(\gamma_B)\cdot \triangledown_t
    \end{equation*}
    From this formula we see that the new variable $X'=\mu(X)$ being Laurent is equivalent to the right hand side of this equation being divisible by $N(X)$. To see this we expand the addition formula for angles as in the derivation of \Cref{eq:mutationAngleSum}
    \begin{align*}
       N(X) L(\alpha_t*\overline{\gamma_B})X' L(\beta_w* \overline{\gamma_C}) =& I_X\cdot N(\gamma_C)\cdot L(\alpha_t*\overline{\gamma_B})\sau^{\ell_1}(L(\theta_u*\gamma_A*\alpha_v*\overline{X}))\\& + O_X\cdot N(\gamma_B)\cdot L(X*\beta_v*\gamma_D*\theta_w)L(\beta_w*\overline{\gamma_C})
    \end{align*}
    Now left and right multiplication by $\tau(L(\alpha_t*\overline{\gamma_B})$ and $\tau( L(\beta_w* \overline{\gamma_C}))$ results $X'$ multiplied by norms which clearly factor out of the right hand side. As in the commutative Laurent phenomena we may assume $N(X), N(\gamma_B), N(\gamma_C)$ are coprime.

    We will now show that after any number of intermediate mutations, the mutation at $X$ still produces a sum of angles which is manifestly divisible by the norm of $X$. 

    There are two nontrivial possibilities for how the node at an admissible mutations can be related to the node $X$. Either it is in a polygon which is adjacent to $X$ or it is connected purely commutatively to $X$. 

    Let us first consider the purely commutative case. We assume without loss of generality that $X$ is connected with $k$ commutative arrows out to $Y$. 

    Now before mutation at $Y$ we have 
    $$N(X)m_{\alpha_t}(\triangleleft_t) =
        I_XN(\gamma_C)\blacktriangle_t + O'_XN(Y)^kN(\gamma_{B}) \triangledown_t$$
        where $O'_XN(Y)^k = O_X$.
    By \Cref{lem:MutationPreservesAngleColors}, mutation at $Y$ preserves the angle skeleton and thus the canonical angles at $X$, $\blacktriangle_t, \triangledown_t$. Assuming there are no cancellations of commutative arrows, we have the following formula for the new mutation at $X$
    \begin{equation}
    N(X)m_{\alpha_t}(\tilde{\triangleleft_t}) =
        I_XN(Y')^kN(\gamma_C)\blacktriangle_t + O'_XN(\Out(Y))^kN(\gamma_{B}) \triangledown_t
    \end{equation}
    where $\Out(Y)$ is the product of all variables connected out from $Y$.
    The magic of this formula is that the cancellation of commutative arrows simply corresponds to removing repeated norm factors from each summand of this equation. Since all norm factors are relatively prime as in the commutative Laurent phenomenon, this does not affect the remaining calculation. Now we will show that we can factor out a $N(X)$ from the right hand side. 
    
    First we notice using the mutation rule that $N(Y') = N(X)^k\Gamma+ N(Y)^{-1}N(\Out(Y))$ where $\Gamma$ is some element. Now our previous equation becomes 
        \begin{align*}
    N(X)m_{\alpha_t}(\tilde{\triangleleft_t}) &=
        I_XN(Y)^{-k}N(\Out(Y))^k N(\gamma_C)\blacktriangle_t +O'_XN(\Out(Y))^{k}N(\gamma_{B}) \triangledown_t + N(X)\Gamma'\\
        &=
        N(Y)^{-k}N(\Out(Y))^k \Big(I_X N(\gamma_C)\blacktriangle_t +O'_XN(Y)^{k}N(\gamma_{B}) \triangledown_t \Big) + N(X)\Gamma'\\
        &= N(Y)^{-k}N(\Out(Y))^k(N(X)m_{\alpha_t}(\triangleleft_t)) +N(X)\Gamma'
    \end{align*}
    where $\Gamma'$ is some other element. 

    Thus we see that $N(Y)^{k}m_{\alpha_t}(\tilde{\triangleleft}) $ and $N(X)m_{\alpha_t}(\tilde{\triangleleft})$ are each Laurent polynomials. This implies that $m_{\alpha_t}(\tilde{\triangleleft})$ is too since commutatively the numerators of $N(X)$ and $N(Y)$ are relatively prime. 

    Let us now consider $Y$  an edge of an adjacent polygon. Now mutation at $Y$ changes the canonical angles at $X$, but the rest of the argument remains the same. Assume without loss of generality that $Y$ lives in the $\gamma_B$ portion of the polygon in \Cref{fig:GeneratlMutation} (otherwise consider the angles based at $v$ instead). Let $k$ be the number of commutative arrows from $X$ to $ Y$, so that in total there are $2k+1$ arrows from $X $ to $ Y$.
    
    Break the paths $\gamma_B = \gamma_{B_1}*Y*\gamma_{B_2}$, $\gamma_A =  \gamma_{A_1}*\gamma_{A_2}$ where the angle indicator of $Y$ ends between $A_1$ and $A_2$ and we write $\gamma_E$, $\gamma_F$ for the two paths starting at the end of the angle indicator for the edge $Y$ in the other polygon it lives in, and we write $\Upsilon_t$ for the canonical angle of this polygon moved to $t$ along $\gamma_{B_2}$, see \Cref{fig:LaurentPhenomenaSetup}.
    \input{figureLaurentPhenomenaSetup}

    Now before mutation at $Y$ we have 
    \begin{equation}\label{eq:x_mutation_laurent}
        N(X)m_{\alpha_t}(\triangleleft_t) =
        I_X\cdot N(\gamma_C)\cdot \blacktriangle_t + O'_XN(Y)^k\cdot N(\gamma_{B_1})N(Y)N(\gamma_{B_2})\cdot \triangledown_t
    \end{equation}
    where $O'_XN(Y)^k = O_X$
    
    We claim that mutation at $Y$ produces $Y'$ and changes the angle $\blacktriangle_t$ to $\tilde{\blacktriangle}_t$ which satisfies that 
    \begin{equation}\label{eq:y_mutation_laurent}
        N(Y)\tilde{\blacktriangle}_t =  O_YN(\gamma_{E})\blacktriangle_t + I_YN(\gamma_{B_2})N(X)N(\gamma_{A_2})\Upsilon_t\,.
    \end{equation} 

    Now we compute the new mutation at $X$. Since the mutation at $Y$ is admissible, we know that we get a new ST-compatible quiver. If there is no cancellation of arrows then we would have the formula
    \begin{equation*}
        N(X)m_{\alpha_t}(\tilde{\triangleleft}) =
        N(Y')^k I_X N(\gamma_C)\tilde{\blacktriangle}_t + O'_XO_Y^{2k+1}N(\gamma_{B_1})^{k+1}N(\gamma_{B_2})N(\gamma_E)^{k+1}N(\gamma_{A_1})^{k}\triangledown_t\,.
    \end{equation*}  
    Again, the magic of this formulation of the mutation rule is that admissible cancellations do not change any of the angles which are being summed, but rather simply remove repeated norm factors between the two summands. Such a cancellation does not affect divisibility by $N(X)$ so we ignore it.

    Once again we only keep careful track of terms which do not obviously contain a factor of $N(X)$. First we notice that $N(Y')^k = (N(X)\Gamma+ N(\gamma_E)N(\gamma_{B_1})N(\gamma_{A_1})O_Y^2N(Y)^{-1})^k$ where $\Gamma$ is some central element so that our equation becomes:
    \begin{align*}
       N(X)&m_{\alpha_t}(\tilde{\triangleleft}) =  N(\gamma_E)^{k} N(\gamma_{B_1})^k N(\gamma_{A_1})^k O_Y^{2k} N(Y)^{-k} I_X N(\gamma_C)\tilde{\blacktriangle}_t \\ & \hspace*{4pc}+ O'_XO_Y^{2k+1}N(\gamma_{B_1})^{k+1}N(\gamma_{B_2})N(\gamma_E)^{k+1}N(\gamma_{A_1})^{k}\triangledown_t + N(X)\Gamma' \\
        =~&
        \left(N(\gamma_E)^{k} N(\gamma_{B_1})^k N(\gamma_{A_1})^k O_Y^{2k}N(Y)^{-k - 1}\right) \\
        &\hspace{1.5pc}\cdot \Big(I_X N(\gamma_C)N(Y)\tilde{\blacktriangle_t} + O'_XO_YN(Y)^{k+1}N(\gamma_E)N(\gamma_{B_1})N(\gamma_{B_2})\triangledown_t\Big) + N(X)\Gamma'\\
        =~&\left(N(\gamma_E)^{k} N(\gamma_{B_1})^k N(\gamma_{A_1})^k O_Y^{2k}N(Y)^{-k - 1}\right) \\
        &\hspace{1.5pc}\Big(I_X N(\gamma_C) O_YN(\gamma_{E})\blacktriangle_t + O'_XO_YN(Y)^{k+1}N(\gamma_E)N(\gamma_{B_1})N(\gamma_{B_2})\triangledown_t\Big) + N(X)\Gamma''\\
        =~& \left(N(\gamma_E)^{k+1} N(\gamma_{B_1})^k N(\gamma_{A_1})^k O_Y^{2k+1}N(Y)^{-k-1}\right)\big(N(X)m_{\alpha_t}(\triangleleft_t)\big) +N(X) \Gamma''
    \end{align*}
    Here the third equality uses \Cref{eq:y_mutation_laurent} and the final equality uses 
    \Cref{eq:x_mutation_laurent}.
  
    Thus we see that $N(Y)^{k+1}m_{\alpha_t}(\tilde{\triangleleft}) $ and $N(X)m_{\alpha_t}(\tilde{\triangleleft})$ are each Laurent polynomials. Again, this implies that $m_{\alpha_t}(\tilde{\triangleleft})$ is too since commutatively the numerators of $N(X)$ and $N(Y)$ are relatively prime. 
    
\end{proof}

\begin{remark}
    The Laurent phenomenon naturally extends to mixed 2-colored quivers. This follows from the commutative Laurent phenomenon for mutation at small nodes. When a small node is the intermediate mutation it operates exactly as a large node attached purely commutatively to $X$. The small nodes satisfy the Laurent phenomenon by the usual proof in the commutative case since the exchange relations are relatively prime. 
\end{remark}

\subsection{Evaluation in Clifford algebras}

Recalling \Cref{thm:Clifford_warmup} from the warm-up, we will now consider ring homomorphisms from our noncommutative algebra to Clifford algebras. We use the same set up as in the warm-up, i.e. $\CL(V,q)$ is the Clifford algebra over $(V,q)$, we fix a vector $e\in V$ with $q(e)\neq 0$ and define two anti-involutions $\sigma,\tau$ on $\CL(V,q)$ as follows 
\begin{equation}
    \tau(x)= x^T \hspace{1cm} \sigma(x) = \frac{1}{q(e)}ex^Te.
\end{equation}
\begin{remark}
    We suppress the dependency of $\sigma$ on $e$ in this notation. Moreover, by a slight abuse of mutation, we use both $\sigma$ and $\tau$ to denote an anti-involution on $\clusterAlgebra_Q$ and $\CL(V,q)$.
\end{remark}

\begin{definition}
    A \keyword{$\CL(V,q,e)$-point} of the algebra $\mathcal{A}_Q$ is an algebra homomorphism $\phi:\mathcal{A}_Q \to \CL(V,q)$ such that 
    \begin{enumerate}
        \item $\phi$ respects the two anti-involutions defined on $\mathcal{A}_Q$ and $\CL(V,q)$, i.e.
            \begin{equation*}
                \phi\circ\sigma=\sigma\circ\phi \hspace{1cm} \phi\circ\tau=\tau\circ\phi\,.
            \end{equation*}
        \item $\phi(X) \in \Gamma(V,q)$ for all cluster variables $X$.
        \item $\phi(\Delta) \in eV$ for all angles $\Delta$ in all clusters.
    \end{enumerate}
\end{definition}
Again we may define positive points of $\mathcal{A}_Q$ in exactly the same way as \Cref{def:positive_points}, i.e. pick $V$ to be an $\R$ vector space, $q$ to be a quadratic form of signature $(1,n-1)$, $e$ to be a vector of norm $+1$, and force that the angles lie in $V^+$, the proper cone of vectors with positive norm which contains $e$.

\begin{theorem}\label{thm:clifford_evaluations_general}
Let $Q$ be an ST-compatible quiver and suppose we have an algebra homomorphism
	\begin{equation*}
		\phi:\mathcal{S}_Q\to\CL(V,q)
	\end{equation*}
	which is such that
	\begin{enumerate}
		\item it respects $\sigma$ and $\tau$,
		\item $\phi(X_i)\in \Gamma(V,q)$ for any $i\in Q_0$,
		\item $\phi(\Delta_k(i))\in eV$ for any $k\in\{0,1\},i\in Q_0$.
	\end{enumerate}
	Then $\phi$ uniquely extends to a $\CL(V,q,e)-$point of $\mathcal{A}_Q$ provided that each new cluster variable is mapped to an invertible element of $\CL(V,q)$. The positive points of $\mathcal{A}_Q$ can be constructed simply by the condition that $\phi(\Delta_k(i))\in eV^+$, the invertibility is implied by positivity. 
\end{theorem} 
\begin{proof}
    The proof of this statement follows exactly from the same reasoning as \Cref{thm:Clifford_warmup}, using the proof of theorem \Cref{thm:seedalgebra_iso} to show that the new cluster variables and angles after mutation are in $\CL(V,q)$ and $eV$, respectively. 
\end{proof}

\begin{corollary}
    $\mathcal{A}_Q$ is a noncommutative ring.
\end{corollary}

\section{Examples}\label{sec:examples}
We will now collect some examples of quivers which are either fruitful or have interesting admissible mutation sequences. 

\subsection{Surface Type}

We obtain the algebra $\mathcal{B}_1(P_n)$ from the warm-up (\Cref{sec:AlgebraB1Pn}) as the polygonal cluster algebra associated to a tiling of an $n-$gon by standard decorated triangles.

More generally, let $\Sigma$ be a a \emph{marked surface}, i.e a topological surface with boundary components and marked points on the boundary. We will assume for simplicity now that $\Sigma$ does not have any punctures, but we will explain in the next section how to add these. 
Lets consider the following example:
\begin{example}
    Let $\Sigma$ be an annulus with one marked point on each boundary component. Then $\Sigma$ can be triangulated with two triangles glued in a square with one set of parallel sides identified. The universal cover of this surface is an infinite strip cut into triangles, (\Cref{fig:AnnulusCover}).
    When we consider the angle space based near any vertex of the angle skeleton, we find that each triangle in the universal cover contributes a different angle, and so the angle space is infinite dimensional. Thus we conclude that while the seed algebra itself is finitely generated, the ideal of centrality conditions is not. 

\begin{figure}[!hb]
    \centering
    \begin{tikzpicture}
        \draw[dashed] (0,0) to (1,0);
        \draw (1,0) to node[below] {$\cdots$} (2,0)  to node[below] {$f_2$} (4,0) to node[below] {$f_2$} (6,0) to[dashed] node[below] {$\cdots$} (7,0);
        \draw[dashed] (7,0) to (8,0);
        \draw[dashed] (0,2) to (1,2);
        \draw (1,2) to node[above] {$\cdots$} (2,2)  to node[above] {$f_1$} (4,2) to node[above] {$f_1$} (6,2) to[dashed] node[above] {$\cdots$} (7,2);
        \draw[dashed] (7,2) to (8,2);
        \draw (2,0) to node[left] {1} (2,2);
        \draw (4,0) to node[left] {1} (4,2);
        \draw (6,0) to node[left] {1} (6,2);

        \draw[dashed] (0,0) to (1,1); 
        \draw (1,1) to node[left] {} (2,2);
        \draw (2,0) to node[left] {2} (4,2);
        \draw (4,0) to node[left] {2} (6,2);
        \draw (6,0) to node[left] {} (7,1);
        \draw[dashed] (7,1) to (8,2);
    \end{tikzpicture}
    \caption{Triangulated universal cover of an annulus}
    \label{fig:AnnulusCover}
\end{figure}

\end{example}

Generally we can see that if $\pi_1(\Sigma)$ is not finite then the angle space is no longer finite dimensional. However, such quivers are clearly fruitful, as there are never any crossing arrows or commutative arrows created after mutation. 

In order to obtain more fruitful quivers, one might try to modify quivers of surface type. However, we will find that that is not really possible.
\begin{lemma}
	If a quiver is weaving-isomorphic to a quiver with only parallel arrows, it is of surface type.
\end{lemma}
\begin{proof}
	By \Cref{thm:GSTRestrictions} the assumption implies that the quiver cannot contain any pseudocycles. This in turn implies that all ST-tiles contain at most 3 elements, and can therefore be realized by a triangle with oriented sides. Gluing these triangles we obtain a triangulation of a surface for which the above construction yields the given seed.
\end{proof}

This prevents us from constructing new examples from surface type by simply recoloring the arrows:

\begin{lemma}
    Suppose $Q$ is an ST-compatible quiver which is a coloring of a surface type quiver with underlying surface $\Sigma$. If $\pi_1(\Sigma)$ is trivial, then $Q$ is weaving isomorphic to a quiver of surface type.
\end{lemma}
\begin{proof}
Since $\pi_1(\Sigma)$ is trivial the dual graph to the triangulation is a tree. Every crossing arrow can be removed by inductively weaving from the leaves.
\end{proof}

Passing to surfaces with nontrivial fundamental group, the quiver cannot always be uncrossed but there is a 2-fold cover which can.
\begin{lemma}
    Let $Q$ be an ST-compatible quiver which is a coloring of a surface type quiver. If the unfolded quiver $\bar{Q}$ has a single connected component it provides a 2-fold cover of $Q$ which is of surface type.
\end{lemma}
This corresponds to a double cover $\Sigma'$ of the underlying surface $\Sigma$, and it means that examples in which the noncommutative surface type cluster algebra on $\Sigma'$ is understood, are not new, since their mutation is equivalent to a noncommutative group mutation on $\bar{Q}$. Nevertheless there are many interesting fruitful seeds.

\subsection{Punctured surfaces and ST-quivers with 2-cycles}\label{sec:puncDigon}
We now show how to treat the situation where $\Sigma$ has punctures. We will see that folded ST-tiles naturally appear in this situation. 

Let $\Sigma$ be a digon with a single puncture. Then $\Sigma$ has an ideal triangulation consisting of two triangles which share two sides shown in \Cref{fig:digonA}. In the usual description of cluster algebras associated to surfaces the quiver associated to this triangulation has no arrows between the two nodes associated to the interior arcs. Mutation at one of the arcs results in a triangulation which can only be described with ``tagged arcs'' \cite{fomin2008clusterTri}. 

In \cite{kaufman2023special} folded quivers were used to give an alternate quiver associated to this triangulation. In this formulation, mutation of one node corresponds to flipping an arc to make a folded triangle shown in \Cref{fig:digonB}. This folded quiver can be equivalently described by the two colored quivers coming from the additional decorations of the tilings in \Cref{fig:digon}. 

\input{figureDigon}

Starting with the first quiver in \Cref{fig:digon} we assign variables $A,B$ to the left and right boundary arcs and variables $X,Y$ to the top and bottom arcs respectively.
Our mutation rule for ST-compatible seeds tells us how to perform the first mutation showed. We call this variable $Z$ and know that $Z$ satisfies
\begin{equation*}
    Z = \sigma(A)X^{-1}\tau(Y) + Y \sigma(X^{-1})B\,.
\end{equation*}
The new canonical angle near the puncture is given by $\Delta= \tau(Y)ZY$. We compute the monodromy around the puncture is given by $\sau$ which produces a new angle $\tau(\Delta)$ which must be compatible with the original, i.e. $b(\Delta,\tau(\Delta))$ is central.

We would now like to describe how to give a mutation rule for the interior arc of the self-folded tile. This mutation should represent the sum of angles $\Delta + \tau(\Delta)$ but also agree with the mutation in the commutative cluster algebra. We define the mutated variable $W$ to satisfy 
$$W = Y^{-1}(Z+\tau(Z))\,.$$ We can see that $W$ has central norm since $N(Y)WY = \Delta + \tau(\Delta)$ and the new canonical angle $\Delta' = WZ\tau(W) $ is symmetric since 
$$N(Y)^4\Delta' = N(Y)WY\left(\tau(Y)ZY\right)\tau(Y)\tau(W)N(Y) = (\Delta+\tau(\Delta))(\Delta)(\Delta+\tau(\Delta))$$

One checks that using the same mutation rule again at $W$ returns the variable $Y$. First notice that switching and weaving at $W$ changes the quiver at \Cref{fig:digonC} back to that of \Cref{fig:digonB} so than we can apply the exact same rule again.  This has the effect of replacing $W$ with $\tau(W)$. Now 
$$ \tau(W)^{-1}(Z+\tau(Z)) = \tau(Y) (Z+\tau(Z))^{-1}(Z+\tau(Z)) = \tau(Y) $$ so that after weaving and switching we get $Y$ back. We can then mutate again at $Z$ and continue on our way. 

In this way we can assign a polygonal cluster algebra to any punctured surface such that the flip of triangulations in a punctured digon is given by the above mutation rules. To do this we pick some ideal triangulation of $\Sigma$ and a decoration so that all the triangles are 3-cycles and such that that the monodromy of the angle skeleton around each puncture is $\sau$. 

\begin{remark}
    Commutatively, the variable $W$ is exactly twice the variable which is usually associated to the tagged arc from the bottom of the digon to the puncture. However, we can see that in our noncommutative situation, this factor of $2Z$ is replaced with $Z+\tau(Z)$. 
\end{remark}

\subsection{Examples of fruitful quivers}\label{subsec:Fruitful}
In general it is a hard problem to decide whether a given quiver is fruitful. For example, experimentally quivers corresponding to a decorated tiling with no extra commutative arrows are fruitful. However we are able to prove several large classes of quivers are fruitful. We will also see several classes of quivers with infinite families of admissible mutations. In particular the mutation sequences which realize the Somos-4 and 5 sequences are admissible. 

\subsubsection{Pseudocycle}
The easiest new example is the pseudocycle. So far, every mutable node has two nonempty sides. However, when one side of a node $X$ is empty we can impose the relation
\begin{equation*}
    X=\sigma(X)
\end{equation*}
to obtain a mutable node. Applying this to the pseudocycle, we obtain an example of a fruitful quiver. The mutation class only contains two distinct quivers up to recoloring and isomorphism: The initial pseudocycle, and the quiver obtain after mutation at the node at which it is based. Thus it is simple to check every mutation is admissible. 

This should not come as a surprise because the unfolded quiver is an orientation of an affine type $A$ Dynkin diagram, which is of finite mutation type.

\subsubsection{Quivers with three nodes}
We can expand the previous example:
\begin{proposition}
	Any ST-compatible quiver with only 3 nodes which corresponds to a decorated tiling of a disk is fruitful.
\end{proposition}

\begin{proof}
	The pruned quiver is either a 3-cycle, a pseudocycle or an orientation and coloring of the $A_3$ Dynkin diagram. In all of these cases we can simply check that any mutation is admissible, even after the addition of an arbitrary number of commutative arrows.
\end{proof}

\subsubsection{Quivers with finite forkless part}
In his dissertation, Warkentin introduced a class of quivers which he called forks \cite{warkentin2014exchange}. Recall that a quiver with a node $A$ is a \emph{fork with point of return} $A$ if
\begin{itemize}
    \item $Q$ is abundant, i.e. every arrows has weight at least 2
    \item $\In(A)$ and $\Out(A)$ are acyclic
    \item for every $B\in\In(A)$, $C\in\Out(A)$, we have $b_{BA},b_{AC}<b_{CB}$.
\end{itemize}
Here we used $\In(A)$ and $\Out(A)$ to denote the full subquivers of $Q$ containing all nodes which are connected to $A$ by arrows which come in, respectively go out, at $A$. The numbers $b_{BA}$ etc. are the entries of the exchange matrix of the quiver.

Such quivers are interesting since they are the roots of trees in the exchange graph, and Warkentin proposed a generalization of finite mutation-type quivers could be provided by quivers with \emph{finite forkless part}. It is reasonable to assume that this structure makes checking fruitfulness easier since this ultimately boils down to controlling cancellation between commutative and noncommutative arrows. In fact, we have the following.

\begin{lemma}\label{lemma:FFPadmissibleMut}
    Suppose that $Q$ is an ST-compatible quiver, which is also a fork with point of return $A$. Mutation at any node $B\neq A$ is admissible.
\end{lemma}

\begin{proof}
    Let $Q$ be a fork with point of return $A$. We know that $Q$ is an abundant quiver, that any 3-cycle goes through $A$, and finally that for every directed path $B\to A\to C$, there is an arrow $C\to B$.
    
    The only potential problem in the mutation procedure is the cancellation of commutative against noncommutative arrows, which does not always lead to an admissible coloring. Thus, there are two cases to consider: We introduce a commutative arrow, which cancels an already-existing arrow, or we introduce a noncommutative arrow which cancels against a commutative arrow. For contradiction, assume $B\neq A$ and that mutation at $B$ is not admissible.
    
    \input{figureForkProof}
    
    In the first case, we may assume that we have a 3-pseudocycle with nodes $X, Y, Z$, based at $X$, and that mutation at $B$ introduces enough commutative arrows to flip the direction of the arrows between $X$ and $Y$, which would lead to a pseudocycle with incorrect coloring. This situation is illustrated in \Cref{fig:ForkProofA}, where orange arrows have arbitrary weight. Observe that the arrow between $Z$ and $B$ has to be oriented as shown in the figure for this to be a non-admissible mutation. Since $X,Y,B$ forms a 3-cycle and $B\neq A$ by assumption, we know that $A\in\{X,Y\}$, which is a contradiction since neither of these nodes satisfies the third condition given above.

    In the second case, two nodes $X$ and $Y$ are only connected by commutative arrows, and they form a cycle with the node $B$ s.t. mutation at $B$ introduces a noncommutative arrow from $X$ to $Y$, as shown in \Cref{fig:ForkProofB}. This can be non-admissible as soon as more nodes are involved. However, by our assumption, $B\neq A$, and since $B,X,Y$ form a cycle, $A\in\{X,Y\}$. Together with the assumption on weights of arrows in the definition of forks, this implies that after mutation arrows go from $X$ to $Y$ (in the commutative setting), meaning that all commutative arrows from $Y$ to $X$ get canceled against commutative arrows introduced by mutation, which is always admissible.
\end{proof}

This means that for quivers with finite forkless part which are also ST-compatible, we only need to check the finite forkless part, to conclude that the quiver is fruitful. This allows us to constructively prove:

\begin{proposition}
    For every $n\in\mathbb{N}$ there is an ST-compatible quiver with $n$ nodes which has empty forkless part, i.e. every quiver in the mutation class is a coloring of a fork.
\end{proposition}

\Cref{fig:FFPex} shows the quivers which we construct in the for $n=4$ and $n=5$. The quivers are already forks with point of return $A$, and mutation at $A$ results in a quiver which is again a fork with point of return $A$, since all the arrows from $C_i$ to $B$ get reversed. Thus the forkless part is empty and to check fruitfulness one only has to check that the mutation at $A$ is admissible, which is clear.

\input{figureFFPex}

\begin{proof}
	For any number $n$ of nodes we can pick the quiver with 2-colored exchange matrix the $n\times n$-matrix:
	\begin{align*}
		C=\left(
		\begin{array}{c|c}
			\begin{array}{cc}
				0 & -(2+\varepsilon)\\
				2+\varepsilon & 0 
			\end{array}&
			\begin{array}{ccccc}
				& \cdots & 3+3\varepsilon & \cdots & \\
				& \cdots & -(4+5\varepsilon) & \cdots &
			\end{array}\\
			\hline
			\begin{array}{cc}
				\vdots & \vdots\\
				-(3+3\varepsilon)  & 4+5\varepsilon\\
				\vdots & \vdots
			\end{array}&
			\begin{array}{ccccc}
				0 & -(1+2\varepsilon) & -(1+2\varepsilon) & -(1+2\varepsilon) & \cdots\\
				1+2\varepsilon & 0 & -(1+2\varepsilon) & -(1+2\varepsilon) & \cdots\\
				1+2\varepsilon & 1+2\varepsilon & 0 & -(1+2\varepsilon) & \cdots\\
				&& \vdots && \hphantom{1+2\varepsilon}
			\end{array}
		\end{array}\right)\,.
	\end{align*}
	Then the same argument as above applies.
\end{proof}

\begin{corollary}\label{thm:large_fuitefull_forks}
    There are arbitrarily large fruitful ST-compatible quivers.
\end{corollary}

\subsection{Decorated tilings}\label{sec:tiledSurfaceCAs}
We would like to use the machinery of decorated tilings to define a large class of fruitful seeds for polygonal cluster algebras. However, this turns out to be a hard combinatorial problem, as we will discuss below. Thus, we only have the following conjecture for now.
\begin{definition}
    We call a decorated tiling \keyword{purely noncommutative} if it contains no commutative arrows.
\end{definition}

\begin{conjecture}
	\leavevmode
    \begin{enumerate}
        \item Any purely noncommutative decorated tiling of a polygon provides a fruitful seed.
        \item In fact, this is even true if we choose a decorated tiling of an arbitrary surface with marked points.
    \end{enumerate}
\end{conjecture}
\begin{remark}
    We know that the conjecture is false if one removes the assumption of the decorated tiling being purely noncommutative, as the following innocent-looking example shows:

    \input{figureNonExample}

    Starting with three decorated polygons which are glued and an extra choice of two commutative arrows, after only two flips, we find that the configuration of commutative arrows is no longer admissible, due to the commutative arrow in the bottom pentagon. Cancelling this commutative arrow would break the cyclic order of the pentagon and thus mutation at $B$ is not admissible.
\end{remark}

It turns out that in order to check whether the seeds from the conjecture are fruitful, it is sufficient to prove a statement about the corresponding commutative 2-colored quivers, which are significantly easier to deal with and are also more easily accessible via computers since the exchange matrix over the split complex numbers can be implemented. Namely, we have the following:

\begin{proposition}
    To prove part (1) of the conjecture it is sufficient to prove that the corresponding commutative quivers have no quiver in their mutation class for which the difference in the number of parallel and crossing arrows between two nodes is bigger than 1.
\end{proposition}

\begin{proof}
    The critical step in the mutation is the cancellation of commutative against noncommu\-tative arrows, which might not always be possible due to the order. Passing to the commutative version of the quiver, we can simply carry out the mutation and obtain a new quiver, which is not correctly colored anymore, i.e. it contains a pseudocycle with an even number of crossing arrows. One further mutation then allows us to have a pair of arrows of the same color between two nodes.
\end{proof}

\begin{remark}
    Observe that this is not a sufficient criterion for arbitrary quivers, since even in the surface case double arrows can occur due to the topology of the underlying surface.    
\end{remark}

\subsection{Del Pezzo Quivers}\label{subsec:delpezzz}
The Del Pezzo quivers are a family of quivers associated to the five toric Del Pezzo surfaces, $\CP^2, \Bl_1 \CP^2, \CP^1\times \CP^1, \Bl_1 (\CP^1\times \CP^1), \Bl_2(\CP^1\times \CP^1)$ (\Cref{fig:Del Pezzo Quivers}).  
\begin{figure}[ht]
    \centering
    \begin{subfigure}{.19\textwidth}
        \centering
        \begin{tikzcd}[cramped]
            & 1 \arrow[ddl]\arrow[ddl,shift left=2]\arrow[ddl,shift right=2]& \\
            & & \\
            3\arrow[rr]\arrow[rr,shift left = 2]\arrow[rr,shift right = 2] & & 2 \arrow[uul]\arrow[uul,shift left = 2]\arrow[uul,shift right = 2]
        \end{tikzcd}
        \caption{$Q_{\CP^2}$}
    \end{subfigure}
    \begin{subfigure}{.19\textwidth}
        \centering
        \begin{tikzcd}[cramped]
            & 1 \arrow[dl]\arrow[dr]& \\
            4 \arrow[rr,shift left = 1]\arrow[rr,shift right = 1]& & 2\arrow[dl]\arrow[dl,shift left = 2]\arrow[dl,shift right = 2]\\
             & 3\arrow[ul] \arrow[uu,shift left = 1]\arrow[uu,shift right = 1]& 
        \end{tikzcd}
        \caption{$Q_{\Bl_1\CP^2}$}
    \end{subfigure}
    \begin{subfigure}{.19\textwidth}
        \centering
        \begin{tikzcd}[cramped]
            & 1 \arrow[dl,shift left = 1]\arrow[dl,shift right = 1]& \\
            2 \arrow[dr,shift left = 1]\arrow[dr,shift right = 1]& & 3\arrow[ul,shift left = 2]\arrow[ul,shift right = 1]\\
             & 4 \arrow[ur,shift left = 1]\arrow[ur,shift right = 1]& 
        \end{tikzcd}
        \caption{$Q_{\CP^1\times \CP^1}$}
    \end{subfigure}
    \begin{subfigure}{.19\textwidth}
        \centering
        \begin{tikzcd}[cramped,column sep=.6pc]
            & 1 \arrow[dl]\arrow[rr]& & 5\arrow[llld]\arrow[ldd]& \\
            2 \arrow[drr,shift left = 1]\arrow[drr,shift right = 1]& & & & 4\arrow[ul]\arrow[ulll]\\
             & & 3 \arrow[urr,shift left = 1]\arrow[urr,shift right = 1]\arrow[uul]& &
        \end{tikzcd}
        \caption{$Q_{\Bl_1(\CP^1\times \CP^1)}$}
    \end{subfigure}
    \begin{subfigure}{.19\textwidth}
        \centering
        \begin{tikzcd}[cramped,column sep=2pc, row sep=1pc,crossing over clearance=1.0ex]
            & 1\ar[dl]\ar[from=dr] &\\
            2\ar[d] \ar[ddr] & & 6\ar[ll]\ar[from=ldd]\\
            3 \ar[rr,crossing over]\ar[from=uur,crossing over]\ar[rd]& & 5\ar[u]\ar[uul,crossing over]\\
            & 4 \ar[ru]& 
        \end{tikzcd}
        \caption{$Q_{\Bl_2(\CP^1\times \CP^1)}$}
    \end{subfigure}
    \caption{Del Pezzo Quivers}
    \label{fig:Del Pezzo Quivers}
\end{figure}
These quivers have connected cluster algebras with exotic Lagrangian tori \cite{pascaleff_wall-crossing_2020},\cite{albers_floer_2023} and $q-$Painleve equations \cite{bershtein_cluster_2018} and have been known in the cases of $Q_{\Bl_1\CP^2}$ and $Q_{\Bl_1(\CP^1\times \CP^1)}$ to generate the Somos-4 and Somos-5 sequences, see \cite{hone_integrality_2008}. Here we realize each Del Pezzo quiver as a polygonal cluster algebra. The associated surfaces in some cases can be related via a ``blow up'' operation analogous to the original Del Pezzo surfaces.
\begin{example}
    The quiver associated to $\CP^2$ is an oriented cycle of triple arrows. This is realized by a tiling of a twice punctured sphere with three digons with all angle indicators meeting at the same puncture. We then add a cycle of commutative edges bringing the weight to 3. See \Cref{fig:DelPezzoTiling1}.\\
    \begin{figure}[!hb]
        \centering
        \begin{tikzpicture}[very thick]
            \node[circle, minimum size=4cm, draw,dashed] (p) at (0,0) {};
			\draw (p.center) -- (p.north);
            \draw (p.center) -- (p.south west);
            \draw (p.center) -- (p.south east);
			\coordinate (m1) at ($(p.center)!0.5!(p.north)$);
   			\coordinate (m2) at ($(p.center)!0.5!(p.south west)$);
			\coordinate (m3) at ($(p.center)!0.5!(p.south east)$);

			\draw[angleindicator] (m1) -- (p.north west);
   			\draw[angleindicator] (m1) -- (p.north east);
      		\draw[angleindicator] (m2) -- (p.north west);
   			\draw[angleindicator] (m2) -- (p.south);
            \draw[angleindicator] (m3) -- (p.south);
   			\draw[angleindicator] (m3) -- (p.north east);
            \draw[carrow] (m1) to (m2);
            \draw[carrow] (m2) to (m3);
            \draw[carrow] (m3) to (m1);
            \centerarc[crossangle](p.north west)(250:380:0.5)
            \centerarc[crossangle](p.south)(20:160:0.5)
            \centerarc[crossangle](p.north east)(150:300:0.5)
        \end{tikzpicture}
        \caption{Decorated tiling for $Q_{\CP^2}$.}
        \label{fig:DelPezzoTiling1}
    \end{figure}
    The three angles corresponding to each digon at the vertex of the angle skeleton near the center of \Cref{fig:DelPezzoTiling1} where each edge is oriented toward the center are \[\canonAngle{1} = \tau(a_1)a_2\hspace{2pc} \canonAngle{2} = \tau(a_2)a_3 \hspace{2pc} \canonAngle{3}= \tau(a_3)a_1\]
    \begin{remark}
        If we calculate the exchange relation at $a_1$ we find $$a_1' = \sigma(a_2)N(a_2)\sigma(a_1^{-1})+(a_1^{-1})a_3N(a_3).$$ Equivalently, the expression
        \begin{equation*}
            \sigma(a_2)^{-1}a_1'a_3^{-1} = N(a_2)\sigma(a_1)^{-1}a_3^{-1}+\sigma(a_2)^{-1}a_1^{-1}N(a_3)
        \end{equation*}
        represents the exchange relation as a sum of angles around the puncture. If we compute the similar expression for the other initial exchanges, we find that the sum
        \begin{equation*}
            F = N(a_2)\sigma(a_1)^{-1}a_3^{-1}+N(a_3)\sigma(a_2)^{-1}a_1^{-1} + N(a_1)\sigma(a_3)^{-1}a_2^{-1}
        \end{equation*}
        represents the ``total angle sum'' around the puncture. This is a non commutative deformation of the potential described in \cite{albers_floer_2023}.
    \end{remark}
\end{example}
\begin{example}
    The next Del Pezzo surface we consider is $\Bl_1(\CP^2)$. The associated quiver is the quiver corresponding to the Somos-4 sequence. It can be obtained from tiling a disk with one marked point on the boundary and one puncture using three digons and a triangle (\Cref{fig:DelPezzoTilingSomos4}). We take all arcs oriented toward the center with boundary oriented clockwise.\\
    This surface has a nontrivial fundamental group generated by the loop around the boundary component. However the monodromy around this loop is $\sigma\tau$ and thus order 2. So the angle space is eight dimensional, with two generators related by $\sigma\tau$ from each of the four polygons.  Explicitly the angle space is generated by 
    \[\begin{aligned}
        \tau(a_1)a_2 & \hspace*{2pc}& \tau(a_2)a_3 &\hspace*{2pc} &\tau(a_3)a_4 &\hspace*{2pc} &\tau(a_4)f_1 \sigma\tau(a_1)\\
        \sigma(a_1)\sigma\tau(a_2) & & \sigma(a_2)\sigma\tau(a_3) &&\sigma(a_3)\sigma\tau(a_4) & &\sigma(a_4)\sigma\tau(f_1)a_1
    \end{aligned}  \]
    and so there are $\binom{8}{2}$ centrality conditions in addition to those that come from the initial tiling.
    \begin{figure}
        \centering
        \begin{tikzpicture}[very thick]
            \node[circle, minimum size=4cm, draw,dashed] (p) at (0,0) {};
            \node[circle, minimum size=1cm, draw] (b) at (0,0) {};
            \coordinate (mid) at ($(b.south)!0.5!(p.north)$);
			\node[circle,minimum size = 2.5cm, draw] (arc) at (mid) {};
            \node[circle, draw,fill] (punc) at (b.south) {};
            \draw (b.south) -- (p.south west);
            \draw (b.south) -- (p.south east);
            \coordinate (m1a) at (arc.west);
            \coordinate (m1b) at (arc.east);
   			\coordinate (m2) at ($(b.south)!0.5!(p.south west)$);
			\coordinate (m3) at ($(b.south)!0.5!(p.south east)$);
            \node[label=1,left] (v1) at (m1a) {};
            \node[label=2,left] (v1) at (m2) {};
            \node[label=3,right] (v1) at (m3) {};
            \node[label=4,right] (v1) at (m1b) {};
            
			\draw[angleindicator] (m1a) -- (p.west);
            \draw[angleindicator, relative=false] (m1a) to [out=60, in=150] ($(p.north east)!0.5!(b.north)$) to[out=330, in = 0]  (b.south);
   			\draw[angleindicator] (m1b) -- (p.east);
            \draw[angleindicator, relative=false] (m1b) to[out=120, in=30] ($(p.north west)!0.5!(b.north)$) to[out=210, in = 180]  (b.south);
            \draw[angleindicator] (b.north)  -- (p.north);
      		\draw[angleindicator] (m2) -- (p.west);
   			\draw[angleindicator] (m2) -- (p.south);
            \draw[angleindicator] (m3) -- (p.south);
   			\draw[angleindicator] (m3) -- (p.east);
            \draw[carrow] (m1b) to (m2);
            \draw[carrow] (m2) to (m3);
            \draw[carrow] (m3) to (m1a);
            \centerarc[crossangle](p.west)(290:430:0.5)
            \centerarc[crossangle](p.south)(20:160:0.5)
            \centerarc[crossangle](p.east)(110:250:0.5)
        \end{tikzpicture}
        \begin{tikzpicture}[very thick,relative=false]
            \node[circle, minimum size=4cm, draw] (p) at (0,0) {};
            \node[circle, draw,fill] (punc) at (0,0) {};

            \coordinate (m1) at ($(p.south)!0.80!(p.north east)$);
            \coordinate (m2) at ($(p.south)!0.33!(p.north east)$);
   			\coordinate (m3) at ($(p.south)!0.33!(p.north west)$);
			\coordinate (m4) at ($(p.south)!0.80!(p.north west)$);

            \node[label=1,left] at (m1) {};
            \node[label=2,right] at (m2) {};
            \node[label=3,left] at (m3) {};
            \node[label=4,right] at (m4) {};
            
            \draw[] (p.south) to[out=130,in=270] (m3) to[out=90,in=240] (punc);
            \draw[] (p.south) to[out=50,in=270] (m2) to[out=90,in=300] (punc);
            
            \draw[] (p.south) to[out=150,in=210] (m4) to[out=30,in=120] (punc.north);         
            \draw[] (p.south) to[out=30,in=330] (m1) to[out=150,in=60] (punc.north);

            \centerarc[crossangle](p.center)(300:420:0.5)
            \centerarc[crossangle](p.center)(240:300:0.5)
            \centerarc[crossangle](p.center)(120:240:0.5)

            \draw[angleindicator] (m2)  to[out=180,in=270] (punc.south);
            \draw[angleindicator] (m2)  to[out=0,in=0] (punc.east);
            \draw[angleindicator] (m3)  to[out=0, in=270] (punc.south);
            \draw[angleindicator] (m3)  to[out=180, in=180] (punc.west);
            \draw[angleindicator] (m1)  to[] (punc.east);
            \draw[angleindicator] (m4)  to[] (punc.west);
            \draw[angleindicator] (m4)  to[out=120, in=130] ($(punc.north)!0.89!(p.north east)$) to[out=310, in=0] (p.south);
            \draw[angleindicator] (m1)  to[out=60, in=50] ($(punc.north)!0.89!(p.north west)$) to[out=230, in=180] (p.south);
            \draw[angleindicator] (p.north) to (punc.north);

            \draw[carrow] (m3) to (m1);
            \draw[carrow] (m2) to (m3);
            \draw[carrow] (m4) to (m2);
        \end{tikzpicture}

        \caption{Decorated tiling for $Q_{\Bl_1\CP^2}$ (Somos-4).}
        \label{fig:DelPezzoTilingSomos4}
    \end{figure}
    
\begin{remark}
    The tiling for $\Bl_1(\CP^2)$ can be seen by ``blowing up'' the tiling for $\CP^2$ at one of the punctures mirroring the operation of surfaces. To blow up a decorated tiling choose a puncture $p$ and side $e$ connected to the puncture. The new tiling replaces $p$ with a boundary component $b$ with one marked point and the side $e$ with two sides $e',e''$ connected around $b$ homotopic to $e$ if the boundary was shrunk to a point. We see the first tiling in \Cref{fig:DelPezzoTilingSomos4} is obtained from this operation at the middle of \Cref{fig:DelPezzoTiling1}. 
\end{remark}
\end{example}

\begin{example}
    The quiver for $\CP^1\times \CP^1$ is an oriented 4 cycle of double edges. This is realized as a tiling of twice punctured torus with two squares (\Cref{fig:DelPezzoTilingCP1xCP1}). We remark that as a quiver this retains a mixed color two cycle between the opposite corners. As such mutation at any node will have to follow rules analogous to the rules for a punctured digon given in \Cref{sec:puncDigon}.
    \begin{figure}[hb]
        \centering
        \scalebox{.75}{
        \begin{tikzpicture}[very thick,relative=false]
            \node[circle, draw,fill] (t1f) at (0,4) {};
            \node[circle, draw] (t2o) at (4,4) {};
            \node[circle, draw,fill] (t3f) at (8,4) {};
            \node[circle, draw] (b1o) at (0,0) {};
            \node[circle, draw,fill] (b2f) at (4,0) {};
            \node[circle, draw] (b3o) at (8,0) {};

            \coordinate (m1) at ($(t2o)!0.50!(b2f)$);
            \coordinate (m2b) at ($(b1o)!0.50!(b2f)$);
            \coordinate (m2t) at ($(t2o)!0.50!(t3f)$);
            \coordinate (m3l) at ($(b1o)!0.50!(t1f)$);
            \coordinate (m3r) at ($(t3f)!0.50!(b3o)$);
            \coordinate (m4t) at ($(t1f)!0.50!(t2o)$);
            \coordinate (m4b) at ($(b2f)!0.50!(b3o)$);
            
            \draw (t1f) to node[above] {4} (t2o)  to node[above] {2} (t3f);   
            \draw (b1o) to node[below] {2} (b2f)  to node[below] {4} (b3o);   
            \draw (t1f) to node[left] {3} (b1o);
            \draw (t2o) to node[left] {1} (b2f);
            \draw (t3f) to node[right] {3} (b3o);

            \draw[angleindicator] (m1) to (t1f);
            \draw[angleindicator] (m2b) to (t1f);
            \draw[angleindicator] (m3l) to (t2o);
            \draw[angleindicator] (m4t) to (b1o);
            
            \draw[angleindicator] (m1) to (t3f);
            \draw[angleindicator] (m2t) to (b2f);
            \draw[angleindicator] (m3r) to (t2o);
            \draw[angleindicator] (m4b) to (t2o);

            \centerarc[crossangle](b2f)(90:180:0.5)
            \centerarc[crossangle](t2o)(270:360:0.5)
        \end{tikzpicture}
        \hspace{.5cm}
        \begin{tikzpicture}[]
            \pic[name=A,rotate=180] at (4,2) {NCnode};
            \node[vertex] () at (4,2) [label ={[xshift=3mm,yshift=3mm]$1$}]{};
            
            \pic[name=Bb,rotate=90] at (2,0) {NCnodeHalfR};
            \node[vertex] () at (2,0) [label ={[xshift=-5mm, yshift=-5mm]$2$}]{};
            
            \pic[name=Cl,rotate=0] at (0,2) {NCnodeHalfR};
            \node[vertex] () at (0,2) [label ={[xshift=-2mm, yshift=3mm]$3$}]{};
            
            \pic[name=Dt,rotate=270] at (2,4) {NCnodeHalfR};
            \node[vertex] () at (2,4) [label ={[xshift=-5mm,yshift=0mm]$4$}]{};
            
            \pic[name=Bt,rotate=90] at (6,4) {NCnodeHalfL};
            \node[vertex] () at (6,4) [label ={[xshift=5mm, yshift=0mm]$2$}]{};
            
            \pic[name=Cr,rotate=0] at (8,2) {NCnodeHalfL};
            \node[vertex] () at (8,2) [label ={[xshift=2mm, yshift=3mm]$3$}]{};
            
            \pic[name=Db,rotate=270] at (6,0) {NCnodeHalfL};
            \node[vertex] () at (6,0) [label ={[xshift=5mm,yshift=-5mm]$4$}]{};
            
            \draw[barrow] (Cl-circle.30) to[out=30,in=210] (Dt-circle.210);
            \draw[barrow] (Dt-circle.300) to[out=300,in=135] (A-circle.135);
            \draw[xarrow] (A-circle.210) to[out=210,in=60] (Bb-circle.60);
            \draw[barrow] (Bb-circle.150) to[out=150,in=330] (Cl-circle.330);
            \draw[barrow] (Dt-circle.330) to[in=30,out=330] (Bb-circle.30);
            \draw[barrow] (A-circle.240) to[in=300,out=240] (Cl-circle.300);
            
            \draw[xarrow] (A-circle.30) to[out=30,in=240] (Bt-circle.240);
            \draw[barrow] (Bt-circle.300) to[out=300,in=150] (Cr-circle.150);
            \draw[barrow] (Cr-circle.210) to[out=210,in=60] (Db-circle.60);
            \draw[barrow] (Db-circle.120) to[out=120,in=330] (A-circle.330);
            \draw[xarrow] (Cr-circle.240) to[out=240,in=300] (A-circle.300);
            \draw[xarrow] (Bt-circle.330) to[out=330,in=30] (Db-circle.30);

        \end{tikzpicture}}
        
        \caption{Decorated tiling for $Q_{\CP^1\times\CP^1}$.}
        \label{fig:DelPezzoTilingCP1xCP1}
    \end{figure}
\end{example}

\begin{example}
    The next example is $\Bl_1(\CP^1 \times \CP^1)$ which also corresponds to the Somos-5 quiver. In \Cref{fig:DelPezzoTilingSomos5} we see the tiled surface constructed by ``blowing up'' the tiling for  $\CP^1\times \CP^1$ at the open puncture along 1. 
    \begin{figure}[!hb]
        \centering
        \scalebox{0.8}{
        \begin{tikzpicture}[very thick,relative=false]
            \node[circle, draw,fill] (t1f) at (0,2) {};
            \node[circle, draw] (t2lo) at (2,4) {};
            \node[circle, draw] (t2ro) at (6,4) {};
            \node[circle, draw,fill] (t3f) at (8,2) {};
            \node[circle, draw] (b1o) at (2,0) {};
            \node[circle, draw,fill] (b2f) at (4,2) {};
            \node[circle, draw] (b3o) at (6,0) {};

            \coordinate (m1) at ($(t2ro)!0.50!(b2f)$);
            \coordinate (m2b) at ($(b1o)!0.50!(b2f)$);
            \coordinate (m2t) at ($(t2ro)!0.50!(t3f)$);
            \coordinate (m3l) at ($(b1o)!0.50!(t1f)$);
            \coordinate (m3r) at ($(t3f)!0.50!(b3o)$);
            \coordinate (m4t) at ($(t1f)!0.50!(t2lo)$);
            \coordinate (m4b) at ($(b2f)!0.50!(b3o)$);
            \coordinate (m5) at ($(t2lo)!0.50!(b2f)$);
            \coordinate (mf) at ($(t2lo)!0.50!(t2ro)$);
            
            \draw (t1f) to node[above] {4} (t2lo) to node[above] {f} (t2ro) to node[above] {2} (t3f);   
            \draw (b1o) to node[below] {2} (b2f)  to node[below] {4} (b3o);   
            \draw (t1f) to node[left] {3} (b1o);
            \draw (t2ro) to node[left] {1} (b2f);
            \draw (t2lo) to node[right] {5} (b2f);
            \draw (t3f) to node[right] {3} (b3o);

            \draw[angleindicator] (m5) to (t1f);
            \draw[angleindicator] (m2b) to (t1f);
            \draw[angleindicator] (m3l) to (t2lo);
            \draw[angleindicator] (m4t) to (b1o);
            
            \draw[angleindicator] (m1) to (t3f);
            \draw[angleindicator] (m2t) to (b2f);
            \draw[angleindicator] (m3r) to (t2ro);
            \draw[angleindicator] (m4b) to (t2ro);

            \draw[angleindicator] (m1) to (t2lo);
            \draw[angleindicator] (mf) to (b2f);
            \draw[angleindicator] (m5) to (t2ro);

            \centerarc[crossangle](b2f)(135:225:0.5)
            \centerarc[crossangle](t2ro)(225:315:0.5)
        \end{tikzpicture}
        \hspace{1cm}
        \begin{tikzpicture}[very thick, relative=false]
            \node[circle, draw,fill] (A) at (0,0) {};
            \node[circle, draw,fill] (B) at (0,4) {};
            \node[circle, draw,fill] (C) at (4,4) {};
            \node[circle, draw,fill] (D) at (4,0) {};
            \node[circle, draw,minimum size=1cm] (bound) at (2,2) {};
            \node[circle, draw,fill=white] (E) at (bound.south west) {};

            \draw[dashed] (A) to (B) to (C) to (D) to (A);
            
            \coordinate (m1) at (2,3.25);
            \coordinate (m2) at (0.45,3);
            \coordinate (m3) at ($(A)!0.50!(E)$);
            \coordinate (m4) at (3,0.3);
            \coordinate (m5) at (3.25,2);
            \coordinate (mf) at (bound.north east);
            \draw (C) to[out=210,in=30] node[below] {1} (m1) to[out=180+30,in=180] (E);
            \draw (C) to[out=240,in=60] node[left] {5} (m5) to[out=180+60,in=270] (E);
            \draw (A) to node[left] {3} (E);
            \draw (B) to[out=310,in=110] node[right] {2} (m2) to[out=290,in=210] (E);
            \draw (D) to[out=160,in=330] node[above] {4} (m4) to[out=150,in=240] (E);

            \draw[angleindicator] (mf) to (C);
            \draw[angleindicator] (m1) to[out=0,in=90] (3,3) to[out=270, in=300] (E.south east);
            \draw[angleindicator] (m1) to[out=180, in=330] (B);
            \draw[angleindicator] (m2) to[out=30,in=190] (C);
            \draw[angleindicator] (m2) to[] (A);
            \draw[angleindicator] (m3) to[out=135,in=0] ($(A)!0.25!(B)$); \draw[angleindicator] ($(D)!0.25!(C)$) to[out=180,in=260] (E);
            \draw[angleindicator] (m3) to[out=360-45,in=90] ($(A)!0.25!(D)$); \draw[angleindicator] ($(B)!0.25!(C)$) to[out=270,in=190] (E);
            \draw[angleindicator] (m4) to[out=45,in=180] ($(C)!0.5!(D)$); \draw[angleindicator] ($(B)!0.5!(A)$) to[out=0,in=217] (E);
             \draw[angleindicator] (m4) to[out=180+45,in=90] ($(A)!0.4!(D)$); \draw[angleindicator] ($(B)!0.4!(C)$) to[out=270,in=180] (E);
            \draw[angleindicator] (m5) to[out=90,in=0] (3,3) to[out=180, in=120] (E.north west);
            \draw[angleindicator] (m5) to[out=270,in=150] (D);

            \centerarc[crossangle](B)(270:300:0.8);\centerarc[crossangle](C)(250:270:1);
            \centerarc[crossangle](E)(105:150:1)
        \end{tikzpicture}}
        
        \caption{Decorated tiling for $Q_{\Bl_1(\CP^1\times\CP^1)}$ (Somos-5).}
        \label{fig:DelPezzoTilingSomos5}
    \end{figure}
    
\end{example}

\begin{example}
    The final Del Pezzo quiver, $\Bl_2(\CP^1\times \CP^1)$, is a tiling of a torus with two boundary components by two squares and two triangles (\Cref{fig:DelPezzoTilingB2CP1xCP1}). It can be obtained by blowing up the remaining puncture of the tiling for Somos-5 along the edge labeled 3 in \Cref{fig:DelPezzoTilingSomos5}. This leaves a mixed color two-cycle between edges 2 and 5. Nevertheless mutation at the other antipodal pairs $(1,4)$ or $(3,6)$ incident to the boundary edges returns to an isomorphic tiling. There is the same subtly as in Somos-5 where a two cycle of commuting arrows must be canceled ``across the topology'' to return the angle indicators in the squares exactly. 
    \begin{figure}[!htb]
        \centering
        \scalebox{.8}{
        \begin{tikzpicture}[very thick, relative=false]
            \node[circle, draw] (t1) at (0,4) {};
            \node[circle, draw] (t2) at (4,4) {};
            \node[circle, draw,fill] (t3) at (6,4) {};
            \node[circle, draw,fill] (t4) at (10,4) {};
            \node[circle, draw,] (t5) at (12,4) {};
            \node[circle, draw,fill] (b1) at (2,0) {};
            \node[circle, draw] (b2) at (8,0) {};
            \node[circle, draw,fill] (b3) at (14,0) {};

            \coordinate (mf) at ($(t1)!0.50!(t2)$);
            \coordinate (m2t) at ($(t2)!0.50!(t3)$);
            \coordinate (mg) at ($(t3)!0.50!(t4)$);
            \coordinate (m5t) at ($(t4)!0.50!(t5)$);
            \coordinate (m6l) at ($(t1)!0.50!(b1)$);
            \coordinate (m1) at ($(t2)!0.50!(b1)$);
            \coordinate (m3) at ($(t3)!0.50!(b2)$);
            \coordinate (m4) at ($(t4)!0.50!(b2)$);
            \coordinate (m6r) at ($(t5)!0.50!(b3)$);
            \coordinate (m5b) at ($(b1)!0.50!(b2)$);
            \coordinate (m2b) at ($(b2)!0.50!(b3)$);
            \draw (t1) to node[above] {f} (t2) to node[above] {2} (t3) to node[above] {g} (t4) to node[above] {5} (t5);
            \draw (t1) to node[left] {6} (b1);
            \draw (t2) to node[right] {1} (b1);
            \draw (t3) to node[left] {3} (b2);
            \draw (t4) to node[right] {4} (b2);
            \draw (t5) to node[right] {6} (b3);
            \draw (b1) to node[below] {5} (b2) to node[below] {2} (b3);
            \draw[angleindicator] (mf) to (b1);
            \draw[angleindicator] (m6l) to (t2);
            \draw[angleindicator] (m1) to (t1);
            \draw[angleindicator] (m1) to (t3);
            \draw[angleindicator] (m2t) to (b1);
            \draw[angleindicator] (m5b) to (t2);
            \draw[angleindicator] (m3) to (t2);
            \draw[angleindicator] (m3) to (t4);
            \draw[angleindicator] (mg) to (b2);
            \draw[angleindicator] (m4) to (t3);
            \draw[angleindicator] (m4) to (t5);
            \draw[angleindicator] (m2b) to (t4);
            \draw[angleindicator] (m5t) to (b2);
            \draw[angleindicator] (m6r) to (t4);
            \centerarc[crossangle](b3)(116:180:0.5);
            \centerarc[crossangle](t2)(244:360:0.5);
        \end{tikzpicture}}
        \caption{Decorated tiling for $Q_{\Bl_2(\CP^!\times \CP^1)}$.}
        \label{fig:DelPezzoTilingB2CP1xCP1}
    \end{figure}
\end{example}

\subsection{Somos Sequences}\label{sec:somos}
The classical cluster algebras associated to Del Pezzo quivers associated to $\Bl_1(\CP^2)$ and $\Bl_1(\CP^1\times\CP^1)$ can produce the Somos-4 and Somos-5 sequences. This is seen by writing the sequence of cluster variables found along the mutation sequence obtained by mutating nodes $1,2,3,\dots$ and evaluating the initial cluster variables at 1. 

The Somos-4 sequence is generated by $\Bl_1(\CP^2)$. We see that mutation at the right (or left) most arc of the tiling obtains the same tiling after relabeling the arcs. In this way we obtain a noncommutative deformation of the Somos-4 recurrence
    \begin{equation*}
        a_5 = a_4 a_1^{-1} a_2 + a_3 \tau(a_3) f_1 \sigma(a_1)^{-1}\,.
    \end{equation*}
Although we do not know this quiver is fruitful, the mutations needed to generate the Somos-4 sequence are admissible. For example we chose initial values in the space of $2\times 2$-matrices. Here $\sigma$ is matrix transpose and $\tau$ is the adjoint map. The usual Somos sequence is recovered by taking all initial and frozen cluster variables to be the identity matrix. Then every cluster variables along the mutation sequence has the form $a_i=\lambda_i I$ for $\lambda_i$ the $i^{\text{th}}$ term of the usual Somos sequence. However there are other possible initial integer matrices. For example
    \[ a_1= \begin{pmatrix}
 1 & 0 \\
 0 & 1 \\
\end{pmatrix}
\hspace{1pc} 
a_2= 
\begin{pmatrix}
 5 & -2 \\
 -2 & 1 \\
\end{pmatrix}
\hspace{1pc}
a_3=
\begin{pmatrix}
 0 & -1 \\
 1 & 3 \\
\end{pmatrix}
\hspace{1pc}
a_4= 
\begin{pmatrix}
 2 & -1 \\
 -1 & 1 \\
\end{pmatrix}
 \]
    One can easily check that all the angles are symmetric and positive definite in this case. Checking centrality is trivial as every nonzero element has central norm in this case. The sequence then continues:
    \[a_5=
\begin{pmatrix}
 13 & -5 \\
 -7 & 5 \\
\end{pmatrix}
\hspace{1pc}
a_6 = 
\begin{pmatrix}
 2 & 2 \\
 8 & 22 \\
\end{pmatrix}
\hspace{1pc}
a_7=
\begin{pmatrix}
 57 & -19 \\
 13 & 6 \\
\end{pmatrix}
\hspace{1pc}
a_8= 
\begin{pmatrix}
 3889 & -86 \\
 100 & 61 \\
\end{pmatrix}
    \]
One can even find families of initial matrices in a free parameter $t$. When these are chosen so that specializing $t$ to $0$ recovers four identity matrices, we consider these families to be noncommutative deformations of the classic Somos sequence. This is analogous to the noncommutative deformations of Markov numbers described in \cite{GreenbergEtAl2024MathrmSL_2} One example of such a family is
\begin{equation}\label{eqn:Somos4SymFamily}
a_1=
 \begin{pmatrix}
  1 & 0 \\
 0 & 1 \\
\end{pmatrix}
 \hspace{1pc}
a_2= 
 \begin{pmatrix}
  1 & 0 \\
 0 & 1 \\
\end{pmatrix}
 \hspace{1pc}
 a_3=
 \begin{pmatrix}
  1 & -t \\
 -t & t^2+1 \\
\end{pmatrix}
\hspace{1pc}
 a_4= 
 \begin{pmatrix}
  t^2+1 & t \\
 t & 1 \\
\end{pmatrix}
\end{equation}
 which continues 
 \[
 \begin{pmatrix}
  t^2+2 & t \\
 t & 2 \\
\end{pmatrix},
 \begin{pmatrix}
  3 & -t \\
 -t & t^2+3 \\
\end{pmatrix},
\begin{pmatrix}
  3 t^4+12 t^2+7 & 3 t ^3+9t \\
 3 t^3+9t & 3 t^2+7 \\
\end{pmatrix},
\begin{pmatrix}
  3 t^4+17 t^2+23 & t^3+2t \\
 t^3+2t & 2 t^4+15 t^2+23 \\
\end{pmatrix} ,\dots
\]

The Somos-5 quiver is given by $\Bl_1(\CP^1 \times \CP^1)$. Here mutation at 1 obtains a quiver which is weaving isomorphic to the original quiver after weaving at 2 and then cyclically reindexing the arcs. There is small subtly due to the topology of the surface here. The resulting quiver is actually isomorphic as ordinary quivers, but not as ST-compatible quivers. The problem is that the arrows connecting side $3$ and $5$ (before relabeling) belong to the ``wrong'' tiles. This can be fixed by adding a a two cycle of commutative arrows between $3$ and $5$ which when taken to live in separate square tiles modify the angle indicators to be exactly the same before and after mutation.\\
    The mutation at $1$ corresponds the following noncommutative recurrence generalizing Somos-5 where we assume every arc is oriented from $\bullet$ to $\circ$ and $f$ is oriented from 1 to 5:
    \begin{equation*}
        (a_1,a_2,a_3,a_4,a_5; ~f) \mapsto (\sau(a_2),a_3,a_4,a_5, \sau(a_2)a_1^{-1}a_5 + \sau(a_3)\sigma(a_4)\sigma(a_1)^{-1}f;~f)
    \end{equation*}
As with Somos-4 we obtain noncommutative deformations in a free parameter $t$. One such family mirrors \Cref{eqn:Somos4SymFamily} as follows: 
\[a_1=a_2=a_3= \begin{pmatrix}
    1 &0 \\ 0 & 1
\end{pmatrix} \hspace{1pc} a_4 = \begin{pmatrix}
    1 &-t \\ -t & 1+t^2
\end{pmatrix} \hspace{1pc} a_5 = \begin{pmatrix}
    1+t^2 & t \\ t & 1
\end{pmatrix}\]
This continues 
\begin{align*}
a_6=&\begin{pmatrix}
  t^2+2 & 0 \\
 0 & t^2+2 \\
\end{pmatrix}
 \hspace{1pc}
a_7=\begin{pmatrix}
  t^4+4 t^2+3 & t^3+2 t \\
 t^3+2 t & 2 t^2+3 \\
\end{pmatrix}
 \hspace{1pc}
a_8=\begin{pmatrix}
  t^6+6 t^4+10 t^2+5 & t^5+4 t^3+3 t \\
 t^5+4 t^3+3 t & 2 t^4+7 t^2+5 \\
\end{pmatrix}
\end{align*}
To illustrate the breadth of possible initial families we give another deformation with nonsymmetric matrices:
\begin{align*}
    a_1=& \begin{pmatrix}1 &0 \\ 0 & 1 \end{pmatrix} \hspace{1pc}
    a_2= \begin{pmatrix}1 &0 \\ 0 & 1 \end{pmatrix} \hspace{1pc}
    a_3= \begin{pmatrix}1 &0 \\ -t & 1 \end{pmatrix} \hspace{1pc}
    a_4= \begin{pmatrix}1 &t \\ 0 & 1 \end{pmatrix} \hspace{1pc}
    a_5= \begin{pmatrix}1 &0 \\ 0 & 1 \end{pmatrix} 
\end{align*}
\begin{align*}
    a_6=& \begin{pmatrix}t^2+2 &t \\t  & 2 \end{pmatrix} \hspace{1pc}
    a_7= \begin{pmatrix}2t^2+3 & 3t\\0  & 3 \end{pmatrix} \hspace{1pc}
    a_8= \begin{pmatrix}2t^2+5 &t \\-2t^3-4t  &5  \end{pmatrix} \hspace{1pc}
    a_9= \begin{pmatrix}3t^2+11 & t^3+4t\\ t^3+4t& t^4+7t^2+11 \end{pmatrix} \hspace{1pc}
\end{align*}
Note that at $t=0$ both sequence recover the classic Somos-5 sequence $1,1,1,1,1,2,3,5,11,\dots$ multiplied by the identity matrix as needed.
\section{Polygonal cluster algebras in higher Teichm\"uller theory}\label{sec:TM}

We explained in \Cref{sec:introTM} that this work is motivated by the study of decorated representations into groups of $\Theta$-type $B_n$. Consider again \Cref{fig:FGQuivers} which shows the quivers that are inscribed into the triangles of an ideal triangulation of a punctured surface $\Sigma$ in order to describe the cluster structure on the space of decorated representations of $\pi_1(\Sigma)$ into $G=\Spin(2,3),\Spin(3,4)$, i.e. the split real forms of type $B_2$ and $B_3$. These give the cluster structure on the configuration space of transverse triples of flags in $\mathcal{F}=G/U$, where $U$ is the unipotent radical of a minimal parabolic subgroup of $G$. Our first goal is to endow these quivers with the structure of ordered mixed quivers. This will describe the cluster structure of the configuration space of transverse triples of flags in $\mathcal{F}_\Theta$, the flag variety corresponding to the $\Theta$-positive structure of type $B_2$ and $B_3$. It should also underlie the cluster structure for decorated $\Theta$-positive representations into the corresponding groups.

\subsection{Ordered quivers}
Let us redraw the quivers in \Cref{fig:FGQuivers} in such a way that they are manifestly ST-compatible ordered mixed quivers. This is shown in \Cref{fig:orderedFGQuivers}.

\input{figureOrderedFGQuivers}

We have drawn arrows between big (noncommutative) and small (commutative) nodes in green for better visibility and to emphasize their commutative nature. These quivers are clearly ST-compatible. In fact, the pruned quivers correspond to the decorated tilings of the triangle / disk with three marked points which are given in \Cref{fig:FGTilings}.

\input{figureFGTilings}

Of course these are the elements $Q_2,Q_3$ of an infinite family $(Q_n)_{n\in\mathbb{N}}$ of quivers describing the cluster structure for decorated representations into the split real groups $\Spin(n,n+1)$ as shown in \cite{le2019cluster}. The tilings associated to the ordered structure on all the quivers are obtained in the same way, by gluing $n-1$ digons to the regular decorated triangle as shown in \Cref{fig:FGTilings}.

We have not specified the orientation of the arcs here, i.e. we have only described the quivers up to switching isomorphism. Additionally, any weaving isomorphic quiver would also define the structure of an ST-compatible quiver on these quivers.

\begin{remark}
    Observe that these pruned quivers are initially of surface type since we could inflate the digons to regular decorated triangles. We will see that this is not stable under mutation for the second quiver, and in particular not once we start gluing these quivers.
\end{remark}

\subsection{Fruitfulness}
The remarkable thing about the quiver in \Cref{fig:orderedFGQuiverB2} is that it is of finite type (the nodes on the edges of the underlying triangle are frozen), and thus one can check explicitly that all mutations are admissible.
\begin{proposition}
    The ST-compatible quiver $Q_2$ in \Cref{fig:orderedFGQuiverB2} is fruitful.
\end{proposition}

\begin{proof}
    This can be checked by hand, see \Cref{sec:B2quiver}. In fact, the corresponding polygonal cluster algebra is of finite type and we give an explicit noncommutative Laurent expression for every cluster variable.
\end{proof}

We expect this to hold much more generally:
\begin{conjecture}
	The following quivers provide fruitful seeds:
	\begin{enumerate}
		\item the quivers $Q_n$ constructed in \cite{le2019cluster} with the structure of an ordered quiver described above,
		\item any quiver obtained by gluing these according to an ideal triangulation of a bordered surface.
	\end{enumerate}
\end{conjecture}

The associated polygonal cluster algebras should capture the cluster structure for configurations of transverse triples of flags in $\mathcal{F}_\Theta$ and decorated representations for groups of $\Theta$-type $B_n$.

\subsection{Big polygons}\label{sec:bigPolygons}
In general the quivers $Q_n$ and those quivers obtained by gluing them according to a triangulation will be equivalent via admissible mutations to quivers which contain more complicated ST-tiles than just bigons and regular triangles. We will illustrate this on two examples:

First, consider again the ordered quiver $Q_3$ in \Cref{fig:orderedFGQuiverB3}. Mutation at the numbered nodes in the indicated order immediately gives us a quiver which contains a pseudocycle (consisting of the node labeled 1 and the two adjacent big nodes) in its pruned part. This shows that the pseudocycle shows up naturally in this construction.

Let us also consider what happens when we glue copies of the quiver $Q_2$ according to a triangulation. In \Cref{fig:bigPolygonSetup}, we see three copies of $Q_2$ glued into a disk with five marked points in the boundary after mutation at the big nodes in the left and right triangle.

\input{figureBigPolygonSetup}

Further mutation at numbered nodes in the indicated order yields the quiver in \Cref{fig:bigPolygonEx} which contains an ST-tile, indicated by the colored nodes, that corresponds to a decorated pentagon.

\input{figureBigPolygonEx}

\begin{remark}
    This process does not immediately yield larger ST-tiles. However, we expect the ST-compatible quiver obtained by gluing copies of $Q_2$ into a disk with a sufficiently large number of marked points in the boundary to be equivalent via admissible mutations to a quiver which obtains an ST-tile that corresponds to a decorated $n$-gon for arbitrarily large $n$.
\end{remark}

\appendix

\section{Clifford algebras and Clifford groups}\label{app:Clifford}
Let $\mathbb{F}$ be a field, and $V$ a $\mathbb{F}$-vector space endowed with a nondegenerate symmetric bilinear form $b$ with induced quadratic form $q$.
\begin{definition}
	The \keyword{Clifford algebra} $\CL(V,q)$ is the unital $\mathbb{F}$-algebra generated by $V$ subject to the relation $v^2=q(v)\cdot 1$ for any $v\in V$.
\end{definition}
The original vector space $V$ embeds as a subspace into the Clifford algebra and it is customary to denote the image of $v\in V$ under the inclusion by $v$ again. We collect some basic facts about the Clifford algebra:
\begin{proposition}\
	\begin{itemize}
		\item $vw+wv=2b(v,w)$ for any $v,w\in V$.
		\item $\dim\CL(V,q)=2^{\,\dim V}$.
		\item $\CL(V,q)$ carries a $\mathbb{Z}/2\mathbb{Z}$-grading, i.e. $\CL(V,q)=\CL(V,q)^{(0)}\oplus\CL(V,q)^{(1)}$.
	\end{itemize}
\end{proposition}
The second and third property follow from the following observation: If $\{e_1,\dots,e_n\}$ is a basis of $V$, then $\left\{\prod_{i=1}^n e_i^{k_i}\big|k_i\in\{0,1\}\right\}$ is a basis of $\CL(V,q)$, and the grading comes from distinguishing between odd and even products of vectors.

There are a some interesting involutions on $\CL(V,q)$. We define them on the basis and extend them linearly:
\begin{itemize}
	\item $\alpha(e_{i_1}\cdots e_{i_k}):=(-1)^k e_{i_1}\cdots e_{i_k}$,
	\item $\tau(e_{i_1}\cdots e_{i_k}):=(e_{i_1}\cdots e_{i_k})^T = e_{i_k}\cdots e_{i_1}$.
\end{itemize}
Notice that the map $\alpha$ gives us the grading. Namely, $\CL(V,q)^{(0)}$ is the $+1$-eigenspace of $\alpha$, while $\CL(V,q)^{(1)}$ is the $-1$-eigenspace.

Next, we define a special subgroup of the group of units:
\begin{definition}
	The \keyword{Clifford group} $\Gamma(V,q)$ is the subgroup of $\CL(V,q)^\times$ whose elements $x$ satisfy
	\begin{equation*}
		\alpha(x)vx^{-1}\in V
	\end{equation*}
	for any $v\in V$.
\end{definition}

Its elements act as linear transformations of $V$ which means that we have a homomorphism $\Gamma(V,q)\to\mathrm{GL}(V)$ which one can show to have kernel $\R^\times$. A consequence of this is that $x\tau(x)\in\R^\times$, and we define:
\begin{definition}
	The \keyword{norm map} is the homomorphism
	\begin{align*}
		N: \Gamma(V,q)&\to\R^\times\,,\\
		x&\mapsto x\tau(x).
	\end{align*}
\end{definition}
Finally, the following is a standard fact:
\begin{proposition}
	Any element of the Clifford group can be written as a product of vectors with nonzero norm, and consequently lies either in the odd or even part of the Clifford algebra.
\end{proposition}

Let $e \in V$ have $q(e)\neq 0$. We define another involution on $\CL(V,q)$ by 
\begin{equation*}
	\sigma_e(x):=e\tau(x)e\,.
\end{equation*}
There is an injection $ V \to \CL(V,q)^{(0)}$ given by $v \to ev$ and we write $eV$ for this image. This image is clearly invariant under $\sigma_e$. There is an action of the Clifford group on the vector space $eV$ by $$ x\cdot \alpha = \sigma_e(x) \alpha x\,. $$ This action is equivalent to first acting on the vector space $V$ in the usual way and then embedding into $eV$.

From now on assume that $V$ is a real vector space and that the quadratic form $q$ on $V$ has signature $(1,n)$, i.e. $V=\R^{1,n}$. We denote the resulting Clifford algebra and group by $\CL(1,n)$ and $\Gamma(1,n)$, respectively to emphasize the special case. Fix an orthonormal basis $\{e,e_1,\dots, e_{n}\}$ for which $q(e)=1=-q(e_i)$ for any $1\leq i\leq n$. We denote by $V^+$ the proper cone of vectors with positive norm which contains $e$. 
\begin{proposition}
    The cone $eV^+$ is invariant by the action of elements $x \in \CL(1,n)$ with positive norm. 
\end{proposition}

The algebra  $M_2(\C) $ can be realized as the Clifford algebra $\CL(1,2)$. Under this isomorphisms we find that the even Clifford group is $\GL_2(\R)$, that the automorphism $\tau$ is given by the complex conjugate of the adjoint matrix, the element $e$ can be chosen to be the matrix $\left(\begin{smallmatrix}
  0 & i \\
 -i & 0 \\
\end{smallmatrix}\right)$, the automorphism $\sigma_e$ is given by matrix transpose, and the vector space $eV^+$ is given by the cone of positive definite symmetric matrices.

\section{A polygonal cluster algebra of finite type}\label{sec:B2quiver}
Consider again the quiver in \Cref{fig:FGQuivers}(a). Its mutable portion is an orientation of a $B_2$ Dynkin diagram, and it therefore provides a seed of finite type. Moreover, the resulting cluster algebra only contains 6 clusters, making it very accessible by an explicit calculation. The quivers in the mutation class are depicted in \Cref{fig:B2mutation}. Observe that another mutation at 5 in the last quiver gives the initial quiver with the color of all arrows connected to 2 switched. Thus the two quivers are weaving isomorphic and their seed algebras are canonically identified.

\input{figureB2mutation.tex}

Let us also compute the cluster variables. For this, we associate variables
\begin{equation*}
	\{X_1,X_2,X_3,x_4,x_5,x_6,x_7,X_8\}
\end{equation*}
to the initial quiver to obtain a seed. Note that the variables $x_i$ commute with every other variable, whereas the variables $X_i$ come with 3 additional variables, namely $\sigma(X_i),\tau(X_i),\sigma\tau(X_i)$, which we are suppressing here. The initial quiver uses the same decoration as \Cref{fig:orderedFGQuiverB2}. 

We start with the angles
\begin{align*}
	\Delta_1(X_1)&=\sigma(X_1)\tau(X_2)\\
	\Delta_1(X_2)&=\tau(X_1)\sigma(X_2)\\
	\Delta_2(X_2)&=\sau(X_8)X_2\sau(X_3)\\
	\Delta_2(X_3)&=\sau(X_2)X_3\sau(X_8)\\
	\Delta_2(X_8)&=\sau(X_3)X_8\sau(X_2)\,.
\end{align*}

Therefore mutation at 2 yields the new variable
\begin{equation*}
	X_2'=x_6\cdot\sigma(X_1)X_2^{-1}\sigma(X_8)+x_5\cdot\sigma(X_2)^{-1}X_3\,,
\end{equation*}
and the new angles
\begin{align*}
	\Delta_1(X_1)&=\tau(X_2')\sigma(X_1)\sau(X_3)\\
	\Delta_1(X_2')&=\sau(X_3)\sigma(X_2')\tau(X_1)\\
	\Delta_2(X_2')&=X_2'\sau(X_8)\\
	\Delta_2(X_3)&=\tau(X_1)X_3\tau(X_2')\\
	\Delta_2(X_8)&=\sau(X_2')X_8\,.
\end{align*}

If we only focus on the full subquiver of the big nodes, this was simply a surface type mutation, following Berenstein-Retakh's mutation rule. Next, we carry out a mutation at 5, a small node, to obtain the new variable
\begin{equation*}
	x_5'=x_5^{-1}\big(x_4x_6N(X_8)+x_7N(X_2')\big)\,.
\end{equation*}
The angles are unchanged except for
\begin{align*}
	\Delta_2(X_2')&=X_8X_2'\\
	\Delta_2(X_8)&=X_8X_2'\,.
\end{align*}

Now another mutation at 2 gives
\begin{align*}
	X_2''&=x_4\cdot X_3X_2'^{-1}\tau(X_8)+x_5'\cdot X_1\sigma(X_2')^{-1}\\
	&=x_4x_5^{-1}\cdot\sigma(X_2)\tau(X_8)+x_5^{-1}x_6x_7\cdot N(X_1)\sigma\tau(X_2)^{-1}\tau(X_8)+x_7\tau(X_2)^{-1}\sigma\tau(X_3)
\end{align*}
with angles
\begin{align*}
	\Delta_1(X_1)&=X_8\sigma(X_1)\sau(X_2'')\\
	\Delta_1(X_2'')&=\sigma(X_2'')\sau(X_3)\\
	\Delta_2(X_2'')&=\tau(X_1)X_2''X_8\\
	\Delta_2(X_3)&=\tau(X_2'')X_3\\
	\Delta_2(X_8)&=X_2''X_8\sigma(X_1))\,.
\end{align*}

Performing one more mutation at 5, we obtain
\begin{align*}
	x_5''&=x_5'^{-1}\big(x_4x_7N(X_3)+x_6N(X_2'')\big)\\
	&=x_5^{-1}x_6x_7\cdot N(X_1)+x_4x_5^{-1}\cdot N(X_2)\,,
\end{align*}
and the angles are unchanged except for
\begin{align*}
	\Delta_1(X_2'')&=X_3\sigma(X_2'')\\
	\Delta_2(X_3)&=X_3\sigma(X_2'')\,.
\end{align*}

One last mutation at 2 leads to the variable
\begin{align*}
	X_2'''=x_7\cdot\sigma\tau(X_3)X_2''^{-1}X_1+x_5''\cdot\sigma(X_2'')^{-1}\sigma\tau(X_8)=\tau(X_2)\,.
\end{align*}
The corresponding angles after this mutation read
\begin{align*}
	\Delta_1(X_1)&=\sau(X_2''')\sigma(X_1)\\
	\Delta_1(X_2''')&=X_8\sigma(X_2''')X_3\\
	\Delta_2(X_2''')&=X_2'''\tau(X_1)\\
	\Delta_2(X_3)&=\sigma(X_2''')X_3\sau(X_8)\\
	\Delta_2(X_8)&=\sau(X_3)X_8\sigma(X_2''')\,.
\end{align*}

Finally, through another mutation at 5 we obtain
\begin{equation*}
	x_5'''=x_5''^{-1}\big(x_6x_7N(X_1)+x_4N(X_2''')\big)=x_5\,
\end{equation*}
with angles
\begin{align*}
	\Delta_1(X_1)&=\sigma(X_1)X_2'''\\	\Delta_1(X_2''')&=X_8\sigma(X_2''')X_3=X_8\sau(X_2)X_3=\sau(\Delta_2(X_2))\\	\Delta_2(X_2''')&=\sigma(X_1)X_2'''=\sigma(X_1)\tau(X_2)=\sau(\Delta_1(X_2))\\
	\Delta_2(X_3)&=\sigma(X_2''')X_3\sau(X_8)\\
	\Delta_2(X_8)&=\sau(X_3)X_8\sigma(X_2''')\,.
\end{align*}

So, we have arrived at a seed which can be identified with the original one by weaving at 2 as needed.

\bibliographystyle{alpha}
\bibliography{main.bib}

\end{document}

%% file: figurePentagonTriangulationQuiver.tex
\begin{figure}[ht]
    \centering
    \begin{tikzpicture}
        \begin{scope}[color=gray,very thick]
            \node[regular polygon,draw,regular polygon sides=5, minimum size=3cm, very thick] (p) at (0,0) {};
            \draw[] (p.corner 1) -- (p.corner 3);
            \draw[] (p.corner 1) -- (p.corner 4);
        \end{scope}
        \node[mutable]	(1)	at (p.side 1)	[]	{};
        \node[mutable]	(2)	at (p.side 2)	[]	{};
        \node[mutable]	(3)	at (p.side 3)	[]	{};
        \node[mutable]	(4)	at (p.side 4)	[]	{};
        \node[mutable]	(5)	at (p.side 5)	[]	{};
        \node[mutable]  (6) at ($0.5*(p.corner 1)+0.5*(p.corner 3)$)  []  {};
        \node[mutable]  (7) at ($0.5*(p.corner 1)+0.5*(p.corner 4)$)  []  {};

        \draw[arrow]    (1) to (6);
        \draw[arrow]    (6) to (2);
        \draw[arrow]    (2) to (1);

        \draw[arrow]    (6) to (7);
        \draw[arrow]    (7) to (3);
        \draw[arrow]    (3) to (6);

        \draw[arrow]    (7) to (5);
        \draw[arrow]    (5) to (4);
        \draw[arrow]    (4) to (7);
    \end{tikzpicture}
    \caption{Fan triangulation of a disk with 5 marked points on the boundary and associated quiver.}
    \label{fig:PentagonTriangulationQuiver}
\end{figure}

%% file: figureFlipSurfaceCase.tex
\begin{figure}[htb]
	\begin{center}
		\begin{tikzpicture}[scale=0.8]
			\begin{scope}[xscale=0.577]
				\node[vertex]	(l) at (0,0)						[label = above: $\ell$]		{};
				\node[vertex]	(i) at (-2*\distL,-2*\distL)		[label = left: $i$]			{};
				\node[vertex]	(j) at (2*\distL,-2*\distL)		    [label = right: $j$]		{};
				\node[vertex]	(k) at (4*\distL,0)				    [label = above: $k$]		{};
				\draw[qarrow, postaction = decorate]		(i) to (j);
				\draw[qarrow, postaction = decorate]		(j) to (k);
				\draw[qarrow, postaction = decorate]		(i) to (k);
				\draw[qarrow, postaction = decorate]		(i) to (l);
				\draw[qarrow, postaction = decorate]		(l) to (k);
				\node[vertex]	(ij) at ($0.4*(i)+0.6*(j)+(0,-0.1)$)	[label = below: $x_{ij}$]	{};
				\node[vertex]	(jk) at ($0.5*(j)+0.5*(k)+(0.1,-0.1)$)	[label = right: $x_{jk}$]	{};
				\node[vertex]	(ik) at ($0.4*(i)+0.6*(k)+(0,-0.2)$)	[label = below: $x_{ik}$]	{};
				\node[vertex]	(il) at ($0.5*(i)+0.5*(l)+(-0.1,0)$)	[label = left: $x_{i\ell}$]	{};
				\node[vertex]	(lk) at ($0.5*(l)+0.5*(k)+(0,0.1)$)		[label = above: $x_{\ell k}$]	{};
			\end{scope}
			\begin{scope}[very thick,shift={(3*\distL,-\distL)}]
				\draw [<->] (0,0)--(\distL,0)
				node [above,midway]{\small flip};
			\end{scope}
			\begin{scope}[xscale=0.577,shift={(10*\distL,0)}]
				\node[vertex]	(l) at (0,0)					[label = above: $\ell$]		{};
				\node[vertex]	(i) at (-2*\distL,-2*\distL)	[label = left: $i$]			{};
				\node[vertex]	(j) at (2*\distL,-2*\distL)		[label = right: $j$]		{};
				\node[vertex]	(k) at (4*\distL,0)				[label = above: $k$]		{};
				\draw[qarrow, postaction = decorate]		(i) to (j);
				\draw[qarrow, postaction = decorate]		(j) to (k);
				\draw[qarrow, postaction = decorate]		(l) to (j);
				\draw[qarrow, postaction = decorate]		(i) to (l);
				\draw[qarrow, postaction = decorate]		(l) to (k);
				\node[vertex]	(ij) at ($0.4*(i)+0.6*(j)+(0,-0.1)$)	[label = below: $x_{ij}$]	{};
				\node[vertex]	(jk) at ($0.5*(j)+0.5*(k)+(0.1,-0.1)$)	[label = right: $x_{jk}$]	{};
				\node[vertex]	(ik) at ($0.5*(l)+0.5*(j)+(0.2,0)$)		[label = right: $x_{\ell j}$]	{};
				\node[vertex]	(il) at ($0.5*(i)+0.5*(l)+(-0.1,0)$)	[label = left: $x_{i\ell}$]	{};
				\node[vertex]	(lk) at ($0.5*(l)+0.5*(k)+(0,0.1)$)		[label = above: $x_{\ell k}$]	{};
			\end{scope}
		\end{tikzpicture}
	\end{center}
	\caption{Flips are related to exchange relations.}
	\label{fig:FlipSurfaceCase}
\end{figure}

%% file: figureFGQuivers.tex
\begin{figure}[htb]
	\centering
		\begin{subfigure}[b]{0.48\textwidth}
			\centering
			\scalebox{0.8}{
			\begin{tikzpicture}[node distance = 50,xscale=0.577,scale=1.5]
				\node[vertex]		(c1)	at (0,0)							[label = above:$2$]{};
				\node[vertex]		(c2)	at (-3*\distance,-3*\distance)	[label = left:$1$]	{};
		        \node[vertex]		(c3)	at (3*\distance,-3*\distance)	[label = right:$0$]	{};
				
				\draw[line] (c1) to (c2);
				\draw[line] (c2) to (c3);
				\draw[line] (c3) to (c1);
					
				\node[frozenBig]		(1)		at (	-\distance,-\distance)		[]	{};
				\node[mutableBig]	(2)		at (0,-\distance)				[]	{};
				\node[frozenBig]		(3)		at (\distance,-\distance)		[]	{};
				\node[frozen]		(4)		at (-2*\distance,-2*\distance)	[]	{};
				\node[mutable]		(5)		at (0,-2*\distance)				[]	{};
			    \node[frozen]		(6)		at (2*\distance,-2*\distance)	[]	{};
				\node[frozen]		(7)		at (-\distance,-3*\distance)		[]	{};
				\node[frozenBig]		(8)		at (\distance,-3*\distance)		[]	{};
				
				\draw[barrow] (1) to (2);
				\draw[barrow] (2) to (3);
				
				\draw[arrow] (4) to (5);
				\draw[arrow] (5) to (6);
				
				\draw[arrow] (5) to (1);
				\draw[arrow] (6) to (2);
				\draw[arrow] (2) to (5);
				
				\draw[barrow] (3) to (8);
				\draw[barrow] (8) to (2);
				
				\draw[arrow] (5) to (7);
				\draw[arrow] (7) to (4);
				
				\draw[halfarrow, bend left] (3) to (6);
				\draw[halfarrow, bend angle = 12, bend right] (7) to (8);
				\draw[halfarrow, bend right] (1) to (4);
			\end{tikzpicture}}
			\caption{}
		\end{subfigure}
		\hfill
		\begin{subfigure}[b]{0.48\textwidth}
			\centering
			\scalebox{0.8}{
			\begin{tikzpicture}[node distance = 50,xscale=0.577]
				\node[vertex]		(c1)	at (0,0)							[label = above:$2$]{};
				\node[vertex]		(c2)	at (-4.5*\distance,-4.5*\distance)	[label = left:$1$]	{};
		        \node[vertex]		(c3)	at (4.5*\distance,-4.5*\distance)	[label = right:$0$]	{};
				
				\draw[line] (c1) to (c2);
				\draw[line] (c2) to (c3);
				\draw[line] (c3) to (c1);
				
				\node[frozenBig]		(1)		at (-1.5*\distance,-1.5*\distance)	[] {};
				\node[mutableBig]	(2)		at (-0.5*\distance,-1.5*\distance)	[] {};
				\node[mutableBig]	(3)		at (0.5*\distance,-1.5*\distance)	[] {};
				\node[frozenBig]		(4)		at (1.5*\distance,-1.5*\distance)	[] {};
				\node[frozen]		(5)		at (-2.5*\distance,-2.5*\distance)	[] {};
				\node[mutable]		(6)		at (-0.5*\distance,-2.5*\distance)	[] {};
				\node[mutable]		(7)		at (0.5*\distance,-2.5*\distance)	[] {};
				\node[frozen]		(8)		at (2.5*\distance,-2.5*\distance)	[] {};
				\node[frozen]		(9)		at (-3.5*\distance,-3.5*\distance)	[] {};
				\node[mutable]		(10)		at (-0.5*\distance,-3.5*\distance)	[] {};
				\node[mutable]		(11)		at (0.5*\distance,-3.5*\distance)	[] {};
				\node[frozen]		(12)		at (3.5*\distance,-3.5*\distance)	[] {};
				\node[frozen]		(13)		at (-2.5*\distance,-4.5*\distance)	[] {};
				\node[frozen]		(14)		at (0*\distance,-4.5*\distance)		[] {};
				\node[frozenBig]		(15)		at (2.5*\distance,-4.5*\distance)	[] {};
				
				\draw[barrow] (1) to (2);
				\draw[barrow] (2) to (3);
				\draw[barrow] (3) to (4);
				
				\draw[arrow] (5) to (6);
				\draw[arrow] (6) to (7);
				\draw[arrow] (7) to (8);
				
				\draw[arrow] (9) to(10);
				\draw[arrow] (10)to(11);
				\draw[arrow] (11)to(12);
				
				\draw[arrow] (2) to (6);
				\draw[arrow] (3) to (7);
				\draw[arrow] (6) to (1);
				\draw[arrow] (7) to (2);
				\draw[arrow] (8) to (3);
				
				\draw[arrow] (6) to(10);
				\draw[arrow] (7) to(11);
				\draw[arrow] (10)to (5);
				\draw[arrow] (11)to (6);
				\draw[arrow] (12)to (7);
				
				\draw[arrow] (10)to(13);
				\draw[arrow] (13)to (9);
				\draw[arrow] (11)to(14);
				\draw[arrow] (14)to(10);
				\draw[barrow] (4) to(15);
				\draw[barrow] (15)to (3);
				
				\draw[halfarrow, bend right] (1) to (5);
				\draw[halfarrow, bend right] (5) to (9);
				\draw[halfarrow, bend left] (4) to (8);
				\draw[halfarrow, bend left] (8) to(12);
				\draw[halfarrow, bend angle = 14, bend right] (13) to (14);
				\draw[halfarrow, bend angle = 14, bend right] (14) to (15);
			\end{tikzpicture}}
			\caption{}
		\end{subfigure}
	\caption{Quivers associated to the cluster structure for decorated representations into split real forms of (A) $B_2$ and (B) $B_3$.}
	\label{fig:FGQuivers}
\end{figure}

%% file: figureColoredMutation.tex
\begin{figure}[h]
	\centering
			\begin{tikzpicture}[scale=0.55]
				\begin{scope}
					\node[mutableBig]	(B) at (0,0)					[label = above: $2$]	{};
					\node[mutableBig]	(A) at (-\distL,-\distL)		[label = left: $1$]		{};
					\node[mutableBig]	(C) at (\distL,-\distL)		[label = right: $3$]	{};
					
					\draw[barrow] (A) to (B);
					\draw[xarrow] (B) to (C);
				\end{scope}
				\begin{scope}[shift={(1.8*\distL,0)}]
					\node[vertex] (v1) at (0,-0.5) []{};
					\node[vertex] (v2) at (1.5*\distL,-0.5) []{};
					\draw[<->] (v1) -- (v2)
					node[above,midway,text centered, text width=2cm]{$\mu_2$};
				\end{scope}
				\begin{scope}[shift={(5*\distL,0)}]
					\node[mutableBig]	(B) at (0,0)					[label = above: $2$]	{};
					\node[mutableBig]	(A) at (-\distL,-\distL)		[label = left: $1$]		{};
					\node[mutableBig]	(C) at (\distL,-\distL)		[label = right: $3$]	{};
					
					\draw[barrow] (B) to (A);
					\draw[xarrow] (C) to (B);
					\draw[xarrow] (A) to (C);
				\end{scope}
			\end{tikzpicture}
%
%
	\caption{Example of a mutation on a 2-colored quiver.}
	\label{fig:colMutEx}
\end{figure}

%% file: figureWeavingIso.tex
\begin{figure}[h]
	\centering
		\begin{tikzpicture}[scale=0.55]
			\begin{scope}
				\node[mutableBig]	(B) at (0,0)					[label = above: $2$]	{};
				\node[mutableBig]	(A) at (-\distL,-\distL)		[label = left: $1$]		{};
				\node[mutableBig]	(C) at (\distL,-\distL)		[label = below: $3$]	{};
				\node[mutableBig]	(D) at (3*\distL,-\distL)	[label = below: $4$]	{};
				
				\draw[barrow] (A) to (B);
				\draw[xarrow] (B) to (C);
				\draw[xarrow] (C) to (A);
				\draw[barrow] (C) to (D);
			\end{scope}
			\begin{scope}[shift={(4.8*\distL,-0.7*\distL)}]
				\node[vertex] (s) at (0,0) [label = above: $\scalebox{1.5}{\weaviso{}}$]{};
			\end{scope}
			\begin{scope}[shift={(8*\distL,0)}]
				\node[mutableBig]	(B) at (0,0)					[label = above: $2$]	{};
				\node[mutableBig]	(A) at (-\distL,-\distL)		[label = left: $1$]		{};
				\node[mutableBig]	(C) at (\distL,-\distL)		[label = below: $3$]	{};
				\node[mutableBig]	(D) at (3*\distL,-\distL)	[label = below: $4$]	{};
				
				\draw[barrow] (A) to (B);
				\draw[barrow] (B) to (C);
				\draw[barrow] (C) to (A);
				\draw[xarrow] (C) to (D);
			\end{scope}
		\end{tikzpicture}
		\caption{Weaving isomorphic 2-colored quivers that are not isomorphic.}
	\label{fig:weavingIso}
\end{figure}

%% file: figureMarkedColoredMutation.tex
\begin{figure}[h]
	\centering
		\begin{tikzpicture}[scale=0.55]
			\begin{scope}
				\node[mutable]		(B) at (0,0)					[label = above: $2$]	{};
				\node[mutableBig]	(A) at (-\distL,-\distL)		[label = left: $1$]		{};
				\node[mutableBig]	(C) at (\distL,-\distL)		[label = right: $3$]	{};
				
				\draw[arrow] (A) to (B);
				\draw[arrow] (B) to (C);
			\end{scope}
			\begin{scope}[shift={(1.8*\distL,0)}]
				\node[vertex] (v1) at (0,-0.5) []{};
				\node[vertex] (v2) at (1.5*\distL,-0.5) []{};
				\draw[<->] (v1) -- (v2)
				node[above,midway,text centered, text width=2cm]{$\mu_2$};
			\end{scope}
			\begin{scope}[shift={(5*\distL,0)}]
				\node[mutable]		(B) at (0,0)					[label = above: $2$]	{};
				\node[mutableBig]	(A) at (-\distL,-\distL)		[label = left: $1$]		{};
				\node[mutableBig]	(C) at (\distL,-\distL)		[label = right: $3$]	{};
				
				\draw[arrow] (B) to (A);
				\draw[arrow] (C) to (B);
				\draw[barrow] ($(A)+(0.3,0.15)$) to ($(C)+(-0.3,0.15)$);
				\draw[xarrow] ($(A)+(0.3,-0.15)$) to ($(C)+(-0.3,-0.15)$);
			\end{scope}
			\end{tikzpicture}
	\caption{Example of a mutation at a small node in a mixed 2-colored quiver.}
	\label{fig:mixedColMutEx}
\end{figure}

%% file: figureAngleEx.tex
\begin{figure}[h]
	\centering
	\scalebox{.8}{
	\begin{tikzpicture}
		\begin{scope}[scale=0.55]
			\pic[name=A,rotate=-225] at (-5,3) {NCnodeHalfL};
			\node[vertex] () at (-5,3) [label ={[yshift=4mm]$1$}]{};
			
			\pic[name=B,rotate=-135] at (-4,-4) {NCnodeHalfL};
			\node[vertex] () at (-6,-2) [label ={[xshift=-6mm, yshift=-3mm]$2$}]{};
			
			\pic[name=C,rotate=-135] at (-6,-2) {NCnodeHalfL};
			\node[vertex] () at (-4,-4) [label ={[xshift=-6mm, yshift=-3mm]$3$}]{};
			
			\pic[name=D,rotate=225] at (5,3) {NCnodeHalfR};
			\node[vertex] () at (5,3) [label ={[yshift=4mm]$4$}]{};
			
			\pic[name=F,rotate=135] at (5,-3) {NCnodeHalfR};
			\node[vertex] () at (5,-3) [label ={[xshift=6mm, yshift=-3mm]$5$}]{};
			
			\pic[name=E] at (0,0) {NCnode};
			\node[vertex] () at (0,0) [label ={[yshift=7mm]$6$}]{};
			
			\draw[barrow] (A-circle.-20) to [out=-20,in=140] (E-circle.140);
			\draw[barrow] (A-circle.10) to [out=10,in=-30] (B-circle.-30);
			
			\draw[barrow] (B-circle.90) to [out=90,in=20]  (C-circle.20);
			
			\draw[barrow] (C-circle.90) to [out=90,in=-90] (A-circle.-90);
			
			\draw[xarrow] (E-circle.210) to [out=210,in=0] (B-circle.0);
			\draw[barrow] (E-circle.250) to [out=250,in=-10] (C-circle.-10);
			\draw[barrow] (E-circle.45) to [out=45,in=-180] (D-circle.-180);
			
			\draw[barrow] (D-circle.-90) to [out=-90,in=90] (F-circle.90);
			
			\draw[barrow] (F-circle.180) to [out=180,in=-45] (E-circle.-45);

           \draw[carrow] (D-circle.center) to[out=160,in=20] (A-circle.center);
		\end{scope}
	\end{tikzpicture}}
	\caption{An ordered quiver.}
	\label{fig:angleEx}
\end{figure}

%% file: figurePseudoCycles.tex
\begin{figure}[h!]
	\centering
	\begin{subfigure}[b]{0.24\textwidth}
		\centering
		\begin{tikzpicture}[yscale=0.577]
			\node[mutableBig]	(X1) at (0,0)				[label = left: $1$]	{};
			\node[mutableBig]	(X2) at (\distL,\distL)		[label = right: $2$]	{};
			\node[mutableBig]	(X3) at (\distL,-\distL)		[label = right: $3$]	{};
			
			\draw[barrow] (X1) to (X2);
			\draw[barrow] (X2) to (X3);
			\draw[barrow] (X3) to (X1);
		\end{tikzpicture}
		\caption{}
	\end{subfigure}
	\begin{subfigure}[b]{0.24\textwidth}
		\centering
		\begin{tikzpicture}[yscale=0.577]
			\node[mutableBig]	(X1) at (0,0)				[label = left: $1$]	{};
			\node[mutableBig]	(X2) at (\distL,\distL)		[label = right: $2$]	{};
			\node[mutableBig]	(X3) at (\distL,-\distL)		[label = right: $3$]	{};
			
			\draw[xarrow] (X1) to (X2);
			\draw[barrow] (X2) to (X3);
			\draw[xarrow] (X3) to (X1);
		\end{tikzpicture}
		\caption{}
	\end{subfigure}
	\begin{subfigure}[b]{0.24\textwidth}
		\centering
		\begin{tikzpicture}[yscale=0.577]
			\node[mutableBig]	(X1) at (0,0)				[label = left: $1$]	{};
			\node[mutableBig]	(X2) at (\distL,\distL)		[label = right: $2$]	{};
			\node[mutableBig]	(X3) at (\distL,-\distL)		[label = right: $3$]	{};
			
			\draw[xarrow] (X1) to (X2);
			\draw[xarrow] (X2) to (X3);
			\draw[barrow] (X3) to (X1);
		\end{tikzpicture}
		\caption{}
	\end{subfigure}
	\begin{subfigure}[b]{0.24\textwidth}
		\centering
		\begin{tikzpicture}[yscale=0.577]
			\node[mutableBig]	(X1) at (0,0)				[label = left: $1$]	{};
			\node[mutableBig]	(X2) at (\distL,\distL)		[label = right: $2$]	{};
			\node[mutableBig]	(X3) at (\distL,-\distL)		[label = right: $3$]	{};
			
			\draw[barrow] (X1) to (X2);
			\draw[xarrow] (X2) to (X3);
			\draw[xarrow] (X3) to (X1);
		\end{tikzpicture}
		\caption{}
	\end{subfigure}
	\caption{All possible types of 3-cycles.}
	\label{fig:3cycles}
\end{figure}

\begin{figure}[h!]
	\centering
	\begin{subfigure}[b]{0.24\textwidth}
		\centering
		\begin{tikzpicture}[yscale=0.577]
			\node[mutableBig]	(X1) at (0,0)				[label = left: $1$]	{};
			\node[mutableBig]	(X2) at (\distL,\distL)		[label = right: $2$]	{};
			\node[mutableBig]	(X3) at (\distL,-\distL)		[label = right: $3$]	{};
			
			\draw[barrow] (X1) to (X2);
			\draw[xarrow] (X3) to (X2);
			\draw[barrow] (X3) to (X1);
		\end{tikzpicture}
		\caption{}
	\end{subfigure}
	\begin{subfigure}[b]{0.24\textwidth}
		\centering
		\begin{tikzpicture}[yscale=0.577]
			\node[mutableBig]	(X1) at (0,0)				[label = left: $1$]	{};
			\node[mutableBig]	(X2) at (\distL,\distL)		[label = right: $2$]	{};
			\node[mutableBig]	(X3) at (\distL,-\distL)		[label = right: $3$]	{};
			
			\draw[xarrow] (X1) to (X2);
			\draw[xarrow] (X3) to (X2);
			\draw[xarrow] (X3) to (X1);
		\end{tikzpicture}
		\caption{}
	\end{subfigure}
	\begin{subfigure}[b]{0.24\textwidth}
		\centering
		\begin{tikzpicture}[yscale=0.577]
			\node[mutableBig]	(X1) at (0,0)				[label = left: $1$]	{};
			\node[mutableBig]	(X2) at (\distL,\distL)		[label = right: $2$]	{};
			\node[mutableBig]	(X3) at (\distL,-\distL)		[label = right: $3$]	{};
			
			\draw[xarrow] (X1) to (X2);
			\draw[barrow] (X3) to (X2);
			\draw[barrow] (X3) to (X1);
		\end{tikzpicture}
		\caption{}
	\end{subfigure}
	\begin{subfigure}[b]{0.24\textwidth}
		\centering
		\begin{tikzpicture}[yscale=0.577]
			\node[mutableBig]	(X1) at (0,0)				[label = left: $1$]	{};
			\node[mutableBig]	(X2) at (\distL,\distL)		[label = right: $2$]	{};
			\node[mutableBig]	(X3) at (\distL,-\distL)		[label = right: $3$]	{};
			
			\draw[barrow] (X1) to (X2);
			\draw[barrow] (X3) to (X2);
			\draw[xarrow] (X3) to (X1);
		\end{tikzpicture}
		\caption{}
	\end{subfigure}
	\caption{All possible types of 3-pseudocycles based at $1$.}
	\label{fig:3pseudocycles}
\end{figure}

%% file: figureLTTcolor.tex
\begin{figure}[ht]
\centering
    \begin{subfigure}[b]{0.3\textwidth}
        \centering
        \scalebox{0.9}{
        \begin{tikzpicture}[thick]
            \node[circle] (j) at (-2,2) {$j$};
            \node[circle] (k) at (0,3.5) {$k$};
            \node[circle] (l) at (2,2) {$l$};
            \node[circle] (i) at (0,0) {$i$};

            \draw[->] (j) edge node[fill=white, anchor=center, inner sep=2pt] {$t_1$} (i);
            \draw[->] (k) edge node[fill=white, pos=0.75, inner sep=2pt] {$t_2$} (i);
            \draw[->] (i) edge node[fill=white, anchor=center, inner sep=2pt] {$t_3$} (l);

            \draw[->] (j) edge node[fill=white, anchor=center, inner sep=2pt] {$t_1+t_2+1$} (k);
            \draw[->] (k) edge node[fill=white, anchor=center, inner sep=2pt] {$t_2+t_3+1$} (l);
            \draw[->] (l) edge node[fill=white, anchor=center, inner sep=2pt] {$t_1+t_3$} (j);
        \end{tikzpicture}}
    \end{subfigure}
    \hfill
    \begin{subfigure}[b]{0.3\textwidth}
        \centering
        \scalebox{0.9}{
        \begin{tikzpicture}[thick]
            \node[circle] (j) at (-2,2) {$j$};
            \node[circle] (k) at (0,3.5) {$k$};
            \node[circle] (l) at (2,2) {$l$};
            \node[circle] (i) at (0,0) {$i$};

            \draw[->] (j) edge node[fill=white, anchor=center, inner sep=2pt] {$t_1$} (i);
            \draw[->] (k) edge node[fill=white, pos=0.75, inner sep=2pt] {$t_2$} (i);
            \draw[->] (i) edge node[fill=white, anchor=center, inner sep=2pt] {$t_3$} (l);

            \draw[->] (j) edge node[fill=white, anchor=center, inner sep=2pt] {$t_1+t_2+1$} (k);
            \draw[->] (k) edge node[fill=white, anchor=center, inner sep=2pt] {$t_2+t_3+1$} (l);
            \draw[->] (j) edge node[fill=white, anchor=center, inner sep=2pt] {$t_1+t_3+1$} (l);
        \end{tikzpicture}}
    \end{subfigure}
    \hfill
    \begin{subfigure}[b]{0.3\textwidth}
        \centering
        \scalebox{0.9}{
        \begin{tikzpicture}[thick]
            \node[circle] (j) at (-2,2) {$j$};
            \node[circle] (k) at (0,3.5) {$k$};
            \node[circle] (l) at (2,2) {$l$};
            \node[circle] (i) at (0,0) {$i$};

            \draw[->] (j) edge node[fill=white, anchor=center, inner sep=2pt] {$t_1$} (i);
            \draw[->] (k) edge node[fill=white, pos=0.75, inner sep=2pt] {$t_2$} (i);
            \draw[->] (l) edge node[fill=white, anchor=center, inner sep=2pt] {$t_3$} (i);

            \draw[->] (j) edge node[fill=white, anchor=center, inner sep=2pt] {$t_1+t_2+1$} (k);
            \draw[->] (k) edge node[fill=white, anchor=center, inner sep=2pt] {$t_2+t_3+1$} (l);
            \draw[->] (j) edge node[fill=white, anchor=center, inner sep=2pt] {$t_1+t_3+1$} (l);
        \end{tikzpicture}}
    \end{subfigure}
\caption{Coloring of a locally transitive tournament determined by $i$.}
\label{fig:LTTcoloring}
\end{figure}

%% file: figureMutationExample.tex
\begin{figure}[ht]
	\centering
	\scalebox{.8}{
	\begin{tikzpicture}
		\begin{scope}[scale=0.55]
			\pic[name=A,rotate=-225] at (-5,3) {NCnodeHalfL};
			\node[vertex] () at (-5,3) [label ={[yshift=4mm]$1$}]{};
			
			\pic[name=B,rotate=-135] at (-4,-4) {NCnodeHalfL};
            \node[vertex] () at (-4,-4) [label ={[xshift=-6mm, yshift=-3mm]$3$}]{};
			
			\pic[name=C,rotate=-135] at (-6,-2) {NCnodeHalfL};
			\node[vertex] () at (-6,-2) [label ={[xshift=-6mm, yshift=-3mm]$2$}]{};
			
			\pic[name=D,rotate=225] at (5,3) {NCnodeHalfR};
			\node[vertex] () at (5,3) [label ={[yshift=4mm]$4$}]{};
			
			\pic[name=F,rotate=135] at (5,-3) {NCnodeHalfR};
			\node[vertex] () at (5,-3) [label ={[xshift=6mm, yshift=-3mm]$5$}]{};
			
			\pic[name=E] at (0,0) {NCnode};
			\node[vertex] () at (0,0) [label ={[yshift=7mm]$6$}]{};
			
			\draw[barrow] (A-circle.-20) to [out=-20,in=140] (E-circle.140);
			\draw[barrow] (A-circle.10) to [out=10,in=-30] (B-circle.-30);
			
			\draw[barrow] (B-circle.90) to [out=90,in=20]  (C-circle.20);
			
			\draw[barrow] (C-circle.90) to [out=90,in=-90] (A-circle.-90);
			
			\draw[xarrow] (E-circle.210) to [out=210,in=0] (B-circle.0);
			\draw[barrow] (E-circle.250) to [out=250,in=-10] (C-circle.-10);
			\draw[barrow] (E-circle.45) to [out=45,in=-180] (D-circle.-180);
			
			\draw[barrow] (D-circle.-90) to [out=-90,in=90] (F-circle.90);
			
			\draw[barrow] (F-circle.180) to [out=180,in=-45] (E-circle.-45);

            \draw[carrow] (D-circle.center) to[out=160,in=20] (A-circle.center);
		\end{scope}
		\begin{scope}[ultra thick, shift={(3.5,0)}]
			\draw[->] (0,0)--(2,0)
			node [above,midway]{$\mu_6$};
		\end{scope}
		\begin{scope}[scale=0.55,shift={(18,0)}]
			\pic[name=A,rotate=-225] at (-5,3) {NCnodeHalfL};
			\node[vertex] () at (-5,3) [label ={[yshift=4mm]$1$}]{};
			
			\pic[name=B,rotate=-135] at (-4,-4) {NCnodeHalfL};
			\node[vertex] () at (-4,-4) [label ={[xshift=-6mm, yshift=-3mm]$3$}]{};
			
			\pic[name=C,rotate=-135] at (-6,-2) {NCnodeHalfL};
            \node[vertex] () at (-6,-2) [label ={[xshift=-6mm, yshift=-3mm]$2$}]{};
			
			\pic[name=D,rotate=225] at (5,3) {NCnodeHalfR};
			\node[vertex] () at (5,3) [label ={[yshift=4mm]$4$}]{};
			
			\pic[name=F,rotate=135] at (5,-3) {NCnodeHalfR};
			\node[vertex] () at (5,-3) [label ={[xshift=6mm, yshift=-3mm]$5$}]{};
			
			\pic[name=E,rotate=-90] at (0,0) {NCnode};
			\node[vertex] () at (0,0) [label ={[yshift=7mm]$6$}]{};
			
			\draw[barrow] (E-circle.140) to [out=140,in=-70] (A-circle.-70);
            \draw[xarrow] (D-circle.-80) to[out=-80,in=0] ($(E-circle.center)+(90:2.8)$) to[out=180,in=-100] (A-circle.-100);
            \draw[carrow] (A-circle.center) to [out=20,in=-30] (B-circle.center);
			
			\draw[xarrow] (B-circle.100) to [out=100,in=210] (E-circle.210);
			\draw[barrow] (B-circle.80) to [out=80,in=20]  (C-circle.20);
			\draw[xarrow] (F-circle.160) to [out=160,in=20] (B-circle.20);
			
			\draw[barrow] (C-circle.90) to [out=90,in=250] (E-circle.250);
			\draw[barrow] (F-circle.180) to [out=180,in=-10] (C-circle.-10);
			
			\draw[barrow] (D-circle.-110) to [out=-110,in=45] (E-circle.45);
			
			\draw[barrow] (E-circle.-45) to [out=-45,in=90] (F-circle.90);
		\end{scope}
	\end{tikzpicture}}
	\caption{An ST-compatible ordered quiver and its mutation.}
	\label{fig:mutationEx}
\end{figure}

%% file: figureSmallNodeMutNonAdmissible.tex
\begin{figure}[h]
	\centering
		\begin{tikzpicture}[xscale=0.8*0.71,yscale=0.8]
            \begin{scope}[shift={(-6.5*\distL,0)}]
				\node[mutableBig]   (A) at (0,0)                [label = above: $1$]    {};
				\node[mutableBig]	(B) at (\distL,-\distL)		[label = right: $2$]	{};
				\node[mutableBig]	(C) at (-\distL,-\distL)	[label = left: $3$]		{};
				\node[mutable]	    (d) at (-2*\distL,0)		[label = above: $4$]	{};
                \node[mutable]	    (e) at (2*\distL,0)  		[label = above: $5$]	{};
				
				\draw[barrow] (A) to (B);
				\draw[barrow] (B) to (C);
                \draw[xarrow] (A) to (C);
                \draw[arrow] (d) to (A);
                \draw[arrow] (C) to (d);
                \draw[arrow] (B) to (e);
                \draw[arrow] (e) to (A);
			\end{scope}
            \begin{scope}[shift={(-4*\distL,0)}]
				\node[vertex] (v1) at (0,-0.5*\distL) []{};
				\node[vertex] (v2) at (1.5*\distL,-0.5*\distL) []{};
				\draw[<->] (v1) -- (v2)
				node[above,midway,text centered, text width=2cm]{$\mu_4$};
			\end{scope}
			\begin{scope}
                \node[mutableBig]   (A) at (0,0)                [label = above: $1$]    {};
				\node[mutableBig]	(B) at (\distL,-\distL)		[label = right: $2$]	{};
				\node[mutableBig]	(C) at (-\distL,-\distL)	[label = left: $3$]		{};
				\node[mutable]	    (d) at (-2*\distL,0)		[label = above: $4$]	{};
                \node[mutable]	    (e) at (2*\distL,0)  		[label = above: $5$]	{};
				
				\draw[barrow] (A) to (B);
				\draw[barrow] (B) to (C);
                \draw[barrow] (C) to (A);
                \draw[arrow] (A) to (d);
                \draw[arrow] (d) to (C);
                \draw[arrow] (B) to (e);
                \draw[arrow] (e) to (A);
			\end{scope}
			\begin{scope}[shift={(2.5*\distL,0)}]
				\node[vertex] (v1) at (0,-0.5*\distL) []{};
				\node[vertex] (v2) at (1.5*\distL,-0.5*\distL) []{};
				\draw[<->] (v1) -- (v2)
				node[above,midway,text centered, text width=2cm]{$\mu_5$};
			\end{scope}
			\begin{scope}[shift={(6.5*\distL,0)}]
				\node[mutableBig]   (A) at (0,0)                [label = above: $1$]    {};
				\node[mutableBig]	(B) at (\distL,-\distL)		[label = right: $2$]	{};
				\node[mutableBig]	(C) at (-\distL,-\distL)	[label = left: $3$]		{};
				\node[mutable]	    (d) at (-2*\distL,0)		[label = above: $4$]	{};
                \node[mutable]	    (e) at (2*\distL,0)  		[label = above: $5$]	{};
				
				\draw[xarrow] (B) to (A);
				\draw[barrow] (B) to (C);
                \draw[barrow] (C) to (A);
                \draw[arrow] (A) to (d);
                \draw[arrow] (d) to (C);
                \draw[arrow] (A) to (e);
                \draw[arrow] (e) to (B);
			\end{scope}
			\end{tikzpicture}
	\caption{Mutations at small nodes in a mixed 2-colored quiver.}
	\label{fig:mixedColMutEx2}
\end{figure}

%% file: figureDecPoly.tex
\begin{figure}[h]
\centering
\begin{subfigure}[b]{0.18\textwidth}
	\centering
	\begin{tikzpicture}
 	\node[regular polygon,draw,regular polygon sides=3, minimum size=2.5cm, very thick] (p) at (0,0) {};
  	\draw[angleindicator] (p.side 1) -- (p.corner 3);
  	\draw[angleindicator] (p.side 2) -- (p.corner 1);
  	\draw[angleindicator] (p.side 3) -- (p.corner 2);
    \centerarc[straightangle](p.corner 1)(240:300:0.5)
    \centerarc[straightangle](p.corner 2)(0:60:0.5)
    \centerarc[straightangle](p.corner 3)(120:180:0.5)
    \node[regular polygon,draw,regular polygon sides=3, minimum size=2.5cm, thick] (p) at (0,0) {};
	\end{tikzpicture}
\end{subfigure}
\hfill
\begin{subfigure}[b]{0.18\textwidth}
	\centering
	\begin{tikzpicture}
 	\node[regular polygon,draw,regular polygon sides=3, minimum size=2.5cm, very thick] (p) at (0,0) {};
  	\draw[angleindicator] (p.side 2) -- (p.corner 1);
  	\draw[angleindicator] (p.side 3) to[out=150,in=270] (p.corner 1);
  	\draw[angleindicator] (p.side 1) to[out=30,in=270] (p.corner 1);
    \centerarc[straightangle](p.corner 1)(240:300:0.5)
    \centerarc[straightangle](p.corner 2)(0:60:0.5)
    \centerarc[straightangle](p.corner 3)(120:180:0.5)
    \node[regular polygon,draw,regular polygon sides=3, minimum size=2.5cm, thick] (p) at (0,0) {};
	\end{tikzpicture}
\end{subfigure}
\hfill
\begin{subfigure}[b]{0.18\textwidth}
	\centering
	\begin{tikzpicture}
 	\node[regular polygon,draw,regular polygon sides=4, minimum size=2.5cm, very thick] (p) at (0,0) {};
 	\draw[angleindicator] (p.side 1) -- (p.corner 4);
  	\draw[angleindicator] (p.side 2) -- (p.corner 1);
	\draw[angleindicator] (p.side 3) -- (p.corner 1);
	\draw[angleindicator] (p.side 4) -- (p.corner 2);
    \centerarc[straightangle](p.corner 1)(180:270:0.5)
    \centerarc[crossangle](p.corner 2)(270:360:0.5)
    \centerarc[crossangle](p.corner 3)(0:90:0.5)
    \centerarc[crossangle](p.corner 4)(90:180:0.5)
    \node[regular polygon,draw,regular polygon sides=4, minimum size=2.5cm, thick] (p) at (0,0) {};
	\end{tikzpicture}
\end{subfigure}
\hfill
\begin{subfigure}[b]{0.18\textwidth}
	\centering
	\begin{tikzpicture}
 	\node[regular polygon,draw,regular polygon sides=5, minimum size=2.3cm, very thick] (p) at (0,0) {};
 	\draw[angleindicator] (p.side 1) -- (p.corner 4);
  	\draw[angleindicator] (p.side 2) -- (p.corner 1);
	\draw[angleindicator] (p.side 3) -- (p.corner 1);
	\draw[angleindicator] (p.side 4) -- (p.corner 2);
	\draw[angleindicator] (p.side 5) -- (p.corner 2);
    \centerarc[crossangle](p.corner 1)(216:324:0.4)
    \centerarc[straightangle](p.corner 2)(288:396:0.4)
    \centerarc[crossangle](p.corner 3)(0:108:0.4)
    \centerarc[straightangle](p.corner 4)(72:180:0.4)
    \centerarc[straightangle](p.corner 5)(144:252:0.4)
    \node[regular polygon,draw,regular polygon sides=5, minimum size=2.3cm, thick] (p) at (0,0) {};
	\end{tikzpicture}
\end{subfigure}
\hfill
\begin{subfigure}[b]{0.18\textwidth}
	\centering
	\begin{tikzpicture}
 	\node[regular polygon,draw,regular polygon sides=8, minimum size=2.3cm, very thick] (p) at (0,0) {};
 	\draw[angleindicator] (p.side 1) -- (p.corner 5);
  	\draw[angleindicator] (p.side 2) -- (p.corner 7);
	\draw[angleindicator] (p.side 3) -- (p.corner 1);
	\draw[angleindicator] (p.side 4) -- (p.corner 1);
	\draw[angleindicator] (p.side 5) -- (p.corner 2);
	\draw[angleindicator] (p.side 6) -- (p.corner 2);
	\draw[angleindicator] (p.side 7) -- (p.corner 3);
	\draw[angleindicator] (p.side 8) -- (p.corner 3);
    \centerarc[crossangle](p.corner 1)(180:315:0.3)
    \centerarc[straightangle](p.corner 2)(225:360:0.3)
    \centerarc[straightangle](p.corner 3)(-90:45:0.3)
    \centerarc[crossangle](p.corner 4)(-45:90:0.3)
    \centerarc[crossangle](p.corner 5)(0:135:0.3)
    \centerarc[straightangle](p.corner 6)(45:180:0.3)
    \centerarc[straightangle](p.corner 7)(90:225:0.3)
    \centerarc[straightangle](p.corner 8)(135:270:0.3)
    \node[regular polygon,draw,regular polygon sides=8, minimum size=2.3cm, thick] (p) at (0,0) {};
	\end{tikzpicture}
\end{subfigure}
\caption{Some examples of decorated polygons, except for the second all are regular.}
\label{fig:decPolyEx}
\end{figure}

%% file: figureDecPolyToQuiver.tex
\begin{figure}[!hb]
\centering
\begin{subfigure}[b]{0.18\textwidth}
	\centering
	\begin{tikzpicture}
     	\node[regular polygon,draw,regular polygon sides=3, minimum size=2.5cm, very thick] (p) at (0,0) {};
      	\draw[angleindicator] (p.side 1) -- (p.corner 3);
      	\draw[angleindicator] (p.side 2) -- (p.corner 1);
      	\draw[angleindicator] (p.side 3) -- (p.corner 2);
        \orientedarc (p.corner 1:p.corner 2)
        \orientedarc (p.corner 2:p.corner 3)
        \orientedarc (p.corner 3:p.corner 1)
        \centerarc[straightangle](p.corner 1)(240:300:0.5)
        \centerarc[crossangle](p.corner 2)(0:60:0.5)
        \centerarc[crossangle](p.corner 3)(120:180:0.5)
        \node[regular polygon,draw,regular polygon sides=3, minimum size=2.5cm, thick] (p) at (0,0) {};
	\end{tikzpicture}
\end{subfigure}
\hfill
\begin{subfigure}[b]{0.18\textwidth}
    \centering
	\begin{tikzpicture}
 	\node[regular polygon,draw,regular polygon sides=3, minimum size=2.5cm, thick] (p) at (0,0) {};
  	\draw[angleindicator] (p.side 2) -- (p.corner 1);
  	\draw[angleindicator] (p.side 3) to[out=150,in=270] (p.corner 1);
  	\draw[angleindicator] (p.side 1) to[out=30,in=270] (p.corner 1);
    \orientedarc (p.corner 1:p.corner 2)
    \orientedarc (p.corner 2:p.corner 3)
    \orientedarc (p.corner 3:p.corner 1)
    \centerarc[straightangle](p.corner 1)(240:300:0.5)
    \centerarc[straightangle](p.corner 2)(0:60:0.5)
    \centerarc[straightangle](p.corner 3)(120:180:0.5)
    \node[regular polygon,draw,regular polygon sides=3, minimum size=2.5cm, thick] (p) at (0,0) {};
	\end{tikzpicture}
\end{subfigure}
\hfill
\begin{subfigure}[b]{0.18\textwidth}
	\centering
	\begin{tikzpicture}
     	\node[regular polygon,draw,regular polygon sides=4, minimum size=2.5cm, thick] (p) at (0,0) {};
     	\draw[angleindicator] (p.side 1) -- (p.corner 4);
      	\draw[angleindicator] (p.side 2) -- (p.corner 1);
    	\draw[angleindicator] (p.side 3) -- (p.corner 1);
    	\draw[angleindicator] (p.side 4) -- (p.corner 2);
        \orientedarc (p.corner 1:p.corner 2)
        \orientedarc (p.corner 3:p.corner 2)
        \orientedarc (p.corner 3:p.corner 4)
        \orientedarc (p.corner 4:p.corner 1)
        \centerarc[straightangle](p.corner 1)(180:270:0.5)
        \centerarc[straightangle](p.corner 2)(270:360:0.5)
        \centerarc[straightangle](p.corner 3)(0:90:0.5)
        \centerarc[crossangle](p.corner 4)(90:180:0.5)
        \node[regular polygon,draw,regular polygon sides=4, minimum size=2.5cm, thick] (p) at (0,0) {};
	\end{tikzpicture}
\end{subfigure}
\hfill
\begin{subfigure}[b]{0.18\textwidth}
	\centering
	\begin{tikzpicture}
     	\node[regular polygon,draw,regular polygon sides=5, minimum size=2.3cm, thick] (p) at (0,0) {};
     	\draw[angleindicator] (p.side 1) -- (p.corner 4);
      	\draw[angleindicator] (p.side 2) -- (p.corner 5);
    	\draw[angleindicator] (p.side 3) -- (p.corner 1);
    	\draw[angleindicator] (p.side 4) -- (p.corner 2);
    	\draw[angleindicator] (p.side 5) -- (p.corner 3);
        \orientedarc (p.corner 1:p.corner 2)
        \orientedarc (p.corner 2:p.corner 3)
        \orientedarc (p.corner 3:p.corner 4)
        \orientedarc (p.corner 4:p.corner 5)
        \orientedarc (p.corner 5:p.corner 1)
        \centerarc[straightangle](p.corner 1)(216:324:0.4)
        \centerarc[straightangle](p.corner 2)(288:396:0.4)
        \centerarc[straightangle](p.corner 3)(0:108:0.4)
        \centerarc[straightangle](p.corner 4)(72:180:0.4)
        \centerarc[straightangle](p.corner 5)(144:252:0.4)
        \node[regular polygon,draw,regular polygon sides=5, minimum size=2.3cm, thick] (p) at (0,0) {};
	\end{tikzpicture}
\end{subfigure}\\\vspace{.5cm}
\begin{subfigure}[b]{0.18\textwidth}
    \centering 
     \scalebox{0.5}{
    \begin{tikzpicture}
        \node[regular polygon,draw,regular polygon sides=3, minimum size=2*2.5cm, thick,dotted] (p) at (0,0) {};
        \pic[name=A,rotate=150] at (p.side 1) {NCnode};            
        \pic[name=B,rotate=270] at (p.side 2) {NCnode};
        \pic[name=C,rotate=30] at (p.side 3) {NCnode};
        
        \draw[barrow] (A-circle.0) to[out=0,in=180] (C-circle.180);
        \draw[xarrow] (C-circle.240) to[out=240,in=60] (B-circle.60);
        \draw[xarrow] (B-circle.120) to[out=120,in=300] (A-circle.300);
    \end{tikzpicture}
    }
\end{subfigure}\hfill
\begin{subfigure}[b]{0.18\textwidth}
    \centering
     \scalebox{0.5}{
    \begin{tikzpicture}
        \node[regular polygon,draw,regular polygon sides=3, minimum size=2*2.5cm, very thick,dotted] (p) at (0,0) {};
        \pic[name=A,rotate=150] at (p.side 1) {NCnode};     
        \pic[name=B,rotate=270] at (p.side 2) {NCnode};        
        \pic[name=C,rotate=30] at (p.side 3) {NCnode};
        
        \draw[xarrow] (C-circle.270) to[out=270,in=270] (A-circle.270);
        \draw[barrow] (C-circle.240) to[out=240,in=60] (B-circle.60);
        \draw[barrow] (B-circle.120) to[out=120,in=300] (A-circle.300);
    \end{tikzpicture}
    }
\end{subfigure}\hfill
\begin{subfigure}[b]{0.18\textwidth}
    \centering
     \scalebox{0.5}{
    \begin{tikzpicture}
        \node[regular polygon,draw,regular polygon sides=4, minimum size=2*2.5cm, very thick,dotted] (p) at (0,0) {};
        \pic[name=A,rotate=90] at (p.side 1) {NCnode};  
        \pic[name=B,rotate=0] at (p.side 2) {NCnode};        
        \pic[name=C,rotate=270] at (p.side 3) {NCnode};
        \pic[name=D,rotate=0] at (p.side 4) {NCnode};
        
        \draw[barrow] (B-circle.30) to[out=30,in=240] (A-circle.240);
        \draw[barrow] (A-circle.300) to[out=300,in=150] (D-circle.150);
        \draw[xarrow] (D-circle.210) to[out=210,in=60] (C-circle.60);
        \draw[barrow] (C-circle.120) to[out=120,in=330] (B-circle.330);
        \draw[xarrow] (C-circle.150) to[out=150,in=210] (A-circle.210);
        \draw[barrow] (D-circle.240) to[out=240,in=300] (B-circle.300);
    \end{tikzpicture}
    }
\end{subfigure}\hfill
\begin{subfigure}[b]{0.18\textwidth}
 \scalebox{0.5}{
        \begin{tikzpicture}
        \node[regular polygon,draw,regular polygon sides=5, minimum size=2*2.3cm, very thick,dotted] (p) at (0,0) {};        
        \pic[name=A,rotate=270-2*72] at (p.side 1) {NCnode};  
        \pic[name=B,rotate=270-72] at (p.side 2) {NCnode};        
        \pic[name=C,rotate=270] at (p.side 3) {NCnode};
        \pic[name=D,rotate=270+72] at (p.side 4) {NCnode};
        \pic[name=E,rotate=270+2*72] at (p.side 5) {NCnode};

        \draw[barrow] (C-circle.120) to[out=120,in=360-12] (B-circle.360-12);
        \draw[barrow] (B-circle.48) to[out=48,in=276] (A-circle.276);
        \draw[barrow] (A-circle.336) to[out=336,in=204] (E-circle.204);
        \draw[barrow] (E-circle.264) to[out=264,in=132] (D-circle.132);
        \draw[barrow] (D-circle.192) to[out=192,in=60] (C-circle.60);

        \draw[xarrow] (C-circle.90+60) to[out=90+60,in=306-60] (A-circle.306-60);
        \draw[xarrow] (B-circle.18+60) to[out=18+60,in=234-60] (E-circle.234-60);
        \draw[xarrow] (A-circle.6) to[out=306+60,in=162-60] (D-circle.162-60);
        \draw[xarrow] (E-circle.234+60) to[out=234+60,in=90-60] (C-circle.90-60);
        \draw[xarrow] (D-circle.162+60) to[out=162+60,in=18-60] (B-circle.18-60);
    \end{tikzpicture}
            }
\end{subfigure}

\caption{Translation from decorated polygons to quivers}
\label{fig:decPolyToQuiverEx}
\end{figure}

%% file: figureDecPolyTournProof.tex
\begin{figure}[h]
\centering
	\begin{minipage}{0.48\textwidth}
    \centering
	\scalebox{0.6}{
	\begin{tikzpicture}
		\begin{scope}[line width=0.8mm]
		 	\node[regular polygon,draw,regular polygon sides=4, minimum size=4.5cm] (p) at (0,0) {};
		 	\draw[angleindicator,line width=0.8mm] (p.side 1) -- (p.corner 4);
		  	\draw[angleindicator,line width=0.8mm] (p.side 2) -- (p.corner 1);
			\draw[angleindicator,line width=0.8mm] (p.side 3) -- (p.corner 1);
			\draw[angleindicator,line width=0.8mm] (p.side 4) -- (p.corner 2);
		\end{scope}
		\begin{scope}[shift={(2.1*\distL,0)}]
			\node[vertex] at (0,-0.2)	[label = \Huge$\leftrightarrow$]		{};
		\end{scope}
		\begin{scope}[shift={(4*\distL,0)}]
			\node[mutableBig](4) at (-\distL,0) {};
			\node[mutableBig](2) at (\distL,0) {};
			\node[mutableBig](3) at (0,-\distL) {};
			\node[mutableBig](1) at (0,\distL) {};
			\draw[barrow] (1) to (2);
			\draw[barrow] (2) to (3);
			\draw[barrow] (2) to (4);
			\draw[barrow] (3) to (1);
			\draw[barrow] (3) to (4);
			\draw[barrow] (4) to (1);
		\end{scope}
	\end{tikzpicture}}
	\caption{Decorated polygon with corresponding tournament.}
	\label{fig:decPolyToQuiver}
	\end{minipage}
	\hfill
	\begin{minipage}{0.48\textwidth}
    \centering
	\scalebox{0.6}{
	\begin{tikzpicture}
		\begin{scope}[line width=0.8mm]
			\node[regular polygon,draw,regular polygon sides=8, minimum size=4cm] (p) at (0,0) {};
		 	\draw[angleindicator,line width=0.8mm] (p.side 1) -- (p.corner 5);
		  	\draw[angleindicator,line width=0.8mm] (p.side 2) -- (p.corner 7);
			\draw[angleindicator,line width=0.8mm] (p.side 3) -- (p.corner 1);
			\draw[angleindicator,line width=0.8mm] (p.side 4) -- (p.corner 1);
			\draw[angleindicator,line width=0.8mm] (p.side 5) -- (p.corner 2);
			\draw[angleindicator,line width=0.8mm] (p.side 6) -- (p.corner 2);
			\draw[angleindicator,line width=0.8mm] (p.side 7) -- (p.corner 3);
			\draw[angleindicator,line width=0.8mm] (p.side 8) -- (p.corner 3);
            \draw[color=Red,line width=1.1mm] (p.corner 5) -- (p.corner 6);
		\end{scope}
		\begin{scope}[shift={(2.1*\distL,0)}]
			\node[vertex] at (0,-0.2)	[label = \Huge$\rightarrow$]		{};
		\end{scope}
		\begin{scope}[shift={(4*\distL,0)}, line width=0.8mm]
			\node[regular polygon,draw,regular polygon sides=4, minimum size=4.5cm] (p) at (0,0) {};
		 	\draw[angleindicator,line width=0.8mm] (p.side 1) to (p.corner 3);
		  	\draw[angleindicator,line width=0.8mm] (p.side 2) to[out=-70,in=60] (p.corner 3);
			\draw[angleindicator,line width=0.8mm] (p.side 3) to[out=160,in=30] (p.corner 3);
			\draw[angleindicator,line width=0.8mm] (p.side 4) to (p.corner 3);
			\node at (p.corner 3) [circle,draw=Red,fill=Red,inner sep=0pt,minimum size=3mm] {};
		\end{scope}
	\end{tikzpicture}}
	\caption{Collapsing to the incoming nodes of the marked node.}
	\label{fig:decPolyCollapse}
	\end{minipage}
\end{figure}

%% file: figureTournamentProof.tex
\begin{figure}[h]
\centering
\begin{tikzpicture}[thick]
    \node[circle] (u) at (0,0.5) {$u$};
    \node[circle] (v') at (1,0) {$v'$};
    \node[circle] (u') at (2,-0.5) {$u'$};
  
    \node[circle] (v) at (3,-2) {$v$};
  
    \node[circle] (w) at (4,-0.5) {$w$};
    \node[circle] (w') at (6,0.5) {$w'$};
  
    \draw[->] (u) -- (v');
    \draw[->] (v') -- (u');
  
    \draw[->] (u) -- (v.160);
    \draw[->] (v') -- (v);
    \draw[->] (u') -- (v.110);
  
    \draw[->] (v) -- (w);
    \draw[->] (v) -- (w'.225);

    \draw[->] (w'.190) -- (v');
    \draw[->] (v') -- (w.160);
  
    \draw[->,color=Orange] (w) -- (w');
    \draw[->,color=Orange] (u') -- (w);
    \draw[->,color=Orange] (w') -- (u);
\end{tikzpicture}
\caption{A locally transitive tournament, schematically.}
\label{fig:tournProof}
\end{figure}

%% file: figureCorrespondenceColorProof.tex
\begin{figure}[h]
    \centering
    \begin{subfigure}[b]{0.2\textwidth}
	    \centering
        \begin{tikzpicture}
            \def\size{1cm}
            
            \draw[very thick] (0,0) circle (\size);

            \path (0,0) ++(90:\size) coordinate (i);
            \node[above] at (i) {i};
            
            \path (0,0) ++(190:\size) coordinate (j);
            \node[below left] at (j) {j};
            
            \path (0,0) ++(350:\size) coordinate (k);
            \node[below right] at (k) {k};

            \path (0,0) ++(270:\size) coordinate (ie);
            \path (0,0) ++(135:\size) coordinate (ke);
            \path (0,0) ++(45:\size) coordinate (je);

            \draw[angleindicator] (i) -- (ie);
            \draw[angleindicator] (j) -- (je);
            \draw[angleindicator] (k) -- (ke);
        \end{tikzpicture}
    \end{subfigure}
    \hspace{.5cm}
    \begin{subfigure}[b]{0.2\textwidth}
	    \centering
        \begin{tikzpicture}
            \def\size{1cm}
            
            \draw[very thick] (0,0) circle (\size);

            \path (0,0) ++(90:\size) coordinate (i);
            \node[above] at (i) {i};
            
            \path (0,0) ++(190:\size) coordinate (j);
            \node[below left] at (j) {j};
            
            \path (0,0) ++(350:\size) coordinate (k);
            \node[below right] at (k) {k};

            \path (0,0) ++(270:\size) coordinate (ie);
            \path (0,0) ++(250:\size) coordinate (ke);
            \path (0,0) ++(290:\size) coordinate (je);

            \draw[angleindicator] (i) -- (ie);
            \draw[angleindicator] (j) -- (je);
            \draw[angleindicator] (k) -- (ke);
        \end{tikzpicture}
    \end{subfigure}
    \caption{Possible configurations of three sides of a decorated polygon.}
    \label{fig:CorrespondenceColorProof}
\end{figure}

%% file: figureStandardDecPoly.tex
\begin{figure}[h]
\centering
\begin{subfigure}[b]{0.18\textwidth}
	\centering
	\begin{tikzpicture}
 	\node[regular polygon,draw,regular polygon sides=3, minimum size=2.5cm, very thick] (p) at (0,0) {};
  	\draw[angleindicator] (p.side 1) -- (p.corner 3);
  	\draw[angleindicator] (p.side 2) -- (p.corner 1);
  	\draw[angleindicator] (p.side 3) -- (p.corner 2);
    \centerarc[straightangle](p.corner 1)(240:300:0.5)
    \centerarc[straightangle](p.corner 2)(0:60:0.5)
    \centerarc[straightangle](p.corner 3)(120:180:0.5)
    \node[regular polygon,draw,regular polygon sides=3, minimum size=2.5cm, thick] (p) at (0,0) {};
	\end{tikzpicture}
\end{subfigure}
\hspace{0.7cm}
\begin{subfigure}[b]{0.18\textwidth}
	\centering
	\begin{tikzpicture}
 	\node[regular polygon,draw,regular polygon sides=5, minimum size=2.3cm, very thick] (p) at (0,0) {};
 	\draw[angleindicator] (p.side 1) -- (p.corner 4);
  	\draw[angleindicator] (p.side 2) -- (p.corner 5);
	\draw[angleindicator] (p.side 3) -- (p.corner 1);
	\draw[angleindicator] (p.side 4) -- (p.corner 2);
	\draw[angleindicator] (p.side 5) -- (p.corner 3);
    \centerarc[straightangle](p.corner 1)(216:324:0.4)
    \centerarc[straightangle](p.corner 2)(288:396:0.4)
    \centerarc[straightangle](p.corner 3)(0:108:0.4)
    \centerarc[straightangle](p.corner 4)(72:180:0.4)
    \centerarc[straightangle](p.corner 5)(144:252:0.4)
    \node[regular polygon,draw,regular polygon sides=5, minimum size=2.3cm, thick] (p) at (0,0) {};
	\end{tikzpicture}
\end{subfigure}
\hspace{0.7cm}
\begin{subfigure}[b]{0.18\textwidth}
	\centering
	\begin{tikzpicture}
 	\node[regular polygon,draw,regular polygon sides=7, minimum size=2.2cm, very thick] (p) at (0,0) {};
 	\draw[angleindicator] (p.side 1) -- (p.corner 5);
  	\draw[angleindicator] (p.side 2) -- (p.corner 6);
	\draw[angleindicator] (p.side 3) -- (p.corner 7);
	\draw[angleindicator] (p.side 4) -- (p.corner 1);
    \draw[angleindicator] (p.side 5) -- (p.corner 2);
	\draw[angleindicator] (p.side 6) -- (p.corner 3);
	\draw[angleindicator] (p.side 7) -- (p.corner 4);
    \centerarc[straightangle](p.corner 1)(1440/7:2340/7:0.3)
    \centerarc[straightangle](p.corner 2)(1800/7:2700/7:0.3)
    \centerarc[straightangle](p.corner 3)(2160/7:3060/7:0.3)
    \centerarc[straightangle](p.corner 4)(0:900/7:0.3)
    \centerarc[straightangle](p.corner 5)(360/7:1260/7:0.3)
    \centerarc[straightangle](p.corner 6)(720/7:1620/7:0.3)
    \centerarc[straightangle](p.corner 7)(1080/7:1980/7:0.3)
    \node[regular polygon,draw,regular polygon sides=7, minimum size=2.2cm, thick] (p) at (0,0) {};
	\end{tikzpicture}
\end{subfigure}
\caption{Standard decorated polygons.}
\label{fig:stdDecPolyEx}
\end{figure}

%% file: figureCommArrowRules.tex
\begin{figure}[h]
\centering
\begin{subfigure}[b]{0.3\textwidth}
  \centering
  \begin{tikzpicture}
  \begin{scope}[very thick]
    \def\sidelength{1cm}
    \coordinate (A) at (0,\sidelength);
    \coordinate (M) at (0.5*\sidelength,1.1*\sidelength);
    \coordinate (B) at (\sidelength,1.2*\sidelength);
    \coordinate (C) at (2*\sidelength,1.2*\sidelength);
    \coordinate (D) at (2*\sidelength,-1.2*\sidelength);
    \coordinate (E) at (\sidelength,-1.2*\sidelength);
    \coordinate (N) at (0.5*\sidelength,-1.1*\sidelength);
    \coordinate (F) at (0,-\sidelength);
    \coordinate (G) at (2.5*\sidelength,0);
    \coordinate (H) at (-0.7*\sidelength,0);
    \coordinate (I) at (-0.5*\sidelength,0.2*\sidelength);
    \coordinate (J) at (-0.5*\sidelength,-0.2*\sidelength);
    
    \draw (A) -- (B);
    \draw (C)--(D);
    \draw (F)--(E);
    \begin{scope}[angleindicator]
	    \draw (H)--(G);
	    \draw (M)--(J);
	    \draw (N)--(I);
	\end{scope}
   	\begin{scope}[dotted]
   		\draw (B)--(C);
   		\draw (E)--(D);
   		\draw (A)--+(-0.5*\sidelength,-0.2\sidelength);
   		\draw (F)--+(-0.5*\sidelength,0.2\sidelength);
   	\end{scope}
   	
   \node [above] (l1) at (M) {$B$};
   \node [below] (l2) at (N) {$C$};
   \node [right] (l3) at (C) {$A$};
  \end{scope}
  
  \begin{scope}[very thick,shift={(0,-2)}]
  	\def\sidelength{1cm}
  	\draw [->] (1*\sidelength,0) -- (1*\sidelength,-1*\sidelength)
  	node [right,midway]{\small flip at $A$};
  \end{scope}
  
  \begin{scope}[very thick,shift={(0,-5)}]
    \def\sidelength{1cm}
    \coordinate (A) at (0,\sidelength);
    \coordinate (M) at (0.5*\sidelength,1.1*\sidelength);
    \coordinate (B) at (\sidelength,1.2*\sidelength);
    \coordinate (C) at (2*\sidelength,1.2*\sidelength);
    \coordinate (D) at (2*\sidelength,-1.2*\sidelength);
    \coordinate (E) at (\sidelength,-1.2*\sidelength);
    \coordinate (N) at (0.5*\sidelength,-1.1*\sidelength);
    \coordinate (F) at (0,-\sidelength);
    \coordinate (G) at (2.5*\sidelength,0);
    \coordinate (H) at (-0.7*\sidelength,0);
    \coordinate (I) at (-0.5*\sidelength,0.2*\sidelength);
    \coordinate (J) at (-0.5*\sidelength,-0.2*\sidelength);
    
    \draw (A) -- (B);
    \draw (H)--(G);
    
    \draw (F)--(E);
    \begin{scope}[angleindicator]
    	\draw (C)--(D);
	    \draw (M)--(G);
	    \draw (N)--(G);
	\end{scope}
    \draw[->>,color=Green] (M)--(N);
   	\begin{scope}[dotted]
   		\draw (B)--(C);
   		\draw (E)--(D);
   		\draw (A)--+(-0.5*\sidelength,-0.2\sidelength);
   		\draw (F)--+(-0.5*\sidelength,0.2\sidelength);
   	\end{scope}
   	
   \node [above] (l1) at (M) {$B$};
   \node [below] (l2) at (N) {$C$};
   \node [above right] (l3) at (H) {$A'$};
  \end{scope}
  \end{tikzpicture}
  \caption{}
\end{subfigure}
\hfill
\begin{subfigure}[b]{0.3\textwidth}
  \centering
  \begin{tikzpicture}
  \begin{scope}[very thick]
    \def\sidelength{1cm}
    \coordinate (A) at (0,\sidelength);
    \coordinate (M) at (0.5*\sidelength,1.1*\sidelength);
    \coordinate (B) at (\sidelength,1.2*\sidelength);
    \coordinate (C) at (2*\sidelength,1.2*\sidelength);
    \coordinate (D) at (2*\sidelength,-1.2*\sidelength);
    \coordinate (E) at (\sidelength,-1.2*\sidelength);
    \coordinate (N) at (0.5*\sidelength,-1.1*\sidelength);
    \coordinate (F) at (0,-\sidelength);
    \coordinate (G) at (2.5*\sidelength,0);
    \coordinate (H) at (-0.7*\sidelength,0);
    \coordinate (I) at ($0.5*(B)+0.5*(C)$);
    \coordinate (J) at ($0.5*(E)+0.5*(D)$);
    
    \draw (A) -- (B);
    \draw (C)--(D);
    \draw (F)--(E);
    \begin{scope}[angleindicator]
	    \draw (H)--(G);
	    \draw (M)--(J);
	    \draw (N)--(I);
	\end{scope}
   	\begin{scope}[dotted]
   		\draw (B)--(C);
   		\draw (E)--(D);
   		\draw (A)--+(-0.5*\sidelength,-0.2\sidelength);
   		\draw (F)--+(-0.5*\sidelength,0.2\sidelength);
   	\end{scope}
   	
   \node [above] (l1) at (M) {$B$};
   \node [below] (l2) at (N) {$C$};
   \node [right] (l3) at (C) {$A$};
  \end{scope}
  
  \begin{scope}[very thick,shift={(0,-2)}]
  	\def\sidelength{1cm}
  	\draw [->] (1*\sidelength,0) -- (1*\sidelength,-1*\sidelength)
  	node [right,midway]{\small flip at $A$};
  \end{scope}
  
  \begin{scope}[very thick,shift={(0,-5)}]
    \def\sidelength{1cm}
    \coordinate (A) at (0,\sidelength);
    \coordinate (M) at (0.5*\sidelength,1.1*\sidelength);
    \coordinate (B) at (\sidelength,1.2*\sidelength);
    \coordinate (C) at (2*\sidelength,1.2*\sidelength);
    \coordinate (D) at (2*\sidelength,-1.2*\sidelength);
    \coordinate (E) at (\sidelength,-1.2*\sidelength);
    \coordinate (N) at (0.5*\sidelength,-1.1*\sidelength);
    \coordinate (F) at (0,-\sidelength);
    \coordinate (G) at (2.5*\sidelength,0);
    \coordinate (H) at (-0.7*\sidelength,0);
    \coordinate (I) at (-0.5*\sidelength,0.2*\sidelength);
    \coordinate (J) at (-0.5*\sidelength,-0.2*\sidelength);
    
    \draw (A) -- (B);
    \draw (H)--(G);
    
    \draw (F)--(E);
    \begin{scope}[angleindicator]
    	\draw (C)--(D);
	    \draw (M)--(G);
	    \draw (N)--(G);
	\end{scope}
   	\begin{scope}[dotted]
   		\draw (B)--(C);
   		\draw (E)--(D);
   		\draw (A)--+(-0.5*\sidelength,-0.2\sidelength);
   		\draw (F)--+(-0.5*\sidelength,0.2\sidelength);
   	\end{scope}
   	
   \node [above] (l1) at (M) {$B$};
   \node [below] (l2) at (N) {$C$};
   \node [above right] (l3) at (H) {$A'$};
  \end{scope}
  \end{tikzpicture}
  \caption{}
\end{subfigure}
\hfill
\begin{subfigure}[b]{0.3\textwidth}
  \centering
  \begin{tikzpicture}
  \begin{scope}[very thick]
    \def\sidelength{1cm}
    \coordinate (A) at (0,\sidelength);
    \coordinate (M) at (\sidelength,\sidelength);
    \coordinate (B) at (2*\sidelength,\sidelength);
    \coordinate (C) at (2*\sidelength,1.2*\sidelength);
    \coordinate (D) at (2*\sidelength,-1.2*\sidelength);
    \coordinate (E) at (2*\sidelength,-1.2*\sidelength);
    \coordinate (N) at (\sidelength,-1.2*\sidelength);
    \coordinate (F) at (0,-1.2*\sidelength);
    
    \draw (A) -- (B);
    \draw (F)--(E);
    \begin{scope}[angleindicator]
	    \draw (M)--(F);
	    \draw (N)--(A);
	\end{scope}
    \draw[->>,color=Green] (N)--(M);
   	\begin{scope}[dotted]
   		\draw (B)--+(0.5*\sidelength,-0.2\sidelength);
   		\draw (E)--+(0.5*\sidelength,0.2\sidelength);
   		\draw (A)--+(-0.5*\sidelength,-0.2\sidelength);
   		\draw (F)--+(-0.5*\sidelength,0.2\sidelength);
   	\end{scope}
   	
   \node [above] (l1) at (M) {$B$};
   \node [below] (l2) at (N) {$C$};
  \end{scope}
  
  \begin{scope}[very thick,shift={(0,-2)}]
  	\def\sidelength{1cm}
  	\draw [->] (1*\sidelength,0) -- (1*\sidelength,-1*\sidelength)
  	node [right,midway]{\small cancellation};
  \end{scope}
  
  \begin{scope}[very thick,shift={(0,-5)}]
     \def\sidelength{1cm}
    \coordinate (A) at (0,\sidelength);
    \coordinate (M) at (\sidelength,\sidelength);
    \coordinate (B) at (2*\sidelength,\sidelength);
    \coordinate (C) at (2*\sidelength,1.2*\sidelength);
    \coordinate (D) at (2*\sidelength,-1.2*\sidelength);
    \coordinate (E) at (2*\sidelength,-1.2*\sidelength);
    \coordinate (N) at (\sidelength,-1.2*\sidelength);
    \coordinate (F) at (0,-1.2*\sidelength);
    
    \draw (A) -- (B);
    \draw (F)--(E);
    \begin{scope}[angleindicator]
	    \draw (M)--(E);
	    \draw (N)--(B);
	\end{scope}
   	\begin{scope}[dotted]
   		\draw (B)--+(0.5*\sidelength,-0.2\sidelength);
   		\draw (E)--+(0.5*\sidelength,0.2\sidelength);
   		\draw (A)--+(-0.5*\sidelength,-0.2\sidelength);
   		\draw (F)--+(-0.5*\sidelength,0.2\sidelength);
   	\end{scope}
   	
   \node [above] (l1) at (M) {$B$};
   \node [below] (l2) at (N) {$C$};
  \end{scope}
  \end{tikzpicture}
  \caption{}
\end{subfigure}
\caption{Rules for introduction and cancellation of commutative arrows.}
\label{fig:commArrowRules}
\end{figure}

%% file: figureDecPentaFlip.tex
\begin{figure}[ht]
\centering
	\begin{tikzpicture}[very thick]
		\begin{scope}
			\node[regular polygon,draw,regular polygon sides=5, minimum size=4cm] (p) at (0,0) {};
			\draw (p.corner 1) -- (p.corner 4);
			\coordinate (m1) at ($(p.corner 1)!0.5!(p.corner 4)$);
			\draw[angleindicator] (m1) -- (p.corner 2);
			\draw[angleindicator] (m1) -- (p.corner 5);
			
		 	\draw[angleindicator] (p.side 1) -- (p.corner 3);
		 	\draw[angleindicator] (p.side 2) -- (p.corner 1);
		 	\draw[angleindicator] (p.side 3) -- (p.corner 2);
		 	\draw[angleindicator] (p.side 4) -- (p.corner 1);
		 	\draw[angleindicator] (p.side 5) -- (p.corner 4);

            \orientedarc (p.corner 1:p.corner 2)
            \orientedarc (p.corner 2:p.corner 3)
            \orientedarc (p.corner 3:p.corner 4)
            \orientedarc (p.corner 5:p.corner 4)
            \orientedarc (p.corner 1:p.corner 5)
            \orientedarc (p.corner 4:p.corner 1)
            
            \centerarc[straightangle](p.corner 1)(216:288:0.5)
            \centerarc[straightangle](p.corner 1)(288:324:0.5)
            \centerarc[straightangle](p.corner 2)(288:396:0.5)
            \centerarc[straightangle](p.corner 3)(0:108:0.5)
            \centerarc[crossangle](p.corner 4)(108:180:0.5)
            \centerarc[straightangle](p.corner 4)(72:108:0.5)
            \centerarc[straightangle](p.corner 5)(144:252:0.5)

            \draw[carrow] (p.side 5) to (p.side 1);
		\end{scope}
		\begin{scope}[shift={(2.25*\distL,0)}]
			\draw[<->] (-0.5,0) -- (0.5,0);
		\end{scope}
		\begin{scope}[shift={(4.5*\distL,0)}]
			\node[regular polygon,draw,regular polygon sides=5, minimum size=4cm] (p) at (0,0) {};
			\draw (p.corner 2) -- (p.corner 5);
			\coordinate (m1) at ($(p.corner 2)!0.5!(p.corner 5)$);
			\draw[angleindicator] (m1) -- (p.corner 1);
			\draw[angleindicator] (m1) -- (p.corner 4);
			
		 	\draw[angleindicator] (p.side 1) to[out=6,in=236] (p.corner 1);
		 	\draw[angleindicator] (p.side 2) -- (p.corner 5);
		 	\draw[angleindicator] (p.side 3) -- (p.corner 5);
		 	\draw[angleindicator] (p.side 4) -- (p.corner 2);
		 	\draw[angleindicator] (p.side 5) to[out=174,in=-56] (p.corner 1);

            \orientedarc (p.corner 1:p.corner 2)
            \orientedarc (p.corner 2:p.corner 3)
            \orientedarc (p.corner 3:p.corner 4)
            \orientedarc (p.corner 5:p.corner 4)
            \orientedarc (p.corner 1:p.corner 5)
            \orientedarc (p.corner 2:p.corner 5)
            
            \centerarc[straightangle](p.corner 1)(216:324:0.5)
            \centerarc[straightangle](p.corner 2)(288:396:0.5)
            \centerarc[straightangle](p.corner 3)(0:108:0.5)
            \centerarc[crossangle](p.corner 4)(72:180:0.5)
            \centerarc[straightangle](p.corner 5)(144:252:0.5)
			
			\draw[carrow] (p.side 1) to (p.side 3);
		\end{scope}
	\end{tikzpicture}	
\caption{Performing a flip on a decorated tiling of a pentagon.}
\label{fig:decPentaFlip}
\end{figure}

%% file: figureAngleSkeleton.tex
\begin{figure}[h]
\centering
	\begin{tikzpicture}[very thick]
		\begin{scope}
			\node[regular polygon,draw,regular polygon sides=5, minimum size=4cm] (p) at (0,0) {};
			\draw (p.corner 1) -- (p.corner 4);
			\coordinate (m1) at ($(p.corner 1)!0.5!(p.corner 4)$);
			\draw[angleindicator] (m1) -- (p.corner 2);
			\draw[angleindicator] (m1) -- (p.corner 5);
			
		 	\draw[angleindicator] (p.side 1) -- (p.corner 3);
		 	\draw[angleindicator] (p.side 2) -- (p.corner 1);
		 	\draw[angleindicator] (p.side 3) -- (p.corner 2);
		 	\draw[angleindicator] (p.side 4) -- (p.corner 1);
		 	\draw[angleindicator] (p.side 5) -- (p.corner 4);

            \centerarc[crossangle](p.corner 1)(216:288:0.4)
            \centerarc[straightangle](p.corner 1)(288:324:0.4)
            \centerarc[straightangle](p.corner 2)(288:396:0.4)
            \centerarc[straightangle](p.corner 3)(0:108:0.4)
            \centerarc[straightangle](p.corner 4)(108:180:0.4)
            \centerarc[crossangle](p.corner 4)(72:108:0.4)
            \centerarc[crossangle](p.corner 5)(144:252:0.4)

            \orientedarc(p.corner 5:p.corner 1)
            \orientedarc(p.corner 4:p.corner 5)
            \orientedarc(p.corner 1:p.corner 4)
            \orientedarc(p.corner 2:p.corner 1)
            \orientedarc(p.corner 3:p.corner 2)
            \orientedarc(p.corner 4:p.corner 3)

            \node[] () at ($(p.corner 1)!0.5!(p.corner 2)$) [label=above:$A_1$]{};
            \node[] () at ($(p.corner 1)!0.5!(p.corner 4)$) [label=left:$A_2$]{};
            \node[] () at ($(p.corner 3)!0.5!(p.corner 4)$) [label=below:$A_3$]{};
            \node[] () at ($(p.corner 2)!0.5!(p.corner 3)$) [label=left:$A_4$]{};
            \node[] () at ($(p.corner 4)!0.5!(p.corner 5)$) [label=right:$A_6$]{};
            \node[] () at ($(p.corner 1)!0.5!(p.corner 5)$) [label=above:$A_5$]{};
		\end{scope}
		\begin{scope}[shift={(2.25*\distL,0)}]
			\draw[->] (-0.5,0) -- (0.5,0);
		\end{scope}
		\begin{scope}[shift={(4.5*\distL,0)}]
			\node[regular polygon,draw,regular polygon sides=4, minimum size=4cm, color=gray] (p) at (0,0) {};
            \coordinate (q) at ($(p.corner 1)!0.5!(p.corner 4)+(2.45,0)$);
            \draw[color=gray] (p.corner 1) -- (q);
            \draw[color=gray] (p.corner 4) -- (q);

            \centerarc[crossangle](p.corner 1)(180:270:0.71)
            \centerarc[straightangle](p.corner 2)(270:360:0.71)
            \centerarc[straightangle](p.corner 3)(0:90:0.71)
            \centerarc[straightangle](p.corner 4)(90:180:0.71)
            \centerarc[crossangle](p.corner 4)(30:90:0.71)
            \centerarc[crossangle](q)(150:210:0.71)
            \centerarc[straightangle](p.corner 1)(270:330:0.71)
            
            \node[mutable]	(a) at ($(p.corner 1)!0.25!(p.corner 2)$)	[label=above:$v_1$] {};
            \node[mutable]	(b) at ($(p.corner 1)!0.75!(p.corner 2)$)	[label=above:$v_0$] {};
            \node[mutable]	(c) at ($(p.corner 2)!0.25!(p.corner 3)$)	[] {};
            \node[mutable]	(d) at ($(p.corner 2)!0.75!(p.corner 3)$)	[] {};
            \node[mutable]	(e) at ($(p.corner 3)!0.25!(p.corner 4)$)	[] {};
            \node[mutable]	(f) at ($(p.corner 3)!0.75!(p.corner 4)$)	[label=below:$v_2$] {};
            \node[mutable]	(g) at ($(p.corner 4)!0.25!(p.corner 1)$)	[] {};
            \node[mutable]	(h) at ($(p.corner 4)!0.25!(q)$)	        [] {};
            \node[mutable]	(i) at ($(p.corner 4)!0.75!(q)$)	        [] {};
            \node[mutable]	(j) at ($(q)!0.25!(p.corner 1)$)            [] {};
            \node[mutable]	(k) at ($(q)!0.75!(p.corner 1)$)	        [] {};
            \node[mutable]	(l) at ($(p.corner 4)!0.75!(p.corner 1)$)	[] {};

            \orientedarc(b:a)
            \orientedarc(d:c)
            \orientedarc(f:e)
            \orientedarc(l:g)
            \orientedarc(j:k)
            \orientedarc(h:i)

            \node[] () at ($(p.corner 1)!0.5!(p.corner 2)$) [label=above:$A_1$]{};
            \node[] () at ($(p.corner 1)!0.5!(p.corner 4)$) [label=left:$A_2$]{};
            \node[] () at ($(p.corner 3)!0.5!(p.corner 4)$) [label=below:$A_3$]{};
            \node[] () at ($(p.corner 2)!0.5!(p.corner 3)$) [label=left:$A_4$]{};
            \node[] () at ($(p.corner 4)!0.5!(q)$) [label=below:$A_6$]{};
            \node[] () at ($(p.corner 1)!0.5!(q)$) [label=above:$A_5$]{};
		\end{scope}
	\end{tikzpicture}	
\caption{A decorated tiling of the disk and the associated angle skeleton.}
\label{fig:AngleSkeleton}
\end{figure}

%% file: figureAngleVsCanonicalAngle.tex
\begin{figure}[h]
    \centering
	\begin{tikzpicture}[very thick]
 	    \node[regular polygon,draw,regular polygon sides=5, minimum size=3cm, very thick] (p) at (0,0) {};
        \draw[angleindicator] (p.side 3) -- (p.corner 1);

        \centerarc[straightangle](p.corner 1)(216:324:0.4)
        \centerarc[straightangle](p.corner 2)(288:396:0.4)
        \centerarc[straightangle](p.corner 3)(0:108:0.4)
        \centerarc[straightangle](p.corner 4)(72:180:0.4)
        \centerarc[straightangle](p.corner 5)(144:252:0.4)

        \node[regular polygon,draw,regular polygon sides=5, minimum size=3cm, thick] (p) at (0,0) {};

        \draw[color=Purple, ultra thick] ($(p.corner 1)!0.25!(p.corner 5)$) -- ($(p.corner 1)!0.79!(p.corner 5)$);
        \centerarc[color=Purple, ultra thick](p.corner 5)(144:252:0.4)
        \draw[->, color=Purple, ultra thick,>=Triangle] ($(p.corner 5)!0.21!(p.corner 4)$) -- ($(p.corner 5)!0.6!(p.corner 4)$);
        \draw[-, color=Purple, ultra thick] ($(p.corner 5)!0.21!(p.corner 4)$) -- ($(p.corner 5)!0.79!(p.corner 4)$);
        \centerarc[color=Purple, ultra thick](p.corner 4)(72:180:0.4)
        \draw[color=Purple, ultra thick] ($(p.corner 4)!0.21!(p.corner 3)$) -- ($(p.corner 4)!0.79!(p.corner 3)$);

        \node[mutable]	(v) at ($(p.corner 1)!0.25!(p.corner 5)$)	[label=above:$v$] {};
        \node[mutable]	(w) at ($(p.corner 3)!0.25!(p.corner 4)$)	[] {};

        \node[vertex]   (l) at (p.side 4) [label=below right:\textcolor{Purple}{$\gamma$}] {};
        \node[vertex]   (l) at (p.side 3) [label=below:$k$] {};
	\end{tikzpicture}
\caption{Relating canonical angles to angles for ST-compatible quivers.}
\label{fig:angleVsCanonicalAngle}
\end{figure}

%% file: figureMutationVariables.tex
\begin{figure}[!hb]
    \centering
    \begin{tikzpicture}
        \begin{scope}
            \coordinate (w) at (2,0);
            \coordinate (t) at (0,2);
            \coordinate (v) at (4,2);
            \coordinate (u) at (2,4);
            \node[mutable,label=below:$w_1$] (w1) at ($(w)!0.25!(t)$) {};
            \node[mutable,label=below:$w_2$] (w2) at ($(w)!0.25!(v)$) {};
            \node[mutable,label=$u_1$] (u1) at ($(v)!0.75!(u)$) {};
            \node[mutable,label=$u_2$] (u2) at ($(t)!0.75!(u)$) {};
            \node[mutable,label=175:$t$] (tM) at ($(t)!0.25!(v)$) {};
            \node[mutable,label=5:$v$] (vM) at ($(t)!0.75!(v)$) {};
            \node[mutable,label=left:$t_1$] (t1) at ($(w)!0.75!(t)$) {};
            \node[mutable,label=left:$t_2$] (t2) at ($(t)!0.25!(u)$) {};
            \node[mutable,label=right:$v_1$] (v1) at ($(v)!0.25!(u)$) {};
            \node[mutable,label=right:$v_2$] (v2) at ($(w)!0.75!(v)$) {};
            \orientedarc (t:w)
            \orientedarc (u:t)
            \orientedarc (v:w)
            \orientedarc (u:v)
            \orientedarc (t:v)
            \draw[] (tM) to[out=90,in=0] node[right] {$\alpha_t$} (t2);
            \draw[] (tM) to[out=270,in=0] node[right] {$\beta_t$} (t1);
            \draw[] (vM) to[out=90,in=180] node[left] {$\alpha_v$} (v1);
            \draw[] (vM) to[out=270,in=180] node[left] {$\beta_v$} (v2);
            \draw[] (u1) to[out=225,in=315] node[below left] {$\theta_u$} (u2);
            \draw[] (w1) to[out=45,in=135] node[above right] {$\theta_w$} (w2);

		    \draw[angleindicator] ($(t)!0.5!(v)$) -- (u);
  		    \draw[angleindicator] ($(t)!0.5!(v)$) -- (w);
            \node[label=below right:$\gamma_D$] (D) at ($(w)!0.5!(v)$) {};
            \node[label=above right:$\gamma_A$] (A) at ($(v)!0.5!(u)$) {};
            \node[label=above left:$\gamma_B$] (B) at ($(t)!0.5!(u)$) {};
            \node[label=below left:$\gamma_C$] (C) at ($(w)!0.5!(t)$) {};
            \node[label=$X$] (X) at ($(t)!0.5!(v)$) {};
        \end{scope}
        \begin{scope}[shift={(3.7*\distL,0)}]
			\draw[->, very thick] (-0.5,2) -- (0.5,2);
		\end{scope}
		\begin{scope}[shift={(4.5*\distL,0)}]
            \coordinate (w) at (2,0);
            \coordinate (t) at (0,2);
            \coordinate (v) at (4,2);
            \coordinate (u) at (2,4);
            \node[mutable,label=below:$w_1$] (w1) at ($(w)!0.25!(t)$) {};
            \node[mutable,label=below:$w_2$] (w2) at ($(w)!0.25!(v)$) {};
            \node[mutable,label=$u_1$] (u1) at ($(v)!0.75!(u)$) {};
            \node[mutable,label=$u_2$] (u2) at ($(t)!0.75!(u)$) {};
            \node[mutable,label=120:$u$] (uM) at ($(u)!0.25!(w)$) {};
            \node[mutable,label=290:$w$] (wM) at ($(u)!0.75!(w)$) {};
            \node[mutable,label=left:$t_1$] (t1) at ($(w)!0.75!(t)$) {};
            \node[mutable,label=left:$t_2$] (t2) at ($(t)!0.25!(u)$) {};
            \node[mutable,label=right:$v_1$] (v1) at ($(v)!0.25!(u)$) {};
            \node[mutable,label=right:$v_2$] (v2) at ($(w)!0.75!(v)$) {};
            \orientedarc (t:w)
            \orientedarc (u:t)
            \orientedarc (v:w)
            \orientedarc (u:v)
            \orientedarc (u:w)
            \draw[] (uM) to[out=0,in=270] node[below] {$\alpha_u$} (u1);
            \draw[] (uM) to[out=180,in=270] node[below] {$\beta_u$} (u2);
            \draw[] (wM) to[out=180,in=90] node[above] {$\beta_w$} (w1);
            \draw[] (wM) to[out=0,in=90] node[above] {$\alpha_w$} (w2);
            \draw[] (t1) to[out=45,in=315] node[below right] {$\theta_t$} (t2);
            \draw[] (v1) to[out=225,in=135] node[above left] {$\theta_v$} (v2);

		    \draw[angleindicator] ($(u)!0.5!(w)$) -- (t);
  		    \draw[angleindicator] ($(u)!0.5!(w)$) -- (v);
            \node[label=below right:$\gamma_D$] (D) at ($(w)!0.5!(v)$) {};
            \node[label=above right:$\gamma_A$] (A) at ($(v)!0.5!(u)$) {};
            \node[label=above left:$\gamma_B$] (B) at ($(t)!0.5!(u)$) {};
            \node[label=below left:$\gamma_C$] (C) at ($(w)!0.5!(t)$) {};
            \node[label=left:$X'$] (X) at ($(u)!0.5!(w)$) {};
        \end{scope}
    \end{tikzpicture}
    \caption{General Mutation}
    \label{fig:GeneratlMutation}
\end{figure}

%% file: figureLaurentPhenomenaSetup.tex
\begin{figure}[!hb]
    \centering
    \scalebox{0.75}{
    \begin{tikzpicture}
        \begin{scope}
            \coordinate (w) at (4.5,0);
            \coordinate (t) at (2,4);
            \coordinate (v) at (7,4);
            \coordinate (u) at (4.5,8);
            \coordinate (s) at (0,7);
            \coordinate (A) at (7,6);
            \coordinate (B1) at ($(s)!0.5!(u)$);
            \coordinate (B2) at ($(s)!0.5!(t)$);
            
            \node[mutable,label=below:$w_1$] (w1) at ($(w)!0.25!(t)$) {};
            \node[mutable,label=below:$w_2$] (w2) at ($(w)!0.25!(v)$) {};
            \node[mutable,label=$u_1$] (u1) at ($(A)!0.75!(u)$) {};
            \node[mutable,label=$u_2$] (u2) at ($(B1)!0.75!(u)$) {};
            \node[mutable,label=175:$t$] (tM) at ($(t)!0.25!(v)$) {};
            \node[mutable,label=5:$v$] (vM) at ($(t)!0.75!(v)$) {};
            \node[mutable,label=left:$t_1$] (t1) at ($(w)!0.75!(t)$) {};
            \node[mutable,label=left:$t_2$] (t2) at ($(t)!0.25!(B2)$) {};
            \node[mutable,label=right:$v_1$] (v1) at ($(v)!0.25!(A)$) {};
            \node[mutable,label=right:$v_2$] (v2) at ($(w)!0.75!(v)$) {};
            \node[mutable] (y1) at ($(B1)!0.75!(B2)$) {};
            \node[mutable] (y2) at ($(B1)!0.25!(B2)$) {};
            \orientedarc (t:w)
            \orientedarc (u:B1)
            \orientedarc (B1:B2)
            \orientedarc (B2:t)
            \orientedarc (v:w)
            \orientedarc (u:A)
            \orientedarc (A:v)
            \orientedarc (t:v)
            \orientedarc (s:B1)
            \orientedarc (s:B2)
            \draw[] (tM) to[out=90,in=15] node[above] {$\alpha_t$} (t2);
            \draw[] (tM) to[out=270,in=0] node[right] {$\beta_t$} (t1);
            \draw[] (vM) to[out=90,in=180] node[above] {$\alpha_v$} (v1);
            \draw[] (vM) to[out=270,in=180] node[left] {$\beta_v$} (v2);
            \draw[] (u1) to[out=225,in=315] node[below left] {$\theta_u$} (u2);
            \draw[] (w1) to[out=45,in=135] node[above right] {$\theta_w$} (w2);

		    \draw[angleindicator] ($(t)!0.5!(v)$) -- (u);
  		    \draw[angleindicator] ($(t)!0.5!(v)$) -- (w);
            \draw[angleindicator] ($(B1)!0.5!(B2)$) -- (s);
  		    \draw[angleindicator] ($(B1)!0.5!(B2)$) -- (A);
            \node[label=below right:$\gamma_D$] (D) at ($(w)!0.5!(v)$) {};
            \node[label=above right:$\gamma_{A_1}$] (A1) at ($(u)!0.5!(A)$) {};
            \node[label=above right:$\gamma_{A_2}$] (A2) at ($(A)!0.5!(v)$) {};
            \node[label=above left:$\gamma_{B_1}$] (gB1) at ($(u)!0.5!(B1)$) {};
            \node[label=above left:$Y$] (Y) at ($(B1)!0.5!(B2)$) {};
            \node[label=above left:$\gamma_{B_2}$] (gB2) at ($(B2)!0.5!(t)$) {};
            \node[label=below left:$\gamma_C$] (C) at ($(w)!0.5!(t)$) {};
            \node[label=$X$] (X) at ($(t)!0.5!(v)$) {};
            \node[label=above:$\gamma_F$] (F) at ($(B1)!0.5!(s)$) {};
            \node[label=below:$\gamma_E$] (E) at ($(B2)!0.5!(s)$) {};

            \draw[carrow] ($(t)!0.5!(v)$) to node[right] {$k$} ($(B1)!0.5!(B2)$);
        \end{scope}

    \end{tikzpicture}
    }
    \scalebox{0.75}{
    \begin{tikzpicture}
        \begin{scope}[]
            \coordinate (w) at (4.5,0);
            \coordinate (t) at (2,4);
            \coordinate (v) at (7,4);
            \coordinate (u) at (4.5,8);
            \coordinate (s) at (0,7);
            \coordinate (A) at (7,6);
            \coordinate (B1) at ($(s)!0.5!(u)$);
            \coordinate (B2) at ($(s)!0.5!(t)$);
            
            \node[mutable,label=below:$w_1$] (w1) at ($(w)!0.25!(t)$) {};
            \node[mutable,label=below:$w_2$] (w2) at ($(w)!0.25!(v)$) {};
            \node[mutable,label=$u_1$] (u1) at ($(A)!0.75!(u)$) {};
            \node[mutable,label=$u_2$] (u2) at ($(B1)!0.75!(u)$) {};
            \node[mutable,label=175:$t$] (tM) at ($(t)!0.25!(v)$) {};
            \node[mutable,label=5:$v$] (vM) at ($(t)!0.75!(v)$) {};
            \node[mutable,label=left:$t_1$] (t1) at ($(w)!0.75!(t)$) {};
            \node[mutable,label=left:$t_2$] (t2) at ($(t)!0.25!(B2)$) {};
            \node[mutable,label=right:$v_1$] (v1) at ($(v)!0.25!(A)$) {};
            \node[mutable,label=right:$v_2$] (v2) at ($(w)!0.75!(v)$) {};
            \node[mutable] (y1) at ($(s)!0.75!(A)$) {};
            \node[mutable] (y2) at ($(s)!0.25!(A)$) {};
            \orientedarc (t:w)
            \orientedarc (u:B1)
            \orientedarc (A:s)
            \orientedarc (B2:t)
            \orientedarc (v:w)
            \orientedarc (u:A)
            \orientedarc (A:v)
            \orientedarc (t:v)
            \orientedarc (s:B1)
            \orientedarc (s:B2)
            \draw[] (tM) to[out=90,in=15] node[above] {$\alpha_t$} (t2);
            \draw[] (tM) to[out=270,in=0] node[right] {$\beta_t$} (t1);
            \draw[] (vM) to[out=90,in=180] node[above] {$\alpha_v$} (v1);
            \draw[] (vM) to[out=270,in=180] node[left] {$\beta_v$} (v2);
            \draw[] (u1) to[out=225,in=315] node[below left] {$\theta_u$} (u2);
            \draw[] (w1) to[out=45,in=135] node[above right] {$\theta_w$} (w2);

		    \draw[angleindicator] ($(t)!0.5!(v)$) -- (s);
  		    \draw[angleindicator] ($(t)!0.5!(v)$) -- (w);
            \draw[angleindicator] ($(s)!0.5!(A)$) -- (B1);
  		    \draw[angleindicator] ($(s)!0.5!(A)$) -- (B2);
            \node[label=below right:$\gamma_D$] (D) at ($(w)!0.5!(v)$) {};
            \node[label=above right:$\gamma_{A_1}$] (A1) at ($(u)!0.5!(A)$) {};
            \node[label=above right:$\gamma_{A_2}$] (A2) at ($(A)!0.5!(v)$) {};
            \node[label=above left:$\gamma_{B_1}$] (gB1) at ($(u)!0.5!(B1)$) {};
            \node[label=above:$Y'$] (Y) at ($(s)!0.5!(A)$) {};
            \node[label=above left:$\gamma_{B_2}$] (gB2) at ($(B2)!0.5!(t)$) {};
            \node[label=below left:$\gamma_C$] (C) at ($(w)!0.5!(t)$) {};
            \node[label=$X$] (X) at ($(t)!0.5!(v)$) {};
            \node[label=above:$\gamma_F$] (F) at ($(B1)!0.5!(s)$) {};
            \node[label=below:$\gamma_E$] (E) at ($(B2)!0.5!(s)$) {};

            \draw[carrow] ($(s)!0.5!(A)$) to[out=300,in=45] node[right] {$k$} ($(t)!0.5!(v)$);
            \draw[carrow] ($(t)!0.5!(v)$) to[out=180-15,in=180] node[right] {$k$} ($(u)!0.5!(A)$);
            \draw[carrow] ($(t)!0.5!(v)$) to[out=180-15,in=0] node[left] {$k$} ($(B2)!0.6!(s)$);
            \draw[carrow] ($(t)!0.5!(v)$) to[out=180-15,in=270] node[left,pos=0.6] {$k+1$} ($(u)!0.5!(B1)$);
        \end{scope}
    \end{tikzpicture}
    }
    \caption{Setup for proving Laurent Phenomena.}
    \label{fig:LaurentPhenomenaSetup}
\end{figure}

%% file: figureDigon.tex
\begin{figure}
\centering

\begin{subfigure}[b]{0.3\textwidth}
  \centering
  \scalebox{0.85}{
  \begin{tikzpicture}[thick]
    \def\circleRadius{2cm}
    \coordinate (p1) at (0,\circleRadius);
    \coordinate (p2) at (0,0);
    \coordinate (p3) at (0,-\circleRadius);

    \coordinate (q1) at (-\circleRadius,0);
    \coordinate (q2) at (\circleRadius,0);

    \draw[angleindicator] ($(p2)!0.5!(p1)$) to[out=0,in=30] (p3);
    \draw[angleindicator] ($(p2)!0.5!(p1)$) to[out=180,in=150] (p3);
    \draw[angleindicator] ($(p2)!0.5!(p3)$) to[out=0,in=-30] (p1);
    \draw[angleindicator] ($(p2)!0.5!(p3)$) to[out=180,in=-150] (p1);
    \draw[angleindicator] (p2) to (q1);
    \draw[angleindicator] (p2) to (q2);
    
    \draw (p2) circle [radius=\circleRadius];
    \fill (p1) circle [radius=2pt];
    \fill (p2) circle [radius=2pt];
    \fill (p3) circle [radius=2pt];

    \draw[->, thick,>=Triangle] (p2) -- ($(p2)!0.7!(p1)$);
    \draw[-,thick] (p2) -- (p1);
    \node[vertex] () at ($(p2)!0.7!(p1)$) [label = left: $X$]{};
    \draw[->, thick,>=Triangle] (p3) -- ($(p3)!0.4!(p2)$);
    \draw[-,thick] (p3) -- (p2);
    \node[vertex] () at ($(p3)!0.4!(p2)$) [label = left: $Y$]{};

    \draw[->, thick,>=Triangle] ($(q1)+(0,0.1)$) -- ($(q1)+(0,-0.1)$);
    \node[vertex] () at (q1) [label = left: $A$]{};

    \draw[->, thick,>=Triangle] ($(q2)+(0,0.1)$) -- ($(q2)+(0,-0.1)$);
    \node[vertex] () at (q2) [label = right: $B$]{};

    \centerarc[crossangle](p2)(270:450:0.5)
    \centerarc[crossangle](p3)(10:90:0.5)
  \end{tikzpicture}}
\end{subfigure}
\hfill
\begin{subfigure}[b]{0.3\textwidth}
  \centering
  \scalebox{0.85}{
  \begin{tikzpicture}[thick]
    \def\circleRadius{2cm}
    \coordinate (p1) at (0,\circleRadius);
    \coordinate (p2) at (0,0);
    \coordinate (p3) at (0,-\circleRadius);

    \coordinate (q1) at (-\circleRadius,0);
    \coordinate (q2) at (\circleRadius,0);

    \draw[angleindicator] (q1) to[out=60,in=170] ($(p2)!0.75!(p1)$) to[out=-10,in=90]
    ($(p2)!0.75!(q2)$) to[out=270,in=10] (p3);
    \draw[angleindicator] (q2) to[out=120,in=10] ($(p2)!0.75!(p1)$) to[out=190,in=90]
    ($(p2)!0.75!(q1)$) to[out=270,in=170] (p3);
    \draw[angleindicator] ($(p2)!0.5!(p1)$) to (p1);
    
    \draw[angleindicator] ($(p2)!0.5!(p1)$) to (p2);

    \draw (p2) circle [radius=\circleRadius];
    \fill (p1) circle [radius=2pt];
    \fill (p2) circle [radius=2pt];
    \fill (p3) circle [radius=2pt];
    
    \draw[->, thick,>=Triangle] (p3) -- ($(p3)!0.5!(p2)$);
    \draw[-,thick] (p3) -- (p2);
    \node[vertex] () at ($(p3)!0.5!(p2)$) [label = left: $Y$]{};
    
    \draw ($(p2)!0.5!(p1)$) to[out=0,in=30] (p3);
    \draw ($(p2)!0.5!(p1)$) to[out=180,in=150] (p3);
    \draw[->, thick,>=Triangle] ($(p2)!0.5!(p1)+(-0.1,0)$) -- ($(p2)!0.5!(p1)+(0.1,0)$);
    \node[vertex] () at ($(p2)!0.5!(p1)$) [label=below: $Z$]{};

    \draw[->, thick,>=Triangle] ($(q1)+(0,0.1)$) -- ($(q1)+(0,-0.1)$);
    \node[vertex] () at (q1) [label = left: $A$]{};

    \draw[->, thick,>=Triangle] ($(q2)+(0,0.1)$) -- ($(q2)+(0,-0.1)$);
    \node[vertex] () at (q2) [label = right: $B$]{};

    \centerarc[crossangle](p2)(-90:270:0.3)
    \centerarc[crossangle](p3)(42:90:0.5)
  \end{tikzpicture}}
\end{subfigure}
\hfill
\begin{subfigure}[b]{0.3\textwidth}
  \centering
  \scalebox{0.85}{
  \begin{tikzpicture}[thick]
    \def\circleRadius{2cm}
    \coordinate (p1) at (0,\circleRadius);
    \coordinate (p2) at (0,0);
    \coordinate (p3) at (0,-\circleRadius);

    \coordinate (q1) at (-\circleRadius,0);
    \coordinate (q2) at (\circleRadius,0);

    \draw[angleindicator] (q1) to[out=60,in=170] ($(p2)!0.75!(p1)$) to[out=-10,in=90]
    ($(p2)!0.75!(q2)$) to[out=270,in=10] (p3);
    \draw[angleindicator] (q2) to[out=120,in=10] ($(p2)!0.75!(p1)$) to[out=190,in=90]
    ($(p2)!0.75!(q1)$) to[out=270,in=170] (p3);
    \draw[angleindicator] ($(p2)!0.5!(p1)$) to (p1);
    
    \draw[angleindicator] ($(p2)!0.5!(p1)$) to (p2);

    \draw (p2) circle [radius=\circleRadius];
    \fill (p1) circle [radius=2pt];
    \fill (p2) circle [radius=2pt];
    \fill (p3) circle [radius=2pt];
    
    \draw[->, thick,>=Triangle] (p2) -- ($(p3)!0.5!(p2)$);
    \draw[-,thick] (p3) -- (p2);
    \node[vertex] () at ($(p3)!0.5!(p2)$) [label = left: $W$]{};
    
    \draw ($(p2)!0.5!(p1)$) to[out=0,in=30] (p3);
    \draw ($(p2)!0.5!(p1)$) to[out=180,in=150] (p3);
    \draw[->, thick,>=Triangle] ($(p2)!0.5!(p1)+(-0.1,0)$) -- ($(p2)!0.5!(p1)+(0.1,0)$);
    \node[vertex] () at ($(p2)!0.5!(p1)$) [label=below: $Z$]{};

    \draw[->, thick,>=Triangle] ($(q1)+(0,0.1)$) -- ($(q1)+(0,-0.1)$);
    \node[vertex] () at (q1) [label = left: $A$]{};

    \draw[->, thick,>=Triangle] ($(q2)+(0,0.1)$) -- ($(q2)+(0,-0.1)$);
    \node[vertex] () at (q2) [label = right: $B$]{};

    \centerarc[crossangle](p2)(-90:270:0.3)
    \centerarc[crossangle](p3)(90:138:0.5)
  \end{tikzpicture}}
\end{subfigure}

\medskip
\vspace{0.4cm}

\begin{subfigure}[b]{0.3\textwidth}
  \centering
  \scalebox{0.7}{
  \begin{tikzpicture}
    \pic[name=A,rotate=180] at (-2,0) {NCnodeHalfL};
    \node[vertex] () at (-2,0) [label ={[xshift=-6mm,yshift=-3mm]$A$}]{};
    
    \pic[name=B] at (0,2) {NCnode};
    \node[vertex] () at (0,2) [label ={[yshift=7mm]$X$}]{};
    
    \pic[name=C] at (0,-2) {NCnode};
    \node[vertex] () at (0,-2) [label ={[yshift=-12mm]$Y$}]{};
    
    \pic[name=D,rotate=180] at (2,0) {NCnodeHalfR};
    \node[vertex] () at (2,0) [label ={[xshift=6mm,yshift=-3mm]$B$}]{};
    
    \draw[arrow] (A-circle.45) to [out=45,in=135] (B-circle.135);
    \draw[arrow] (B-circle.45) to [out=45,in=135] (D-circle.135);
    \draw[xarrow] (D-circle.225) to [out=225,in=-45] (C-circle.-45);
    \draw[arrow] (C-circle.225) to [out=225,in=-45] (A-circle.-45);
    \draw[arrow] (B-circle.225) to [out=225,in=135] (C-circle.135);
    \draw[xarrow] (C-circle.45) to [out=45,in=-45] (B-circle.-45);
  \end{tikzpicture}}
  \caption{}
  \label{fig:digonA}
\end{subfigure}
\hfill
\begin{subfigure}[b]{0.3\textwidth}
  \centering
  \scalebox{0.7}{
  \begin{tikzpicture}
    \pic[name=A,rotate=180] at (-2,0) {NCnodeHalfL};
    \node[vertex] () at (-2,0) [label ={[xshift=-7mm,yshift=-3mm]$A$}]{};
    
    \pic[name=B,rotate=270] at (0,2) {NCnode};
    \node[vertex] () at (0,2) [label ={[yshift=7mm]$Z$}]{};
    
    \pic[name=C] at (0,-2) {NCnode};
    \node[vertex] () at (0,-2) [label ={[yshift=-12mm]$Y$}]{};
    
    \pic[name=D,rotate=180] at (2,0) {NCnodeHalfR};
    \node[vertex] () at (2,0) [label ={[xshift=7mm,yshift=-3mm]$B$}]{};
    
    \draw[arrow] (A-circle.45) to [out=45,in=135] (D-circle.135);
    \draw[arrow] (B-circle.135) to [out=135,in=-45] (A-circle.-45);
    \draw[arrow] (D-circle.225) to [out=225,in=45] (B-circle.45);
    \draw[arrow] (C-circle.225) to [out=225,in=225] (B-circle.225);
    \draw[xarrow] (B-circle.-45) to [out=-45,in=-45] (C-circle.-45);
    \draw[xarrow] (C-circle.45) .. controls +(45:1.1) and +(135:1) .. (C-circle.135);
  \end{tikzpicture}}
  \caption{}
  \label{fig:digonB}
\end{subfigure}
\hfill
\begin{subfigure}[b]{0.3\textwidth}
  \centering
  \scalebox{0.7}{
  \begin{tikzpicture}
    \pic[name=A,rotate=180] at (-2,0) {NCnodeHalfL};
    \node[vertex] () at (-2,0) [label ={[xshift=-7mm,yshift=-3mm]$A$}]{};
    
    \pic[name=B,rotate=270] at (0,2) {NCnode};
    \node[vertex] () at (0,2) [label ={[yshift=7mm]$Z$}]{};
    
    \pic[name=C,rotate=180] at (0,-2) {NCnode};
    \node[vertex] () at (0,-2) [label ={[yshift=-12mm]$W$}]{};
    
    \pic[name=D,rotate=180] at (2,0) {NCnodeHalfR};
    \node[vertex] () at (2,0) [label ={[xshift=7mm,yshift=-3mm]$B$}]{};
    
    \draw[arrow] (A-circle.45) to [out=45,in=135] (D-circle.135);
    \draw[arrow] (B-circle.135) to [out=135,in=-45] (A-circle.-45);
    \draw[arrow] (D-circle.225) to [out=225,in=45] (B-circle.45);
    \draw[xarrow] (C-circle.225) to [out=225,in=225] (B-circle.225);
    \draw[arrow] (B-circle.-45) to [out=-45,in=-45] (C-circle.-45);
    \draw[xarrow] (C-circle.135) .. controls +(135:1.1) and +(45:1) .. (C-circle.45);
  \end{tikzpicture}}
  \caption{}
  \label{fig:digonC}
\end{subfigure}
\caption{The punctured digon with mutations and corresponding 2-coloured quivers.}
\label{fig:digon}
\end{figure}

%% file: figureForkProof.tex
\begin{figure}[htb]
	\centering
		\begin{subfigure}[b]{0.4\textwidth}
			\centering
			\begin{tikzpicture}
				\node[mutableBig]	(B) at (0,0)					[label = above: $X$]	{};
				\node[mutableBig]	(A) at (-\distL,-\distL)		[label = left: $Z$]		{};
				\node[mutableBig]	(C) at (\distL,-\distL)		[label = right: $Y$]	{};
				\node[mutableBig]	(D) at (0,-2*\distL)			[label = below: $B$]	{};
				
				\draw[barrow] ($(A)+(0.05,-0.05)$) to ($(A)+(0.05,-0.05)+0.92*(B)-0.92*(A)$);
				\draw[carrow] ($(A)+(-0.05,0.05)$) to ($(A)+(-0.05,+0.05)+0.92*(B)-0.92*(A)$);
				\node[left, Green] (l1) at ($(A)+0.7*(B)-0.7*(A)+(-0.1,0.1)$) {\small $j$};
				
				\draw[barrow] ($(B)+(-0.05,-0.05)$) to ($0.91*(C)+0.91*(-0.05,-0.05)$);
				\draw[carrow] ($0.85*(0.05,0.05)$) to ($0.91*(C)+0.91*(0.05,0.05)$);
				
				\draw[xarrow] ($(A)+(0,0.03)$) to ($(C)+(-0.2,0.03)$);
				\draw[carrow] ($(A)+(0,-0.1)$) to ($(C)+(-0.2,-0.1)$);
				\node[above, Green] (l2) at ($(C)+0.25*(A)-0.25*(C)+(0,-0.55)$) {\small $k$};
				
				\draw[carrow] (D) to (B);
				\node[right, Green] (l3) at ($(B)+(0,-0.7)$) {\small $l$};
				
				\draw[barrow,color=Orange] (C) to (D);
				\draw[barrow,color=Orange] (A) to (D);
				
				\node[mutableBig]	(A2) at (A)	[]	{};
				\node[mutableBig]	(B2) at (B)	[]	{};
				\node[mutableBig]	(D2) at (D)	[]	{};
			\end{tikzpicture}
			\caption{}
            \label{fig:ForkProofA}
		\end{subfigure}
		\begin{subfigure}[b]{0.4\textwidth}
			\centering
			\begin{tikzpicture}
				\node[mutableBig]	(B) at (0,0)					[label = above: $B$]	{};
				\node[mutableBig]	(A) at (-\distL,-\distL)		[label = left: $X$]		{};
				\node[mutableBig]	(C) at (\distL,-\distL)		    [label = right: $Y$]	{};
				
				\draw[barrow] ($(A)+(0.05,-0.05)$) to ($(A)+(0.05,-0.05)+0.92*(B)-0.92*(A)$);
				\draw[barrow] ($(B)+(-0.05,-0.05)$) to ($0.91*(C)+0.91*(-0.05,-0.05)$);
				\draw[carrow] ($(A)+(-0.05,0.05)$) to ($(A)+(-0.05,+0.05)+0.92*(B)-0.92*(A)$);
				\draw[carrow] ($0.85*(0.05,0.05)$) to ($0.91*(C)+0.91*(0.05,0.05)$);
				\draw[carrow] (C) to (A);
				
				\node[above, Green] (l1) at ($(B)+(-0.7,-0.5)$) {\small $j$};
				\node[above, Green] (l2) at ($(C)+(-0.3,0.55)$) {\small $k$};
				\node[above, Green] (l3) at ($(A)+(0.7,-0.5)$) {\small $l$};
				
				\node[mutableBig]	(A2) at (A)	[]	{};
				\node[mutableBig]	(B2) at (B)	[]	{};
			\end{tikzpicture}
			\vspace{8mm}
			\caption{}
            \label{fig:ForkProofB}
		\end{subfigure}
	\caption{Cases that could lead to 2-cycles in an admissible fork.}
	\label{fig:ForkProof}
\end{figure}

%% file: figureFFPex.tex
\begin{figure}[htb]
	\centering
		\begin{subfigure}[b]{0.3\textwidth}
			\centering
			\begin{tikzpicture}
				\node[mutableBig]	(B) at (0,0)					[label = above: $B$]	{};
				\node[mutableBig]	(A) at (-\distL,-\distL)		[label = left: $C_1$]	{};
				\node[mutableBig]	(C) at (\distL,-\distL)		[label = right: $A$]	{};
				\node[mutableBig]	(D) at (0,-2*\distL)			[label = below: $C_2$]	{};
				
				\draw[xarrow] ($(A)+(0.05,-0.05)$) to ($(A)+(0.05,-0.05)+0.92*(B)-0.92*(A)$);
				\draw[carrow] ($(A)+(-0.05,0.05)$) to ($(A)+(-0.05,+0.05)+0.92*(B)-0.92*(A)$);
				\node[left, Green] (l1) at ($(A)+0.7*(B)-0.7*(A)+(-0.1,0.1)$) {\small 4};
				
				\draw[barrow] ($(B)+(-0.05,-0.05)$) to ($0.91*(C)+0.91*(-0.05,-0.05)$);
				\draw[carrow] ($0.85*(0.05,0.05)$) to ($0.91*(C)+0.91*(0.05,0.05)$);
				
				\draw[carrow] (C) to (A);
				\node[above, Green] (l2) at ($(C)+0.7*(A)-0.7*(C)$) {\small 3};
				
				\draw[xarrow] ($(D)+(-0.07,0.25)$) to ($(B)+(-0.07,-0.35)$);
				\draw[carrow] ($(D)+(0.07,0.25)$) to ($(B)+(0.07,-0.35)$);
				\node[right, Green] (l3) at ($(B)+(0.07,-0.9)$) {\small 4};
				
				\draw[carrow] (C) to (D);
				\node[right, Green] (l4) at ($(C)+0.7*(D)-0.7*(C)+(0.1,0)$) {\small 3};
				
				\draw[xarrow] ($(D)+(0.05,0.05)$) to ($0.09*(D)+0.09*(0.05,0.05)+0.91*(A)+0.91*(0.05,0.05)$);
				\draw[carrow] ($(D)+(-0.05,-0.05)$) to ($0.09*(D)+0.09*(-0.05,-0.05)+0.91*(A)+0.91*(-0.05,-0.05)$);
				
				\node[mutableBig]	(A2) at (A)	[]	{};
				\node[mutableBig]	(B2) at (B)	[]	{};
				\node[mutableBig]	(D2) at (D)	[]	{};
			\end{tikzpicture}
			\caption{}
		\end{subfigure}
		\hspace{15mm}
		\begin{subfigure}[b]{0.3\textwidth}
			\centering
			\begin{tikzpicture}[scale=1.3]
				\node[mutableBig]	(B) at (0,0)							[label = above: $B$]		{}; 
				\node[mutableBig]	(A) at (-\distL,-\distL)				[label = left: $C_1$]		{}; 
				\node[mutableBig]	(C) at (\distL,-\distL)				[label = right: $A$]		{}; 
				\node[mutableBig]	(D) at (0,-2*\distL)					[label = below: $C_3$]		{}; 
				\node[mutableBig]	(E) at (-1*\distL,-1.8*\distL)		[label = below left: $C_2$]	{}; 
				
				\draw[carrow] ($(A)+(0,0.1)$) to ($(B)+(-0.15,0.1)$);
				\draw[xarrow] ($(A)$) to ($(B)+(-0.15,0)$);
				\node[left, Green] (l1) at ($(A)+0.7*(B)-0.7*(A)+(-0.15,0.1)$) {\small 4};
				
				\draw[carrow] ($(E)+(-0.01,0.1)$) to ($(B)+(-0.1,-0.1)$);
				\draw[xarrow] ($(E)$) to ($(B)+(-0.05,-0.15)$);
				\node[left, Green] (l2) at ($(B)+(-0.5,-0.8)$) {\small 4};
				
				\draw[carrow] ($(D)+(0.07,0)$) to ($(B)+(0.07,-0.15)$);
				\draw[xarrow] ($(D)$) to ($(B)+(-0.0,-0.15)$);
				\node[left, Green] (l3) at ($(B)+(0.35,-0.8)$) {\small 4};
				
				\draw[carrow] ($(B)+(0.05,0.05)$) to ($(C)+(-0.05,0.15)$);
				\draw[barrow] ($(B)$) to ($(C)+(-0.1,0.1)$);
				
				\draw[carrow] ($(D)+(-0.05,0.13)$) to ($(A)+(0.12,-0.06)$);
				\draw[xarrow] ($(D)$) to ($(A)+(0.08,-0.12)$);
				
				\draw[carrow] ($(D)+(0,-0.08)$) to ($(E)+(0.14,-0.08)$);
				\draw[xarrow] ($(D)$) to ($(E)+(0.14,0)$);
				
				\draw[carrow] ($(E)+(-0.08,0)$) to ($(A)+(-0.08,-0.13)$);
				\draw[xarrow] ($(E)$) to ($(A)+(0,-0.13)$);
				
				\draw[carrow] ($(C)$) to ($(A)+(0.12,0)$);
				\draw[carrow] ($(C)$) to ($(D)+(0.12,0.06)$);
				\draw[carrow] ($(C)$) to ($(E)+(0.12,0.08)$);
				\node[left, Green] (l4) at ($(A)+(0.6,0.18)$) {\small 3};
				\node[left, Green] (l5) at ($(E)+(0.65,0.4)$) {\small 3};
				\node[left, Green] (l6) at ($(D)+(1,0.6)$) {\small 3};

				\node[mutableBig]	(A2) at (A)	[]	{};
				\node[mutableBig]	(B2) at (B)	[]	{};
				\node[mutableBig]	(C2) at (C)	[]	{};
				\node[mutableBig]	(D2) at (D)	[]	{};
				\node[mutableBig]	(E2) at (E)	[]	{};
			\end{tikzpicture}
			\caption{}
		\end{subfigure}

	\caption{Quivers with finite forkless part which provide fruitful seeds.}
	\label{fig:FFPex}
\end{figure}

%% file: figureNonExample.tex
\begin{figure}[ht]
\centering
	\begin{tikzpicture}[very thick]
		\begin{scope}
			\node[regular polygon,draw,regular polygon sides=7, minimum size=4cm] (p) at (0,0) {};
			\draw (p.corner 1) -- (p.corner 5);
			\draw (p.corner 2) -- (p.corner 4);
			\coordinate (m1) at ($(p.corner 1)!0.5!(p.corner 5)$);
			\coordinate (m2) at ($(p.corner 2)!0.5!(p.corner 4)$);
			\draw[angleindicator] (m1) -- (p.corner 2);
			\draw[angleindicator] (m1) -- (p.corner 7);
			\draw[angleindicator] (m2) -- (p.corner 1);
			\draw[angleindicator] (m2) -- (p.corner 3);
			
		 	\draw[angleindicator] (p.side 1) -- (p.corner 4);
		 	\draw[angleindicator] (p.side 2) -- (p.corner 4);
		 	\draw[angleindicator] (p.side 3) -- (p.corner 2);
		 	\draw[angleindicator] (p.side 4) -- (p.corner 2);
		 	\draw[angleindicator] (p.side 5) -- (p.corner 7);
		 	\draw[angleindicator] (p.side 6) -- (p.corner 1);
		 	\draw[angleindicator] (p.side 7) -- (p.corner 6);
		 	
		 	\draw[carrow] (m2) to (p.side 3);
		 	\draw[carrow] (p.side 3) to (m1);
		 	
		 	\node [right] (l1) at (m1) {$A$};
		 	\node [right] (l1) at (m2) {$B$};
		\end{scope}
		\begin{scope}[shift={(1.8*\distL,0)}]
			\coordinate (a1) at (0,0);
			\coordinate (a2) at (1.5*\distL,0);
			\draw [->] (a1)--(a2)
			node[above,midway,text centered, text width=2cm]{\small Flip at $A$, then at $B$};
		\end{scope}
		\begin{scope}[shift={(5*\distL,0)}]
			\node[regular polygon,draw,regular polygon sides=7, minimum size=4cm] (p) at (0,0) {};
			\draw (p.corner 2) -- (p.corner 7);
			\draw (p.corner 3) -- (p.corner 7);
			\coordinate (m1) at ($(p.corner 2)!0.5!(p.corner 7)$);
			\coordinate (m2) at ($(p.corner 3)!0.5!(p.corner 7)$);
			\draw[angleindicator] (m1) -- (p.corner 1);
			\draw[angleindicator] (m1) -- (p.corner 3);
			\draw[angleindicator] (m2) -- (p.corner 2);
			\draw[angleindicator] (m2) -- (p.corner 4);
			
			\draw[angleindicator] (p.side 1) -- (p.corner 7);
			\draw[angleindicator] (p.side 2) -- (p.corner 7);
			\draw[angleindicator] (p.side 3) -- (p.corner 7);
			\draw[angleindicator] (p.side 4) -- (p.corner 3);
			\draw[angleindicator] (p.side 5) -- (p.corner 3);
			\draw[angleindicator] (p.side 6) -- (p.corner 3);
			\draw[angleindicator] (p.side 7) -- (p.corner 2);
			
			\draw[carrow] (p.side 1) to ($0.01*(p.side 1)+0.99*(p.side 4)$);
			\draw[carrow] (p.side 3) to ($0.01*(p.side 3)+0.99*(p.side 7)$);
			\draw[carrow] ($(p.side 3)+(0.08,-0.02)$) to ($0.04*(p.side 3)+0.96*(p.side 4)$);
			\draw[carrow] ($(p.side 4)+(0.06,0)$) to ($(m1)+(0.04,0)$);
			\draw[carrow] (p.side 5) to ($0.01*(p.side 5)+0.99*(p.side 7)$);
			\draw[carrow] (m1) to ($0.03*(m1)+0.97*(p.side 3)+(-0.08,0.02)$);
			
		\end{scope}
	\end{tikzpicture}	
\caption{Two flips can lead to a non-admissible configuration of commutative arrows.}
\label{fig:NonEx}
\end{figure}

%% file: figureOrderedFGQuivers.tex
\begin{figure}[htb]
	\centering
		\begin{subfigure}[b]{0.48\textwidth}
			\centering
			\scalebox{0.7}{
			\begin{tikzpicture}[]
                \node[regular polygon,draw,regular polygon sides=3, minimum size=10cm, very thick, color=gray] (p) at (0,0) {};
				\pic[name=A,rotate=-30] at ($0.5*(p.corner 1)+0.5*(p.corner 2)$) {NCnode};
                \pic[name=B,rotate=0] at ($0.5*(p.corner 1)+0.5*(p.side 2)$) {NCnode};
                \pic[name=C,rotate=30] at ($0.5*(p.corner 1)+0.5*(p.corner 3)$) {NCnode};
                \pic[name=D,rotate=90] at ($0.65*(p.corner 3)+0.35*(p.corner 2)$) {NCnode};

                \draw[barrow] (A-circle.40) to[out=40,in=135] (B-circle.135);
                \draw[barrow] (B-circle.45) to[out=45,in=140] (C-circle.140);
                \draw[barrow] (C-circle.260) to[out=260,in=45] (D-circle.45);
                \draw[barrow] (D-circle.135) to[out=135,in=310] (B-circle.310);
                
                \node[mutable] (a) at ($0.2*(p.corner 1)+0.8*(p.corner 2)$) {};
                \node[mutable] (c) at ($0.2*(p.corner 1)+0.8*(p.corner 3)$) {};
                \node[mutable] (b) at ($0.5*(a)+0.5*(c)$) {};
                \node[mutable] (d) at ($0.65*(p.corner 2)+0.35*(p.corner 3)$) {};

                \draw[arrow] (a) to (b);
                \draw[arrow] (b) to (c);
                \draw[arrow] (b) to (d);
                \draw[arrow] (d) to (a);

                \draw[arrow, color=Green] (B-circle.center) to (b);
                \draw[arrow, color=Green] (b) to (A-circle.center);
                \draw[arrow, color=Green] (c) to (B-circle.center);

                \draw[halfarrow, bend right, color=Green] (A-circle.center) to (a);
                \draw[halfarrow, bend left, color=Green] (C-circle.center) to (c);
				\draw[halfarrow, bend angle = 25, bend right, color=Green] (d) to (D-circle.center);
			\end{tikzpicture}}
			\caption{}
           \label{fig:orderedFGQuiverB2}
		\end{subfigure}
		\hfill
		\begin{subfigure}[b]{0.48\textwidth}
			\centering
			\scalebox{0.7}{
			\begin{tikzpicture}[]
				\node[regular polygon,draw,regular polygon sides=3, minimum size=10cm, very thick, color=gray] (p) at (0,0) {};

                \node[] (h1) at ($0.67*(p.corner 2)+0.33*(p.corner 3)$) {};
                \node[] (h2) at ($0.67*(p.corner 3)+0.33*(p.corner 2)$) {};
    
				\pic[name=A,rotate=-30] at ($0.5*(p.corner 1)+0.5*(p.corner 2)$) {NCnode};
                \pic[name=B1,rotate=0] at ($0.5*(p.corner 1)+0.5*(h1)$) {NCnode};
                \pic[name=B2,rotate=0] at ($0.5*(p.corner 1)+0.5*(h2)$) {NCnode};
                \node[vertex] () at ($0.5*(p.corner 1)+0.5*(h2)$) [label ={[xshift=-1.5mm, yshift=4mm]\color{Purple}{1}}]{};
                \pic[name=C,rotate=30] at ($0.5*(p.corner 1)+0.5*(p.corner 3)$) {NCnode};
                \pic[name=D,rotate=90] at ($0.75*(p.corner 3)+0.25*(p.corner 2)$) {NCnode};

                \draw[barrow] (A-circle.40) to[out=40,in=135] (B1-circle.135);
                \draw[barrow] (B1-circle.45) to[out=45,in=135] (B2-circle.135);
                \draw[barrow] (B2-circle.45) to[out=45,in=140] (C-circle.140);
                \draw[barrow] (C-circle.260) to[out=260,in=45] (D-circle.45);
                \draw[barrow] (D-circle.135) to[out=135,in=310] (B2-circle.310);

                \node[mutable] (a1) at ($0.3*(p.corner 1)+0.7*(p.corner 2)$) {};
                \node[mutable] (c1) at ($0.3*(p.corner 1)+0.7*(p.corner 3)$) {};
                \node[mutable] (b11) at ($0.3*(p.corner 1)+0.7*(h1)$) {};
                \node[vertex] () at ($0.3*(p.corner 1)+0.7*(h1)$) [label ={[xshift=-2mm, yshift=-4.5mm]\color{Purple}{3}}]{};
                \node[mutable] (b21) at ($0.3*(p.corner 1)+0.7*(h2)$) {};
                \node[vertex] () at ($0.3*(p.corner 1)+0.7*(h2)$) [label ={[xshift=-1.5mm, yshift=-4.5mm]\color{Purple}{2}}]{};
                \node[mutable] (d1) at ($0.75*(p.corner 2)+0.25*(p.corner 3)$) {};
                
                \node[mutable] (a2) at ($0.15*(p.corner 1)+0.85*(p.corner 2)$) {};
                \node[mutable] (c2) at ($0.15*(p.corner 1)+0.85*(p.corner 3)$) {};
                \node[mutable] (b12) at ($0.15*(p.corner 1)+0.85*(h1)$) {};
                \node[mutable] (b22) at ($0.15*(p.corner 1)+0.85*(h2)$) {};
                \node[mutable] (d2) at ($0.5*(p.corner 2)+0.5*(p.corner 3)$) {};

                \draw[arrow] (a1) to (b11);
                \draw[arrow] (b11) to (b21);
                \draw[arrow] (b21) to (c1);

                \draw[arrow] (a2) to (b12);
                \draw[arrow] (b12) to (b22);
                \draw[arrow] (b22) to (c2);
                
                \draw[arrow] (b12) to (d1);
                \draw[arrow] (d1) to (a2);

                \draw[arrow] (b22) to (d2);
                \draw[arrow] (d2) to (b12);

                \draw[arrow] (b11) to (b12);
                \draw[arrow] (b21) to (b22);
                
                \draw[arrow] (b12) to (a1);
                \draw[arrow] (b22) to (b11);
                \draw[arrow] (c2) to (b21);

                \draw[arrow, color=Green] (B1-circle.center) to (b11);
                \draw[arrow, color=Green] (B2-circle.center) to (b21);
                \draw[arrow, color=Green] (b11) to (A-circle.center);
                \draw[arrow, color=Green] (b21) to (B1-circle.center);
                \draw[arrow, color=Green] (c1) to (B2-circle.center);

                \draw[halfarrow, bend right, color=Green] (A-circle.center) to (a1);
                \draw[halfarrow, bend right] (a1) to (a2);
                \draw[halfarrow, bend left, color=Green] (C-circle.center) to (c1);
                \draw[halfarrow, bend left] (c1) to (c2);
				\draw[halfarrow, bend angle = 25, bend right, color=Green] (d2) to (D-circle.center);
                \draw[halfarrow, bend angle = 25, bend right] (d1) to (d2);
			\end{tikzpicture}}
			\caption{}
           \label{fig:orderedFGQuiverB3}
		\end{subfigure}
	\caption{Structure of an ordered quiver on the quivers $Q_2,Q_3$ from \Cref{fig:FGQuivers}.}
	\label{fig:orderedFGQuivers}
\end{figure}
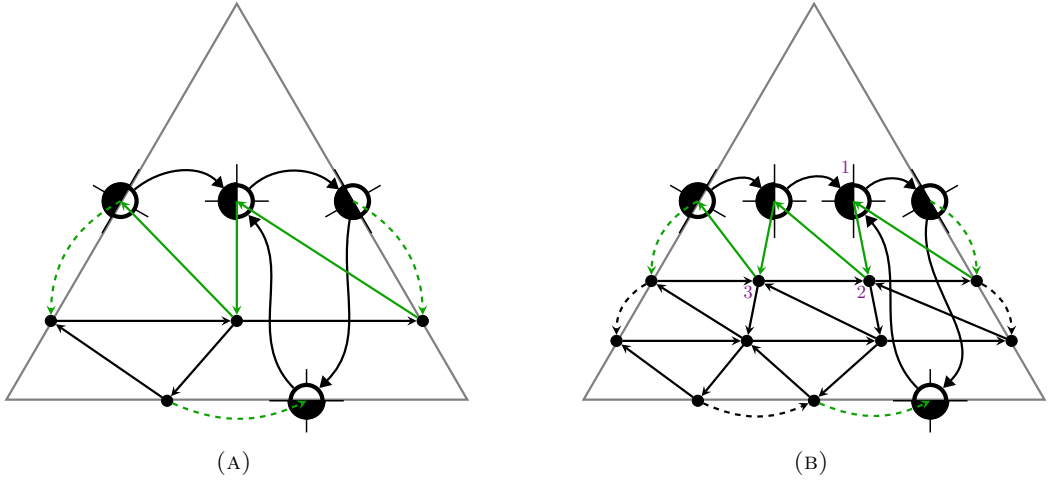


%% file: figureFGTilings.tex
\begin{figure}[ht]
\centering
\begin{subfigure}[b]{0.25\textwidth}
	\centering
	\begin{tikzpicture}[very thick]
 	\node[regular polygon,draw,regular polygon sides=3, minimum size=3.5cm, very thick] (p) at (0,0) {};
    \node[vertex] (h) at ($0.8*(p.side 1)+0.2*(p.corner 3)$) {};
    \draw (p.corner 1) to[out=260,in=60] (h) to[out=240,in=40] (p.corner 2);
    
  	\draw[angleindicator] (h) -- (p.corner 3);
  	\draw[angleindicator] (p.side 2) to[out=80,in=280] (p.corner 1);
  	\draw[angleindicator] (p.side 3) to[out=220,in=20] (p.corner 2);

    \draw[angleindicator] (p.side 1) to[out=260,in=40] (p.corner 2);
    \draw[angleindicator] (h) to[out=210,in=40] (p.corner 2);
   
    \centerarc[straightangle](p.corner 1)(240:300:0.6)
    \centerarc[straightangle](p.corner 2)(0:30:1.1)
    \centerarc[crossangle](p.corner 2)(30:60:1.1)
    \centerarc[straightangle](p.corner 3)(120:180:0.6)
    
    \node[regular polygon,draw,regular polygon sides=3, minimum size=3.5cm, thick] (p) at (0,0) {};
    \draw (p.corner 1) to[out=260,in=60] (h) to[out=240,in=40] (p.corner 2);
	\end{tikzpicture}
\end{subfigure}
\hspace{0.7cm}
\begin{subfigure}[b]{0.25\textwidth}
	\centering
	\begin{tikzpicture}[very thick]
 	\node[regular polygon,draw,regular polygon sides=3, minimum size=3.5cm, very thick] (p) at (0,0) {};
    \node[vertex] (h) at ($0.8*(p.side 1)+0.2*(p.corner 3)$) {};
    \node[vertex] (h2) at ($1.2*(p.side 1)-0.2*(p.corner 3)$) {};
    \draw (p.corner 1) to[out=260,in=60] (h) to[out=240,in=40] (p.corner 2);
    \draw (p.corner 1) to[out=220,in=60] (h2) to[out=240,in=80] (p.corner 2);
    
  	\draw[angleindicator] (h) -- (p.corner 3);
  	\draw[angleindicator] (p.side 2) to[out=80,in=280] (p.corner 1);
  	\draw[angleindicator] (p.side 3) to[out=220,in=20] (p.corner 2);

    \draw[angleindicator] (p.side 1) to[out=260,in=40] (p.corner 2);
    \draw[angleindicator] (h) to[out=210,in=40] (p.corner 2);

    \draw[angleindicator] (p.side 1) to[out=220,in=80] (p.corner 2);
    \draw[angleindicator] (h2) to[out=270,in=80] (p.corner 2);
   
    \centerarc[straightangle](p.corner 1)(220:300:0.6)
    \centerarc[straightangle](p.corner 2)(0:40:1.1)
    \centerarc[crossangle](p.corner 2)(40:80:1.1)
    \centerarc[straightangle](p.corner 3)(120:180:0.6)
    
    \node[regular polygon,draw,regular polygon sides=3, minimum size=3.5cm, thick] (p) at (0,0) {};
    \draw (p.corner 1) to[out=260,in=60] (h) to[out=240,in=40] (p.corner 2);
    \draw (p.corner 1) to[out=220,in=60] (h2) to[out=240,in=80] (p.corner 2);
	\end{tikzpicture}
\end{subfigure}
\caption{Tilings associated to the pruned quivers.}
\label{fig:FGTilings}
\end{figure}
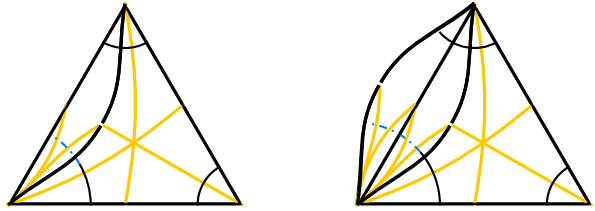

%% file: figureBigPolygonSetup.tex
\begin{figure}[h]
    \centering
    \begin{tikzpicture}
        \node[regular polygon,draw,regular polygon sides=3, minimum size=5cm, color=gray,rotate=-60] (p) at (0,0) {};
        \node[regular polygon,draw,regular polygon sides=3, minimum size=5cm, color=gray,xshift=2.17cm,yshift=-1.25cm, rotate=-120] (q) at (0,0) {};
        \node[regular polygon,draw,regular polygon sides=3, minimum size=5cm, color=gray,xshift=-2.17cm,yshift=-1.25cm, rotate=-120] (r) at (0,0) {};

        \node[mutableBig] (p1) at ($0.6*(p.corner 1)+0.4*(p.corner 2)$)     []{};
        \node[mutableBig] (p2) at ($0.6*(p.corner 1)+0.4*(p.side 2)$)       []{};
        \node[mutableBig] (p3) at ($0.6*(p.corner 1)+0.4*(p.corner 3)$)     [label=below:\color{Purple}{2}]{};
        \node[mutableBig] (p4) at ($0.6*(p.corner 3)+0.4*(p.corner 2)$)     [label=below:\color{Purple}{4}]{};

        \node[mutable] (p5) at ($0.3*(p.corner 1)+0.7*(p.corner 2)$)     []{};
        \node[mutable] (p6) at ($0.3*(p.corner 1)+0.7*(p.side 2)$)       []{};
        \node[mutable] (p7) at ($0.3*(p.corner 1)+0.7*(p.corner 3)$)     [label=below:\color{Purple}{1}]{};
        \node[mutable] (p8) at ($0.3*(p.corner 3)+0.7*(p.corner 2)$)     [label=above:\color{Purple}{3}]{};

        \draw[barrow] (p1) to (p2);
        \draw[barrow] (p2) to (p3);
        \draw[barrow] (p3) to (p4);
        \draw[barrow] (p4) to (p2);

        \draw[arrow] (p5) to (p6);
        \draw[arrow] (p6) to (p7);
        \draw[arrow] (p6) to (p8);
        \draw[arrow] (p8) to (p5);

        \draw[arrow] (p2) to (p6);
        \draw[arrow] (p7) to (p2);
        \draw[arrow] (p6) to (p1);

        \draw[arrow] (p8) to (p4);

        \draw[halfarrow, bend right] (p1) to (p5);

        \node[mutableBig] (q1) at ($0.6*(q.corner 1)+0.4*(q.corner 2)$)     []{};
        \node[mutableBig] (q2) at ($0.6*(q.corner 1)+0.4*(q.side 2)$)       []{};
        \node[mutableBig] (q3) at ($0.6*(q.corner 1)+0.4*(q.corner 3)$)     []{};
        \node[mutableBig] (q4) at ($0.4*(q.corner 3)+0.6*(q.corner 2)$)     []{};

        \node[mutable] (q5) at ($0.3*(q.corner 1)+0.7*(q.corner 2)$)     []{};
        \node[mutable] (q6) at ($0.3*(q.corner 1)+0.7*(q.side 2)+(-0.1,-0.2)$)       []{};
        \node[mutable] (q7) at ($0.3*(q.corner 1)+0.7*(q.corner 3)$)     []{};
        \node[mutable] (q8) at ($0.7*(q.corner 3)+0.3*(q.corner 2)$)     []{};

        \draw[barrow] (q1) to (q2);
        \draw[barrow] (q2) to (q3);
        \draw[barrow] (q2) to (q4);
        \draw[barrow] (q4) to (q1);

        \draw[arrow] (q5) to (q6);
        \draw[arrow] (q6) to (q7);
        \draw[arrow] (q7) to (q8);
        \draw[arrow] (q8) to (q6);

        \draw[arrow] (q6) to (q2);
        \draw[arrow] (q3) to (q6);
        \draw[arrow] (q2) to (q5);

        \draw[arrow] (q4) to (q8);

        \draw[halfarrow, bend right] (q7) to (q3);
        \draw[halfarrow, bend left] (q5) to (q1);

        \node[mutableBig] (r1) at ($0.6*(r.corner 1)+0.4*(r.corner 2)$)     []{};
        \node[mutableBig] (r2) at ($0.6*(r.corner 1)+0.4*(r.side 2)$)       []{};
        \node[mutableBig] (r3) at ($0.6*(r.corner 1)+0.4*(r.corner 3)$)     []{};
        \node[mutableBig] (r4) at ($0.4*(r.corner 3)+0.6*(r.corner 2)$)     []{};

        \node[mutable] (r5) at ($0.3*(r.corner 1)+0.7*(r.corner 2)$)     []{};
        \node[mutable] (r6) at ($0.3*(r.corner 1)+0.7*(r.side 2)+(-0.1,-0.2)$)       []{};
        \node[mutable] (r7) at ($0.3*(r.corner 1)+0.7*(r.corner 3)$)     []{};
        \node[mutable] (r8) at ($0.7*(r.corner 3)+0.3*(r.corner 2)$)     []{};

        \draw[barrow] (r1) to (r2);
        \draw[barrow] (r2) to (r3);
        \draw[barrow] (r2) to (r4);
        \draw[barrow] (r4) to (r1);

        \draw[arrow] (r5) to (r6);
        \draw[arrow] (r6) to (r7);
        \draw[arrow] (r7) to (r8);
        \draw[arrow] (r8) to (r6);

        \draw[arrow] (r6) to (r2);
        \draw[arrow] (r3) to (r6);
        \draw[arrow] (r2) to (r5);

        \draw[halfarrow, bend right] (r4) to (r8);
        \draw[halfarrow, bend right] (r7) to (r3);
    \end{tikzpicture}
    \caption{Gluing three copies of $Q_2$.}
    \label{fig:bigPolygonSetup}
\end{figure}
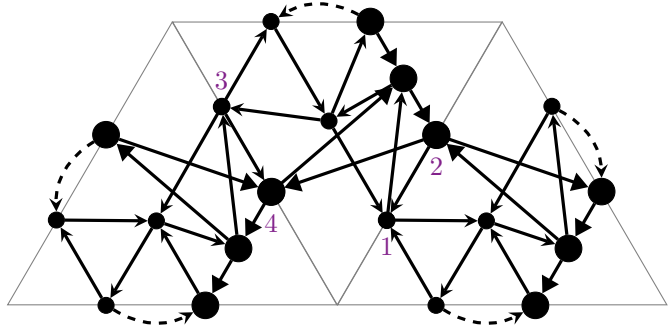

%% file: figureBigPolygonEx.tex
\begin{figure}[ht]
    \centering
    \begin{tikzpicture}
        \node[regular polygon,draw,regular polygon sides=3, minimum size=5cm, color=gray,rotate=-60] (p) at (0,0) {};
        \node[regular polygon,draw,regular polygon sides=3, minimum size=5cm, color=gray,xshift=2.17cm,yshift=-1.25cm, rotate=-120] (q) at (0,0) {};
        \node[regular polygon,draw,regular polygon sides=3, minimum size=5cm, color=gray,xshift=-2.17cm,yshift=-1.25cm, rotate=-120] (r) at (0,0) {};

        \node[mutableBig] (p1) at ($0.6*(p.corner 1)+0.4*(p.corner 2)$)     []{};
        \node[mutableBig,fill=Pink] (p2) at ($0.6*(p.corner 1)+0.4*(p.side 2)$)       []{};
        \node[mutableBig,fill=Pink] (p3) at ($0.6*(p.corner 1)+0.4*(p.corner 3)$)     []{};
        \node[mutableBig,fill=Pink] (p4) at ($0.6*(p.corner 3)+0.4*(p.corner 2)$)     [label ={[xshift=5.4mm, yshift=1.2mm] \footnotesize{2}}]{};

        \node[mutable] (p5) at ($0.3*(p.corner 1)+0.7*(p.corner 2)$)     []{};
        \node[mutable] (p6) at ($0.3*(p.corner 1)+0.7*(p.side 2)+(-0.2,0.3)$)       []{};
        \node[mutable] (p7) at ($0.3*(p.corner 1)+0.7*(p.corner 3)$)     []{};
        \node[mutable] (p8) at ($0.3*(p.corner 3)+0.7*(p.corner 2)$)     []{};

        \node[mutableBig] (q1) at ($0.6*(q.corner 1)+0.4*(q.corner 2)$)     []{};
        \node[mutableBig] (q2) at ($0.6*(q.corner 1)+0.4*(q.side 2)$)       []{};
        \node[mutableBig] (q3) at ($0.6*(q.corner 1)+0.4*(q.corner 3)$)     []{};
        \node[mutableBig,fill=Pink] (q4) at ($0.4*(q.corner 3)+0.6*(q.corner 2)$)     []{};

        \node[mutable] (q5) at ($0.3*(q.corner 1)+0.7*(q.corner 2)$)     []{};
        \node[mutable] (q6) at ($0.3*(q.corner 1)+0.7*(q.side 2)+(-0.0,-0.4)$)       []{};
        \node[mutable] (q7) at ($0.3*(q.corner 1)+0.7*(q.corner 3)$)     []{};
        \node[mutable] (q8) at ($0.7*(q.corner 3)+0.3*(q.corner 2)$)     []{};

        \node[mutableBig,fill=Pink] (r1) at ($0.6*(r.corner 1)+0.4*(r.corner 2)$)     []{};
        \node[mutableBig,fill=Pink] (r2) at ($0.6*(r.corner 1)+0.4*(r.side 2)$)       []{};
        \node[mutableBig] (r3) at ($0.6*(r.corner 1)+0.4*(r.corner 3)$)     []{};
        \node[mutableBig,fill=Pink] (r4) at ($0.4*(r.corner 3)+0.6*(r.corner 2)$)     []{};

        \node[mutable] (r5) at ($0.3*(r.corner 1)+0.7*(r.corner 2)$)     []{};
        \node[mutable] (r6) at ($0.3*(r.corner 1)+0.7*(r.side 2)+(-0.3,-0.4)$)       []{};
        \node[mutable] (r7) at ($0.3*(r.corner 1)+0.7*(r.corner 3)$)     []{};
        \node[mutable] (r8) at ($0.7*(r.corner 3)+0.3*(r.corner 2)$)     []{};

        \draw[barrow] (p1) to (p2);

        \draw[barrow] (p2) to [out=210,in=60] (p4);
        \draw[xarrow] (p2) to (p3);

        \draw[barrow] (p3) to (p4);

        \draw[barrow] (p4) to (r4);
        \draw[xarrow] (p4) to (r2);
        \draw[barrow] (p4) to [out=0,in=130] (q2);
        \draw[arrow] (p4) to (p7);

        \draw[arrow] (p5) to (p8);

        \draw[arrow] (p6) to (p1);
        \draw[arrow] (p6) to (p2);
        \draw[arrow] (p6) to [out=-10,in=150] (p3);
        \draw[arrow] (p6) to (q6);
        \draw[arrow] (p6) to [out=230,in=20] (r6);
        \draw[arrow] (p6) to [out=170,in=30] (p8);

        \draw[arrow] (p7) to (p8);
        \draw[arrow] (p7) to (p6);
        \draw[arrow] (p7) to (p2);
        \draw[arrow] (p7) to (q1);
        \draw[arrow] (p7) to (q6);
        \draw[arrow] (p7) to (q7);

        \draw[arrow] (p8) to (p4);

        \draw[barrow] (q1) to (p3);

        \draw[arrow] (q2) to (q5);
        \draw[carrow] (q2) to [out=110,in=-10] (p2);
        \draw[arrow] (q2) to [out=150,in=-5] (p8);
        \draw[barrow] (q2) to (q3);

        \draw[arrow] (q3) to (q6);

        \draw[arrow] (q5) to (q6);

        \draw[arrow] (q6) to (p3);

        \draw[arrow] (q7) to [out=55,in=-30] (p2);

        \draw[barrow] (r2) to (r3);
        \draw[barrow] (r2) to [out=160,in=-60] (r4);
        \draw[xarrow] (r2) to [out=30,in=230] (p2);
        \draw[xarrow] (r2) to [out=10,in=230] (p3);
        \draw[arrow] (r2) to (p5);

        \draw[arrow] (r3) to (r6);

        \draw[barrow] (r4) to [out=40,in=170] (p2);
        \draw[barrow] (r4) to (p3);
        \draw[arrow] (r4) to (p8);

        \draw[arrow] (r6) to (p8);
        \draw[arrow] (r6) to (r7);

        \draw[arrow] (r7) to (r8);

        \draw[arrow] (r8) to (r6);

        \draw[halfarrow, bend right] (p1) to (p5);
        \draw[halfarrow, bend right] (q7) to (q3);
        \draw[halfarrow, bend left] (q5) to (q1);
        \draw[halfarrow, bend right] (r4) to (r8);
        \draw[halfarrow, bend right] (r7) to (r3);
    \end{tikzpicture}
    \caption{Quiver from \Cref{fig:bigPolygonSetup} after mutations.}
    \label{fig:bigPolygonEx}
\end{figure}
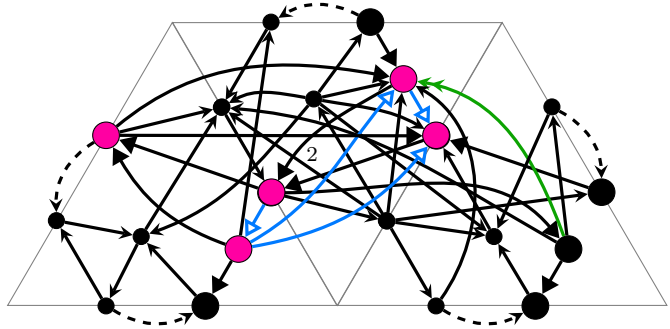

%% file: figureB2mutation.tex
\begin{figure}[htb]
	\centering
		\begin{subfigure}[c]{0.45\textwidth}
			\centering
			\scalebox{0.8}{
			\begin{tikzpicture}[node distance = 50,xscale=0.577,scale=1.5]
				\node[vertex]		(c1)	at (0,0)							[]{};
				\node[vertex]		(c2)	at (-3*\distance,-3*\distance)	[]	{};
		        \node[vertex]		(c3)	at (3*\distance,-3*\distance)	[]	{};
				
				\draw[line] (c1) to (c2);
				\draw[line] (c2) to (c3);
				\draw[line] (c3) to (c1);
					
				\node[frozenBig]		(1)		at (	-\distance,-\distance)	[label = above left: 1]	{};
				\node[mutableBig]	(2)		at (0,-\distance)			[label = above: 2]{};
				\node[frozenBig]		(3)		at (\distance,-\distance)	[label = above right: 3]{};
				\node[frozen]		(4)		at (-2*\distance,-2*\distance)	[label = left: 4]	{};
				\node[mutable]		(5)		at (0,-2*\distance)				[label = below right: 5]	{};
			    \node[frozen]		(6)		at (2*\distance,-2*\distance)	[label = right: 6]{};
				\node[frozen]		(7)		at (-\distance,-3*\distance)		[label = below: 7]{};
				\node[frozenBig]		(8)		at (\distance,-3*\distance)		[label = below: 8]{};
				
				\draw[barrow] (1) to (2);
				\draw[barrow] (2) to (3);
				\draw[arrow] (4) to (5);
				\draw[arrow] (5) to (6);
				\draw[arrow] (5) to (1);
				\draw[arrow] (6) to (2);
				\draw[arrow] (2) to (5);
				\draw[barrow] (3) to (8);
				\draw[barrow] (8) to (2);
				\draw[arrow] (5) to (7);
				\draw[arrow] (7) to (4);
				
				\draw[halfarrow, bend left] (3) to (6);
				\draw[halfarrow, bend angle = 12, bend right] (7) to (8);
				\draw[halfarrow, bend right] (1) to (4);
			\end{tikzpicture}}
		\end{subfigure}
		\hfill
		\begin{subfigure}[c]{0.08\textwidth}
			\centering
			\begin{tikzpicture}
				\draw[->,very thick] (-0.7,0) to (0.7,0)
				node[above,midway,text centered]{$\mu_2$};
			\end{tikzpicture}
		\end{subfigure}
		\hfill
		\begin{subfigure}[c]{0.45\textwidth}
			\centering
			\scalebox{0.8}{
			\begin{tikzpicture}[node distance = 50,xscale=0.577,scale=1.5]
				\node[vertex]		(c1)	at (0,0)							[]{};
				\node[vertex]		(c2)	at (-3*\distance,-3*\distance)	[]	{};
		        \node[vertex]		(c3)	at (3*\distance,-3*\distance)	[]	{};
				
				\draw[line] (c1) to (c2);
				\draw[line] (c2) to (c3);
				\draw[line] (c3) to (c1);
					
				\node[frozenBig]		(1)		at (	-\distance,-\distance)	[label = above left: 1]	{};
				\node[mutableBig]	(2)		at (0,-1.2*\distance)			[label = above: 2]{};
				\node[frozenBig]		(3)		at (\distance,-\distance)	[label = above right: 3]{};
				\node[frozen]		(4)		at (-2*\distance,-2*\distance)	[label = left: 4]	{};
				\node[mutable]		(5)		at (0,-2*\distance)				[label = above left: 5]	{};
			    \node[frozen]		(6)		at (2*\distance,-2*\distance)	[label = right: 6]{};
				\node[frozen]		(7)		at (-\distance,-3*\distance)		[label = below: 7]{};
				\node[frozenBig]		(8)		at (\distance,-3*\distance)		[label = below: 8]{};
				
				\draw[barrow] (2) to (1);
				\draw[barrow] (3) to (2);
				\draw[barrow, bend angle=15, bend left] (1) to (3);
				\draw[arrow] (4) to (5);
				\draw[arrow] (6) to (5);
				\draw[arrow] (2) to (6);
				\draw[arrow] (5) to (2);
				\draw[arrow] (8) to (5);
				\draw[barrow, bend left] (2) to (8);
				\draw[arrow] (5) to (7);
				\draw[arrow] (7) to (4);
				
				\draw[halfarrow, bend right] (6) to (3);
				\draw[halfarrow, bend angle = 12, bend right] (7) to (8);
				\draw[halfarrow, bend right] (1) to (4);
			\end{tikzpicture}}
		\end{subfigure}\\
		\begin{subfigure}[c]{0.45\textwidth}
			\centering
			\begin{tikzpicture}
			\end{tikzpicture}
		\end{subfigure}
		\hfill
		\begin{subfigure}[c]{0.08\textwidth}
			
		\end{subfigure}
		\hfill
		\begin{subfigure}[c]{0.45\textwidth}
			\centering
			\vspace{5mm}
			\hspace{7mm}
			\begin{tikzpicture}
				\draw[->, very thick] (0,0) to (0,-1.4);
				\node[vertex] () at (0,-0.7) [label = right:$\mu_5$]{};
			\end{tikzpicture}
			\vspace{5mm}
		\end{subfigure}\\
		\begin{subfigure}[c]{0.45\textwidth}
			\centering
			\scalebox{0.8}{
			\begin{tikzpicture}[node distance = 50,xscale=0.577,scale=1.5]
				\node[vertex]		(c1)	at (0,0)							[]{};
				\node[vertex]		(c2)	at (-3*\distance,-3*\distance)	[]	{};
		        \node[vertex]		(c3)	at (3*\distance,-3*\distance)	[]	{};
				
				\draw[line] (c1) to (c2);
				\draw[line] (c2) to (c3);
				\draw[line] (c3) to (c1);
					
				\node[frozenBig]		(1)		at (	-\distance,-\distance)	[label = above left: 1]	{};
				\node[mutableBig]	(2)		at (0,-\distance)			[label = above: 2]{};
				\node[frozenBig]		(3)		at (\distance,-\distance)	[label = above right: 3]{};
				\node[frozen]		(4)		at (-2*\distance,-2*\distance)	[label = left: 4]	{};
				\node[mutable]		(5)		at (0,-2*\distance)				[label = below right: 5]	{};
			    \node[frozen]		(6)		at (2*\distance,-2*\distance)	[label = right: 6]{};
				\node[frozen]		(7)		at (-\distance,-3*\distance)		[label = below: 7]{};
				\node[frozenBig]		(8)		at (\distance,-3*\distance)		[label = below: 8]{};
				
				\draw[barrow] (2) to (1);
				\draw[xarrow] (2) to (3);
				\draw[arrow] (5) to (4);
				\draw[arrow] (6) to (5);
				\draw[arrow, bend angle = 15, bend right] (4) to (6);
				\draw[arrow] (5) to (2);
				\draw[arrow] (7) to (5);
				\draw[xarrow, bend angle = 12, bend right] (8) to (2);
				\draw[barrow] (3) to (8);
				\draw[arrow] (1) to (5);
				\draw[arrow,bend angle = 17, bend right] (2) to (7);
				
				\draw[halfarrow, bend left] (3) to (6);
				\draw[halfarrow, bend angle = 12, bend right] (7) to (8);
				\draw[halfarrow, bend left] (4) to (1);
			\end{tikzpicture}}
		\end{subfigure}
		\hfill
		\begin{subfigure}[c]{0.08\textwidth}
			\centering
			\begin{tikzpicture}
				
			\end{tikzpicture}
		\end{subfigure}
		\hfill
		\begin{subfigure}[c]{0.45\textwidth}
			\centering
			\scalebox{0.8}{
			\begin{tikzpicture}[node distance = 50,xscale=0.577,scale=1.5]
				\node[vertex]		(c1)	at (0,0)							[]{};
				\node[vertex]		(c2)	at (-3*\distance,-3*\distance)	[]	{};
		        \node[vertex]		(c3)	at (3*\distance,-3*\distance)	[]	{};
				
				\draw[line] (c1) to (c2);
				\draw[line] (c2) to (c3);
				\draw[line] (c3) to (c1);
					
				\node[frozenBig]		(1)		at (	-\distance,-\distance)	[label = above left: 1]	{};
				\node[mutableBig]	(2)		at (0,-1.2*\distance)			[label = above: 2]{};
				\node[frozenBig]		(3)		at (\distance,-\distance)	[label = above right: 3]{};
				\node[frozen]		(4)		at (-2*\distance,-2*\distance)	[label = left: 4]	{};
				\node[mutable]		(5)		at (0,-2*\distance)				[label = above left: 5]	{};
			    \node[frozen]		(6)		at (2*\distance,-2*\distance)	[label = right: 6]{};
				\node[frozen]		(7)		at (-\distance,-3*\distance)		[label = below: 7]{};
				\node[frozenBig]		(8)		at (\distance,-3*\distance)		[label = below: 8]{};
				
				\draw[barrow] (2) to (1);
				\draw[barrow] (3) to (2);
				\draw[barrow, bend angle=15, bend left] (1) to (3);
				\draw[arrow] (5) to (4);
				\draw[arrow] (5) to (6);
				\draw[arrow] (6) to (7);
				\draw[arrow] (2) to (5);
				\draw[arrow] (5) to (8);
				\draw[xarrow, bend right] (8) to (2);
				\draw[arrow] (7) to (5);
				\draw[arrow] (4) to (2);
				
				\draw[halfarrow, bend right] (6) to (3);
				\draw[halfarrow, bend angle = 12, bend left] (8) to (7);
				\draw[halfarrow, bend right] (1) to (4);
			\end{tikzpicture}}
		\end{subfigure}
		\begin{subfigure}[c]{0.45\textwidth}
			\centering
			\vspace{5mm}
			\begin{tikzpicture}
				\draw[<-, very thick] (0,0) to (0,-1.4);
				\node[vertex] () at (0,-0.7) [label = left:$\mu_2$]{};
			\end{tikzpicture}
			\hspace{4mm}
			\vspace{5mm}
		\end{subfigure}
		\hfill
		\begin{subfigure}[c]{0.08\textwidth}
			
		\end{subfigure}
		\hfill
		\begin{subfigure}[c]{0.45\textwidth}
			\centering
			\vspace{5mm}
			\hspace{7mm}
			\begin{tikzpicture}
				\draw[->, very thick] (0,0) to (0,-1.4);
				\node[vertex] () at (0,-0.7) [label = right:$\mu_2$]{};
			\end{tikzpicture}
			\vspace{5mm}
		\end{subfigure}\\
		\begin{subfigure}[c]{0.45\textwidth}
			\centering
			\scalebox{0.8}{
			\begin{tikzpicture}[node distance = 50,xscale=0.577,scale=1.5]
				\node[vertex]		(c1)	at (0,0)							[]{};
				\node[vertex]		(c2)	at (-3*\distance,-3*\distance)	[]	{};
		        \node[vertex]		(c3)	at (3*\distance,-3*\distance)	[]	{};
				
				\draw[line] (c1) to (c2);
				\draw[line] (c2) to (c3);
				\draw[line] (c3) to (c1);
					
				\node[frozenBig]		(1)		at (	-\distance,-\distance)	[label = above left: 1]	{};
				\node[mutableBig]	(2)		at (0,-\distance)			[label = above: 2]{};
				\node[frozenBig]		(3)		at (\distance,-\distance)	[label = above right: 3]{};
				\node[frozen]		(4)		at (-2*\distance,-2*\distance)	[label = left: 4]	{};
				\node[mutable]		(5)		at (0,-2*\distance)				[label = above left: 5]	{};
			    \node[frozen]		(6)		at (2*\distance,-2*\distance)	[label = right: 6]{};
				\node[frozen]		(7)		at (-\distance,-3*\distance)		[label = below: 7]{};
				\node[frozenBig]		(8)		at (\distance,-3*\distance)		[label = below: 8]{};
				
				\draw[barrow] (1) to (2);
				\draw[xarrow] (3) to (2);
				\draw[arrow] (5) to (4);
				\draw[arrow] (6) to (5);
				\draw[arrow, bend angle = 15, bend right] (4) to (6);
				\draw[arrow] (2) to (5);
				\draw[arrow] (5) to (7);
				\draw[arrow] (5) to (3);
				\draw[xarrow, bend left] (2) to (8);
				\draw[xarrow, bend left] (8) to (1);
				\draw[arrow,bend angle = 15, bend left] (7) to (2);
				
				\draw[halfarrow, bend left] (3) to (6);
				\draw[halfarrow, bend angle = 12, bend left] (8) to (7);
				\draw[halfarrow, bend left] (4) to (1);
			\end{tikzpicture}}
		\end{subfigure}
		\hfill
		\begin{subfigure}[c]{0.08\textwidth}
			\centering
			\begin{tikzpicture}
				\draw[<-,very thick] (-0.7,0) to (0.7,0)
				node[above,midway,text centered]{$\mu_5$};
			\end{tikzpicture}
		\end{subfigure}
		\hfill
		\begin{subfigure}[c]{0.45\textwidth}
			\centering
			\scalebox{0.8}{
			\begin{tikzpicture}[node distance = 50,xscale=0.577,scale=1.5]
				\node[vertex]		(c1)	at (0,0)							[]{};
				\node[vertex]		(c2)	at (-3*\distance,-3*\distance)	[]	{};
		        \node[vertex]		(c3)	at (3*\distance,-3*\distance)	[]	{};
				
				\draw[line] (c1) to (c2);
				\draw[line] (c2) to (c3);
				\draw[line] (c3) to (c1);
					
				\node[frozenBig]		(1)		at (	-\distance,-\distance)	[label = above left: 1]	{};
				\node[mutableBig]	(2)		at (0,-\distance)			[label = above: 2]{};
				\node[frozenBig]		(3)		at (\distance,-\distance)	[label = above right: 3]{};
				\node[frozen]		(4)		at (-2*\distance,-2*\distance)	[label = left: 4]	{};
				\node[mutable]		(5)		at (0,-2*\distance)				[label = above left: 5]	{};
			    \node[frozen]		(6)		at (2*\distance,-2*\distance)	[label = right: 6]{};
				\node[frozen]		(7)		at (-\distance,-3*\distance)		[label = below: 7]{};
				\node[frozenBig]		(8)		at (\distance,-3*\distance)		[label = below: 8]{};
				
				\draw[barrow] (1) to (2);
				\draw[barrow] (2) to (3);
				\draw[arrow] (4) to (5);
				\draw[arrow] (5) to (6);
				\draw[arrow] (6) to (7);
				\draw[arrow] (5) to (2);
				\draw[arrow] (3) to (5);
				\draw[xarrow, bend left] (2) to (8);
				\draw[xarrow, bend left] (8) to (1);
				\draw[arrow] (7) to (5);
				\draw[arrow] (2) to (4);
				
				\draw[halfarrow, bend right] (6) to (3);
				\draw[halfarrow, bend angle = 12, bend left] (8) to (7);
				\draw[halfarrow, bend left] (4) to (1);
			\end{tikzpicture}}
		\end{subfigure}
	\caption{Mutation class of the quiver in \Cref{fig:FGQuivers}(A).}
	\label{fig:B2mutation}
\end{figure}
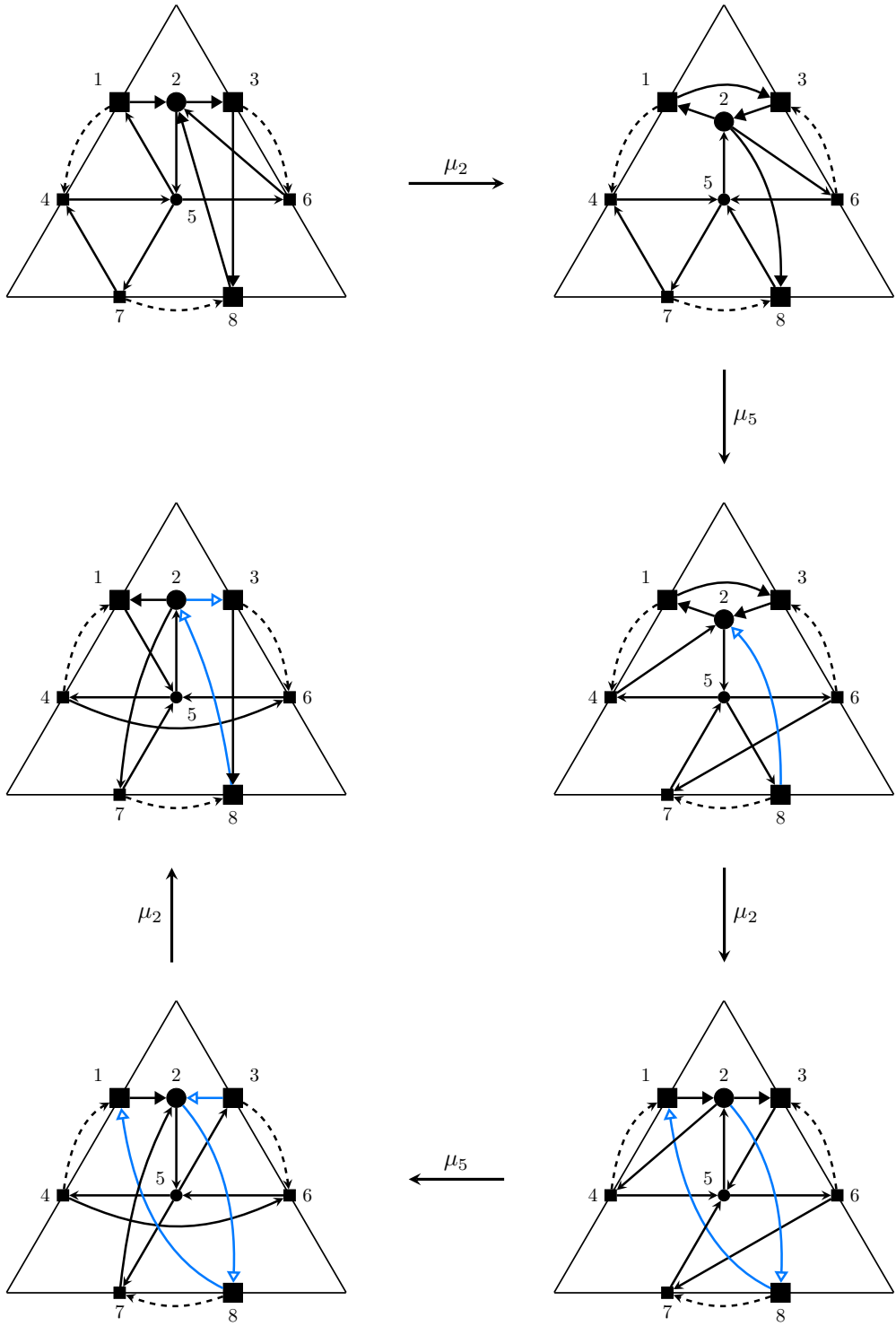

%% file: main.bbl
\newcommand{\etalchar}[1]{$^{#1}$}
\begin{thebibliography}{AGRW19}

\bibitem[ABR{\etalchar{+}}22]{alessandrini2022-SP2Asigma}
Daniele Alessandrini, Arkady Berenstein, Vladimir Retakh, Eugen Rogozinnikov, and Anna Wienhard.
\newblock Symplectic groups over noncommutative algebras.
\newblock {\em Selecta Math. (N.S.)}, 28(4):Paper No. 82, 119, 2022.

\bibitem[ABR23]{albers_floer_2023}
Peter Albers, Maria Bertozzi, and Markus Reineke.
\newblock Floer potentials, cluster algebras and quiver representations, September 2023.
\newblock arXiv:2309.16009 [math] version: 1.

\bibitem[AGRW19]{alessandrini2019noncommutative}
Daniele Alessandrini, Olivier Guichard, Evgenii Rogozinnikov, and Anna Wienhard.
\newblock Noncommutative coordinates for symplectic representations.
\newblock {\em arXiv preprint arXiv:1911.08014}, 2019.

\bibitem[BGL{\etalchar{+}}24]{beyrer2024}
Jonas Beyrer, Olivier Guichard, Fran{\c{c}}ois Labourie, Beatrice Pozzetti, and Anna Wienhard.
\newblock Positivity, cross-ratios and the collar lemma.
\newblock {\em arXiv preprint arXiv:2409.06294}, 2024.

\bibitem[BGM18]{bershtein_cluster_2018}
M.~Bershtein, P.~Gavrylenko, and A.~Marshakov.
\newblock Cluster integrable systems, q-{Painlevé} equations and their quantization.
\newblock {\em Journal of High Energy Physics}, 2018(2):77, February 2018.

\bibitem[BP21]{beyrer2021}
Jonas Beyrer and Beatrice Pozzetti.
\newblock Positive surface group representations in po (p, q).
\newblock {\em arXiv preprint arXiv:2106.14725}, 2021.

\bibitem[BR18]{berenstein2018noncommutative}
Arkady Berenstein and Vladimir Retakh.
\newblock Noncommutative marked surfaces.
\newblock {\em Advances in Mathematics}, 328:1010--1087, 2018.

\bibitem[Bro80]{brouwer1980enumeration}
Andries~E Brouwer.
\newblock The enumeration of locally transitive tournaments.
\newblock {\em Stichting Mathematisch Centrum. Zuivere Wiskunde}, 138(ZW 138/80), 1980.

\bibitem[FG06]{fock2006moduli}
Vladimir Fock and Alexander Goncharov.
\newblock Moduli spaces of local systems and higher {T}eichm{\"u}ller theory.
\newblock {\em Publications Math{\'e}matiques de l'IH{\'E}S}, 103:1--211, 2006.

\bibitem[FST08]{fomin2008clusterTri}
Sergey Fomin, Michael Shapiro, and Dylan Thurston.
\newblock Cluster algebras and triangulated surfaces. {P}art {I}: Cluster complexes.
\newblock {\em Acta Mathematica}, 201, 2008.

\bibitem[FST12]{felikson2012cluster}
Anna Felikson, Michael Shapiro, and Pavel Tumarkin.
\newblock Cluster algebras of finite mutation type via unfoldings.
\newblock {\em International Mathematics Research Notices}, 2012(8):1768--1804, 2012.

\bibitem[FWZ17]{fomin2017introduction}
Sergey Fomin, Lauren Williams, and Andrei Zelevinsky.
\newblock Introduction to cluster algebras, chapters {I}--{III}.
\newblock {\em arXiv preprint arXiv:1608.05735}, 2017.

\bibitem[FZ02]{fomin2002cluster}
Sergey Fomin and Andrei Zelevinsky.
\newblock Cluster algebras {I}: Foundations.
\newblock {\em Journal of the American mathematical society}, 15(2):497--529, 2002.

\bibitem[Gil21]{gilles2021fock}
S~Gilles.
\newblock {\em Fock--Goncharov coordinates for semisimple Lie groups}.
\newblock PhD thesis, University of Maryland, College Park, 2021.

\bibitem[GK21]{goncharov2021spectral}
Alexander Goncharov and Maxim Kontsevich.
\newblock Spectral description of non-commutative local systems on surfaces and non-commutative cluster varieties.
\newblock {\em arXiv preprint arXiv:2108.04168}, 2021.

\bibitem[GKW24]{GreenbergEtAl2024MathrmSL_2}
Zachary Greenberg, Dani Kaufman, and Anna Wienhard.
\newblock $\mathrm{SL}_2$-like properties of matrices over noncommutative rings and generalizations of markov numbers.
\newblock {\em arXiv preprint arXiv:2402.19300}, 2024.

\bibitem[GLW21]{guichard2021}
Olivier Guichard, Fran{\c{c}}ois Labourie, and Anna Wienhard.
\newblock Positivity and representations of surface groups.
\newblock {\em arXiv preprint arXiv:2106.14584}, 2021.

\bibitem[GS19]{goncharov2019quantum}
Alexander Goncharov and Linhui Shen.
\newblock Quantum geometry of moduli spaces of local systems and representation theory.
\newblock {\em arXiv preprint arXiv:1904.10491}, 2019.

\bibitem[GW22]{guichard2022generalizing}
Olivier Guichard and Anna Wienhard.
\newblock Generalizing {L}usztig's total positivity.
\newblock {\em arXiv preprint arXiv:2208.10114}, 2022.

\bibitem[HS08]{hone_integrality_2008}
Andrew N.~W. Hone and Christine Swart.
\newblock Integrality and the {Laurent} phenomenon for {Somos} 4 and {Somos} 5 sequences.
\newblock {\em Mathematical Proceedings of the Cambridge Philosophical Society}, 145(1):65--85, July 2008.

\bibitem[Kau23]{kaufman2023special}
Dani Kaufman.
\newblock Special folding of quivers and cluster algebras.
\newblock {\em arXiv preprint arXiv:2304.07510}, 2023.

\bibitem[KR22]{kineider2022-FramedLocalSystems}
Clarence Kineider and Eugen Rogozinnikov.
\newblock On partial abelianization of framed local systems, 2022.

\bibitem[Le19]{le2019cluster}
Ian Le.
\newblock Cluster structures on higher teichmuller spaces for classical groups.
\newblock In {\em Forum of Mathematics, Sigma}, volume~7, page e13. Cambridge University Press, 2019.

\bibitem[Mar13]{marsh2013lectureOnClusterAlgebras}
Robert~J. Marsh.
\newblock {\em Lecture notes on cluster algebras}.
\newblock Zurich {Lectures} in {Advanced} {Mathematics}. European Mathematical Society (EMS), Zürich, 2013.

\bibitem[OS19]{os-ClusterAlgebrasWithGrassmanVariables}
Valentin Ovsienko and Michael Shapiro.
\newblock Cluster algebras with {G}rassmann variables.
\newblock {\em Electron. Res. Announc. Math. Sci.}, 26:1--15, 2019.

\bibitem[Pen87]{penner1987decorated}
Robert~C Penner.
\newblock The decorated {T}eichm{\"u}ller space of punctured surfaces.
\newblock {\em Communications in Mathematical Physics}, 113:299--339, 1987.

\bibitem[PT20]{pascaleff_wall-crossing_2020}
James Pascaleff and Dmitry Tonkonog.
\newblock The wall-crossing formula and {Lagrangian} mutations.
\newblock {\em Advances in Mathematics}, 361, February 2020.

\bibitem[Rog20]{rogozinnikov2020symplectic}
Evgenii Rogozinnikov.
\newblock {\em Symplectic groups over noncommutative rings and maximal representations}.
\newblock PhD thesis, Universit\"at Heidelberg, 2020.

\bibitem[She23]{shemyakova-superClusterAlgebras}
Ekaterina Shemyakova.
\newblock On super cluster algebras based on super {P}l\"ucker and super {P}tolemy relations.
\newblock {\em J. Geom. Phys.}, 188:Paper No. 104776, 16, 2023.

\bibitem[War14]{warkentin2014exchange}
Matthias Warkentin.
\newblock {\em Exchange graphs via quiver mutation}.
\newblock PhD thesis, Technische Universit\"at Chemintz, 2014.

\bibitem[Zic20]{zickert2020RankTwoLieGroups}
Christian~K. Zickert.
\newblock Fock-{Goncharov} coordinates for rank two {Lie} groups.
\newblock {\em Mathematische Zeitschrift}, 294(1-2):251--286, 2020.

\end{thebibliography}
